\documentclass[11pt]{article}
\usepackage{latexsym,amssymb,amsfonts}
\usepackage{mathrsfs}
\usepackage{psfrag}
\usepackage{amsmath, amsbsy}
\usepackage{amsopn, amstext}
\def\sqr#1#2{{\vcenter{\vbox{\hrule height.#2pt
              \hbox{\vrule width.#2pt height#1pt \kern#1pt \vrule width.#2pt}
              \hrule height.#2pt}}}}
\def\signed #1{{\unskip\nobreak\hfil\penalty50
              \hskip2em\hbox{}\nobreak\hfil#1
              \parfillskip=0pt \finalhyphendemerits=0 \par}}
\def\endpf{\signed {$\sqr69$}}

\def\dbR{{\mathop{\rm l\negthinspace R}}}
\def\dbC{{\mathop{\rm l\negthinspace\negthinspace\negthinspace C}}}

\def\1n{\negthinspace }
\def\2n{\negthinspace \negthinspace }
\def\3n{\negthinspace \negthinspace \negthinspace }

\def\ns{\noalign{\ss}}

\def\dbC{{\mathbb{C}}}

\def\dbE{{\mathbb{E}}}
\def\dbF{{\mathbb{F}}}

\def\dbH{{\mathbb{H}}}

\def\dbN{{\mathbb{N}}}

\def\dbP{{\mathbb{P}}}

\def\dbR{{\mathbb{R}}}
\def\dbS{{\mathbb{S}}}

\def\dbm{\mathbf{m}}

\def\limsup{\mathop{\overline{\rm lim}}}
\def\liminf{\mathop{\underline{\rm lim}}}

\def\lan{\langle}
\def\ran{\rangle}

\def\div{\mathop{\nabla\cd\neg}}
\def\esssup{\mathop{\rm esssup}}

\def\pa{\partial}
\def\h{\widehat}
\def\wt{\widetilde}
\def\cd{\cdot}
\def\cds{\cdots}

\def\dim{\hbox{\rm dim$\,$}}

\def\ae{\hbox{\rm a.e.{ }}}
\def\as{\hbox{\rm a.s.{ }}}

\def\supp{\hbox{\rm supp$\,$}}

\def\cl{\overline}

\def\deq{\mathop{\buildrel\D\over=}}

\def\({\Big (}
\def\){\Big )}
\def\[{\Big[}
\def\]{\Big]}

\def\={\buildrel \triangle \over =}
\def\wh{\widehat}

\def\resp{{\it resp. }}

\def\-{\mbox{-}}
\def\={\buildrel \triangle \over =}

\def\resp{{\it resp. }}
\def\ds{\displaystyle}
\def\mE{\mathbb{E}}
\def\a{\alpha}
\def\b{\beta}
\def\g{\gamma}
\def\d{\delta}
\def\e{\varepsilon}

\def\k{\kappa}
\def\l{\lambda}

 \def\n{\nabla}
\def\si{\sigma}
\def\t{\times}
\def\f{\varphi}
\def\th{\theta}
\def\om{\omega}

\def\mR{\mathbb{R}}

\def\G{\Gamma}
\def\D{\Delta}
\def\Th{\Theta}
\def\L{\Lambda}
\def\Si{\Sigma}

\def\Om{\Omega}

\def\cA{{\cal A}}
\def\cB{{\cal B}}
\def\cC{{\cal C}}
\def\cD{{\cal D}}
\def\cE{{\cal E}}
\def\cF{{\cal F}}

\def\cG{{\cal G}}
\def\cH{{\cal H}}
\def\cI{{\cal I}}
\def\cJ{{\cal J}}

\def\cL{{\cal L}}
\def\cM{{\cal M}}
\def\cN{{\cal N}}

\def\cP{{\cal P}}
\def\cQ{{\cal Q}}

\def\cS{{\cal S}}

\def\cU{{\cal U}}
\def\cV{{\cal V}}

\def\cX{{\cal X}}
\def\cY{{\cal Y}}
\def\cZ{{\cal Z}}
\def\cl{{\cal l}}
\def\no{\noindent}

\def\ss{\smallskip}
\def\ms{\medskip}
\def\bs{\bigskip}
\def\q{\quad}
\def\qq{\qquad}
\def\hb{\hbox}

\def\limsup{\mathop{\overline{\rm lim}}}
\def\liminf{\mathop{\underline{\rm lim}}}

\def\esssup{\mathop{\rm esssup}}

\def\pa{\partial}
\def\h{\widehat}
\def\wh{\widehat}
\def\wt{\widetilde}

\def\cd{\cdot}
\def\cds{\cdots}

\def\dim{\hbox{\rm dim$\,$}}

\def\ae{\hbox{\rm a.e.{ }}}
\def\as{\hbox{\rm a.s.}}

\def\supp{\hbox{\rm supp$\,$}}

\def\cl{\overline}

\def\deq{\mathop{\buildrel\D\over=}}

\def\|{\Big |}
\def\({\Big (}
\def\){\Big )}
\def\[{\Big[}
\def\]{\Big]}

\def\Om{\Omega}

\def\limsup{\mathop{\overline{\rm lim}}}
\def\liminf{\mathop{\underline{\rm lim}}}

\def\esssup{\mathop{\rm esssup}}
\def\max{\mathop{\rm max}}

\def\exp{\mathop{\rm exp}}
\def\sup{\mathop{\rm sup}}

\def\pa{\partial}
\def\h{\widehat}
\def\wt{\widetilde}
\def\cd{\cdot}
\def\cds{\cdots}

\def\div{\hbox{\rm div$\,$}}
\def\inf{\hbox{\rm inf$\,$}}
\def\esssup{\hbox{\rm ess$\,$\rm sup$\,$}}

\def\ae{\hbox{\rm a.e.{ }}}
\def\as{\hbox{\rm a.s.{ }}}

\def\supp{\hbox{\rm supp$\,$}}

\def\cl{\overline}

\def\deq{\mathop{\buildrel\D\over=}}

\def\|{\Big |}
\def\({\Big (}
\def\){\Big )}
\def\[{\Big[}
\def\]{\Big]}

\def\be{\begin{equation}}
	\def\bel{\begin{equation}\label}
		\def\ee{\end{equation}}
	\def\bt{\begin{theorem}}
		\def\bcd{\begin{condition}}
			\def\ecd{\end{condition}}
		\def\et{\end{theorem}}
	\def\bc{\begin{corollary}}
		\def\ec{\end{corollary}}
	\def\bde{\begin{definition}}
		\def\ede{\end{definition}}
	\def\bl{\begin{lemma}}
		\def\el{\end{lemma}}
	\def\bp{\begin{proposition}}
		\def\ep{\end{proposition}}
	\def\br{\begin{remark}}
		\def\er{\end{remark}}
	\def\ba{\begin{array}}
		\def\ea{\end{array}}
	\def\ed{\end{document}}
\def\ns{\noalign{\ms}}
\def\ds{\displaystyle}

\newtheorem{lemma}{Lemma}[section]
\newtheorem{remark}{Remark}[section]
\newtheorem{example}{Example}[section]
\newtheorem{theorem}{Theorem}[section]
\newtheorem{corollary}{Corollary}[section]

\newtheorem{definition}{Definition}[section]
\newtheorem{proposition}{Proposition}[section]
\newtheorem{condition}{Condition}[section]

\makeatletter
   
   \@addtoreset{equation}{section}
\makeatother

\allowdisplaybreaks

\begin{document}
	
\title{\bf A Concise Introduction to Control Theory for Stochastic Partial
    Differential Equations\thanks{This is a lecture
    	notes of ``a concise introduction to control theory for
    	stochastic partial differential equations". It was
    	written for the EECI (European Embedded Control
    	Institute) Summer School ``International
    	Graduate School on Control" at Sichuan
    	University, Chengdu, China from July 8 to July
    	27, 2019. We refer to the monograph \cite{LZ3.1} for a more detailed presentation of the whole theory, while the survey paper \cite{LZ7} is mainly for the engineering-oriented readers.}}

\author{Qi L\"u\thanks{School of Mathematics,
Sichuan University, Chengdu 610064, Sichuan
Province, China. The author is partially
supported by NSF of China
under grants 12025105, 11971334 and 11931011, and
by the Chang Jiang Scholars Program
from the Chinese Education Ministry. {\small\it E-mail:} {\small\tt
lu@scu.edu.cn}.} \q and\q Xu Zhang\thanks{School of Mathematics,
Sichuan University, Chengdu 610064, Sichuan
Province, China.  The author is supported by the
NSF of China under grants 11931011 and
11821001. {\small\it E-mail:} {\small\tt zhang$\_$xu@scu.edu.cn}.}}

\date{}

\maketitle

\begin{abstract} 
The aim of this notes is to give a
concise introduction to control theory
for systems governed by stochastic
partial differential equations. We
shall mainly focus on  controllability
and optimal control problems for these systems. For the
first one, we present results for the
exact controllability of stochastic
transport equations, null and approximate
controllability of stochastic parabolic
equations  and lack of exact
controllability of stochastic
hyperbolic equations. For the second one, we first
introduce the stochastic linear
quadratic optimal control problems and
then the Pontryagin type maximum
principle for general optimal control
problems. It deserves mentioning that, in order to solve some difficult problems in this field, one
has to develop new tools, say, the stochastic transposition method introduced in our previous works.
\end{abstract}

\bs

\no{\bf 2010 Mathematics Subject
Classification}.  Primary 93E20, 60H15, 93B05, 93B07.

\bs

\no{\bf Key Words}. Stochastic partial differential equation,
controllability, observability,
optimal control, linear quadratic control
problem, Pontryagin type maximum principle.

\ms

%%%%%%%%%%%%%%%%%%%%%%%%%%%%%%%%%%%%%%%%%%%%%%%%%%%%%%%%

\tableofcontents

%%%%%%%%%%%%%%%%%%%%%%%%%%%%%%%%%%%%%%%%%%%%%%%%%%%%%%%%%

\section{Some preliminaries}\label{s1}

For the readers' convenience, we collect in this section
some basic knowledge of Functional
Analysis, Partial Differential Equations (PDEs for short),  Stochastic Analysis and Stochastic Partial Differential Equations (SPDEs for short), which will be used later. For more details, we refer to
\cite{Prato, Halmos, LZ3.1, Meyer1, Pazy, Yos95, Zeidler}. Unless  otherwise stated, all definitions, theorems and examples in this section are taken from these references.
Throughout this paper, $\dbN$ is the set of positive integers, while
$\dbR$ and $\dbC$ stand for respectively
the fields of real numbers and complex numbers. Also, we shall denote by
$\cC$ a generic positive constant which may change from line to line (unless otherwise stated).

%%%%%%%%%%%%%%%%%%%%%%%%%%%%%%%%%%%%%%%%%%%%%%%

\subsection{Banach spaces, Hilbert spaces and linear operators}
\label{sec-mpr-1.3}

%%%%%%%%%%%%%%%%%%%%%%%%%%%%%%%%%%%%%%%%%%%%%%%

In this notes, all linear spaces
are over the field $\dbR$ or over the field $\dbC$. Clearly, any linear space
over the field $\dbC$ is also a linear space
over the field $\dbR$. For any $c\in\dbC$, denote by $\bar c$ its complex conjugate.

%%%%%%%%%%%%%%%%%%%%%%%%%%%%%%%%%%%%%%%%%%%%%%%

\subsubsection{Banach spaces and Hilbert spaces}

%%%%%%%%%%%%%%%%%%%%%%%%%%%%%%%%%%%%%%%%%%%%%%%

\begin{definition}\label{2.1-def1} Let $\cX$ be a linear space.
    A map $|\cd|_\cX:\cX\to\dbR$ is called a {\it norm}
    on $\cX$ if it satisfies the following:
    \begin{equation}\label{2.1-eq1}
        \left\{
        \begin{array}{ll}\ds
            |x|_\cX \ge0,\q\forall\; x\in
            \cX;\hbox{ and }|x|_\cX=0~\iff~x=0;\\
            \ns\ds|\a x|_\cX=|\a||x|_\cX,\q\forall\;\a\in \dbC,\,x\in
            \cX;\\
            \ns\ds|x+y|_\cX\le|x|_\cX+|y|_\cX,\qq\forall\; x,y\in
            \cX.
        \end{array}
        \right.
    \end{equation}
    A linear space $\cX$ with the above norm $|\cd|_\cX$ is called a {\it normed linear
        space} and denoted by $(\cX, |\cd|_\cX)$ (or simply by $\cX$ if the norm $|\cd|_\cX$ is clear from the context).
\end{definition}
\begin{definition}
    Let  $\cX$ be a normed linear space.

    1) We call $\{x_k\}_{k=1}^\infty\subset \cX$  a {\it Cauchy
        sequence} (in $(\cX, |\cd|_\cX)$) if for any $\e>0$, there is $k_0\in \dbN$ such that
    $|x_k-x_j|_\cX<\e$ for all $k,j\ge k_0$.

    2)
    We call $\{x_k\}_{k=1}^\infty\subset \cX$ converges to some $x_0(\in \cX)$ in $\cX$ if $\ds\lim_{k\to\infty}x_k=x_0$ in $\cX$, i.e., $\ds\lim_{k\to\infty}|x_k-x_0|_\cX=0$.
\end{definition}
\begin{definition}
    Let  $\cX$ be a normed linear space.

    1) A subset $G\subset\cX$ is called {\it bounded}, if there is a constant $\cC$ such that
    $|x|_\cX\le \cC$ for any $x\in G$.

    2) A subset $G$ is said to be {\it dense}  in $\cX$ if for any $x\in \cX$, one can find a sequence $\{x_k\}_{k=1}^\infty\subset G$ such that $\ds\lim_{k\to\infty}x_k=x$ in $\cX$.
\end{definition}

\begin{definition}\label{2.1-def3}
    A normed linear space
    $(\cX,|\cd|_\cX)$ is called a {\it Banach space} if it
    is {\it complete}, i.e., for any Cauchy sequence
    $\{x_k\}_{k=1}^\infty\subset \cX$,
    there exists $x_0\in \cX$ so that $\ds\lim_{k\to\infty}x_k=x_0$ in $\cX$.
\end{definition}

\begin{definition}\label{2.1-def3.1}
    A Banach space $\cX$ is called {\it separable} if there exists a countable dense subset of $\cX$.
\end{definition}

\begin{definition}\label{2.1-def1.1} Let $\cX$ be a linear space.
    A map $\lan\cd,\cd\ran_\cX:\cX\times \cX\to \dbC$ is called an
    {\it inner product} on $\cX$ if it satisfies the
    following:
    \begin{equation}\label{2.1-eq2}
        \left\{
        \begin{array}{ll}\ds
            \lan x,x\ran_\cX\ge0,\q\forall\; x\in \cX;\hbox{ and }\lan x,x\ran_\cX=0~\iff~x=0;
            \\
            \ns\ds\lan x, y\ran_\cX=\overline{\lan y, x\ran_\cX},\qq\forall\; x,y\in
            \cX; \\
            \ns\ds \lan \a x+\b
            y,z\ran_\cX=\a\lan x,y\ran_\cX+\b\lan y,z\ran_\cX, \qq\forall\;\a,\b\in \dbC,x,y,z\in
            \cX.
        \end{array}
        \right.
    \end{equation}
    A linear space $\cX$ with the inner product $\lan\cd,\cd\ran_\cX$ is called an {\it inner product
        space} and denoted by $(\cX, \lan\cd,\cd\ran_\cX)$ (or simply by $\cX$ if the inner product $\lan\cd,\cd\ran_\cX$ is clear from the context).
\end{definition}

The following result gives a
relationship between norm and inner
products.

\begin{proposition}\label{2.1-prop4}
    Let $(\cX, \lan\cd,\cd\ran_\cX)$  be an
    inner product
    space. Then, the map
    $x\mapsto\sqrt{\lan x,x\ran_\cX}$, $\forall\; x\in \cX$, is a norm
    on $\cX$.
\end{proposition}

By Proposition \ref{2.1-prop4}, any inner product space $(\cX, \lan\cd,\cd\ran_\cX)$ can be regarded as a
normed linear space. We call
$|x|_\cX=\sqrt{\lan x,x\ran_\cX}$ the norm {\it induced}
by $\lan\cd,\cd\ran_\cX$.

\begin{definition}\label{2.1-def5}
    An inner product space $\cX$  is called
    a {\it Hilbert space} if it is complete under
    the norm induced by its inner product.
\end{definition}
%

%%%%%%%%%%%%%%%%%%%%%%%%%%%%%%%%%%%%%%%%%%%%%%%

\subsubsection{Bounded linear operator}

%%%%%%%%%%%%%%%%%%%%%%%%%%%%%%%%%%%%%%%%%%%%%%%

In the rest of this section, unless otherwise stated,  $\cX_1$ and
$\cY_1$ are normed linear spaces, and $\cX$ and
$\cY$ are Banach spaces.

\begin{definition}\label{1def6}
    A map $A:\cX_1\to
    \cY_1$ is called a {\it bounded linear operator} if it is linear, i.e.,
    \begin{equation}\label{2.1-eq8}
        A(\a x+\b y)=\a Ax+\b Ay,\qq\forall\; x,y\in \cX_1,\a,\b\in \dbC,
    \end{equation}
    and $A$ maps bounded subsets of $\cX_1$ into bounded
    subsets of $\cY_1$.
\end{definition}

Denote by $\cL(\cX_1;\cY_1)$  the set of all bounded linear
operators from $\cX_1$ to $\cY_1$. We simply write $\cL(\cX_1)$ instead
of $\cL(\cX_1;\cX_1)$. For any
$\a,\b\in \dbR$ and $A,B\in\cL(\cX_1;\cY_1)$, we define $\a
A+\b B$ as follows:
\begin{equation}\label{2.1-eq9}
    (\a A+\b B)(x)=\a Ax+\b Bx,\qq\forall\; x\in \cX_1.
\end{equation}
Then, $\cL(\cX_1;\cY_1)$ is also a linear space. Let
\begin{equation}\label{2.1-eq10}
    |A|_{\cL(\cX_1;\cY_1)}\deq\sup_{x\in\cX_1\setminus\{0\}}{\frac{|Ax|_{\cY_1}}{|x|_{\cX_1}}}.
\end{equation}
One can show that, $|\cd|_{\cL(\cX_1;\cY_1)}$ defined by
\eqref{2.1-eq10} is a norm on $\cL(\cX_1;\cY_1)$, and  $\cL(\cX;\cY)$ is a
Banach space under such a norm.

We shall use the following two theorems.
\begin{theorem}[Inverse Mapping]\label{2.1-th10}
    If
    $A\in \cL(\cX;\cY)$ is bijective, then $A^{-1}\in\cL(\cY;\cX)$.
\end{theorem}

\begin{theorem}[Uniform Boundedness]\label{2.1ubth10}
    If $\L$ is a given index set,
    $A_\l\in \cL(\cX;\cY)$ for each $\l\in\L$ and,
    $$
    \sup_{\l\in \L}|A_\l x|_\cY<\infty,\q\forall\;x\in \cX,
    $$
    then $\ds\sup_{\l\in \L}|A_\l|_{\cL(\cX;\cY)}<\infty$.
\end{theorem}

Let us consider the special case $\cY_1=\dbC$. Any $f\in \cL(\cX_1;\dbC)$ is
called a {\it bounded linear functional} on $\cX_1$.
Hereafter, we write $\cX'_1=\cL(\cX_1;\dbC)$ and call it
the {\it dual} (space) of $\cX_1$. We also
denote
\begin{equation}\label{2.1-eq15}
    f(x)=\langle
    f,x\rangle_{\cX'_1,\cX_1},\qq\forall\;
    x\in \cX_1.
\end{equation}
The symbol
$\langle\cd\,,\cd\rangle_{\cX'_1,\cX_1}$
is referred to as the {\it duality pairing}
between $\cX'_1$ and $\cX_1$. It follows from
\eqref{2.1-eq10} that
\begin{equation}\label{2.1-eq16}
    |f|_{\cX'_1}=\sup_{x\in \cX_1,|x|_{\cX_1}\le1}|f(x)|,\qq\forall\; f\in \cX'_1.
\end{equation}
Clearly, both $\cX'_1$  and $\cX'$ are  Banach spaces. Particularly , $\cX$ is called {\it reflexive} if $\cX''=\cX$.

The following result is quite useful.
\begin{theorem}[Hahn-Banach]\label{Hahn-Banach}
    Let $\cX_0$ be a linear subspace of $\cX_1$ and $f_0\in \cX'_0$. Then, there exists
    $f\in \cX'_1$ with $|f|_{\cX'_1}=|f_0|_{\cX'_0}$, such that
    $$
    f(x)=f_0(x),\quad\forall\;x\in \cX_0.
    $$
\end{theorem}
An immediate consequence of Theorem \ref{Hahn-Banach} is as follows:
\begin{proposition}\label{2.1-cor16}
    For any $x_0\in \cX$, there is
    $f\in \cX'$ with $|f|_{\cX'}=1$ such that
    $
    f(x_0)=|x_0|_{\cX}$.
\end{proposition}

We shall use the following theorem.
\begin{theorem}[Riesz Representation]\label{1t10s}
    If $(\cX, \lan\cd,\cd\ran_\cX)$  is a Hilbert space and $F\in \cX'$, then there is $y\in \cX$ such that
    $$
    F(x)={\lan x,y\ran}_{\cX},\q\forall\;x\in\cX.
    $$
\end{theorem}

For any
$A\in\cL(\cX;\cY)$, we define a map $A^*: \cY'\to
\cX'$ by the following:
\begin{equation}\label{2.1-eq26}
    \langle A^*y',x\rangle_{\cX',\cX}=\langle
    y',Ax\rangle_{\cY',\cY},\qq\forall\; y'\in
    \cY',\;x\in \cX.
\end{equation}
Clearly, $A^*$ is linear and bounded. We call
$A^*$ the {\it adjoint operator} of $A$. It is easy to check that for any
$A,B\in\cL(\cX;\cY)$ and $\a,\b\in \dbR$, $(\a A+\b
B)^*=\a A^*+\b B^*$.

Let $V_1$ and $V_2$ be separable Hilbert spaces,  and $\{e_j\}_{j=1}^\infty$ be an orthonormal basis  of $V_1$. $F\in \cL(V_1;V_2)$ is called a
{\it Hilbert-Schmidt
    operator} if
$\sum_{j=1}^\infty | Fe_j|_{V_2}^2 <\infty. $
Denote by $\cL_2(V_1;V_2)$
the space of all Hilbert-Schmidt operators from
$V_1$ into $V_2$.  One can show that, $\cL_2(V_1;V_2)$ equipped with the
inner product
$$
\langle F,G\rangle_{\cL_2(V_1;V_2)}
= \sum_{j=1}^\infty \langle
Fe_j,Ge_j\rangle_{V_2},\;\, \forall\;F,G\in
\cL_2(V_1;V_2),
$$
is a separable Hilbert space.
%%%%%%%%%%%%%%%%%%%%%%%%%%%%%%%%%%%%%%%%%%%%%%%

\subsubsection{Unbounded linear operator}

%%%%%%%%%%%%%%%%%%%%%%%%%%%%%%%%%%%%%%%%%%%%%%%

In this section, unless otherwise stated, $H$ is a Hilbert space.
\begin{definition}\label{1def8}
    For a linear map $A : D(A) \to
    H$, where $D(A)$ is a linear
    subspace in $ H$, we call $D(A)$ the {\it
        domain} of $A$. The {\it graph} of $A$ is the subset
    of $H \times H$ consisting of all elements of
    the form $(x,Ax)$ with $x \in D(A)$.
    The operator $A$
    is called {\it closed} (\resp {\it densely defined}) if its graph is a closed
    subspace of $H\times H$ (\resp $D(A)$ is dense in $H$).
\end{definition}

Unlike bounded linear operators, the linear operator $A$ in Definition \ref{1def8} might be unbounded, i.e., it may map some bounded sets in $H$
to unbounded sets in $H$.
A typical densely defined,
closed unbounded linear  operator (on $H=L^2(0,1)$) is the
differential operator $\frac{d}{dx}$, with the
domain
$$
\begin{array}{ll}\ds\Big\{y\in L^2(0,1) \;\Big|\;y \mbox{ is absolutely continuous on $[0,1]$, }
    \frac{dy}{dx} \in L^2 (0,1),y(0)=0\Big\}.
\end{array}
$$

In the sequel, we shall assume that the operator $A$
is densely defined.
The domain $D(A^*)$ of the {\it adjoint operator} $A^*$
of $A$ is
defined as the set of all $f \in H$
such that, for some $g_f \in
H$,
$$
\lan Ax,f\ran_H=\lan x,g_f\ran_H,\q\forall\;x
\in D(A).
$$
In this case, we define $A^*f \deq g_f$.

Denote by $I$ the identity operator on $H$. The {\it resolvent set} $\rho(A)$ of  $A$ is defined by
\begin{equation*}
    \begin{array}{ll}
        \rho(A) \deq \big\{\lambda \in \dbC\,\big|\,
        \lambda I- A:\;D(A)\to H\hbox{ is bijectve and } (\lambda I-
        A )^{-1}\in \cL(H)\big\}.
    \end{array}
\end{equation*}
The resolvent set $\rho(A)$ is open in $\dbC$.

%%%%%%%%%%%%%%%%%%%%%%%%%%%%%%%%%%%%%%%%%

\subsection{Measure and integration}\label{sec-mea-3}

%%%%%%%%%%%%%%%%%%%%%%%%%%%%%%%%%%%%%%%%%
In this subsection, we recall some basic definitions and
results on measure and integration.

Let $\Omega$ be a nonempty set. For any subset $E\subset \Omega$,
denote by $\chi_E(\cdot)$ the {\it characteristic
    function} of $E$, defined on $\Omega$, that is,
$$
\chi_E(\om)=\begin{cases}\ds
    1, & \mbox{ if }\om\in E,\\
    \ns\ds 0, & \mbox{ if }\om\in \Om\setminus E.
\end{cases}
$$
\begin{definition}
    Let $\cF$ be a family
    of nonempty subsets of $\Omega$. $\cF$ is
    called a {\it $\si$-field on $\Omega$} if 1)
    $\Omega\in\cF$; 2) $\Omega\setminus E\in \cF$
    for any $E\in\cF$; and 3)
    $\ds\bigcup_{i=1}^\infty E_i\in\cF$ whenever each
    $E_i\in \cF$.
\end{definition}

If $\cF$ is a $\si$-field on $\Omega$,
then $(\Omega,\cF)$ is called a
{\it measurable space}.  Any element $E\in\cF$
is called a
{\it measurable set} on $(\Omega,\cF)$, or
simply a measurable set.

In the sequel, we shall fix a measurable space $(\Omega,\cF)$ .

\begin{definition}
    A set function
    $\mu:\ \cF\to[0,+\infty]$ is called a
    {\it measure} on
    $(\Omega,\cF)$ if $\mu(\emptyset)=0$
    and $\mu$ is
    countably additive, i.e., $\ds\mu(
    \bigcup_{i=1}^\infty
    E_i)=\sum_{i=1}^\infty\mu(E_i)$
    whenever $\{E_i\}_{i=1}^\infty\subset\cF$ are mutually disjoint, i.e., $E_i\cap E_j=\emptyset$ for all $i,j\in\dbN$ with $i\not= j$. The triple $(\Omega,\cF,\mu)$ is
    called a {\it measure
        space}.
\end{definition}

We shall fix below a measure space $(\Omega,\cF,\mu)$.

\begin{definition}
    The measure $\mu$ is called
    finite (\resp $\si$-finite) if $\mu
    (\Omega)<\infty$ (\resp there exists a
    sequence
    $\{E_i\}_{i=1}^\infty\subset\cF$ so
    that $\Omega=\bigcup _{i=1}^\infty E_i$
    and $\mu(E_i)<\infty$ for each
    $i\in\dbN$).
\end{definition}

We call any $E\subset\cF$ a {\it $\mu$-null (measurable)
    set} if $\mu(E)=0$.

\begin{definition}\label{1del17}
    The measure space $(\Omega,\cF,\mu)$ is said to be {\it complete} (and $\mu$ is said to be complete on $\cF$) if
    $$
    \cN=\{\wt E\subset\Omega\;|\;\wt E\subset E\hb{ for
        some} \hb{ $\mu$-null set }E\}\subset\cF.
    $$
\end{definition}
If the measure space $(\Omega,\cF,\mu)$ is not
complete,  then the class $\overline\cF$ of
all sets of the form $(E\setminus
N)\bigcup(N\setminus E)$, with $E\in \cF$ and
$N\in\cN$ (See Definition \ref{1del17} for $\cN$), is a $\si$-field which contains $\cF$
as a proper sub-class, and the set function
$\bar \mu$ defined by $\bar \mu((E\setminus
N)\cup(N\setminus E))=\mu(E)$ is a complete
measure on $\overline\cF$ . The measure
$\bar\mu$ is called the  {\it completion of $\mu$}.

If a certain proposition concerning the points of $(\Omega,\cF,\mu)$ is true for every point $\omega\in\Omega$, with the exception at most of a
set of points which form a $\mu$-null set, it is
customary to say that the proposition is true $\mu$-a.e. (or simply a.e., if the measure $\mu$ is clear from the context). A
function $g: (\Omega,\cF,\mu)\to\dbR$ is called {\it essentially bounded} if it is bounded $\mu$-a.e., i.e.
if there exists a positive, finite constant $c$ such that $\big\{\omega\in\Omega\;\big|\; |g(\omega) | > c\big\}$
is a $\mu$-null set. The infimum of the values of $c$ for which
this statement is true is called the {\it essential supremum} of $|g|$,
abbreviated to $\esssup_{\omega\in\Omega}|g(\omega)|$.

Let $(\Omega_1,\cF_1),\cdots,(\Omega_n,\cF_n)$ ($n\in\dbN$)
be measurable spaces. Denote by
$\cF_1\times\cdots\times \cF_n$ the
$\sigma$-field\footnote{Note that here and
    henceforth the $\sigma$-field
    $\cF_1\times\cdots\times \cF_n$ does not stand
    for the Cartesian product of $\cF_1$, $\cdots$,
    $\cF_n$.} (on the Cartesian product
$\Omega_1\times\cdots\times \Omega_n$) generated
by the subsets of the form
$E_1\times\cdots\times E_n$, where
$E_i\in\cF_i$, $1\leq i\leq n$.
Let $\mu_i$ be a measure  on $(\Omega_i,\cF_i)$.
We call $\mu$ a  {\it product measure} on $(\Omega_1\times\cdots\times
\Omega_n,\cF_1\times\cdots\times \cF_n)$ induced
by $\mu_1,\cdots,\mu_n$ if
$$
\mu(E_1\times\cdots\times E_n)=\prod_{i=1}^n
\mu_i(E_i),\quad \forall\; E_i\in\cF_i.
$$
\begin{theorem}\label{5.14-th2}
    Let $(\Omega_1,\cF_1,\mu_1),\cdots,(\Omega_n,\cF_n,\mu_n)$ be  $\si$-finite
    measure spaces. Then there is
    a unique product measure $\mu$ (denoted
    by
    $\mu_1\times\cdots\times \mu_n$) on
    $(\Omega_1\times\cdots\times
    \Omega_n,\cF_1\times\cdots\times \cF_n)$ induced
    by $\{\mu_i\}_{i=1}^n$.
\end{theorem}

In the rest of this subsection, we fix a Banach space $\cX$.
\begin{definition}
    The smallest $\si$-field
    containing all open sets of $\cX$ is called the
    {\it Borel $\si$-field} of $\cX$ and denoted by
    $\cB(\cX )$\footnote{More generally, when $\cX$ is a topology space, one can define its Borel $\si$-field $\cB(\cX )$ in the same way.}. Any set $E\in \cB(\cX)$
    is called a {\it Borel set} (in $\cX$).
\end{definition}
Let $(\Omega',\cF')$ be another
measurable space, and $f:\ \Omega\to\Omega'$ be
a map\footnote{When $\Omega'=\cX$, we call $f$ an $\cX$-valued function.}.
\begin{definition}\label{12.25-def1}
    The above map $f$ is said to be
    {\it $\cF/\cF'$-measurable} or
    simply {\it $\cF$-measurable} or even {\it measurable} (in
    the case that no confusion would occur) if
    $f^{-1}(\cF')\subset \cF$. Particularly if
    $(\Omega',\cF')=(\cX,\cB(\cX))$, then $f$ is said to
    be an ($\cX$-valued)
    $\cF$-measurable (or simply measurable)
    function.
\end{definition}

Let $\mathfrak{P}$ be a property concerning the above map
$f$ at some elements in $\Omega$. We shall simply denote by
$\{\mathfrak{P}\}$
the subset $ \{\omega\in \Omega\;|\;\mathfrak{P}
\mbox{ holds for } f(\omega)\}$. For example,
when $\Omega'=\dbR$, to simplify the notation we
denote $\{\omega\in \Omega\;|\;f(\omega)\geq
0\}$ by $\{f\geq 0\}$.

For a measurable map $f:\
(\Omega,\cF)\to(\Omega',\cF')$, it is obvious that $f^{-1}(\cF')$ is a sub-$\si$-field of
$\cF$. We call it the {\it $\si$-field generated by
    $f$}, and denote it by $\si(f)$. Further,  for a
given index set $\Lambda$ and a family of
measurable maps $\{f_\l\}_{\l\in \Lambda}$
(defined on $(\Omega,\cF)$,  with possibly
different ranges), we denote by $\si(f_\l;\l\in
\Lambda)$ the $\si$-field generated by
$\bigcup_{\l\in \Lambda}\si(f_\l)$.

\begin{definition}
    Let $\{f_i\}_{i=0}^\infty$ be a sequence of
    $\cX$-valued functions defined on $\Omega$.
    The sequence
    $\{f_k\}_{k=1}^\infty$ is said to
    converge to $f_0$ (denoted by $\ds\lim_{k\to\infty}f_k=f_0$) in $\cX$,
    $\mu$-a.e., if $\{ \ds\lim_{k\to\infty}
    f_k \not=f_0\}\in\cF$, and
    $
    \ds\mu(\{ \lim_{k\to\infty} f_k \not=f_0\})=0.
    $

\end{definition}

In the setting of infinite dimensions, the notion of $\cF$-measurability in Definition \ref{12.25-def1} may not
provide any means for
approximation arguments. Thus, we need
another notion of measurability. We restrict ourselves to Banach space valued functions, although some of the results below can be generalized
to functions with values in metric spaces.

\begin{definition}\label{Ap-def1}
    Let $f:\ \Omega\to\cX$ be an $\cX$-valued function.

    {\rm 1)}  We call
    $f(\cd)$ an {\it $\cF$-simple function} (or simply a {\it simple function} when $\cF$ is clear from the context) if
    \begin{equation}\label{Axzp-def5-eq1}
        f(\cdot) =\sum_{i=1}^k\chi_{E_i}(\cdot)h_i,
    \end{equation}
    for some $k\in\dbN$, $h_i\in \cX$, and mutually
    disjoint sets $E_1,\cdots,E_k$ $\in \cF$ satisfying
    $\ds\bigcup_{i=1}^kE_i=\Omega$;

    {\rm 2)}   The function $f(\cd)$ is
    said to be {\it strongly $\cF$-measurable with respect to (w.r.t. for short) $\mu$} (or
    simply strongly measurable) if there exists a
    sequence of $\cF$-simple functions
    $\{f_k\}_{k=1}^\infty$ converging to
    $f$ in $\cX$,
    $\mu$-a.e. In this case, sometimes we also say that $f: (\Omega,\cF,\mu)\to \cX$ is strongly measurable.
\end{definition}

In order to characterize the strong  measurability of $\cX$-valued
functions, we introduce some terminology.
\begin{definition}
    A function $f :\Omega\to \cX$ is called  {\it $\mu$-separably valued} (or simply  {\it separably valued}) if there exists a separable closed subspace $\cX_0\subset\cX$ such
    that $f(\om)\in \cX_0$ for $\mu$-a.e. $\om\in\Om$.
\end{definition}
\begin{theorem}\label{12.25-th1}
    Let
    $(\Omega,\cF,\mu)$ be a $\si$-finite measure space.
    For a function $f :\Om\to\cX$, the following assertions are equivalent:

    {\rm 1)} $f$ is strongly measurable;

    {\rm 2)} $f$ is separably valued and measurable.
\end{theorem}

By Theorem \ref{12.25-th1},  if the Banach space $\cX$ is separable and the measure space $(\Omega,\cF,\mu)$ is
$\si$-finite, then, $f$ is strongly measurable $\Leftrightarrow$ $f$ is  measurable.

Let us fix below a
$\si$-finite measure space $(\Omega,\cF,\mu)$.

\begin{definition}
    Let $f(\cd)$ be an ($\cX$-valued) simple
    function in the form \eqref{Axzp-def5-eq1}. We
    call $f(\cd)$ {\it Bochner integrable} if
    $\mu(E_i)<\infty$ for each $i=1,\cdots,k$.
    In
    this case, for any $E\in\cF$, the {\it Bochner
        integral} of $f(\cd)$ over $E$ is defined by
    $$\int_E
    f(s)d\mu=\sum_{i=1}^k \mu(E\cap E_i)h_i.
    $$
\end{definition}

In general, we have the following notion.
\begin{definition}\label{Ap-def5}
    A strongly  measurable function $f(\cd): $
    $\Omega\to \cX$ is said to be  {\it Bochner integrable}
    (w.r.t. $\mu$) if there exists a sequence of
    Bochner integrable  simple functions
    $\{f_i(\cd)\}_{i=1}^\infty$ converging
    to $f(\cd)$ in $\cX$, $\mu$-a.e., so that
    $$
    \lim_{i,j\to\infty}\int_\Omega|f_i(s)-
    f_j(s)|_{\cX}d\mu =0.
    $$
    In this case, we also say that $f(\cd):\; (\Omega,\cF,\mu)\to \cX $
    is Bochner integrable. For any $E\in\cF$, the
    {\it Bochner integral} of $f(\cd)$ over $E$ is
    defined by
    \begin{equation}\label{Ap-def5-eq1}
        \int_Ef(s)d\mu =
        \lim_{i\to\infty}\int_\Omega\chi_E(s)f_i(s)d\mu(s)\q
        \mbox{ in } \cX.
    \end{equation}
\end{definition}

It is easy to verify that the limit in the right
hand side of \eqref{Ap-def5-eq1} exists and its
value is independent of the choice of the
sequence $\{f_i(\cd)\}_{i=1}^\infty$. Clearly,
when $\cX=\dbR^n$ (for some $n\in\dbN$), the above
Bochner integral coincides the usual Lebesgue  integral
for $\dbR^n$-valued functions.

The following result reveals the relationship
between the Bochner integral (for vector-valued
functions) and the usual Lebesgue integral (for scalar
functions).

\begin{theorem}\label{Ap-th2}
    Let $f(\cd):\; \Omega\to \cX$ be strongly
    measurable. Then, $f(\cd)$ is
    Bochner integrable (w.r.t. $\mu$) if and only if
    the scalar function $|f(\cd)|_{\cX}:\; \Omega\to
    \dbR$ is integrable (w.r.t. $\mu$).
\end{theorem}

Further properties for Bochner integral are
collected as follows.

\begin{theorem}\label{Ap-tha2}
    Let $f(\cd), \ g(\cd):\; (\Omega,\cF,\mu)\to \cX $
    be Bochner integrable. Then,

    {\rm 1)} For any $a,b\in\dbR$, the function
    $af(\cd)+bg(\cd)$ is Bochner integrable, and  for any $E\in\cF$,
    $$
    \int_E\big(af(s)+bg(s)\big)d\mu=a\int_Ef(s)d\mu+b\int_Eg(s)d\mu.
    $$

    {\rm 2)} For any $E\in\cF$,
    $$
    \Big|\int_Ef(s)d\mu\Big|_{\cX }\leq
    \int_E|f(s)|_{\cX }d\mu.
    $$

    {\rm 3)}  The Bochner integral is
    $\mu$-absolutely continuous, that is,
    $$
    \lim_{E\in\cF,\;\mu(E)\to 0}\int_E f(s)d\mu=0\q
    \mbox{ in }\; \cX .
    $$

    {\rm 4)} If $F\in \cL(\cX;\cY)$,  then $Ff(\cd)$
    is a $\cY$-valued Bochner integrable
    function, and for any $E\in\cF$,
    $$
    \int_E Ff(s)d\mu = F\int_E f(s)d\mu.
    $$
\end{theorem}

The following result, known as {\it Dominated Convergence
    Theorem}, is very useful.
\begin{theorem}\label{Ap-th4}
    Let $f: (\Omega,\cF,\mu)\to \cX $ be strongly
    measurable, and let $g: (\Omega,\cF,\mu)$ $
    \to\dbR$ be a real valued nonnegative integrable
    function. Assume that, $f_i:\;
    (\Omega,\cF,\mu)\to \cX $ is Bochner integrable such
    that $|f_i|_{\cX }\leq g$, $\mu$-a.e. for each
    $i\in\dbN$ and $\ds\lim_{i\to\infty}f_i=f$ in $\cX$,
    $\mu$-a.e. Then, $f$ is Bochner integrable (w.r.t. $\mu$), and
    $$
    \lim_{i\to\infty}\int_E f_i(s)d\mu=\int_E
    f(s)d\mu\;\;\hbox{ in }\cX,\quad\forall\;E\in\cF.
    $$
\end{theorem}

Also, one has the following {\it Fubini Theorem} (on
Bochner integrals).
\begin{theorem}\label{Ap-tash4}
    Let $(\Omega_1,\cF_1,\mu_1)$ and
    $(\Omega_2,\cF_2,\mu_2)$ be $\si$-finite measure
    spaces. If $f(\cdot,\cdot):
    \Omega_1\times\Omega_2\to \cX $ is Bochner
    integrable (w.r.t. $\mu_1\times\mu_2$), then the functions
    $
    y(\cdot)\equiv\int_{\Omega_1}f(t,\cdot)d\mu_1(t)$
    and
    $z(\cdot)\equiv\int_{\Omega_2}f(\cdot,s)d\mu_2(s)
    $ are a.e. defined on
    $\Omega_2$ and $\Omega_1$ and Bochner integrable w.r.t. $\mu_2$ and $\mu_1$, respectively.
    Moreover,
    $$
    \int_{\Omega_1\times\Omega_2}f(t,s)d\mu_1\times\mu_2(t,s)=\int_{\Omega_1}z(t)d\mu_1(t)=\int_{\Omega_2}y(s)d\mu_2(s).
    $$
\end{theorem}

For any $p\in [1,\infty)$, denote by
\index{$L_\cF^p(\Omega;\cX )$}
\index{$L^p(\Omega,\cF,\mu;\cX )$}
$L_\cF^p(\Omega;\cX )\deq L^p(\Omega,\cF,\mu;\cX )$
the set of all (equivalence classes of) strongly
measurable functions $f:\ (\Omega,\cF,\mu)\to \cX $
such that $\int_\Omega|f|^p_{\cX }d\mu<\infty$. It is a Banach space with
the norm
\begin{equation*}\label{1s1e2}
    |f|_{L_\cF^p(\Omega;\cX )}=\left(\int_\Omega
    |f|^p_{\cX } d\mu\right)^{1/p}.
\end{equation*}
When $\cX $ is a Hilbert space, so is
$L_\cF^2(\Omega;\cX )$. Denote by
$L_\cF^\infty(\Omega;\cX )\deq
L^\infty(\Omega,\cF,\mu;\cX )$ the set of all
(equivalence classes of) strongly measurable
($X$-valued) functions $f$ such that
$\ds\esssup_{\omega\in\Omega}|f(\omega)|_{\cX }<\infty$.
It is also a Banach space with the norm
\begin{equation*}\label{1s1e2.1}
    |f|_{L_\cF^\infty(\Omega;\cX )}=\mathop{\hspace*{2pt}\esssup}_{\omega\in\Omega}|f(\omega)|_{\cX }.
\end{equation*}
For $1\leq p\leq\infty$ and  any
non-empty open subset $ G$ of $ \dbR^n$, we shall
simply  denote $L^p(G,$ $\cL,\dbm;\cX )$ by
$L^p(G;\cX )$, where $\cL$
is the family of Lebesgue measurable sets in
$G$, and $\dbm$ is the Lebesgue measure on
$G$. Also, we simply denote
$L_\cF^p(\Omega;\dbR)$  and $L^p(G;\dbR)$ by
$L_\cF^p(\Omega)$  and $L^p(G)$, respectively. Particularly, if $G=(0,T)\subset  \dbR$ for some $T>0$, we simply write $L^p((0,T);\cX )$ and $L^p((0,T))$ respectively by $L^p(0,T;\cX )$ and $L^p(0,T)$.

For any $p\in[1, \infty)$, denote by $p'$ the H\"older conjugate of $p$, i.e., $$p'=\left\{\ba{ll}\frac{p}{p-1},&\hbox{if }p>1,\\
\infty,&\hbox{if }p=1.\ea\right.
$$
The following result characterizes the dual space of
$L^p(\Omega;X)$.

\begin{proposition}\label{2.3-prop5}
    Let $\cX $ be a reflexive
    Banach space and $p\in[1, \infty)$. Then,
    \begin{equation}\label{2.3-eq8}
        L^p_\cF(\Om;\cX )'=L^{p'}_\cF(\Om;\cX').
    \end{equation}
\end{proposition}

Let  $\Phi: (\Omega,\cF)\to
(\Omega',\cF')$ be a measurable map. Then, for the measure $\mu$ on $(\Omega,\cF)$, the map $\Phi$
induces a measure $\mu'$ on $(\Omega',\cF')$ via
\begin{equation}\label{4.21-eq1}
    \mu'(E')\deq \mu(\Phi^{-1}(E')),\qq \forall \;
    E'\in\cF'.
\end{equation}
The following is a change-of-variable formula:
\begin{theorem}\label{4.21-th1}
    A function $f(\cd):\ \Omega'\to \cX $ is Bochner
    integrable w.r.t.  $\mu'$ if and only if
    $f(\Phi(\cd))$ (defined on $(\Omega,\cF)$) is
    Bochner integrable w.r.t.  $\mu$. Furthermore,
    \bel{chang1} \int_{\Omega'}
    f(\omega')d\mu'(\omega')=\int_{\Omega}
    f(\Phi(\omega)) d\mu(\omega). \ee
\end{theorem}

\subsection{Continuity and differentiability of vector-valued functions}\label{subsec0120}

In this subsection, we recall the notions of continuity and differentiability for vector-valued functions. Let $\cX$ and $\cY$ be Banach spaces,  $\cX_0\subset \cX$ and let $F:\,\cX_0\to \cY$ be a
function (not necessarily linear).

\begin{definition}\label{def-con}
    {\rm 1)} We say
    that $F$ is {\it continuous at $x_0\in
        \cX_0$} if
    $$\lim_{x\in \cX_0,x\to
        x_0}|F(x)-F(x_0)|_\cY=0.$$  If $F$ is
    continuous at each point of $\cX_0$, we
    say that $F$ is {\it continuous on $\cX_0$}.

    \ss

    {\rm 2)}  We say that $F$ is  {\it Fr\'echet
        differentiable} at $x_0\in \cX_0$ if there
    exists $F_1\in\cL(\cX;\cY)$ such that
    \begin{equation}\label{2.2-eq37}
        \lim_{x\in \cX_0,x\to
            x_0}\frac{|F(x)-F(x_0)-F_1(x-x_0)|_{\cY}}{|x-x_0|_{\cX}}=0.
    \end{equation}
    We call $F_1$ the {\it Fr\'{e}chet derivative} of $F$ at
    $x_0$, and write $F_x(x_0)=F_1$. If
    $F$ is Fr\'echet
    differentiable at each point of $
    \cX_0$, we say that $F$ is {\it Fr\'echet differentiable on
        $\cX_0$}.  Moreover, when the map $F_x:\cX_0\to \cL(\cX;\cY)$ is continuous, we say that $F$ is {\it continuous Fr\'echet differentiable on
        $\cX_0$}.

    \ss

    {\rm 3)}  We say that $F$ is {\it  second order Fr\'echet
        differentiable} at $x_0\in \cX_0$ if $F:\,\cX_0\to \cY$ is continuous Fr\'echet differentiable and there
    exists $F_2\in\cL(\cX;\cL(\cX;\cY))$ such that
    \begin{equation}\label{2.2-eq37-1}
        \lim_{x\in \cX_0,x\to
            x_0}\frac{|F_x(x)-F_x(x_0)-F_2(x-x_0)|_{\cL(\cX;\cY)}}{|x-x_0|_{\cX}}=0.
    \end{equation}
    We call $F_2$ the {\it second order Fr\'{e}chet derivative} of $F$ at
    $x_0$, and write $F_{xx}(x_0)=F_2$. If
    $F$ is  second order Fr\'echet
    differentiable at each point of $
    \cX_0$, we say that $F$ is {\it  second order Fr\'echet differentiable on
        $\cX_0$}. Moreover, when the map $F_{xx}:\cX_0\to \cL(\cX;\cL(\cX;\cY))$ is continuous, we say that $F$ is {\it second order continuous  Fr\'echet differentiable on
        $\cX_0$}.
\end{definition}

The
set of all continuous (\resp continuous Fr\'echet
differentiable,  second order continuous Fr\'echet differentiable) functions from $\cX_0$ to
$\cY$ is denoted by $C(\cX_0;\cY)$ (\resp $C^1(\cX_0; \cY)$, $C^2(\cX_0; \cY)$).
When
$\cY=\dbR$, we simply denote it by
$C(\cX_0)$ (\resp $C^1(\cX_0)$, $C^2(\cX_0)$).

Next,  we recall the definition of bounded  bilinear operator. Let $\cZ$ be another Banach space.

\begin{definition}\label{def-con-2}
    A mapping $M:\cX\times \cZ\to \cY$ is called a
    {\it bounded bilinear operator} if  $M$ is
    linear in each argument and there is a constant
    $\cC>0$ such that
    $$
    |M(x,z)|_\cY\leq \cC|x|_\cX|z|_\cZ,\qq \forall\; (x,z)\in\cX\times\cZ.
    $$

\end{definition}
Denote by $\cL(\cX,\cZ;\cY)$ the set of all
bounded bilinear operators from $\cX\times \cZ$
to $\cY$. The norm of $M\in \cL(\cX,\cZ;\cY)$ is
defined by
$$
|M|_{\cL(\cX,\cZ;\cY)}=\sup_{x\in\cX\setminus\{0\},z\in\cZ\setminus\{0\}}\frac{|M(x,z)|_{\cY}}{|x|_\cX|z|_\cZ}.
$$
$\cL(\cX,\cZ;\cY)$ is a Banach space w.r.t.  this norm.

Any $L\in \cL(\cX;\cL(\cZ;\cY))$ defines a
bounded bilinear operator $\wt L\in
\cL(\cX,\cZ;\cY)$ as follows:
$$
\wt L(x,z)=\big(L(x)\big)(z),\q \forall\; (x,z)\in\cX\times\cZ.
$$
Conversely, any  $\wt L\in \cL(\cX,\cZ;\cY)$ defines an $L\in \cL(\cX;\cL(\cZ;\cY))$ in the following way:
$$
\big(L(x)\big)(z)=\wt L(x,z),\q \forall\; (x,z)\in\cX\times\cZ.
$$
Hence, any $F\in C(\cX_0;\cL(\cX;\cL(\cZ;\cY)))$ can be regarded as an element $\wt F\in C(\cX_0;$ $\cL(\cX, \cZ;\cY)$. In the rest of this paper,  if there is no confusion,   we identify $F\in C(\cX_0;\cL(\cX;\cL(\cZ;\cY)))$ with $\wt F\in C(\cX_0;\cL(\cX,$ $\cZ;\cY)$.

%

%%%%%%%%%%%%%%%%%%%%%%%%%%%%%%%%%%%%%%%%%%%%%%%%

\subsection{Generalized functions and Sobolev spaces}\label{sec-mpr-6.1}

%%%%%%%%%%%%%%%%%%%%%%%%%%%%%%%%%%%%%%%%%%%%%%%%

We begin with introducing some notations.
For any
$\a=(\a_1,\a_2,\cds,\a_n)$ (with $\a_j\in\dbN\cup\{0\}$, $j=1,\cds,n$), put $\ds |\a|
=\sum_{j=1}^n|\a_j|$ and
$
\pa^\a\equiv{\frac{\pa^{\a_1}}{\pa
        x_1^{\a_1}}}{\frac{\pa^{\a_2}}{\pa
        x_2^{\a_2}}}\cds{\frac{\pa^{\a_n}}{\pa
        x_n^{\a_n}}}.
$
Let $ G\subset \dbR^n$ be a bounded
domain with the boundary $\G $, and
write $\cl G$ for its closure. For any $x^0\in\dbR^n$ and $r>0$, denote by $B\left(x^{0}, r\right)$ the open ball with radius $r$, centered at $x^0$.

\begin{definition}
    We say the boundary $\G $ (of $G$) is $C^{k}$ (for some $k\in\dbN$) if
    for each point $x^{0} \in \G $ there
    exist $r>0$ and  a $C^{k}$ function
    $\gamma: \mathbb{R}^{n-1} \rightarrow
    \mathbb{R}$ such that, upon relabeling
    and reorienting the coordinates axes if
    necessary,
    $$
    G \cap B\left(x^{0}, r\right)=\left\{x \in B\left(x^{0}, r\right) \mid x_{n}>\gamma\left(x_{1}, \ldots, x_{n-1}\right)\right\}.
    $$
    We say that  $\G $ is $C^{\infty}$ if $\G
    $ is $C^{k}$ for any $k=1,2, \cdots$.
\end{definition}

If $\G $ is $C^{1}$, then along $\G $ one can define the
unit outward normal vector field (of $G$):
$$
\nu=\left(\nu^{1}, \ldots, \nu^{n}\right).
$$
The unit outward normal vector at each point $x \in
\G $ is hence
$\nu\left(x\right)=(\nu^{1}(x),
\cdots, $ $\nu^{n}(x))$. For any $u \in
C^{1}(\overline{G}) $, we call
$\frac{\pa u}{\pa\nu}\deq \nabla
u\cd\nu$ the outward normal derivative
of $u$.

Although most of the results in this
notes hold when $\G $ is $C^k$ for some
suitable $k\in\dbN$ (usually, $k=1$ or
$2$), for simplicity, in the sequel we shall always  assume
that $\G $ is $C^\infty$ (unless other stated).

For any $m\in \dbN\cup\{0\}$, denote by $C^m( G)$ and
$C^m(\cl G)$ the sets of all $m$
times continuously differentiable
functions on $ G$ and $\cl G$,
respectively, and by $C^m_0( G)$
the set of all functions $f\in C^m(
G)$, such that $\supp f\,\deq\,\{x\in
G \;|\;f(x)\ne0\}$ is compact in $ G$.
$C^m(\cl G)$
is a Banach space with the norm
\begin{equation}\label{2.6-eq1}
    |f|_{C^m(\cl G)}=\max_{x\in\cl
        G}\sum_{|\a|\le m}|\pa^\a f(x)|,\q  \forall\; f\in C^m(\cl G).
\end{equation}
In the sequel, we shall  write $C(\cl G)$ for $C^0(\cl G)$.

Denote by $C^\infty( G)$ the set of infinitely differentiable
functions on $ G$. Put
$$
C_0^\infty(G)=\big\{f\in C^\infty(G)\;\big|\;\supp f \hb{ is compact in
}G\big\}.
$$

\begin{definition}\label{1d23}
    Let $\{f_k\}_{k=1}^\infty$ be a sequence in $C_0^\infty(G)$. We say that
    $f_k\to 0$ in $C_0^\infty(G)$ as $k\to\infty$, if

    1) For some compact subset $K$ in $G$,
    $$
    \supp f_k\subset K,\q \forall\; k\in\dbN;
    $$

    2) For all multi-index $\a=(\a_1,\cdots,\a_n)$,
    $$
    \sup_{x\in K}|\pa^\a f_k(x)|\to 0,\q\hb{ as }k\to\infty.
    $$
    Generally, for $f\in C_0^\infty(G)$, we say that $f_k\to f$  in $C_0^\infty(G)$ as $k\to\infty$, if $f_k-f\to 0$  in
    $C_0^\infty(G)$ as $k\to\infty$.
\end{definition}

Denote by $\cD(G)$ the linear space $C_0^\infty(G)$ equipped with the sequential convergence given in Definition \ref{1d23}. A linear functional $F: \cD(G)\to\dbC$ is said to be a $\cD'(G)$ generalized function, in symbol $F\in\cD'(G)$, if $F$ is also
continuous in $\cD(G)$, i.e., for any  sequence $\{f_k\}_{k=1}^\infty\subset C_0^\infty(G)$ with
$f_k\to 0$ in $C_0^\infty(G)$ as $k\to\infty$, one has $\ds\lim_{k\to\infty}F(f_k)=0$.

\begin{example}
    Any function $f\in L^1(G)$ can be identified as a $\cD'(G)$
    generalized function in the following way:
    $$
    F_f(\f)=\int_G f(x)\f(x) dx,\qq \forall\;\f\in\cD(G).
    $$
\end{example}

\begin{example}
    Let $a\in G$. The Dirac $\d_a$-function (in $G$) is defined as
    $$
    \d_a(\f)=\f(a),\qq \forall\;\f\in\cD(G).
    $$
    It is easy to show that $\d_a\in\cD'(G)$. Also, one can show that $\d_a$ cannot be identified as a ``usual" function in $G$.
\end{example}

Let $F\in\cD'(G)$. The {\it generalized
    derivative} of $F$ w.r.t.  $x_k$ (for
some $k\in\{1,2,\cdots,n\}$),
$\frac{\pa F}{\pa x_k}$ is defined by
$$
\(\frac{\pa F}{\pa
    x_k}\)(\f)=-F\(\frac{\pa\f}{\pa
    x_k}\),\qq\forall\;\f\in \cD'(G).
$$
One can easily show that $\frac{\pa
    F}{\pa x_k}\in \cD'(G)$.

Let $p\in [1,\infty)$. For $f\in C^m(\cl G)$, define
\begin{equation}\label{2.6-eq5}
    |f|_{W^{m,p}( G)}\deq\Big(\int_
    G\sum_{|\a|\le m}|\pa^\a
    f|^pdx\Big)^{1/p},
\end{equation}
which is a norm on $C^m(\cl G)$. Write $W^{m,p}(
G)$ for the
completion of $C^m(\cl G)$ w.r.t.  the norm
\eqref{2.6-eq5}. Similarly, the completion of
$C^m_0( G)$ under the norm \eqref{2.6-eq5} is
denoted by $W^{m,p}_0 ( G)$.  For $p=2$, we also
write $H^m( G)$ for $W^{m,2}( G)$ and
$H^m_0( G)$ for $W^{m,p}_0(G)$.
Both $H^m(G)$ and $H^m_0(G)$ are
Hilbert spaces with the inner product
\begin{equation}\label{2.6-eq6}
    \lan f,g\ran_{H^m( G)}=\int_
    G\sum_{|\a|\le m}\pa^\a f \pa^\a \bar g dx,\qq\forall\;f,g\in H^m( G).
\end{equation}
It can be proved that a function $y\in
W^{m,p}( G)$, if and only if there
exist functions $f_\a\in L^p( G)$,
$|\a|\le m$, such that
\begin{equation}\label{2.6-eq7}
    \begin{array}{ll}\ds
        \int_ G y\pa^\a\f dx=(-1)^{|\a|}\int_ G
        f_\a\f dx,\ \q\forall\;\f\in C^\infty_0 (
        G).
    \end{array}
\end{equation}
The above function $f_\a$ is referred to as the
$\a$-th {\it generalized derivative}
of $y$.

Noting that $C_0^\infty(G)\subset W_0^{m,p}(G)$, people introduce the following space:

\begin{definition}
    Each element of $(W_0^{m,p}(G))'$, the dual
    space of $W_0^{m,p}(G)$, determines a $\cD'(G)$
    generalized function. All of these generalized
    functions form a subspace of $\cD'(G)$, and we
    denote it by $W^{-m,p'}(G)$, where $p'$ is the
    H\"older conjugate of $p$. Especially, we denote
    $ H^{-m}(G)=W^{-m,2}(G)$.
\end{definition}

$W^{-m,q}(G)$ is a Banach space with  the canonical norm, i.e.,
$$
|F|_{W^{-m,q}(G)}=\sup_{\f\in
    W_0^{m,p}(G)\setminus\{0\}}\frac{F(\f)}{|\f|_{
        W_0^{m,p}(G)}},\qq\forall\;F\in W^{-m,q}(G).
$$
Especially, $H^{-m}(G)$ is a Hilbert space.

\begin{remark}\label{2.6-rmk4} Since $H^m_0( G)$ is
    a Hilbert space, by the classical Riesz
    representation theorem (i.e., Theorem \ref{1t10s}), there exists an
    isomorphism between $H_0^m( G)$ and
    $H^{-m}( G) \equiv(H_0^m( G))'$. But,
    this does not mean that these two
    spaces are the same.  Indeed, the
    elements in $H_0^m( G)$ are functions
    with certain regularities; while the
    elements in $H^{-m}( G)$ need not even
    to be a ``usual" function.
\end{remark}

All the above spaces $W^{m,p}( G)$,
$W^{m,p}_0( G)$ and $W^{-m,p'}(G)$ are called {\it
    Sobolev spaces}.

\subsection{Operator semigroups}\label{ss-semigroup}
Recall that $H$ is a Hilbert
space.

\begin{definition}
    A $C_0$-semigroup   on  $H$ is a
    family of bounded linear operators $\{S(t)\}_{t
        \ge 0}$ (on  $H$) such that $S(t) S(s) =
    S(t+s)$ for any $s,t \ge 0$, $S(0) = I$ and for any $ y\in H$,
    $\ds\lim_{t\to0+} S(t)y=y$ in $ H$. When $|S(t)|_{\cL(H)}\leq 1$, the semigroup is called contractive.
\end{definition}

For any $C_0$-semigroup $\{S(t)\}_{t\geq 0}$ on $H$, define a linear operator $A$ on $H$ as
follows:
\begin{equation}\label{e:defGen}
    \begin{cases}\ds
        D(A)\deq \Big\{x\in H\;\Big|\;\lim_{t \to 0+}
        \frac{S(t)x - x}{t}  \mbox{ exists in
        }H\Big\},\\
        \ns\ds Ax =\lim_{t \to 0+} \frac{S(t)x -
            x}{t},\q\forall\; x\in D(A).
    \end{cases}
\end{equation}
We say that the above operator $A$ is the {\it infinitesimal  generator} of
$\{S(t)\}_{t\geq 0}$, or $A$ {\it generates} the
$C_0$-semigroup $\{S(t)\}_{t\geq 0}$. One can show that, $D(A)$ is dense in
$H$, $A$  is a
closed operator, and $S(t)D(A)\subset D(A)$ for all $t\geq
0$.

The following results hold:

\begin{proposition}\label{2.4-prop13}
    Let $\{S(t)\}_{t\geq 0}$ be a $C_0$-semigroup
    on $H$.  Then

    {\rm 1)}  There exist positive constants
    $M$ and $\a$ such that
    $
    |S(t)|_{\cL(H)}\leq Me^{\a t}$, $\forall\;t\geq
    0$;

    {\rm 2)} $\ds\frac{dS(t)x}{dt}  = AS(t)x = S(t)Ax$ for every $x \in D(A)$
    and $t\ge0$;

    {\rm 3)} For every
    $z \in D(A^*)$ and $x \in H$, the map
    $t \mapsto \lan z,S(t)x\ran_H$ is
    differentiable and
    $\ds\frac{d }{dt}\lan z,S(t)x\ran_H$ $= \lan A^*z,S(t)x
    \ran_H$;

    {\rm 4)} If a function $y\in C^1([0,+\infty); D(A))$ satisfies
    $\ds\frac{dy(t)}{dt} = A y(t)$ for every $t \in [0,+\infty)$,
    then $y(t) = S(t) y(0)$;

    {\rm 5)}
    $\{S(t)^*\}_{t\geq 0}$ is  a $C_0$-semigroup on $H$ with the infinitesimal generator $A^*$.
\end{proposition}

Stimulating from the conclusions 2) and 4) in Proposition \ref{2.4-prop13}, it is natural to consider the following evolution equation:
\bel{1equa4}
\left\{
\ba{ll}
\ds y_t(t) = A y(t),\qq\forall\;t>0,\\[2mm]
\ds y(0)=y_0,
\ea
\right.
\ee
and for any $y_0\in H$, we call $y(\cd)=S(\cd)y_0$ the {\it mild solution} to (\ref{1equa4}).

The following result is quite useful to show that some specific unbounded operator under consideration generates a contraction semigroup.

\begin{theorem}\label{2.4-th12.1}
    Let  $A: \cD(A)\subset H\to H$ be a linear
    densely defined closed operator. If
    $$
    \lan Ax,x \ran_H\leq 0,\q \forall\; x\in D(A)
    $$
    and
    $$
    \lan A^*x',x' \ran_H\leq 0,\q \forall\; x'\in D(A^*),
    $$
    then $A$ generates a contraction semigroup on $H$.
\end{theorem}
%

%Next, for any $C_0$-semigroup $\{S(t)\}_{t\geq
%    0}$, we define
%\begin{equation*}\label{2.4-eq30}
%\!\begin{array}{ll}\ds
%\cP(A)\!=\3n\3n&\ds\Big\{P \;\Big|\;P\hb{ is an closed
%    operator on $H$, }\cD(P)\supset\cD(A);\;
% \exists\;
%K_t< \infty,\;t>0,\\
%\ns&\ds\qq\q\int_0^1K_tdt<\infty,\;|Pe^{At}x|_H\le
%K_t|x|_H, \;\forall\; x\in\cD(A)\;\Big\}.
%\end{array}
%\end{equation*}
%
%Clearly, $\cL(H)\subset\cP(A)$. We have the
%following result for the perturbation of
%generators of $C_0$-semigroups.

%
%\begin{proposition}\label{2.4-prop15}
%    Let $\{S(t)\}_{t\geq 0}$ be a $C_0$-semigroup on
%    $H$. Then, for any $P\in\cP(A)$, the operator
%    $A+P$ generates a $C_0$-semigroup  $H$.
%\end{proposition}

%\begin{example}\label{2.4-ex21}
%    Let $A\in\dbR^{n\times n}$ be an $n\times n$
%    matrix. Then,
%
%    \begin{equation}\label{2.4-eq45}
%    S(t)=\sum_{k=1}^\infty{\frac{A^kt^k}{k!}},\qq
%    t\ge0.
%    \end{equation}
%
%    is a $C_0$-semigroup on $H=\dbR^n$ and  $e^{At}y(0)$ is the  solution of the ordinary
%    differential equation $\frac{dy(t)}{dt}=Ay(t)$.
%\end{example}

\begin{example}\label{2.4-ex22}
    Let $G$ be given as in Subsection \ref{sec-mpr-6.1} and
    convex, and let $O\in\dbR^n$
    satisfy $|O|_{\dbR^n}=1$. Take
    $H=L^2(G)$.  Write $\G ^-\deq \{x\in
    \G \;|\; O\cd\nu(x)\leq 0\}$. Define an
    unbounded linear operator $A$ on $H$ as
    follows:
    $$
    \begin{cases}
        D(A)=\big\{f\in H^1(G)\;\big|\; f=0 \mbox{ on }\G ^-\big\},\\
        \ns\ds Af=-O\cd\nabla f,\qq \forall\; f\in  D(A).
    \end{cases}
    $$
    Then, by Theorem \ref{2.4-th12.1}, $A$ generates a contraction semigroup
    $\{S(t)\}_{t\geq 0}$ on $L^2(G)$. This semigroup
    arises in the study of (homogeneous) transport
    equation:
    \begin{equation}\label{2.4-eq57.1}
        \left\{
        \begin{array}{ll}\ds
            y_t+O\cd\nabla y=0 &\mbox{ in }G\times(0,+\infty),\\
            \ns\ds y=0 &\mbox{ on }\G ^-\times(0,+\infty),\\
            \ns\ds
            y(0)=y_0 &\mbox{ in }G.
        \end{array}
        \right.
    \end{equation}
\end{example}

\begin{example}\label{2.4-ex23}
    Take $H=L^2(G)$, and define an unbounded linear operator $A$ on $H$
    as follows:
    $$
    \begin{cases}
        D(A)=H^2(G)\bigcap H_0^1(G),\\
        \ns\ds Af=\D f =\sum_{k=1}^n\frac{\pa^2 f}{\pa
            x_k^2},\qq \forall\; f\in  D(A).
    \end{cases}
    $$
    Then, by Theorem \ref{2.4-th12.1}, $A$ generates a contraction semigroup
    $\{S(t)\}_{t\geq 0}$ on $L^2(G)$. This semigroup
    arises in the study of (homogeneous) heat
    equation:
    \begin{equation}\label{2.4-eq57}
        \left\{
        \begin{array}{ll}\ds
            y_t-\D y=0 &\mbox{ in }G\times(0,+\infty),\\
            \ns\ds y=0 &\mbox{ on }\G \times(0,+\infty),\\
            \ns\ds
            y(0)=y_0 &\mbox{ in }G.
        \end{array}
        \right.
    \end{equation}
\end{example}
\begin{example}\label{2.4-ex24}
    Let $H=H_0^1(G)\times L^2(G)$ and define an unbounded linear operator $A$ on $H$
    as follows:
    \begin{equation*}\label{2.4-eq58}
        \left\{
        \begin{array}{ll}\ds
            D(A)=\big\{\left(
            y,
            z
            \right)^\top
            \;\big|\; y\in H^{2}(G)\cap H^{1}_0(G),z\in
            H_0^1(G)\big\},\\
            \ns\ds A\left(
            \begin{array}{c}
                y \\
                z \\
            \end{array}
            \right)=\left(
            \begin{array}{c}
                z \\
                \D
                y \\
            \end{array}
            \right)
            \equiv\left(
            \begin{array}{cc}
                0 & I \\
                \D & 0 \\
            \end{array}
            \right)
            \left(
            \begin{array}{c}
                y \\
                z \\
            \end{array}
            \right),\qq \forall\; \left(
            y,
            z
            \right)^\top\in D(A).
        \end{array}
        \right.
    \end{equation*}
    By Theorem \ref{2.4-th12.1},  $A$ generates a contraction semigroup on $H$.  This semigroup arises in the
    study of the following (homogeneous) wave
    equation:
    \begin{equation}\label{2.4-eq63}
        \left\{
        \begin{array}{ll}\ds
            y_{tt}-\D y=0 &\hb{ in }G\times(0,+\infty),\\
            \ns\ds y =0 &\hb{ on }\G \times(0,+\infty),\\
            \ns\ds y(0)=y_0,\q y_t(0)=y_1 &\hb{ in }G.
        \end{array}
        \right.
    \end{equation}
    In fact, if we set $z=y_t$, then,
    \eqref{2.4-eq63} can be transformed into the
    following equation:
    \begin{equation}\label{2.4-eq64}
        \begin{cases}\ds
            w_t=Aw&\hb{ in } (0,+\infty),\\
            \ns\ds  w(0)=w_0,
        \end{cases}
    \end{equation}
    where $w=(y,z)$ and $w_0=(y_0,y_1)$.
\end{example}

As applications of $C_0$-semigroup, people are concerned with the well-posedness of the following
semilinear evolution equation:
\begin{equation}\label{2.5-eq22}
    \left\{
    \begin{array}{ll}\ds
        y_t(t)=Ay(t)+f(t,y(t)),\qq
        t\in(0,T],\cr
        \ns\ds y(0)=y_0,
    \end{array}
    \right.
\end{equation}
where $A:\cD(A)(\subset H)\to H$
generates a $C_0$-semigroup
$\{S(t)\}_{t\geq 0}$ on $H$, $y_0\in H$ and
$f:[0,T]\times H\to H$ satisfies the
following:

\ss

1) For each $x\in H$, $f(\cd,x):[0,T]\to H$ is
strongly measurable; and

\ss

2) There exists a constant $\cC>0$, such that, for a.e. $t\in[0,T]$,
\begin{equation}\label{2.5-eq23}
    \left\{
    \begin{array}{ll}\ds
        |f(t,x_1)-f(t, x_2)|_H\le \cC|x_1-
        x_2|_H,\ \q \forall\; x_1,\,x_2 \in
        H,\\
        \ns\ds|f(t,0)|_H\le \cC.
    \end{array}
    \right.
\end{equation}

\begin{definition}\label{2.5-def1}
    1) We call $y\in C([0,T];H)$ a {\it
        strong solution} to \eqref{2.5-eq22} if
    $y(0)=y_0$, $y$ is
    differentiable for a.e.
    $t\in (0,T)$, $y(t)\in\cD(A)$ for  $\ae
    t\in(0,T)$ and $y_t(t)=Ay(t)+f(t,y(t))$  for a.e.
    $t\in (0,T)$.

    \ss

    2) We call  $y\in C([0,T];H)$ a {\it
        weak solution} to \eqref{2.5-eq22} if
    for any $\f \in\cD(A^*)$ and   $t\in[0,T]$,
    \begin{equation}\label{2.5-eq2}
        \begin{array}{ll}\ds
            \langle y(t),\f \rangle_H =\langle y_0,\f
            \rangle_H+ \int_0^t\big(\langle y(s),A^*\f
            \rangle_H+\langle f(s,y(s)),\f
            \rangle_H\big)ds.
        \end{array}
    \end{equation}

    3) We call  $y\in C([0,T];H)$  a
    {\it mild solution} to \eqref{2.5-eq22}
    if
    \begin{equation}\label{2.5-eq3}
        y(t)=S(t)y_0+\int_0^t S(t-s) f(s,y(s))ds,\qq
        t\in[0,T].
    \end{equation}
\end{definition}

It is not difficult to show the following result.
\begin{proposition}\label{2.5-prop2}
    A function $y\in C([0, T];H)$ is a mild solution
    to \eqref{2.5-eq22} if and only if it
    is a weak solution to the same equation.
\end{proposition}

Hereafter, we will not distinguish
mild  and
weak solutions to \eqref{2.5-eq22}.
The following result is concerned with
the existence and uniqueness of
solutions to \eqref{2.5-eq22}.

\begin{proposition}\label{2.5-prop3} For any
    $y_0\in H$, \eqref{2.5-eq22} admits a
    unique mild solution $y$, and
    \begin{equation*}
        |y(\cd)|_{C([0,T];H)}\le \cC(1+|y_0|_H).
    \end{equation*}
\end{proposition}

\ms

We note that the solution $y$ to
\eqref{2.5-eq22} is not necessarily
differentiable because $A$ is usually
unbounded. This may bring some
trouble in some situations. The
following convergence result sometimes
can help us to remove such an obstacle.

\begin{proposition}\label{2.5-prop4}
    Let $y$ be the mild solution to
    \eqref{2.5-eq22} and let $y^\l$ be the strong
    solution to the following equation:
    \begin{equation}\label{2.5-eq27}
        y^\l(t)=S_\l
        (t)y_0+\int_0^t S_\l(t-s)f(s,y^\l(s))ds,\;\,
        t\in[0,T],
    \end{equation}
    where $\{S_\l(t)\}_{t\geq 0}$ is the $C_0$-semigroup generated by $A_\l=\l A(\l I-A)^{-1}$, where $\l\in\rho(A)$, the resolvent set of $A$. Then,
    \begin{equation}\label{2.5-eq28}
        \lim_{\l\to\infty}\sup_{t\in[0,T]}|y^\l(t)-y(t)|_H=0.
    \end{equation}
\end{proposition}

\subsection{Probability,
    random variables and expectation}

We first recall that a {\it probability
    space} is simply a measure space
$(\Omega,\cF,\dbP)$ for which
$\dbP(\Omega)=1$. In this case, $\dbP$
is called a probability measure.

In the sequel, we fix a
probability space $(\Omega,\cF,\dbP)$ and a Banach space $\cX$. Any
$\omega\in\Omega$ is called a {\it sample
    point}; any $A\in \cF$ is called an
{\it event}, and $\dbP(A)$ represents the
probability of the event $A$. If an event $A\in \cF$ is such that $\dbP(A) = 1$, then we may alternatively
say that $A$ holds, $\dbP$-a.s., or simply $A$ holds a.s. (if the probability $\dbP$ is clear
from the context).

Any $\cX$-valued, strongly
measurable function $f:\
(\Omega,\cF)\to (\cX,\cB(\cX))$ is
called an ($\cX$-valued) {\it random
    variable}.
If $f$ is Bochner integrable w.r.t.  the
measure $\dbP$, then we denote the integral by $\mE f$ and call it the {\it mean} or {\it mathematical expectation} of $f$.

For any $p\in [1,\infty]$, one can define the Banach spaces
$L_{\mathcal F}^p(\Omega;\cX)\equiv
L^p(\Omega,{\mathcal F},\dbP;\cX)$ and $L_{\mathcal F}^p(\Omega)$ as that in Subsection \ref{sec-mea-3}.  For
any $f\in L_{{\mathcal F}}^2(\Omega;\dbR^m)$ (for some $m\in\dbN$), we define
the {\it variance} of $f$ by
$$
\hb{Var}\; f=\dbE[(f-\mE f)(f-\mE f)^\top].
$$

Let $A,B\in\cF$. We say that $A$ and $B$ are
independent if $\dbP(A\cap B)=\dbP(A)\dbP(B)$.
Let $\cJ_1$ and $\cJ_2$ be two subfamilies of $\cF$.
We say that $\cJ_1$ and $\cJ_2$ are independent
if $$\dbP(A\cap B)=\dbP(A)\dbP(B),\qq  \forall\;
(A,B)\in \cJ_1\times\cJ_2.$$ Let $f,\ g:\
(\Omega,\cF)\to(\cX,\cB(\cX))$ be two random
variables. We say that $f$ and $g$ (\resp $f$
and $\cJ_1$) are independent if $\si(f)$ and
$\si(g)$ (\resp $\si(f)$ and $\cJ_1$) are
independent.

Let $X=(X_1,\cds,X_m):\ (\Omega,\cF)\to (\dbR^m,\cB(\dbR^m))$ be a
random variable. We call
$$
F(x)\equiv F(x_1,\cds,x_m)\deq  \dbP\{X_1\le x_1,\cds,X_m\leq x_m\}
$$
the \emph{distribution function} of
$X$. If for some nonnegative function
$f(\cd)$, one has
$$
F(x)\equiv F(x_1,\cds,x_m)=\int_{-\infty}^{x_1}\cds \int_{-\infty}^{x_m}f(\xi_1,\cds,\xi_m)d\xi_1\cds d\xi_m,
$$
then the function $f(\cd)$ is called
the \emph{density} of $X$. If
$$
f(x)=(2\pi)^{-\frac{m}{2}}
|{\rm det}Q|^{-1/2}\exp\Big\{-\frac{1}{2\mu}(x-\l)^\top Q^{-1} (x-\l)\Big\},
$$
where $\l\in\dbR^m$  and $Q$  is a positive definite $m$-dimensional matrix, then $X$ is
called a {\it normally distributed random variable}
(or $X$ has a {\it normal distribution}) and denoted by
$X\sim \cN(\l,Q)$.  When $\l=0$ and $Q$ is the  $m$-dimensional identity matrix $I_m$, we call $X$  a {\it standard normally distributed random variable}.

%%%%%%%%%%%%%%%%%%%%%%%%%%%%%%%%%%%%%%%

\subsection{Stochastic process}

Let $\cI=[0,T]$ with
$T>0$. A
family of $\cX$-valued random variables
$\{X(t)\}_{t\in \cI}$ on
$(\Omega,\cF,\dbP)$  is called a
{\it stochastic process}. In the sequel, we shall interchangeably use
$\{X(t)\}_{t\in \cI}$, $X(\cd)$ or even
$X$ to denote a (stochastic) process.
For any
$\omega\in\Omega$, the map $t\mapsto
X(t,\omega)$ is called a {\it sample path}
(of $X$).
$ X(\cdot)$
is said to be continuous ({\it resp.}
c\'adl\`ag, i.e., right-continuous with
left limits) if there is a $\dbP$-null
set $N \in\cF$, such that for any
$\omega\in\Omega\setminus N$, the
sample path $X(\cdot,\omega)$ is
continuous (\resp c\'adl\`ag) in $\cX$.

\begin{definition}
    Two ($\cX$-valued) processes $X(\cdot)$
    and $\cl X(\cdot)$ are said to be
    {\it stochastically equivalent}
    if
    $$
    \dbP(\{X(t)=\cl X(t)\})=1,\qq \forall\; t\in
    \cI.
    $$
    In this case, one is said to be a
    {\it modification} of the other.
\end{definition}
\begin{definition}
    We call a family of sub-$\si$-fields
    $\{\cF_t\}_{t\in \cI}$ in $\cF$ a
    {\it filtration} if
    $\cF_{t_1}\subset \cF_{t_2}$ for all $
    t_1,t_2\in \cI \mbox{ with } t_1\leq t_2$. For
    any $t\in [0,T)$, we put
    $$
    \cF_{t+}\deq \bigcap_{s\in (t,T]}\cF_s,\qq
    \cF_{t-}\deq \bigcup_{s\in [0,t)}\cF_s.
    $$
    If $\cF_{t+}=\cF_t$ (\resp $\cF_{t-}=\cF_t$),
    then $\{\cF_t\}_{t\in \cI}$ is said to be right
    (\resp left) continuous.
\end{definition}

In the sequel,
we simply denote $\{\cF_t\}_{t\in
    \cI}$ by $\mathbf{F}$ unless we want to emphasize what $\cF_t$
or $\cI$ exactly is. We call
$(\Omega,\cF,\mathbf{F}, \dbP)$ a {\it filtered
    probability space}.

\begin{definition}\label{12.30-def1}
    We say that $(\Omega,\cF,\mathbf{F}, \dbP)$
    satisfies
    {\it the usual condition} if
    $(\Omega,\cF,\dbP)$ is complete, $\cF_0$
    contains all $\dbP$-null sets in $\cF$, and
    $\mathbf{F}$ is right continuous.
\end{definition}

In what follows, unless
otherwise stated, we always assume that $(\Omega,\cF,\mathbf{F}, \dbP)$
satisfies
the usual condition.

\begin{definition}
    Let $X(\cd)$ be an $\cX$-valued
    process.

    {\rm 1)}
    $X(\cd)$ is said to be {\it measurable} if the map
    $(t,\omega)\mapsto X(t,\omega)$ is strongly
    $(\cB(\cI)\t \cF)/\cB(\cX)$-measurable;

    {\rm 2)}
    $X(\cd)$ is said to be {\it $\mathbf{F}$-adapted} if
    it is measurable, and for each $t\in \cI$, the
    map $\omega\mapsto X(t,\omega)$ is strongly
    $\cF_t/\cB(\cX)$-measurable;

    {\rm 3)}  $X(\cd)$ is said to be
    {\it $\mathbf{F}$-progressively measurable} if for
    each $t\in \cI$, the map $(s,\omega)\mapsto
    X(s,\omega)$ from $[0,t]\times\Omega$ to $\cX$
    is strongly $(\cB([0,t])\t
    \cF_t)/\cB(\cX)$-measurable.
\end{definition}
\begin{definition}
    A set $A\in \cI\times \Omega$ is called
    {\it progressively measurable w.r.t.  $\mathbf{F}$} if
    the process $\chi_{A}(\cdot)$ is progressively
    measurable.
\end{definition}
The class of all progressively
measurable sets is a $\si$-field, called the
{\it progressive $\si$-field w.r.t.  $\mathbf{F}$},
denoted by $\dbF$. One can show that, a process
$\f:[0,T]\times\Omega\to \cX$ is
$\mathbf{F}$-progressively measurable if and
only if it is strongly $\dbF$-measurable.

It is clear that if $X(\cd)$ is
$\mathbf{F}$-progressively measurable, it must
be $\mathbf{F}$-adapted. Conversely, it can be
proved that, for any $\mathbf{F}$-adapted
process $X(\cd)$, there is an
$\mathbf{F}$-progressively measurable process
$\wt X(\cd)$ which is stochastically equivalent
to $X(\cd)$. For this reason, in the sequel, by
saying that a process $X(\cd)$ is
$\mathbf{F}$-adapted, we mean that it is
$\mathbf{F}$-progressively measurable.

For any $p,q\in[1,\infty)$, write
$$
\begin{array}{ll}
    \ds L^p_\dbF(\Omega;L^q(0,T;\cX)) \deq  \Big\{ \f:
    (0,T) \times \Omega \to \cX\;\Big|\; \f(\cd)\hb{
        is $\mathbf{F}$-adapted and }\\
    \ns\ds\hspace{6.73cm}\dbE\(\int_0^T |\f(t)|_\cX^qdt\)^{\frac{p}{q}} < \infty\Big\},\\
\end{array}
$$
and
$$
\begin{array}{ll}
    \ds  L^q_\dbF(0,T;L^p(\Omega;\cX)) \deq
    \Big\{ \f :(0,T) \times \Omega \to  \cX\;
    \Big|\; \f(\cd)\hb{ is
        $\mathbf{F}$-adapted and
    }\\
    \ns\ds\hspace{6.73cm}\int_0^T
    \(\dbE|\f(t)|_\cX^p\)^{\frac{q} {p}}dt<
    \infty\Big\}.
\end{array}
$$
Similarly, we may also define (for $1\le
p,q<\infty$)
$$
\left\{
\begin{array}{ll}
    \ds L^\infty_\dbF(\Omega;L^q(0,T;\cX)),\q
    L^p_\dbF(\Omega;L^\infty(0,T;\cX)),\q
    L^\infty_\dbF
    (\Omega;L^\infty(0,T;\cX)),\\
    \ns\ds L^\infty_\dbF(0,T;L^p(\Omega;\cX)),\q
    L^q_\dbF(0,T;L^\infty(\Omega;\cX)),\q
    L^\infty_\dbF(0,T;L^\infty(\Omega;\cX)).
\end{array}
\right.
$$
All these spaces are Banach spaces (with the
canonical norms). In what follows, we shall simply
denote $L^p_\dbF(\Omega;L^p(0,T;\cX))\equiv
L^p_\dbF(0,T;L^p(\Omega;\cX))$ by
$L^p_\dbF(0,T;\cX)$; and further simply write $L^p_\dbF(0,T)$
for $L^p_\dbF(0,T;\dbR)$.

For any $p\in [1,\infty)$, set
$$
\begin{array}{ll}\ds
    L^{p}_{\dbF}(\Omega;C([0,T];\cX))\deq \Big\{\f:[0,T]\times\Omega\to
    \cX \;\Big|\; \f(\cd)\hb{
        is continuous, }\\
    \ns\ds\qq\qq\qq\qq \qq\qq
    \mathbf{F}\mbox{-adapted and
    }\dbE\big(|\f(\cdot)|_{C([0,T];\cX)}^p\big)<\infty\Big\}
\end{array}
$$
and
$$
\begin{array}{ll}\ds
    C_{\dbF}([0,T];L^{p}(\Omega;\cX))\deq \Big\{\f:[0,T]\times\Omega\to
    \cX\;\Big|\;\f(\cd)\hb{
        is $\mathbf{F}$-adapted }\\
    \ns\ds\qq\qq\qq\qq \qq \mbox{ and }\f(\cd):[0,T]
    \to L^p_{\cF_T}(\Omega;\cX)\mbox{ is
        continuous}\Big\}.
\end{array}
$$
One can show that both
$L^{p}_{\dbF}(\Omega;C([0,T];\cX))$ and
$C_{\dbF}([0,T];L^{p}(\Omega;\cX))$ are Banach
spaces with norms
$$|\f(\cd)|_{L^{p}_{\dbF}(\Omega;C([0,T];\cX))} =
\big[\dbE|\f(\cdot)|_{C([0,T];\cX)}^p\big]^{1/p}$$
and
$$|\f(\cd)|_{C_{\dbF}([0,T];L^{p}(\Omega;\cX))}=\max_{t\in
    [0,T]}
\big[\dbE|\f(t)|_\cX^p\big]^{1/p},$$
respectively.
Also, we denote by
$D_{\dbF}([0,T];L^p(\Omega;\cX))$ the Banach
space of all processes which are
right continuous with left limits in $L^p_{\cF_T}(\Omega;\cX)$  w.r.t.
$t\in [0,T]$,  with the norm
$\left[\mathbb{E}|\f(\cdot)|^p_\cX\right]^{1/p}_{L^{\infty}(0,T)}$.

\ms

\begin{definition}\label{def-Brown1}
    A continuous, $\dbR^m$-valued,
    $\mathbf{F}$-adapted process
    $\{W(t)\}_{t\ge 0}$ is called an (standard) {\it
        $m$-dimensional Brownian
        motion}, if

    {\rm 1)} $\dbP(\{W(0)=0\})=1$; and

    {\rm 2)} For any $s,t\in [0,\infty)$
    with $0\le s<t<\infty$, the random
    variable $W(t)-W(s)$ is independent of
    $\cF_s$, and
    $W(t)-W(s)\sim\cN(0,(t-s)I_m)$.
\end{definition}

Similarly, one can define
$\dbR^m$-valued Brownian motions over
any time interval $[a,b]$ or $[a,b)$
with $0\le a<b\le \infty$.

In the rest of this paper, unless otherwise stated,
we fix a $1$-dimensional standard
Brownian motion $W(\cd)$ on
$(\Omega,\cF,\mathbf{F}, \dbP)$. Write
\begin{equation}\label{a2.2}
    \cF_t^W\deq \si(W(s);\; s\in [0,t])\subset\cF_t,\qq\forall\;t\in\cI.
\end{equation}
Generally, the filtration
$\{\cF_t^W\}_{t\in\cI}$ is left-continuous, but
not necessarily\; right-continuous.
Nevertheless, the augmentation $\{\hat
\cF_t^W\}_{t\in\cI}$ of $\{ \cF_t^W\}_{t\in\cI}$
by adding all $\dbP$-null sets is continuous,
and $W(\cd)$ is still a Brownian motion on the
(augmented) filtered probability space
$(\Omega,\cF,\{\hat\cF_t^W\}_{t\in\cI}, \dbP)$.
In the sequel, by saying that $\mathbf{F}$ is
the natural filtration generated by $W(\cd)$, we
mean that $\mathbf{F}$ is generated as in
(\ref{a2.2}) with the above augmentation, and
hence it is continuous.

\subsection{It\^o's integral}

Let $H$ be a Hilbert space.
We now define the {\it It\^o integral}
\bel{a2.3}
\int_0^TX(t)dW(t)
\ee
of an $H$-valued, $\mathbf{F}$-adapted
stochastic process $X(\cd)$ (satisfying suitable
conditions) w.r.t. $W(t)$.
Note that one cannot define (\ref{a2.3}) to be a
Lebesgue-Stieltjes type integral by regarding
$\omega$ as a parameter. Indeed, the map
$t(\in[0,T])\mapsto W(t, \cd)$ is not of bounded
variation,  a.s.

Write $\cL_0$ for  the set of $f\in
L^2_{{\mathbb{F}}}(0,T;H)$ of the form
\bel{6.2}
\displaystyle f(t,\omega)=\sum_{j=0}^n
f_j(\omega)\chi_{[t_j,t_{j+1})}(t), \qq (t,\omega)\in[0,T]\t\Omega,
\ee
where $n\in\dbN$, $0=t_0<t_1<\cdots<t_{n+1}= T$, $f_j$ is
${\mathcal F}_{t_j}$-measurable with
$$\sup\big\{|f_j(\omega)|_H\;\big|\;
j\in\{0,\cdots,n\},\,\omega\in\Omega\big\}<\infty.$$
One can show that $\cL_0$ is dense in
$L_{{\mathbb{F}}}^2(0,T;H)$.

Assume that $f\in\cL_0$ takes the form of
(\ref{6.2}). Then we set
\begin{equation}\label{6.7}
    I(f)(t,\omega)=\sum_{j=0}^n f_j(\omega)[W(t\wedge
    t_{j+1},\omega)-W(t\wedge t_j,\omega)].
\end{equation}
It is easy to show that $I(f )\in L^2_{\cF_t}(\Om;H)$
and the following {\it It\^o isometry} holds:
\begin{equation}\label{6oo7}
    |I(f)|_{L^2_{\cF_t}(\Om;H)}=|f
    |_{L_{{\mathbb{F}}}^2(0,t;H)}.
\end{equation}
Generally, for $f \in
L_{{\mathbb{F}}}^2(0,T;H)$, one can
find a sequence of $\{f
_k\}_{k=1}^\infty\subset \cL_0$ such
that
$$
\lim_{k\to\infty}|f_k-f |_{L_{{\mathbb{F}}}^2(0,T;H)}= 0.
$$
Since
$
|I(f _k)-I(f _j)|_{L^2_{\cF_t}(\Om;H)}=|f _k-f
_j|_{L_{{\mathbb{F}}}^2(0,t;H)},
$
$\{I(f_k)\}_{k=1}^\infty$ is a Cauchy
sequence in $L^2_{\cF_t}(\Om;H)$ and hence, it
converges to a unique element in $
L^2_{\cF_t}(\Om;H)$, which is determined uniquely by
$f$ and is independent of the particular choice
of $\{f _k\}_{k=1}^\infty$. We call this element the {\it It\^o integral}  of $f$ (w.r.t. the Brownian
motion $W(\cd)$) on $[0,t]$ and denote it by $\int_0^tf  dW$.

For $0\leq s <t \leq T$, we call $\int_0^tf  dW-\int_0^sf  dW$ the
{\it It\^o integral} of $f \in
L_{{\mathbb{F}}}^2(0,T;H)$ (w.r.t. the Brownian
motion $W(\cd)$) on $[s,t]$. We shall denote it
by $\int_s^tf (\tau)dW(\tau)$ or simply $\int_s^tf dW$.

For any $p\in (0,\infty)$, write
$L_\dbF^{p,loc}(0,T;H)$ for the set of all $H$-valued,
$\mathbf{F}$-adapted stochastic processes
$f(\cd)$ satisfying only $\int_0^T|f(t)|_{H}^pdt
<\infty$,  a.s. For any $p\in [1,\infty)$, one can also define the
It\^o integral $\int_0^tf  dW$ for $f\in L_\dbF^{p,loc}(0,T;H)$, especially for $f\in
L_\dbF^p(\Om;L^2(0,T;H))$ (See \cite{LZ3.1} for
more details). The It\^o integral has the following properties.
\begin{theorem}
    Let $p\in [1,\infty)$. Let $f, g \in
    L_\dbF^p(\Om;L^2(0,T;H))$ and $a,b\in
    L_{\cF_s}^2(\Omega)$, $0\le s<t\leq T$. Then

    \ms

    {\rm 1)}
    $\ds
    \int_s^t f dW\in
    L^p_\dbF(\Om;C([s,t];H));
    $

    \ms

    {\rm 2)}
    $\ds
    \int_s^t(af +b g )dW=a\int_s^tf dW+b\int_s^t g
    dW,\;  \as\!\!;
    $

    \ms

    {\rm 3)} $\ds\dbE\Big(\int_s^t f
    dW \Big)=0;
    $

    \ms

    {\rm 4)} When $p\geq 2$,
    $\ds\dbE\Big(\Big\langle\int_s^tf dW,\;\int_s^t g
    dW\Big\rangle_H \Big)=\dbE\Big(\int_s^t \lan f (r,\cd),g
    (r,\cd) {\ran}_Hdr \Big)$.

\end{theorem}

The following result, known as the {\it
    Burkholder-Davis-Gundy inequality}, links
It\^o's integral to the Lebesgue/Bochner
integral.

\begin{theorem}\label{BDG} For any $p\in [1,\infty)$, there exists a constant
    $\cC_p>0$ such that for any $T>0$ and $f\in
    L_\dbF^p(\Om;L^2(0,T;H))$,
    \begin{equation*}\label{BDGQ-eq2}
        \ba{ll} \ds\frac{1}{\cC_p}
        \dbE\(\int_0^T|f(s)|_H^2 ds\)^{\frac{p}{2}}\leq
        \dbE\(\sup_{t\in [0,T]}\Big|\int_0^t
        f(s)dW(s)\Big|_H^p\)\leq \cC_p
        \dbE\(\int_0^T|f(s)|_H^2 ds\)^{\frac{p}{2}}. \ea
    \end{equation*}
\end{theorem}

We need the following notion for an important class of stochastic processes.

\begin{definition} An $H$-valued, $\mathbf{F}$-adapted,  continuous
    process $X(\cd)$ is called an {\it It\^o process} if
    there exist two $H$-valued stochastic processes
    $\phi(\cd)\in L_\dbF^{1,loc}(0,T;H)$ and
    $\Phi(\cd)\in L_\dbF^{2,loc}(0,T;H)$ such that
    \begin{equation}\label{5.888}
        X(t)=X(0)+\int_0^t\phi(s)ds+\int_0^t\Phi(s)dW(s),\quad
        \hbox{a.s.}, \;\;\forall\; t\in[0,T].
    \end{equation}
\end{definition}

The following fundamental result is known as
{\it It\^o's formula}.

\begin{theorem}\label{c1t7.1} Let
    $X(\cd)$ be given by \eqref{5.888}. Let
    $F: [0,T]\t H\to \dbR$ be a function
    such that its partial derivatives
    $F_t$, $F_x$ and $F_{xx}$ are uniformly
    continuous on any bounded subset of
    $[0,T]\t H$. Then,
    \begin{equation}\label{7e2}
        \begin{array}{ll}
            F(t,X(t))- F(0,X(0))\\
            \ns\ds=\int_0^t F_x(s,X(s))
            \Phi(s)dW(s)+\int_0^t\Big[F_t(s,X(s))+
            F_x(s,X(s))\phi(s)\\
            \ns \ds\q+\frac{1}{2}\langle
            F_{xx}(s,X(s))\Phi(s),\Phi(s)\rangle_{H}\Big]ds,\quad
            \hbox{a.s.}, \;\;\forall\; t\in[0,T].
        \end{array}
    \end{equation}
\end{theorem}
\begin{remark}
    Usually, people write the formula \eqref{7e2} in the following differential form:
    $$
    \begin{array}{ll}\ds
        dF(t,X(t))\3n&\ds=F_x(t,X(t))
        \Phi(t)dW(t)+F_t(t,X(t))dt +
        F_x(t,X(t))\phi(t) dt\\
        \ns&\ds\q + \frac{1}{2}\langle
        F_{xx}(t,X(t))\Phi(t),\Phi(t)\rangle_{H}dt .
    \end{array}
    $$
\end{remark}

Theorem \ref{c1t7.1} works well for
It\^o processes in the (strong) form
\eqref{5.888}. However, usually this is
too restrictive in the study of
stochastic differential equations in
infinite dimensions. Indeed, in the
infinite dimensional setting sometimes
one has to handle It\^o processes in
a weaker form, to be presented below.

Let $\cV$ be a Hilbert space such that
the embedding $\cV\subset H$ is
continuous and dense. Denote by $\cV^*$
the dual space of $\cV$ w.r.t.
the pivot space $H$. Hence,
$ \cV\subset H =H^*\subset \cV^*, $
continuously and densely and
$$
{\lan z,v\ran}_{\cV,\cV^*} = {\lan
    z,v\ran}_H,\qq \forall\; (v,z)\in H\times
\cV.
$$

We have the following  It\^o's formula
for a weak form of It\^o process.

\begin{theorem}\label{c1-ito1}
    Suppose that $X_0\in
    L^2_{\cF_0}(\Omega;H)$,  $\phi(\cd)\in
    L^2_\dbF(0,T;\cV^*)$, and $\Phi(\cd)\in
    L^p_\dbF(\Omega;L^{2} (0,T;H))$ for
    some $p\geq 1$. Let
    \begin{equation}\label{ctoQ1-eq1}
        X(t) =X_0+ \int^t_0 \phi(s)ds+ \int^t_0
        \Phi(s)dW(s),\q t\in [0,T]. \ee
        If  $X\in L^2_\dbF(0,T;\cV)$, then
        $X(\cd)\in C([0,T]; H)$, a.s., and for
        any $t\in [0,T]$,
        \begin{equation}\label{c1-itoQ1-eq1}
            \begin{array}{ll}\ds
                |X(t)|^2_{H}&\3n\ds = |X_0|^2_H + 2
                \int_0^t {\lan \phi(s),
                    X(s)\ran}_{\cV^*,\cV} ds\\
                \ns&\ds \;\;  + 2\int_0^t {\lan
                    \Phi(s),X(s) \ran}_HdW(s)+ \int_0^t
                |\Phi(s)|_{H}^2 ds,\q \mbox{a.s.}
            \end{array}
        \end{equation}
    \end{theorem}
    %

    %As an immediate corollary of Theorem \ref{c1-ito1}, we have the following result.
    %
    %\begin{corollary}\label{c1-ito1-1}
    %    Suppose that $X_0, Y_0\in
    %    L^2_{\cF_0}(\Omega;H)$,  $\phi_1(\cd),$ $\phi_2(\cd)\!\in\!
    %    L^2_\dbF(0, T;\cV^*)$, and $\Phi_1(\cd),\Phi_2(\cd)\!\in\!
    %    L^p_\dbF(\Omega;L^{2} (0,T;H))$ for
    %    some $p\geq 1$. Let
    %
    %    \begin{equation*}
    %        X(t) =X_0+ \int^t_0 \phi_1(s)ds+ \int^t_0
    %        \Phi_1(s)dW(s),\;\, t\in [0,T]
    %    \end{equation*}
    %
    %    and
    %
    %    \begin{equation*}
    %        Y(t) =Y_0+ \int^t_0 \phi_2(s)ds+ \int^t_0
    %        \Phi_2(s)dW(s),\;\, t\in [0,T].
    %    \end{equation*}
    %    %
    %    If  $X\in L^2_\dbF(0,T;\cV)$ (\resp $Y\in L^2_\dbF(0,T;\cV)$), then
    %    $X(\cd)\in C([0,T]; H)$, a.s.(\resp $Y\in C([0,T]; H)$, a.s.), and for
    %    any $t\in [0,T]$,
    %
    %    \begin{eqnarray}\label{c1-itoQ1-eq1-1}
    %        &&\3n\3n\3n\3n
    %        \lan X(t),Y(t)\ran_{H}-\lan X_0,Y_0\ran_{H}\nonumber\\
    %        &&\3n\3n\3n\3n =
    %        \int_0^t (\lan \phi_2(s),
    %        X(s)\ran+\lan \phi_1(s),
    %        Y(s)\ran_{\cV^*,\cV})_{\cV^*,\cV} ds + \int_0^t
    %        \lan \Phi_1(s),\Phi_2(s)\ran_{H}^2 ds \nonumber\\
    %        && \3n\3n  +  \int_0^t \big(\lan
    %        \Phi_2(s),X(s) \ran_H  + \lan
    %        \Phi_1(s),Y(s) \ran_H)dW(s),\q\,   \mbox{a.s.}
    %    \end{eqnarray}
    %
    %\end{corollary}
    %

    %
    \begin{remark}\label{12.24-rmk1}
        For simplicity, we usually write formula \eqref{c1-itoQ1-eq1} in the following differential form:
        $$
        \begin{array}{ll}\ds
            d|X(t)|_H^2 =2\lan \phi(t),
            X(t)\ran_{\cV^*,\cV}dt+
            |\Phi(t)|_{H}^2dt + 2 \lan
            \Phi(t),X(t) \ran_HdW(t)
        \end{array}
        $$
        and denote $|\Phi(t)|_{H}^2dt$ by $|dX|_{H}^2$ for simplicity.
        %    Similarly, we write formula \eqref{c1-itoQ1-eq1-1} in the following differential form:
        %
        %    $$
        %    \!\begin{array}{ll}\ds
        %    d\lan X(t),Y(t)\ran_H^2\\
        %    \ns\ds = \lan \phi_2(t),\!
        %    X(t)\ran_{\cV^*,\cV}\!+\!  \lan \phi_1(t),\!
        %    Y(t)\ran_{\cV^*,\cV}\!+ \!
        %    \lan\Phi_1(t), \Phi_2(t)\ran_Hdt \\
        %    \ns \ds\q +   \lan
        %    \Phi_2(t),X(t) \ran_HdW(t) + \lan
        %    \Phi_1(t),Y(t) \ran_HdW(t)
        %    \end{array}
        %    $$
        %
        %    and denote $\lan\Phi_1(t), \Phi_2(t)\ran_Hdt$ by $ \lan dX,dY\ran_H$ for simplicity.
    \end{remark}
    %

    %%%%%%%%%%%%%%%%%%%%%%%%%%%%%%%%%%%%%%%%%%%%%%%%%%%%%%%%%%%%%

    \subsection{Stochastic evolution equations}\label{Ch-SEE}

    In the rest of this notes, unless other stated, we shall always assume that $H$
    is a separable Hilbert space, and $A$ is an
    unbounded linear operator (with domain $D(A)$) on
    $H$, which is the infinitesimal generator of a
    $C_0$-semigroup $\{S(t)\}_{t\geq 0}$.

    Consider the following stochastic
    evolution equation (SEE for short):
    \begin{equation}\label{c1-system1}
        \left\{
        \begin{array}{ll}\ds
            dX(t) = \big(A X(t) + F(t,X(t))\big)dt +
            \wt F(t,X(t))dW(t) &\mbox{ in } (0,T],\\
            \ns\ds X(0)=X_0,
        \end{array}
        \right.
    \end{equation}
    where $X_0\in L^p_{\cF_0}(\Omega;H)$
    (for some $p\geq 1$), and $F(\cd,\cd)$
    and $\wt F(\cd,\cd)$ are two given functions from $[0,T]\times\Omega\times
    H$ to $H$.

    First, let us give the definition of {\it strong solutions} to
    \eqref{c1-system1}.

    \begin{definition}\label{def-strong-sol}
        An $H$-valued, $\mathbf{F}$-adapted, continuous stochastic process
        $X(\cd)$ is called a \index{strong solution} strong solution to the equation
        \eqref{c1-system1} if

        {\rm 1)} $X(t)\in D(A)$ for a.e. $(t,\omega)\in
        [0,T]\times\Omega$ and $AX(\cd)\in
        L^1(0,T;H)$ a.s.;

        {\rm 2)} $F(\cd,X(\cd))\in
        L^1(0,T;H)$ a.s., $ \wt F(\cd,X(\cd))\in
        L_\dbF^{2,loc}(0,T;H)$; and

        {\rm 3)} For all $t\in [0,T]$,
        $$
        X(t) = X_0 + \int_0^t \big(A X(s)
        + F(s,X(s))\big) ds  + \int_0^t  \wt F(s,X(s))dW(s),\q\as
        $$
    \end{definition}

    Generally speaking, one needs very restrictive
    conditions to guarantee the existence of
    strong solutions to \eqref{c1-system1}. Thus, people introduce two
    types of ``weak" solutions.

    \begin{definition}\label{c1-def-weak}
        An $H$-valued, $\mathbf{F}$-adapted, continuous stochastic process
        $X(\cd)$ is called a \index{weak solution} weak solution to
        \eqref{c1-system1} if $F(\cd,X(\cd))\in
        L^1(0,T;H)$ a.s., $ \wt F(\cd,X(\cd))\in
        L_\dbF^{2,loc}(0,T;H)$, and for any $t\in [0,T]$ and
        $\xi\in D(A^*)$,
        $$
        \begin{array}{ll}\ds
            \big\langle X(t), \xi \big\rangle_H \3n&\ds=
            \big\langle X_0, \xi \big\rangle_H+ \int_0^t
            \big(\big\langle X(s), A^*\xi \big\rangle_H +
            \big\langle F(s,X(s)), \xi \big\rangle_H \big)
            ds\\
            \ns&\ds \q  + \int_0^t \big\langle \wt
            F(s,X(s)), \xi
            \big\rangle_HdW(s),\qquad \as
        \end{array}
        $$
    \end{definition}

    \begin{definition}\label{c1-def-mild}
        An $H$-valued, $\mathbf{F}$-adapted, continuous stochastic process
        $X(\cd)$ is called a \index{mild solution} mild solution to
        \eqref{c1-system1} if $F(\cd,X(\cd))\in
        L^1(0,T;H)$ a.s., $ \wt F(\cd,X(\cd))\in
        L_\dbF^{2,loc}(0,T;H)$, and for any $t\in [0,T]$,
        $$
        X(t) = S(t)X_0 + \int_0^t
        S(t-s)F(s,X(s))ds + \int_0^t
        S(t-s)\wt
        F(s,X(s))dW(s),\q \mbox{a.s.}
        $$
    \end{definition}

    Clearly, a strong solution to
    \eqref{c1-system1} is a weak/mild
    solution to the same equation. The following result provides a
    sufficient condition for a mild
    solution to be a strong solution to
    \eqref{c1-system1}.

    \begin{proposition}\label{c1-th1}
        A mild solution $X(\cd)$ to
        \eqref{c1-system1} is a strong solution
        (to the same equation) if the following
        three conditions hold  for
        all $x\in H$ and $0\le s\le t\leq T$, a.s.,

        {\rm 1)} $X_0\in
        D(A)$,\;$S(t-s)F(s,x)\in
        D(A)$ and $S(t-s) \wt F(s,x)\in  D(A)$;

        {\rm 2)} $|AS(t-s)F(s,x)|_H \leq
        \a(t-s)|x|_H$ for some real-valued stochastic process $\a(\cd)\in
        L_\dbF^{1,loc}(0,T)$ a.s.; and

        {\rm 3)} $|AS(t-s) \wt F(s,x)|_{H} \leq
        \b(t-s)|x|_H$ for some real-valued stochastic process $\b(\cd)\in
        L_\dbF^{2,loc}(0,T)$ a.s.
    \end{proposition}

    The next result gives the relationship between mild
    and weak solutions to \eqref{c1-system1}.

    \begin{proposition}\label{ch-1-well-rel1}
        Any  weak solution to \eqref{c1-system1} is also
        a mild solution to the same equation and vice versa.
    \end{proposition}
    By Proposition \ref{ch-1-well-rel1}, in what follows, we will not distinguish the
    mild and
    the weak solutions to \eqref{c1-system1}.

    In the rest of this subsection, we assume that both $F(\cd, x)$ and $\wt F(\cd, x)$ are
    $\mathbf{F}$-adapted for each $x\in H$, and there exist two nonnegative (real-valued) functions $L_1(\cd)\in L^1(0, T)$
    and $L_2(\cd) \in L^2(0, T)$ such that for $\ae t\in [0,T]$ and all $y,z\in H$,
    \begin{equation}
        \left\{
        \begin{array}{ll}\ds
            |F(t,y)-F(t,z)|_{H}\leq L_1(t)|y-z|_H,\q\as\!\!,\\
            \ns\ds
            |\wt F(t,y)-\wt F(t,z)|_{H}\leq L_2(t)
            |y-z|_H,\q\as\!\!,\\
            \ns\ds F(\cd,0)\in
            L^p_\dbF(\Omega;L^1(0,T;H)),\q \wt F(\cd,0)\in
            L^p_\dbF(\Omega;L^2(0,T;H)).
        \end{array}
        \right.
    \end{equation}
    We have the
    following result:
    \begin{theorem}\label{ch-1-well-mild}
        Let $p\geq 1$. Then, there is a unique mild
        solution $X(\cd)\in C_{\dbF}([0,T];$
        $L^p(\Omega;H))$ to \eqref{c1-system1}.
        Moreover,
        \begin{equation}\label{ch-1-well-mild-eq1}
            \begin{array}{ll}\ds
                |X(\cd)|_{C_\dbF([0,T];L^p(\Omega;H))}\\
                \ns\ds \leq
                \cC\big(|X_0|_{L^p_{\cF_0}(\Omega;H)} +
                |F(\cd,0)|_{L^p_\dbF(\Omega;L^1(0,T;H))} + |\wt
                F(\cd,0)|_{L^p_\dbF(\Omega;L^2(0,T;H))}\big).
            \end{array}
        \end{equation}
    \end{theorem}

    If the semigroup $\{S(t)\}_{t\geq 0}$ is
    contractive, then one can obtain a better
    regularity for the mild solution w.r.t.
    the time variable as follows:

    \begin{theorem}\label{ch-1-well-mild1}
        If $A$ generates a contraction semigroup and
        $p\ge 1$, then \eqref{c1-system1} admits a
        unique mild solution $X(\cd)\in
        L^p_\dbF(\Omega;C([0,T];H))$. Moreover,
        \begin{equation}\label{ch-1-well-mild1-eq1}
            \begin{array}{ll}\ds
                |X(\cd)|_{L^p_\dbF(\Omega;C([0,T];H))}\\
                \ns\ds\leq \cC\big(|X_0|_{L^p_{\cF_0}(\Omega;H)}
                + |F(\cd,0)|_{L^p_\dbF(\Omega;L^1(0,T;H))}+ |\wt
                F(\cd,0)|_{L^p_\dbF(\Omega;L^2(0,T;H))}\big).
            \end{array}
        \end{equation}
    \end{theorem}

    The following result indicates the regularity of mild solutions to a class of
    stochastic evolution equations.

    \begin{theorem}\label{ch-1-well-mild2}
        Let $p\geq 1$. Assume that $A$ is a
        self-adjoint, negative definite (unbounded
        linear) operator on $H$. Then, the equation
        \eqref{c1-system1} admits a unique mild solution
        $ X(\cd)\in L^p_\dbF(\Omega;C([0,T];H))\cap
        L^p_\dbF(\Omega;L^2(0,T;D((-A)^{\frac{1}{2}}))).
        $
        Moreover,
        \begin{equation}\label{ch-1-well-mild2-eq1}
            \begin{array}{ll}\ds
                |X(\cd)|_{L^p_\dbF(\Omega;C([0,T];H))} + |X(\cd)|_{L^p_\dbF(\Omega;L^2(0,T;D((-A)^{\frac{1}{2}})))} \\
                \ns\ds \leq
                \cC\big(|X_0|_{L^p_{\cF_0}(\Omega;H)}+
                |F(\cd,0)|_{L^p_\dbF(\Omega;L^1(0,T;H))}+ |\wt
                F(\cd,0)|_{L^p_\dbF(\Omega;L^2(0,T;H))}\big).
            \end{array}
        \end{equation}
    \end{theorem}

    Usually, a mild solution does not have ``enough"
    regularities for many situations. For example, when establishing the
    pointwise identity for Carleman estimate on stochastic partial differentia equations
    of second order, we actually
    need the functions in  consideration to be twice
    differentiable in the sense of generalized derivative
    w.r.t.  the spatial  variables.
    Nevertheless, these problems can be solved by
    the following strategy:

    \begin{enumerate}
        \item Introduce some auxiliary equations with
        strong solutions such that the limit of these
        strong solutions is the mild or weak solution to
        the original equation.
        \item  Obtain the desired
        properties for these strong solutions.
        \item Utilize the density argument to establish the desired properties for
        mild/ weak solutions to
        the original equation.
    \end{enumerate}

    There are many methods to implement the above
    three steps in the setting of deterministic PDEs. Roughly
    speaking, any of these methods, which does not
    destroy the adaptedness of the solution  w.r.t
    $\mathbf{F}$, can be applied to SEEs. Here we
    only present one approach.

    Introduce a family of auxiliary equations for
    \eqref{c1-system1} as follows:
    \begin{equation}\label{c1-system3}
        \left\{
        \begin{array}{ll}\ds
            dX^\l(t) = AX^\l(t)dt + R(\l) F(t,X^\l(t))dt +
            R(\l)\wt F(t,X^\l(t))dW(t) &\mbox{ in } (0,T],\\
            \ns\ds X^\l(0)=R(\l)X_0\in D(A).
        \end{array}
        \right.
    \end{equation}
    Here   $\l\in \rho(A)$ and $R(\l)\deq\l(\l I-A)^{-1}$.

    \begin{theorem}\label{ch-2-app1}
        For each $X_0\in L^p_{\cF_0}(\Omega;H)$ with
        $p\geq 2$ and $\l\in\rho(A)$, the equation
        \eqref{c1-system3} admits a unique strong
        solution $X^\l(\cd)\in
        C_\dbF([0,T];L^p(\Omega;D(A)))$. Moreover, as
        $\l\to\infty$, the solution $X^{\l}(\cd)$
        converges to $X(\cd)$ in
        $C_\dbF([0,T];L^p(\Omega;H))$, where $X(\cd)$ is
        the mild solution to \eqref{c1-system1}.
    \end{theorem}

    At last, we give a Burkholder-Davis-Gundy type inequality, which is very useful in the study of mild solutions to SEEs.
    \begin{proposition}\label{BDGQ1}
        For any $p\geq 1$, there exists a constant
        $\cC_p>0$ such that for any $g\!\in\!
        L^p_\dbF(\Omega;L^2(0,\!T;$ $H))$ and $t\in
        [0,T]$,
        \begin{equation}\label{BDG1-eq1}
            \mE\Big|\int_0^t S(t-s)g(s)dW(s)\Big|_H^p \leq
            \cC_p\mE\(\int_0^t
            |g(s)|_{H}^2ds\)^{\frac{p}{2}}.
        \end{equation}
    \end{proposition}

    \ss

    %%%%%%%%%%%%%%%%%%%%%%%%%%%%%%%%%%%%%%%%%%%%%%%%%%%%%%%%%%

    \subsection{Backward stochastic  evolution
        equations}\label{sec-BSEE}

    %%%%%%%%%%%%%%%%%%%%%%%%%%%%%%%%%%%%%%%%%%%%%%%%%%%%%%%%%%%

    Backward stochastic differential equations and
    more generally, backward stochastic  evolution
    equations (BSEEs for short) are by-products in
    the study of stochastic control theory. These equations have independent interest and are applied
    in other places.

    Consider the following BSEE:
    \begin{equation}\label{c1-system4}
        \left\{
        \begin{array}{ll}\ds
            dY(t) = -\big(A^*Y(t) + F(t,Y(t),Z(t))\big)dt -
            Z(t) dW(t) &\mbox{
                in } [0,T),\\
            \ns\ds Y(T)=Y_T.
        \end{array}
        \right.
    \end{equation}
    Here $Y_T\in L^2_{\cF_T}(\Omega;H)$ and $F:[0,T]\times\Omega\times H\times
    H\to H$ is a  given function.

    \br
    The diffusion term ``$Y(t)dW(t)$"  in
    \eqref{c1-system4}  may be replaced by a more general form ``$\big(\wt
    F(t,y(t))+Y(t)\big)dW(t)$" (for a function $\wt
    F:[0,T]\times\Omega\times H\to H$). Clearly, the latter can be reduced to the former one by means of a simple
    transformation:
    $$
    \tilde y(\cd)=y(\cd),\qquad  \wt Y(\cd)=\wt F(\cd,y(\cd))+Y(\cd),
    $$
    under suitable assumptions on the function $\wt F$.
    \er

    Similarly to the case of SEEs, one introduces
    below notions of strong, weak and mild solutions
    to the equation \eqref{c1-system4}.

    \begin{definition}\label{def-strong-sol-b}
        An $H\times H$-valued process
        $(Y(\cd),Z(\cd))$ is called a \index{strong solution} strong solution to
        \eqref{c1-system4} if

        {\rm 1)} $Y(\cd)$ is $\mathbf{F}$-adapted and continuous, $Z(\cd)\in
        L^{2,loc}_\dbF(0,T;H)$;

        {\rm 2)}  $Y(t)\in D(A^*)$ for a.e. $(t,\omega)\in
        [0,T]\times\Omega$, and $A^*Y(\cd)\in
        L^1(0,T;H)$ and $F(\cd,Y(\cd), Z(\cd))\in
        L^1(0,T;H)$ a.s.; and

        {\rm 3)} For any $t\in [0,T]$,
        $$
        Y(t)  \ds=Y_T + \int_t^T \big(A^* Y(s)
        +F(s,Y(s),Z(s)) \big) ds   + \int_t^T
        Z(s) dW(s),\q \as
        $$
    \end{definition}

    \begin{definition}\label{def-weak-sol-b}
        An $H\times H$-valued process
        $(Y(\cd),Z(\cd))$ is called a \index{weak solution} weak solution to
        \eqref{c1-system4} if

        {\rm 1)} $Y(\cd)$ is $\mathbf{F}$-adapted and continuous, $Z(\cd)\in
        L^{2,loc}_\dbF(0,T;H)$, and  $F(\cd,Y(\cd),Z(\cd))\in
        L^1(0,T;H)$ a.s.; and

        {\rm 2}  For any $t\in [0,T]$ and $\eta\in
        D(A)$,
        $$
        \begin{array}{ll}\ds
            {\lan Y(t),\eta\ran}_H\3n&\ds= {\lan Y_T,\eta\ran}_H  +
            \int_t^T
            {\lan Y(s),A\eta\ran}_H ds + \int_t^T
            {\lan F(s,Y(s),Z(s)),\eta\ran}_H ds  \\
            \ns&\ds \q + \int_t^T {\lan Z(s),\eta\ran}_HdW(s),\quad \as
        \end{array}
        $$
    \end{definition}

    \begin{definition}\label{def-mild-sol-b}
        An $H\times H$-valued process
        $(Y(\cd),Z(\cd))$ is called a \index{weak solution} mild solution to
        \eqref{c1-system4} if

        {\rm 1)} $Y(\cd)$ is $\mathbf{F}$-adapted and continuous, $Z(\cd)\in
        L^{2,loc}_\dbF(0,T;H)$, and  $F(\cd,Y(\cd),Z(\cd))\in
        L^1(0,T;H)$ a.s.; and

        {\rm 2}  For any $t\in [0,T]$,
        $$
        \begin{array}{ll}\ds
            Y(t) =  S(T-t)^*Y_T +\int_t^T
            S(s-t)^*F(s,Y(s),Z(s))ds\\
            \ns\ds\qq\q + \int_t^T S(s-t)^* Z(s)
            dW(s),\qq \quad \as
        \end{array}
        $$
    \end{definition}

    Similarly to Proposition \ref{ch-1-well-rel1}, we have
    the following result concerning the  equivalence
    between weak and mild solutions to the equation
    \eqref{c1-system4}.

    \begin{proposition}\label{ch-1-well-brel1}
        An $H\times H$-valued process
        $(Y(\cd),Z(\cd))$ is a weak solution to
        \eqref{c1-system4} if and only if it is a mild
        solution to the same equation.
    \end{proposition}

    According to Proposition \ref{ch-1-well-brel1}, in the rest of this paper, we will not distinguish the
    mild and weak solutions to \eqref{c1-system4}.

    Clearly, if $(Y(\cd),Z(\cd))$ is a strong
    solution to \eqref{c1-system4}, then it is also
    a weak/mild solution to the same equation. On the other hand, starting from Definition \ref{def-weak-sol-b}, it is easy to show
    the following result.

    \begin{proposition}\label{ch-1-well-brel}
        A weak/mild solution $(Y(\cd),Z(\cd))$ to
        \eqref{c1-system4} is a strong solution
        (to the same equation), provided that $Y(t)\in D(A^*)$
        for a.e. $(t,\omega)\in
        [0,T]\times\Omega$, and $A^*Y(\cd)\in
        L^1(0,T;H)$ a.s.
    \end{proposition}

    In the rest of this subsection, we assume that $F(\cd, y,z)$ is $\mathbf{F}$-adapted for each $y,z\in H$, and there exist two nonnegative functions $L_1(\cd) \in L^1(0, T)$ and $L_2(\cd)\in
    L^2(0,T)$ such that, for any for $\ae t\in [0,T]$,
    \begin{equation}\label{c3-Lip}
        \left\{
        \begin{array}{ll}
            \ds |F(t,y_1,z_1)-F(t,y_2,z_2)|_{H} \leq
            L_1(t)|y_1-y_2|_H + L_2(t)|z_1-z_2|_{H},
            \\ \ns\ds \qq\qq\qq\qq\qq\qq\qq\qq\qq\forall\;
            y_1,y_2,z_1,z_2\in {H},  \q \as \\ \ns
            \ds
            F(\cd,0,0)\in L_\dbF^1(0,T;L^2(\Om;H)).
        \end{array}
        \right.
    \end{equation}

    Similarly to Theorem
    \ref{ch-1-well-mild} (but here one
    needs that the filtration $\mathbf{F}$
    is the natural one generated by
    $W(\cd)$), we have the following
    well-posedness result for the equation \eqref{c1-system4} in
    the sense of mild solution.
    \begin{theorem}\label{ch-1-well-bmild}
        Assume that $\mathbf{F}$ is the natural
        filtration generated by $W(\cd)$. Then, for every
        $Y_T\in L^2_{\cF_T}(\Omega;H)$, the
        equation \eqref{c1-system4} admits a unique mild
        solution $(Y(\cd),Z(\cd))\in
        L^2_\dbF(\Omega;C([0,T];H)) \times
        L^2_\dbF(0,T;H)$, and
        \begin{equation}\label{ch-1-well-bmild-eq1}
            |(Y,Z)|_{L^2_\dbF(\Omega;C([0,T];H))\times
                L^2_\dbF(0,T;H)}   \leq
            \cC\big(|Y_T|_{L^2_{\cF_T}(\Omega;H)} +
            |F(\cd,0,0)|_{L_\dbF^1(0,T;L^2(\Om;H))} \big).
        \end{equation}
    \end{theorem}

    Similarly to Theorem \ref{ch-1-well-mild2}, the
    following result describes the smoothing effect
    of mild solutions to a class of backward
    stochastic evolution equations.

    \begin{theorem}\label{ch-1-bwell-mild2}
        Let $\mathbf{F}$ be the natural filtration
        generated by $W(\cd)$, and $A$ be a
        self-adjoint, negative definite (unbounded
        linear) operator on $H$. Then, for any $Y_T\in
        L^2_{\cF_T}(\Omega;H)$, the equation
        \eqref{c1-system4} admits a unique mild solution
        $(Y(\cdot),Z(\cdot))\in
        \big(L^2_\dbF(\Omega;C([0,T]; H))\cap
        L^2_\dbF(0,T;D((-A)^{\frac{1}{2}}))\big)\times
        L^2_\dbF(0,T;H)$.
        Moreover,
        \begin{equation}\label{ch-1-well-bmild2-eq1}
            \begin{array}{ll}\ds
                |Y(\cd)|_{L^2_\dbF(\Omega;C([0,T];H))} + |Y(\cd)|_{L^2_\dbF(0,T;D((-A)^{\frac{1}{2}})))} + |Z(\cd)|_{L^2_\dbF(0,T;H)} \\
                \ns\ds \leq
                \cC\big(|Y_T|_{L^2_{\cF_T}(\Omega;H)} +
                |F(\cd,0,0)|_{L^1_\dbF(0,T;L^2(\Omega;H))}\big).
            \end{array}
        \end{equation}
    \end{theorem}

    Next, for each $\l\in \rho(A^*)$, we introduce the following approximate equation
    of \eqref{c1-system4}:
    \begin{equation}\label{c1-system5}
        \left\{
        \begin{array}{ll}\ds
            dY^\l(t) = -\big(A^*Y^\l(t)+ R^*(\l)
            F(t,Y(t),Z(t))\big)dt \\
            \ns\ds\qq\qq\;\,-
            R^*(\l)Z^\l(t)dW(t) &\mbox{in } [0,T),\\
            \ns\ds Y^\l(T)=R^*(\l)Y_T,
        \end{array}
        \right.
    \end{equation}
    where $R^*(\l)\deq\l(\l I-A^*)^{-1}$, $Y_T\in L^2_{\cF_T}(\Omega;H)$ and $(Y(\cd),Z(\cd))$ is the mild solution to
    \eqref{c1-system4}.
    Similarly to Theorem \ref{ch-2-app1}, we have the
    following result.
    \begin{theorem}\label{app th1}
        Assume that $\mathbf{F}$ is the natural
        filtration generated by $W(\cd)$. Then,
        for each $Y_T\in L^2_{\cF_T}(\Omega;H)$ and
        $\l\in\rho(A^*)$, the equation \eqref{c1-system5}
        admits a unique strong solution
        $(Y^\l(\cd),Z^\l(\cd))\in
        L^2_\dbF(\Omega;C([0,T];D(A^*)))\times
        L^2_\dbF(0,T;D(A^*))$. Moreover, as
        $\l\to\infty$, $(Y^\l(\cd),Z^\l(\cd))$ converges
        to $(Y(\cd),Z(\cd))$ in $
        L^2_\dbF(\Omega;C([0,T];H))\times
        L^2_\dbF(0,T;H)$.
    \end{theorem}

    Note that, in Theorems  \ref{ch-1-well-bmild}--\ref{app th1}, we
    need the filtration $\mathbf{F}$ to be the
    one generated by $W(\cd)$. For the general
    filtration, as we shall see below, we need to
    employ the stochastic transposition method
    (developed in our previous works \cite{LZ, LZ1}) to show the
    well-posedness of the equation
    \eqref{c1-system4}.

    In order to solve the BSEE
    \eqref{c1-system4} in the general
    filtration space, a fundamental idea in the stochastic transposition method (\cite{LZ, LZ1}) is to
    introduce the following test SEE:
    \begin{equation}\label{fsystem1}
        \left\{
        \begin{array}{lll}\ds
            d\f = (A\f + \psi_1)ds +  \psi_2 dW(s) &\mbox{ in }(t,T],\\
            \ns\ds \f(t)=\eta,
        \end{array}
        \right.
    \end{equation}
    where $t\in[0,T)$, $\psi_1\in
    L^1_{\dbF}(t,T;L^2(\Om;H))$, $\psi_2\in
    L^2_{\dbF}(t,T;H)$ and $\eta\in
    L^{2}_{\cF_t}(\Om;H)$.  By Theorem
    \ref{ch-1-well-mild},  the equation
    \eqref{fsystem1} admits a unique solution $\f\in
    C_\dbF([t,T];L^2(\Om;H))$ such that
    \begin{equation}\label{3.25-eq1}
        |\f|_{C_\dbF([t,T];L^2(\Om;H))}\leq
        \cC\left|\big(\psi_1(\cdot),
        \psi_2(\cdot),\eta\big)\right|_{
            L^1_{\dbF}(t,T;L^2(\Om;H))\times
            L^2_{\dbF}(t,T;H)\times L^2_{\cF_t}(\Om;H)}.
    \end{equation}
    \begin{definition}\label{definition1}
        We call $(Y(\cdot), Z(\cdot)) \in
        D_{\dbF}([0,T];L^{2}(\Om;H)) \times
        L^2_{\dbF}(0,T;H)$  a {\it transposition solution} to
        \eqref{c1-system4} if for any $t\in [0,T]$,
        $\psi_1(\cdot)\in L^1_{\dbF}(t,T;L^2(\Om;H))$,
        $\psi_2(\cdot)\in L^2_{\dbF}(t,T; H)$, $\eta\in
        L^2_{\cF_t}(\Om;H)$ and the corresponding
        solution $\f\in C_{\dbF}([t,T];L^2(\Om;H))$ to
        \eqref{fsystem1}, it holds that
        \begin{equation}\label{eq def solzz}
            \begin{array}{ll}\ds
                \dbE \big\langle \f(T),Y_T\big\rangle_{H}
                - \dbE\int_t^T \big\langle \f(s),F(s,Y(s),Z(s) )\big\rangle_Hds\\
                \ns\ds = \dbE \big\langle\eta,Y(t)\big\rangle_H
                + \dbE\int_t^T \big\langle
                \psi_1(s),Y(s)\big\rangle_H ds + \dbE\int_t^T
                \big\langle \psi_2(s),Z(s)\big\rangle_{H} ds.
            \end{array}
        \end{equation}
    \end{definition}
    \begin{remark}
        If \eqref{c1-system4} admits a  weak solution
        $(Y(\cdot), Z(\cdot)) \in
        C_{\dbF}([0,T];L^{2}(\Om;H)) \times
        L^2_{\dbF}(0,T;$ $H)$, then by It\^o's formula,
        people obtain the variational equality \eqref{eq def solzz}.
        This is the very reason for us to introduce Definition
        \ref{definition1}.
    \end{remark}

    We have the following result on the
    well-posedness of the equation \eqref{c1-system4} (See \cite[Theorem 3.1]{LZ1}).
    \begin{theorem}\label{vw-th1}
        The equation \eqref{c1-system4} admits a unique
        transposition solution $(Y(\cdot), Z(\cdot))\in D_{\dbF}([0,T];$ $L^{2}(\Om;H))
        \times L^2_{\dbF}(0,T;H)$. Furthermore,
        \begin{equation}\label{vw-th1-eq1}
            \begin{array}{ll}\ds
                |(Y(\cdot), Z(\cdot))|_{
                    D_{\dbF}([0,T];L^{2}(\Om;H)) \times
                    L^2_{\dbF}(0,T;L^2(\Om;H))}\\
                \ns\ds \leq \cC\big(|Y_T|_{
                    L^2_{\cF_T}(\Om;H)}+|F(\cd,0,0)|_{
                    L^1_{\dbF}(0,T;L^2(\Om;H))}\big).
            \end{array}
        \end{equation}
    \end{theorem}
    Unlike the case that $\mathbf{F}$ is the natural filtration generated by $W(\cd)$, the proof of Theorem \ref{vw-th1} is based on a new stochastic Riesz-Type Representation Theorem (proved in \cite{LYZ, LYZ1}).

    \br
    The stochastic transposition method works well not only for the vector-valued (i.e., $H$-valued) BSEE  \eqref{c1-system4}, but also (and more importantly) for the difficult operator-valued (i.e., $\cL(H)$-valued) BSEEs (See the equations (\ref{5.5-eq6}) and (\ref{op-bsystem3}),  also \cite{LZ1, LZ2, LZ3, LZ5} and particularly \cite{LZ3.1} for more details).
    \er

    %%%%%%%%%%%%%%%%%%%%%%%%%%%%%%%%%%%%%%%%%%%%%%%%%%%

    \section{Control systems governed by stochastic partial differential equations}\label{Ch5-sec1}

    %%%%%%%%%%%%%%%%%%%%%%%%%%%%%%%%%%%%%%%%%%%%%%%%%%%%%%%%%

    It is well-known that Control Theory was founded by N. Wiener in 1948. After that, owing to the great effort of numerous mathematicians and engineers, this theory was greatly extended to various different setting, and widely used in sciences, technologies, engineering and economics, particularly in Artificial Intelligence in recent years.
    Usually, in terms of the
    state-space technique (\cite{Kalman}), people describe the control system under consideration as a suitable state equation.

    Control theory for finite dimensional systems is now relatively mature.
    There exist a huge list of publications on control theory for deterministic distributed parameter
    systems (typically governed by partial differential equations) though it is still quite active; while the same can be
    said for control theory for stochastic (ordinary) differential equations (i.e., stochastic differential equations in finite dimensions).
    By contrast, control theory for stochastic distributed parameter systems (described by
    stochastic differential equations in infinite dimensions, typically by stochastic partial differential equations), is still at its very beginning stage, though it was ``born" in almost the same time as that for deterministic distributed parameter control systems.

    Because of their inherent complexities, many control systems in reality exhibit very complicated
    dynamics, including substantial model uncertainty, actuator and state constraints, and
    high dimensionality (usually infinite). These systems are often best described by stochastic partial differential equations or even more complicated stochastic equations in infinite dimensions (e.g., \cite{Murray}).

    Generally speaking, any PDE can be
    regarded as a  SPDE provided that at least one of its  coefficients,
    forcing terms, initial and boundary conditions
    is random. The terminology SPDE is a little misused, which may
    mean different types of equations in different
    places. The analysis of equations with random
    coefficients differs dramatically from that of equations with stochastic noises.
    In this
    notes, we focus on the latter.

    The study of SPDEs is mainly motivated by
    two aspects:
    \begin{itemize}
        \item  One is due to the rapid development of stochastic analysis recently.  The topic of SPDEs is an interdisciplinary area involved both stochastic processes (random
        fields) and PDEs, which has quickly become a discipline  itself. In the
        last two decades, SPDEs have been one of the most
        active areas in stochastic analysis.

        \item Another is the requirement from some
        branches of sciences.   In many phenomena in
        physics, biology, economics and control theory (including filter theory in particular), stochastic
        effects play a central role. Thus, stochastic
        corrections to the deterministic models are indispensably
        needed. These backgrounds have been influencing
        the development of SPDEs remarkably.
    \end{itemize}

    In this section, we present two typical models described by SPDEs, in which some control actions may be introduced.
    The readers are referred to \cite{Prato, Holden} and so on for more SPDE systems.

    \begin{example}\label{ch-4-ex1}

        {\bf Stochastic parabolic equations}

        \ss

        The following equation was introduced to
        describe the evolution of the density of a
        population (e.g., \cite{Dawson}):
        \begin{equation}\label{ch-4-ex1-eq1}
            \left\{
            \begin{array}{ll}\ds
                dy = \kappa\pa_{xx} y dt + \a\sqrt{y} dW(t)
                &\mbox{ in } (0,T)\times (0,L),\\
                \ns\ds y_x=0 &\mbox{ on } (0,T)\times \{0,L\},
                \\
                \ns\ds y(0) = y_0 &\mbox{ in } (0,L),
            \end{array}
            \right.
        \end{equation}
        where $\kappa>0$, $\a>0$ and $L>0$ are given constants and $y_0\in L^2(0,L)$. The derivation of \eqref{ch-4-ex1-eq1} is as
        follows.

        Suppose that a bacteria population is distributed in
        the interval $[0,L]$. Denote by $y(t,x)$ the density of
        this population at time $t\in[0,T]$ and position $x\in[0,L]$. If there is no large-scale
        migration, the variation of the density is
        governed by
        \bel{200306e11}
        dy(t,x)=\kappa y_{xx}(t,x)dt + d\xi(t,x,y).
        \ee
        In \eqref{200306e11}, $\kappa y_{xx}(t,x)$ describes the population's diffusion from the high density place to the low
        one, while $\xi(t,x,y)$ is a random perturbation caused by lots of small independent
        random disturbances. Suppose that the random perturbation $\xi(t,x,y)$ at
        time $t$ and position $x$  can be approximated
        by a Gaussian stochastic process whose variance ${\rm
            Var }\;\xi$ is monotone w.r.t. $y$. In the study
        of bacteria, $L$ is very small and we may assume that $\xi$ is
        independent of $x$. Further, when $y$ is small,
        the variance of $\xi$ can be approximated by
        $\a^2 y$, where $\a^2$ is the derivative of
        ${\rm Var }\;\xi$ at $y=0$. Under these assumptions, it follows that
        \bel{200306e12}
        d\xi(t,x,y)=\a\sqrt{y}dW(t).
        \ee
        Combining \eqref{200306e11} and \eqref{200306e12}, we arrive at the first equation of \eqref{ch-4-ex1-eq1}. If
        there is no bacteria entering into or leaving $[0,L]$ from
        its boundary, we have
        $y_x(t,0)=y_x(t,1)=0$. As a result, we
        obtain the boundary condition of \eqref{ch-4-ex1-eq1}.

        In order to change the population's density, one can
        put in or draw out some species. Under such actions,
        the equation \eqref{ch-4-ex1-eq1} becomes a
        controlled stochastic parabolic equation:
        \begin{equation}\label{ch-4-ex1-eq1-con}
            \left\{
            \begin{array}{ll}\ds
                dy = (\kappa\pa_{xx} y + u) dt + (\a\sqrt{y}+v)
                dW(t)
                &\mbox{ in } (0,T)\times (0,L),\\
                \ns\ds y_x=0 &\mbox{ on } (0,T)\times \{0,L\},
                \\
                \ns\ds y(0) = y_0 &\mbox{ in } (0,L),
            \end{array}
            \right.
        \end{equation}
        where $u$ and $v$ are the ways of putting in or
        drawing out species.

    \end{example}
    \begin{example}\label{ch-4-ex2}
        {\bf Stochastic wave equations}

        \ss

        To study the vibration of thin string/membrane perturbed by the random force, people
        introduced the following stochastic wave
        equation (e.g. \cite{Funaki}):
        \begin{equation}\label{ch-4-ex2-eq1}
            \left\{
            \begin{array}{ll}\ds
                dy_t =  \D y dt  + \a y  dW(t)
                &\mbox{ in } (0,T)\times G,\\
                \ns\ds y=0 &\mbox{ on } (0,T)\times \G
                ,
                \\
                \ns\ds y(0) = y_0, \; y_t(0)=y_1 &\mbox{ in } G.
            \end{array}
            \right.
        \end{equation}
        Here $G$  is given as in Subsection \ref{sec-mpr-6.1}, $\a(\cd)$ is a
        suitable function, while $(y_0, y_1)$ is an initial datum.
        Let us recall below a derivation
        of \eqref{ch-4-ex2-eq1} for $G=(0,L)$ (for some
        $L>0$) by studying
        the motion of a DNA molecule in a fluid
        (e.g., (\cite{Funaki, MMK})).

        Compared with its length, the diameter of such a DNA
        molecule is very small, and hence, it can be
        viewed as a thin and long elastic string. One can
        describe its position by using an
        $\dbR^3$-valued function $y=(y_1,y_2,y_3)$
        defined on $[0,L]\times [0,+\infty)$. Usually, a DNA molecule  floats in a fluid.
        Thus, it is always struck  by the fluid
        molecules, just as a particle of pollen floating
        in a fluid.

        For simplicity, we assume that the mass of this string per unit
        length is equal to $1$. Then, the acceleration at position $x\in [0,L]$ along
        the string at time $t\in [0,+\infty)$ is $y_{tt}(t,x)$. There
        are mainly four kinds of forces acting on the
        string: the elastic force $F_1(t,x)$, the
        friction $F_2(t,x)$ due to viscosity of the
        fluid, the impact force $F_3(t,x)$ from the
        flowing of the fluid and the random impulse
        $F_4(t,x)$ from the impacts of the fluid
        molecules. By Newton's Second Law, it follows that
        \begin{equation}\label{1.17-eq1}
            y_{tt}(t,x) = F_1(t,x) + F_2(t,x) + F_3(t,x) +
            F_4(t,x).
        \end{equation}
        Similar to the derivation of deterministic
        wave equations, the elastic force
        $F_1(t,x)=y_{xx}(t,x)$. The fiction depends on
        the property of the fluid and the velocity and
        shape of the string. When $y$, $y_x$ and
        $y_t$ are small, $F_2(t,x)$ approximately
        depends on them linearly, i.e., $F_2 =a_1
        y_x+a_2y_t+a_3y$ for some suitable functions $a_1$, $a_2$
        and $a_3$. From the classical theory of
        Statistical Mechanics (e.g., \cite[Chapter
        VI]{Stowe}), the random impulse
        $F_4(t,x)$ at time $t$ and position $x$  can be
        approximated by a Gaussian stochastic process with a
        covariance $k(\cd,y)$,
        which also depends on the property of the fluid, and therefore we may assume that
        $$
        \begin{array}{ll}\ds
            F_4(t,x)  =\int_0^tk(x,y(s,x))dW(s).
        \end{array}
        $$
        Thus, the equation \eqref{1.17-eq1} can be
        rewritten as
        \begin{equation}\label{xx1.17-eq2}
            dy_{t}(t,x) = (y_{xx} + a_1
            y_x+a_2y_t+a_3y)dt +
            F_3(t,x)dt + k(x,y(t,x))dW(t).
        \end{equation}
        When $y$ is small, we may assume that $k(\cd, y)$ is
        linear w.r.t. $y$. In this case, $
        k(x,y(t,x))=a_4 (x)y(t,x)$ for a suitable  function
        $a_4(\cd)$. More generally, one may assume that $a_4$ depends on $t$ (and even on the sample point $\omega\in \Omega$). Thus,  \eqref{xx1.17-eq2} is reduced to the following:
        \begin{equation}\label{1.17-eq2}
            \ds dy_t(t,x) = (y_{xx} + a_1
            y_x+a_2y_t+a_3y)dt +
            F_3(t,x)dt+a_4(t,x)y(t,x)dW(t).
        \end{equation}

        Many biological behaviors are related to the
        motions of DNA molecules. Hence, there is a
        strong motivation for controlling these motions.
        In \eqref{1.17-eq2},  $F_3(\cd,\cd)$ in the
        drift term and some part of the diffusion term
        can be designed as controls acting on the fluid.
        Furthermore, one can introduce some forces on
        the endpoints of the molecule. In this way, we
        arrive at the following controlled stochastic
        wave equation in one space dimension:
        \begin{equation}\label{ch-4-ex2-eq1-con}
            \left\{
            \begin{array}{ll}\ds
                dy_t = (y_{xx}+a_1y_x+a_2y_t+a_3y + u) dt  +
                (a_4y+v) dW(t)
                &\mbox{ in } (0,T)\times (0,L),\\
                \ns\ds y =f_1 &\mbox{ on } (0,T)\times \{0\},\\
                \ns\ds y =f_2 &\mbox{ on } (0,T)\times \{L\},
                \\
                \ns\ds y(0) = y_0,\q y_t(0)=y_1 &\mbox{ in }
                (0,L).
            \end{array}
            \right.
        \end{equation}
        In \eqref{ch-4-ex2-eq1-con}, $(f_1,f_2,u,v)$ are
        controls which belonging to some suitable set, $a_4$
        is a suitable function.

        A natural question is as follows:
        $$
        \ba{ll}\ds
        \hbox{\it Can we drive the DNA
            molecule (governed by \eqref{ch-4-ex2-eq1-con}) from any }(y_0,y_1)\\
        \ds\hbox{\it to any target }(z_0,z_1) \hbox{ \it at time } T \hbox{ \it by choosing suitable controls }
        (f_1,f_2,u,v)\hbox{\it ?}
        \ea
        $$
        Clearly, this is an exact controllability problem for the stochastic wave equation \eqref{ch-4-ex2-eq1-con}.
        Although there are
        four controls in the system
        \eqref{ch-4-ex2-eq1-con}, as we shall show in Section \ref{s-h}, the answer to
        the above question is negative for controls
        belonging to some reasonable set.

        On the other hand, though the system
        \eqref{ch-4-ex2-eq1-con} is not exactly
        controllable, it may well happen that for some
        $(y_0,y_1)$ and $(z_0,z_1)$, there exist some
        controls $(f_1,f_2,u,v)$ (in some Hilbert space
        $H$) so that the corresponding solutions to
        \eqref{ch-4-ex2-eq1-con} verifying $y(T)=z_0$
        and $y_t(T)=z_1$, a.s. Usually, such sort of
        controls are not unique, and therefore, people
        may hope to find an optimal one, among these
        controls, which minimize the following
        functional
        $$|(f_1,f_2,u,v)|_H^2.
        $$
        This is typically an optimal control problem for
        the stochastic wave equation
        \eqref{ch-4-ex2-eq1-con}.
    \end{example}

    Clearly,  it is of fundamental importance to
    ``understand" and then to change the behavior of
    the system under consideration by means of
    suitable ``control" actions in an ``optimal"
    way. As we have seen in the above examples
    (particularly in Example \ref{ch-4-ex2}), this
    leads to the formation of {\it controllability}
    and {\it optimal control} problems, actually two
    fundamental  problems that we will be concerned
    with in the rest of this notes.

    %%%%%%%%%%%%%%%%%%%%%%%%%%%%%%%%%%%%%%%%%%%%%%%%%%%%%%

    \section{Controllability and observability for
        stochastic evolution equations: some simple results}\label{ch6-sec3}

    %%%%%%%%%%%%%%%%%%%%%%%%%%%%%%%%%%%%%%%%%%%%%%%%%%%%%%

    In this section, we shall present some
    results on the controllability and
    observability of SEEs. As the title
    suggests, these results are simple and
    their proofs are easy. However, some of them may not be obvious for beginners. On the other hand, some results in this respect reveal new phenomena for
    controllability problems in the stochastic setting. Unless otherwise stated, the content of this section is taken from \cite[Chapter 7]{LZ3.1}.

    Throughout this section, $U$ is a given Hilbert space.

    %%%%%%%%%%%%%%%%%%%%%%%%%%%%%%%%%%%%%%%%%%%%%%%%%%%%%%

    \subsection{Lack of exact controllability}

    %%%%%%%%%%%%%%%%%%%%%%%%%%%%%%%%%%%%%%%%%%%%%%%%%%%%%%

    We first consider the following control system evolved on $H$:
    \begin{equation}\label{6.18-eq4}
        \left\{
        \begin{array}{ll}\ds
            dX(t)=\big(A+F(t)\big)X(t)dt  +
            Bu(t) dt + K(t)X(t) dW(t) \mbox{ in } (0,T],\\
            \ns\ds X(0)=X_0.
        \end{array}
        \right.
    \end{equation}
    In \eqref{6.18-eq4}, $X_0\in H$, $F,K\in
    L^\infty_\dbF(0,T;\cL(H))$,  $B \in
    \cL(U;H)$ and $u\in L^2_\dbF(0,T;U)$.

    We begin with the notion of exact controllability of \eqref{6.18-eq4}.

    \begin{definition}\label{def-exc-ab}
        The system \eqref{6.18-eq4} is called {\it exactly
            controllable  at time $T$ } if for any $X_0\in H$ and
        $X_1\in L^2_{\cF_T}(\Om;H)$, there is a
        control $u\in L^2_\dbF(0,T;U)$ such
        that the corresponding solution
        $X(\cd)$ to \eqref{6.18-eq4} fulfills that $X(T)=X_1$, a.s.
    \end{definition}

    Define a function $\eta(\cd)$ on $[0,T]$ by
    \begin{equation}\label{6.18-eq13}
        \eta(t) = \left\{
        \begin{array}{ll}
            \ds 1, &\mbox{ for } t\in \ds\Big[\Big(1-\frac{1}{2^{2j}}\Big)T, \Big(1-\frac{1}{2^{2j+1}}\Big)T \Big), \q j=0,1,2,\cds,\\
            \ns\ds  -1,  &\mbox{ otherwise}.
        \end{array}
        \right.
    \end{equation}
    Let us recall the following known  result (\cite[Lemma 2.1]{Peng94}).
    \begin{lemma}\label{6.18-cor1}
        Let $\ds\xi = \int_0^T \eta(t)dW(t)$. It is
        impossible to find $(\varrho_1,\varrho_2)\in
        L^2_\dbF(0,T)\times C_\dbF([0,T];L^2(\Omega))$ and
        $x_0\in \dbR$ such that
        \begin{equation}
            \xi = x_0 + \int_0^T \varrho_1(t)dt + \int_0^T
            \varrho_2(t)dW(t).
        \end{equation}
    \end{lemma}

    We have the following negative controllability result for the
    system \eqref{6.18-eq4}.
    \begin{theorem}\label{6.18-th3}
        Let $K(\cd)\in
        C_\dbF([0,T];L^\infty(\Om;\cL(H)))$. Then the
        system \eqref{6.18-eq4} is not exactly
        controllable  at any time $T>0$.
    \end{theorem}

    {\it Proof}.
        Without loss of generality, we assume that
        $F(\cd)=0$. Let $\f\in D(A^*)$ be such that $|\f|_H=1$.

        We use the contradiction argument. Suppose that \eqref{6.18-eq4} was exactly
        controllable  at some time $T>0$. Then, for $X_0=0$ and $X_1=\(\ds\int_0^T
        \eta(t)dW(t)\)\f$ with $\eta(\cd)$ being given by
        \eqref{6.18-eq13}, one could find a
        control $u\in L^2_\dbF(0,T;U)$ such
        that the corresponding solution
        $X(\cd)$ to \eqref{6.18-eq4} fulfills that $X(T)=X_1$, a.s.
        Hence, by the definition of weak solutions to \eqref{6.18-eq4}, it follows that
        \begin{equation}\label{6.18-eq14}
            \begin{array}{ll}\ds
                \int_0^T \eta(t)dW(t)\3n&\ds = \int_0^T
                \big(\langle X(t), A^*\f\rangle_H  +  \langle B
                u(t), \f\rangle_H \big)dt \\
                \ns&\ds\q  +  \int_0^T \langle K(t)X(t),
                \f\rangle_H dW(t).
            \end{array}
        \end{equation}
        Since $X(\cd)\in C_\dbF([0,T];L^2(\Omega;H))$, we deduce that $
        \langle X(\cd), A^*\f\rangle_H + \langle B
        u(\cd), \f\rangle_H \in L^2_\dbF(0,T) $ and
        $\langle K(\cd)X(\cd), \f\rangle_H \in
        C_\dbF([0,T];L^2(\Omega))$. Clearly,
        \eqref{6.18-eq14} contradicts Lemma
        \ref{6.18-cor1}.
    \endpf

    Theorem \ref{6.18-th3} concludes that an SEE is
    not exactly controllable when the control, which
    is $L^2$ w.r.t. the time variable, is only
    acted in the drift term. Further results for
    controllability problems of such stochastic control systems
    can be found in \cite{DL, Peng94, YZ}.

    %%%%%%%%%%%%%%%%%%%%%%%%%%%%%%%%%%%%%%%%%%%%%%%%%%%

    \subsection{Null and approximate controllability of stochastic evolution equations}
    \label{Ch-abc-na}

    %%%%%%%%%%%%%%%%%%%%%%%%%%%%%%%%%%%%%%%%%%%%%%%%%%%

    To begin with, let us introduce the following concepts.
    \begin{definition}\label{def3.4}
        The system \eqref{6.18-eq4} is called {\it null
            controllable  at time $T$} if for any $X_0\in H$, there is
        a control $u\in L^2_\dbF(0,T;U)$ such
        that the corresponding solution
        $X(\cd)$ to \eqref{6.18-eq4} fulfills that $X(T)=0$, a.s.
    \end{definition}
    \begin{definition}
        The system \eqref{6.18-eq4} is called {\it approximately
            controllable at time $T$} if for any $X_0\in H$, $X_1\in L^2_{\cF_T}(\Om;H)$ and $\e>0$, there is
        a control $u\in L^2_\dbF(0,T;U)$ such
        that the corresponding solution
        $X(\cd)$ to \eqref{6.18-eq4} fulfills that $|X(T)-X_1|_{L^2_{\cF_T}(\Om;H)}<\e$.
    \end{definition}

    It is quite easy to study the null controllability of
    \eqref{6.18-eq4} when the operator $F(\cd)$ is ``deterministic", i.e., $F(\cd) \in
    L^\infty(0,T;\cL(H))$. In this case, let us
    introduce the following deterministic control
    system:
    \begin{equation}\label{6.18-eq5}
        \left\{
        \begin{array}{ll}\ds
            \frac{d\widehat
                X(t)}{dt}=\big(A+F(t)\big)\widehat X(t) +
            B\hat u(t) &\mbox{ in } (0,T],\\
            \ns\ds \widehat X(0)=X_0,
        \end{array}
        \right.
    \end{equation}
    where $X_0\in H$ and $\hat u\in
    L^2(0,T;U)$. Similarly to Definition \ref{def3.4},
    the system \eqref{6.18-eq5} is called {\it null
        controllable  at time $T$} if for any $X_0\in H$, there is
    a control $\hat u\in L^2(0,T;U)$ such
    that the corresponding solution to \eqref{6.18-eq5} satisfies that $\widehat X(T)=0$. We have the following result:

    \begin{theorem}\label{6.18-th1}
        Suppose that $F (\cd) \in
        L^\infty(0,T;\cL(H))$ and $K(\cd) \in
        L^\infty_\dbF(0,T)$, Then, the system \eqref{6.18-eq4} is null
        controllable  at time $T$ if and only if so is the
        system \eqref{6.18-eq5}.
    \end{theorem}
    {\it Proof}.  The ``only if" part. For
        arbitrary $X_0\in H$, let $u(\cd)\in
        L^2_\dbF(0,T;U)$ which steers the
        solution of \eqref{6.18-eq4} to $0$ at
        time $T$. Then, clearly, $\mE X$ is the
        solution to \eqref{6.18-eq5} with the
        control $\hat u=\mE u\in
        L^2(0,T;U)$. Since $X(T)=0$,  a.s.,
        we have that $\mE X(T)=0$, which
        deduces that \eqref{6.18-eq5} is null
        controllable at
        time $T$.

        \vspace{0.15cm}

        The ``if" part. For a given $X_0\in H$, let
        $\hat u(\cd)\in L^2(0,T;U)$ which steers the
        solution of \eqref{6.18-eq5} to $0$ at time
        $T$. Let $$X(t)=e^{\int_0^t
            K(s)dW(s)-\frac{1}{2}\int_0^t K(s)^2ds}\widehat
        X(t).$$ Then,  $X(\cd)$ is the unique solution
        of \eqref{6.18-eq4} with the control
        $$u(\cd)=e^{\int_0^\cd
            K(s)dW(s)-\frac{1}{2}\int_0^\cd K(s)^2ds}\hat
        u(\cd).$$ Since $\widehat X(T)=0$, we have that
        $$X(T)=e^{\int_0^T K(s)dW(s)-\frac{1}{2}\int_0^T
            K(s)^2ds}\widehat X(T)=0.$$ Further,
        $$
        \begin{array}{ll}\ds
            \mE\int_0^T|u(t)|_U^2 dt\3n&\ds = \int_0^T|\hat
            u(t)|_U^2\mE e^{2\int_0^t K(s)dW(s)- \int_0^t
                K(s)^2ds} dt\\
            \ns&\ds \leq \cC\int_0^T|\hat
            u(t)|_U^2dt<\infty.
        \end{array}
        $$
        This completes the proof of Theorem \ref{6.18-th1}.
    \endpf

    \begin{remark}\label{rmk-12-3-1}
        The technique in the proof of Theorem \ref{6.18-th1} cannot be applied to get the
        exact controllability or approximate
        controllability of the system \eqref{6.18-eq4}.
        Indeed, let us assume that \eqref{6.18-eq5} is
        exactly controllable, i.e., for any $X_0, X_1\in H$, there exists a control $\hat u$ such that the solution to \eqref{6.18-eq5} fulfills that $\widehat X(T)= X_1$. Put
        $$
        \begin{array}{ll}\ds
            \cA_T \3n&\ds\deq \big\{X(T)\big|\, X(\cd)
            \mbox{ solves \eqref{6.18-eq4} for some $X_0\in
                H$ and
            }\\
            \ns&\ds \qq\qq\; u(\cd)=e^{\int_0^\cd
                K(s)dW(s)-\frac{1}{2}\int_0^\cd K(s)^2ds}\hat
            u(\cd) \mbox{ for }\hat u(\cd)\in
            L^2(0,T;U)\big\}.
        \end{array}
        $$
        Then, from the exact controllability of \eqref{6.18-eq5}, we see that
        $$
        \cA_T=\big\{e^{\int_0^T
            K(s)dW(s)-\frac{1}{2}\int_0^T
            K(s)^2ds}X_1\big|\, X_1\in H\big\}.
        $$
        Clearly, $\cA_T\neq L^2_{\cF_T}(\Omega;H)$ and
        $\cA_T$ is not dense in $L^2_{\cF_T}(\Omega;H)$. Consequently, the system \eqref{6.18-eq4} is neither exactly controllable nor approximately controllable.
    \end{remark}
    \begin{remark}
        When $A$ generates a $C_0$-group on $H$, the null controllability of \eqref{6.18-eq5} implies the exact/approximate controllability of \eqref{6.18-eq5}. From  Remark \ref{rmk-12-3-1}, it is easy to see that this is incorrect for \eqref{6.18-eq4}, which is a new phenomenon in the stochastic setting.
    \end{remark}

    For deterministic control systems, if we assume
    that the control operator $B$ is
    invertible, then it is easy to get
    controllability results. The following theorem shows that the same
    holds for  stochastic control systems.
    \begin{theorem}\label{6.18-th2}
        If $B \in \cL(U;H)$ is invertible, then the
        system \eqref{6.18-eq4} is both null
        and approximately controllable at any time $T>0$.
    \end{theorem}
    {\it Proof}. Consider the following BSEE:
        \begin{equation}\label{6.18-eq6}
            \left\{
            \begin{array}{ll}\ds
                d Y(t)=-\big(A^* Y(t)+ F(t) Y(t)+ K(t) Z(t)\big)dt + Z(t) dW(t)  &\mbox{ in } [0,T),\\
                \ns\ds Y(T)= Y_T,
            \end{array}
            \right.
        \end{equation}
        where $Y_T\in L^2_{\cF_T}(\Omega;H)$. By \cite[Theorem 7.17]{LZ3.1}, to prove
        that \eqref{6.18-eq4} is null controllable, it suffices to show that
        \begin{equation}\label{6.18-eq7}
            |Y(0)|_H^2 \leq \cC\int_0^T \mE|B^*Y(t)|_U^2dt.
        \end{equation}
        By Theorem \ref{ch-1-well-bmild},  there is
        a constant $\cC>0$ such that
        \begin{equation}\label{6.18-eq8}
            |Y(0)|_H^2 \leq \cC \mE|Y(t)|_H^2,\q\forall\;
            t\in [0,T].
        \end{equation}
        Since $B$ is invertible, there is
        a constant $\cC>0$ such that
        $$
        |y|_H\leq \cC|B^* y|_U,\q\forall\; y\in
        H.
        $$
        This, together with \eqref{6.18-eq8}, implies
        the inequality \eqref{6.18-eq7}.

        To prove that \eqref{6.18-eq4} is approximately
        controllable, it suffices to show that
        \begin{equation}\label{6.18-eq9}
            Y_T = 0, \mbox{ provided that } B^*Y(\cd)=0
            \mbox{ in } L^2_\dbF(0,T;U).
        \end{equation}
        Since $B$ is invertible, we have that $Y(\cd)=0$
        in $L^2_\dbF(0,T;H)$. Noting that $Y(\cd)\in
        C_\dbF([0,T];L^2(\Omega;H))$, we get that
        $Y_T=0$.
    \endpf

    %%%%%%%%%%%%%%%%%%%%%%%%%%%%%%%%%%%%%%%%%%%%%%%%%%%%%%

    \subsection{Observability estimate of stochastic evolution equations}

    %%%%%%%%%%%%%%%%%%%%%%%%%%%%%%%%%%%%%%%%%%%%%%%%%%%%%%

    Consider the following SEE:
    \begin{equation}\label{6.18-eq17}
        \left\{
        \begin{array}{ll}\ds
            dX(t)=\big(A+F(t)\big)X(t)dt   + G(t)X(t) dW(t) & \mbox{ in } (0,T],\\
            \ns\ds X(0)=X_0,
        \end{array}
        \right.
    \end{equation}
    where $X_0\in H$, $F\in L^\infty(0,T;\cL(H))$
    and $G\in L^\infty_\dbF(0,T;\cL(H))$.
    Further, we introduce the following
    deterministic evolution equation:
    \begin{equation}\label{6.18-eq18}
        \left\{
        \begin{array}{ll}\ds
            \frac{d\widehat X(t)}{dt}=\big(A+F(t)\big)\widehat X(t)& \mbox{ in } (0,T],\\
            \ns\ds X(0)=X_0.
        \end{array}
        \right.
    \end{equation}
    Let $\wt U$ be another Hilbert space and $\mathfrak{C}
    \in \cL(H;\wt U)$.
    We say that the equation
    \eqref{6.18-eq17} (\resp \eqref{6.18-eq18}) is
    {\it continuously initially observable on $[0,T]$} for
    the observation operator $\mathfrak{C}$ if
    \begin{equation}\label{4.16-eq1}
        |X_0|_H^2 \leq \cC\int_0^T \mE|\mathfrak{C}
        X(t)|_{\wt U}^2,
    \end{equation}
    (\resp
    \begin{equation}\label{4.16-eq2}
        |X_0|_H^2 \leq \cC\int_0^T |\mathfrak{C}
        \widehat X(t)|_{\wt U}^2.\)
    \end{equation}

    The following
    result reveals the relationship between the observability of \eqref{6.18-eq17} and \eqref{6.18-eq18}.

    \begin{theorem}\label{6.18-th5}
        The equation \eqref{6.18-eq17} is continuously
        initially observable on $[0,T]$ for the
        observation operator $\mathfrak{C}$, provided
        that so is the equation \eqref{6.18-eq18}.
    \end{theorem}
    {\it Proof}. Let $X(\cd)$ be a solution to
        \eqref{6.18-eq17}. Then, $\mE X(\cd)$ is a
        solution to \eqref{6.18-eq18}. If
        \eqref{6.18-eq18} is continuously initially
        observable on $[0,T]$, then, by (\ref{4.16-eq2}), we have
        $$
        |X_0|_H^2 \leq \cC\int_0^T |\mathfrak{C}\mE
        X(t)|_{\wt U}^2.
        $$
        This, together with H\"older's inequality,
        implies \eqref{4.16-eq1}.
    \endpf

    \begin{remark}
        Clearly, if $F$ is stochastic, then one
        cannot use the above argument to solve the observability
        problems for SEEs. For establishing the observability for SPDEs,
        a powerful method is the global Carleman
        estimate (e.g., \cite{Fu, Gao, Liu-Yu1, Luqi3,
            Luqi7, Luqi1, WuCW, Yuan1, Yuan,
            Zhangxu1, Zhangxu2})
    \end{remark}
    \begin{remark}
        In this section, for simplicity, we assume that
        $B\in\cL(U;H)$ and $\mathfrak{C} \in \cL(H;\wt
        U)$. Some results for unbounded $B$ and
        $\mathfrak{C}$ can be found in \cite{Luqi8, LZ3.1}.
    \end{remark}
    %
    %%%%%%%%%%%%%%%%%%%%%%%%%%%%%%%%%%%%%%%%%%%%%%%%%%%%

    \section{Controllability of stochastic partial differential equations I: \linebreak stochastic
        transport equations}\label{s-t}

    %%%%%%%%%%%%%%%%%%%%%%%%%%%%%%%%%%%%%%%%%%%%%%%%%%%%

    Stochastic
    transport equation   can be
    regarded as the simplest SPDE.
    In this section, we consider the exact
    controllability for this equation. By a
    duality argument, the controllability problem is
    reduced to a suitable observability
    estimate for backward stochastic
    transport equations, and we shall
    employ a stochastic version of global
    Carleman estimate to derive such an
    estimate. The content of this section
    is a simplified version of
    \cite{Luqi6}.

    In Sections \ref{s-t}, \ref{s-p} and \ref{s-h},
    we assume that
    $\mathbf{F}=\{\cF_t\}_{t\in[0,T]}$ is the
    natural filtration generated by $W(\cd)$.

    %%%%%%%%%%%%%%%%%%%%%%%%%%%%%%%%%%%%%%%%%%%%%%%%%%%%%%%%%

    \subsection{Formulation of the problem}

    Let $G$ be given as in Subsection \ref{sec-mpr-6.1} and
    convex. For any $T > 0$, put
    \bel{2021-1-16e1}
    Q\equiv (0,T)\times G,\q\Si\equiv
    (0,T)\times \G.
    \ee
    Fix an $ O\in \mR^n$ satisfying
    $|O|_{\dbR^n}=1$. Set
    $$
    \begin{array}{ll}\ds
        \G ^- \deq  \big\{ x \in\G \,\big|\,
        O\cd\nu(x)< 0 \big\},&  \G ^+\deq
        \G \setminus \G ^-,\\
        \ns\ds
        \Si^-\deq (0,T)\times \G ^-, & \Si^+\deq  (0,T)
        \times \G ^+.
    \end{array}
    $$
    Define the Hilbert space $ L^2_O(\G ^-)
    $ as the completion of all $h\in
    C_0^\infty(\G ^-)$ w.r.t. the
    norm
    $$
    |h|_{L^2_O(\G ^-)}\deq \(\ds- \int_{\G^-}O\cd\nu |h|^2 d\G\)^{\frac{1}{2}}.
    $$
    Clearly, $L^2(\G ^-)$ is dense in
    $L^2_O(\G ^-)$.

    Let us first recall the controllability problem
    for the following deterministic transport
    equation:
    \begin{equation}\label{9.19-eq1}
        \begin{cases}\ds
            y_t  +  O \cd \nabla y    =  a_1^0  y
            &\mbox{ in } Q,\\
            \ns\ds y  = u  & \mbox{ on } \Si^-,\\
            \ns\ds y(0) = y_0 &\mbox{ in } G,
        \end{cases}
    \end{equation}
    where $y_0\in L^{2}(G)$ and $a_1^0\in L^{\infty}(0,T;L^{\infty}(G))$.
    In \eqref{9.19-eq1}, $y$ is the {\it state
    } and $u$ is the {\it control}.
    The {\it state}  and  {\it control spaces} are
    chosen to be $L^{2}(G)$ and
    $L^2(0,T;L_O^{2}(\G^-))$, respectively.
    It is
    easy to show that, for any $y_0\in L^{2}(G)$ and
    $u\in L^2(0,T;L_O^{2}(\G^-))$, \eqref{9.19-eq1}
    admits a unique (transposition) solution
    $y(\cd)\in C([0,T];L^2(G))$ (e.g., \cite{Lions2}).

    \begin{definition}
        The system \eqref{9.19-eq1} is said to be
        {\it exactly controllable at time $T$} if  for any
        given $y_0, y_1\in L^2(G)$, one can find a
        control $u\in L^2(0,T;L_O^{2}(\G^-))$ such that
        the corresponding solution $y(\cd)$ to \eqref{9.19-eq1} satisfies
        $y(T)=y_1$.
    \end{definition}

    To study the exact controllability problem of \eqref{9.19-eq1}, people introduce the following dual equation:
    \begin{eqnarray}\label{9.19-eq2}
        \left\{
        \begin{array}{lll}
            \ds z_t  + O\cd\!\nabla z =
            -a_1^0 z  &\mbox{ in } Q,\\
            \ns\ds z  = 0 & \mbox{ on } \Si^+,\\
            \ns\ds z(T) = z_T &\mbox{ in }   G.
        \end{array}
        \right.
    \end{eqnarray}
    By means of the standard duality argument, it is
    easy to show the following result (e.g.,
    \cite{Klibanov-Y}).

    \begin{proposition}\label{20160719pro1}
        The equation \eqref{9.19-eq1} is exactly
        controllable at time $T$ if and only if solutions
        to the equation \eqref{9.19-eq2} satisfy the
        following observability estimate:
        \begin{equation}\label{9.22-eq4}
            |z_T|_{L^2(G)}\le
            \cC|z|_{L^2(0,T;L_O^{2}(\G^-))},\qq\forall\;z_T\in L^2(G).
        \end{equation}
    \end{proposition}

    Controllability
    problems for deterministic
    transport equations are now
    well understood. Indeed, one can use
    the global Carleman estimate to prove
    the observability inequality
    \eqref{9.22-eq4} for $T>2\max_{x\in \G} |x|_{\mR^n}$ (c.f.
    \cite{Klibanov-Y}). In the rest of this
    section, we shall see that things are
    quite different in the stochastic
    setting.

    \ss

    Now, we consider the following controlled
    stochastic transport equation:
    \begin{equation}\label{ch-7-csystem1}
        \left\{
        \begin{array}{ll}\ds
            dy  +  O \cd \nabla y dt   =
            \big(a_1  y   +a_3  v\big)dt
            + \big(a_2 y + v\big) dW(t)   &\mbox{ in } Q,\\
            \ns\ds y  = u  & \mbox{ on } \Si^-,\\
            \ns\ds y(0) = y_0 &\mbox{ in } G.
        \end{array}
        \right.
    \end{equation}
    Here $y_0\!\in\! L^2(G)$, and $a_1,a_2,a_3\!\in\!
    L^\infty_{\dbF}(0,T;L^\infty(G))$. In
    \eqref{ch-7-csystem1}, $y$ is the {\it state}, while
    $u\!\in\! L^2_{\dbF}(0,T;L^2_O(\G^-))$ and $v\in
    L^2_\dbF(0,T;L^2(G))$ are two {\it controls}.

    \begin{remark}
        In practice, very likely, the control $v$ in the diffusion term
        affects also the system through its drift
        term (say, in the form $a_3v$, as in (\ref{ch-7-csystem1})). This means
        one cannot simply put a control in the
        diffusion term only. The cost is that it
        will effect the drift term in an
        undesired way, i.e.,   $a_3v$ does not
        work as a control but as a disturbance.
    \end{remark}

    A solution  to (\ref{ch-7-csystem1}) is also
    understood in the transposition sense.
    For this, we introduce the following ``reference" equation which is a
    backward stochastic transport equation:
    \begin{equation}\label{ch-7-csystem2}
        \left\{
        \begin{array}{ll}\ds dz  +
            O\cd\nabla zdt  =
            - (a_1 z  +  a_2 Z ) dt   +  Z dW(t)  &\mbox{in } Q_\tau,\\
            \ns\ds z  = 0 & \mbox{on } \Si_\tau^+,\\
            \ns\ds z(\tau) = z_\tau &\mbox{in }   G,
        \end{array}
        \right.
    \end{equation}
    where $\tau\in (0,T]$,
    $Q_\tau\deq(0,\tau)\times G$,
    $\Si_\tau^+\deq(0,\tau)\times \G^+$ and $
    z_\tau\in
    L^2_{\cF_\tau}(\Omega;L^2(G))$.

    Recalling the operator $A$ given in Example \ref{2.4-ex22}, as an easy consequence of Theorem
    \ref{ch-1-well-bmild}, we have the
    following well-posedness result for
    \eqref{ch-7-csystem2}.
    \begin{proposition}\label{well posed1}
        For any $\tau\in (0,T]$ and $z_\tau\in
        L^2_{\cF_\tau}(\Omega;L^2(G))$, the
        equation \eqref{ch-7-csystem2} admits a
        unique mild solution $(z,Z)\in
        L^2_\dbF(\Omega;C([0,\tau];L^2(G)))\times
        L^2_\dbF(0,\tau;L^2(G))$ such that
        \begin{equation}\label{best1}
            |z|_{L^2_\dbF(\Omega;C([0,\tau];L^2(G)))}
            + |Z|_{L^2_\dbF(0,\tau;L^2(G))} \leq
            \cC
            |z_\tau|_{L^2_{\cF_\tau}(\Omega;L^2(G))},
        \end{equation}
        where $\cC$ is a constant, independent of $\tau$
        and $z_\tau$.
    \end{proposition}

    We need the following further regularity result for solutions to \eqref{ch-7-csystem2}.

    \begin{proposition}\label{hidden}
        For any $\tau\in
        (0,T]$, solutions to  the equation \eqref{ch-7-csystem2}
        satisfy that
        $$
        |z|^2_{L^2_\dbF(0,\tau;L^2_O(\G^-))} \leq
        \cC\mE|z_\tau|^2_{L^2(G)},\q\forall\,z_\tau\in L^2_{\cF_\tau}(\Omega;L^2(G)).
        $$
    \end{proposition}

    {\it Proof}. For simplicity, we only consider
        the case of $\tau=T$.  Let $A$ be the operator introduced in Example \ref{2.4-ex22} and $H=L^2(G)$. Then the equation \eqref{ch-7-csystem2} can be written in the following form:
        \begin{equation}\label{ch-7-csystem2-a1}
            \left\{
            \begin{array}{ll}\ds dz  +
                A^* zdt  = fdt + ZdW(t)dW(t)  &\mbox{in } Q_\tau,\\
                \ns\ds z(T) = z_T,
            \end{array}
            \right.
        \end{equation}
        where $f=- (a_1 z +  a_2 Z )$.
        Let $\{(z_\l,Z_\l)\}_{\l\in \rho(A^*)}$ be the approximation of $(z,Z)$ introduced in \eqref{c1-system5}, i.e.,  $(z_\l,Z_\l)$ solves
        \begin{equation}\label{ch-7-csystem2-a2}
            \left\{
            \begin{array}{ll}\ds dz_\l  +
                A^* z_\l dt  = R^*(\l)fdt +R^*(\l) Z_\l dW(t)  &\mbox{in } [0,T),\\
                \ns\ds z_\l(T) = R^*(\l)z_T,
            \end{array}
            \right.
        \end{equation}
        where $R^*(\l)\deq\l(\l I-A^*)^{-1}$.
        By Theorem
        \ref{ch-1-well-bmild} and Proposition
        \ref{ch-1-well-brel}, one can show that the
        solution to \eqref{ch-7-csystem2-a2} satisfies that
        $ (z_\l,Z_\l)\in
        \big(L^2_\dbF(\Omega;C([0,T];L^2(G)))\cap
        L^2_\dbF(0,T;D(A^*))\big)\t
        L^2_\dbF(0,T;L^2(G)).
        $
        It follows from
        It\^{o}'s formula that
        $$
        \begin{array}{ll}\ds
            \mE| R(\l)z_T|^2_{L^2(G)}- |z_\l(0)|^2_{L^2(G)}\\
            \ns\ds = -2\mE \int_0^T \int_G  z_\l O \cd \nabla z_\l
            dxdt + \mE \int_0^T  \int_G \big(2z_\l R^*(\l)f +|R^*(\l) Z_\l|^2 \big)dxdt.
        \end{array}
        $$
        Therefore,
        \begin{equation}\label{hidden2}
            \begin{array}{ll}\ds
                -\mE\int_0^T\int_{\G^-} O\cd\nu
                z_\l^2d\G^- dt \\
                \ns\ds\le \mE|R(\l)z_T|^2_{L^2(G)} -
                \mE \int_0^T  \int_G \big(2z_\l R^*(\l)f +| R^*(\l)Z_\l|^2 \big)dxdt.
            \end{array}
        \end{equation}
        By letting $\l$ tend to $\infty$ in \eqref{hidden2}, noting Theorem \ref{app th1},  we obtain that
        $$
        \begin{array}{ll}\ds
            -\mE\int_0^T\int_{\G^-} O\cd\nu
            z^2d\G^- dt \\
            \ns\ds \leq  \mE|z_{T}|^2_{L^2(G)} -
            \mE \int_0^T  \int_G \big(2z f +| Z |^2 \big)dxdt\leq \cC \mE|z_{T}|^2_{L^2(G)}.
        \end{array}
        $$
        This complete the proof of Proposition \ref{hidden}.
    \endpf
    \begin{remark}
        By Proposition \ref{well
            posed1}, the first component $z$ (of the
        solution to \eqref{ch-7-csystem2})
        belongs to $
        L^2_{\dbF}(\Omega;$ $C([0,\tau];L^2(G)))$.
        This does not yield the regularity $z|_{\G^-}\in
        L^2_\dbF(0,\tau;$ $L^2_O(\G^-))$, guaranteed by Proposition
        \ref{hidden}.  Hence, the latter is sometimes called a
        {\it hidden regularity property}.
    \end{remark}

    With the aid of Proposition \ref{hidden}, we can give the following definition.
    \begin{definition}
        A stochastic process $ y\in
        C_{\dbF}([0,T];L^2(\Om;L^2(G)))$ is called a
        {\it transposition solution} to \eqref{ch-7-csystem1}
        if for any $\tau\in (0,T]$ and $z_{\tau}\in
        L^2_{\cF_\tau}(\Om;L^2(G))$, it   holds that
        \begin{equation}\label{deftraoid}
            \begin{array}{ll}
                \ds \mE\langle y(\tau), z_{\tau}\rangle_{L^2(G)}
                - \langle y_0, z(0)\rangle_{L^2(G)} \\
                \ns\ds  = \mE\int_0^\tau \langle
                v,a_3z+Z\rangle_{L^2(G)} dt-\mE\int_0^\tau
                \int_{\G^-}O\cdot\nu uzd\G^- dt.
            \end{array}
        \end{equation}
        Here $(z,Z)$ solves \eqref{ch-7-csystem2}.
    \end{definition}
    \begin{remark}\label{rem20210121}
        In Subsection \ref{sec-BSEE}, we have introduced
        the notion of {\it transposition solution} for BSEEs with general filtration, to overcome the difficulty that
        there is no Martingale
        Representation Theorem. Here and in Section  \ref{s-h}, we use the notion of
        {\it transposition solution} for controlled SPDEs
        to overcome the difficulty arisen from the boundary control. Note however that the
        idea to solve these two different problems is very similar, that
        is, to study the well-posedness of a ``less-understood" equation
        by means of another ``well-understood" one. To keep clarity,
        it seems that we would better use two different notions for them but we prefer keeping the present presentation because the readers should be able to detect the differences easily from the context.
    \end{remark}

    We have the following well-posedness result for
    \eqref{ch-7-csystem1} (See Subsection
    \ref{ch-7-se2} for its proof).

    \begin{proposition}\label{well posed}
        For any $y_0\in L^2(G)$, $u\in
        L^2_{\dbF}(0,T;L^2_O(\G^-))$ and $ v\in
        L^2_{\dbF}(0,T;$ $L^2(G))$, the equation
        \eqref{ch-7-csystem1} admits a unique
        (transposition) solution $y\in C_\dbF([0,T];$
        $L^2(\Omega;L^2(G)))$ such that
        \begin{equation}\label{well posed est-chapt7}
            \begin{array}{ll}\ds
                |y|_{C_\dbF([0,T];L^2(\Omega;L^2(G)))}\leq \cC
                \big(|y_0|_{L^2(G)}+
                |u|_{L^2_{\dbF}(0,T;L^2_O(\G^-))} +
                |v|_{L^2_{\dbF}(0,T; L^2(G))}\big).
            \end{array}
        \end{equation}
    \end{proposition}

    Now we can introduce the notion of exact controllability of \eqref{ch-7-csystem1}.
    \begin{definition}
        The system
        \eqref{ch-7-csystem1} is called {\it exactly
            controllable  at time $T$} if for any
        $y_0\in L^2(G)$ and $y_T\in
        L^2_{\cF_T}(\Om;L^2(G))$, one can find
        a couple of controls $(u,v)\in
        L^2_{\dbF}(0,T;L^2_O(\G^-))\times
        L^2_\dbF(0,T;L^2(G))$ such that the
        corresponding (transposition) solution
        $y$ to \eqref{ch-7-csystem1} satisfies
        that $y(T) = y_T$ in $L^2(G)$, a.s.
    \end{definition}

    The main result in this section is the
    following exact controllability result
    for \eqref{ch-7-csystem1}.

    \begin{theorem}\label{exact th}
        If $T > 2\max_{x\in \G} |x|_{\mR^n}$, then the system
        \eqref{ch-7-csystem1} is exactly
        controllable  at time $T$.
    \end{theorem}
    \begin{remark}
        We introduce two controls into the
        system \eqref{ch-7-csystem1}. Moreover,
        the control $v$ acts on the whole
        domain. Compared with the deterministic
        transport equations, it seems that our
        choice of controls is too restrictive.
        One may consider the following three
        weaker cases:

        \ss

        1) Only one control is acted on the system, that
        is, $u=0$ or $v=0$ in \eqref{ch-7-csystem1}.

        2) Neither $u$ nor $v$ is zero. But $v=0$ in
        $(0,T)\times G_0$, where $G_0$ is a nonempty
        open subset of $G$.

        3) Two controls are imposed on the system. But
        both of them are in the drift term.

        \ss

        For the above three cases,
        according to the controllability result
        for deterministic transport equations,
        it seems that the corresponding control
        system should be exactly controllable.
        However, similarly to the proof of
        Theorem \ref{6.18-th3}, one can show
        that it is not the truth (c.f.
        \cite{Luqi6}).
    \end{remark}

    Using the standard duality argument,
    in order to prove Theorem \ref{exact th}, it
    suffices to establish a suitable
    observability estimate for
    \eqref{ch-7-csystem2} with $\tau=T$.
    The latter will be done in  Subsection
    \ref{ch-7-se3}, by means of a
    stochastic version of global Carleman
    estimate.

    %%%%%%%%%%%%%%%%%%%%%%%%%%%%%%%%%%%%%%%%%%%%%%%%%%%%%

    \subsection{Well-posedness of the control system and a weighted
        identity}\label{ch-7-se2}

    %%%%%%%%%%%%%%%%%%%%%%%%%%%%%%%%%%%%%%%%%%%%%%%%%%%%%

    This subsection is addressed to present some
    preliminary results.

    To begin with, let us prove the well-posedness result (i.e., Proposition \ref{well
        posed}) for the controlled
    stochastic transport equation
    \eqref{ch-7-csystem1}.

    {\it Proof}.[Proof of Proposition \ref{well
            posed}] We first prove the existence
        of transposition solutions to \eqref{ch-7-csystem1}. Since $u\in
        L^2_{\dbF}(0,T; L^2_O(\G^-))$, there
        exists a sequence
        $\{u_m\}_{m=1}^\infty\subset
        C_\dbF^2([0,T];$ $H^{3/2}_0(\G^-))$
        with $u_m(0)=0$ for all $m\in\dbN$ such
        that
        \begin{equation}\label{8.31-eq1}
            \lim_{m\to\infty} u_m = u \q \mbox{ in
            }L^2_{\dbF}(0,T;L^2_O(\G^-)).
        \end{equation}
        For each $m\in\dbN$, we can find a $\tilde u_m
        \in C_\dbF^2([0,T];H^2(G))$ such that $\tilde
        u_m|_{\G^-}=u_m$ and $\tilde u_m(0)=0$.

        Consider the following stochastic transport equation:
        \begin{eqnarray}\label{system2n-cha7}
            \left\{\!\!
            \begin{array}{ll}
                \ds d\tilde y_m \! + \! O\!\cd \!\nabla \tilde y_m dt\! =\!
                \big(a_1  \tilde y_m \!  +\!a_3 v\!+\!\zeta_m\big)dt
                \! +\! \big[a_2 (\tilde y_m \!+\! \tilde u_m)\!+\!v\big] dW(t) &\mbox{in }Q, \\
                \ns\ds
                \tilde y_m=0&\mbox{on }\Si,\\
                \ns\ds \tilde y_m(0)= y_0&\mbox{in }G,
            \end{array}
            \right.
        \end{eqnarray}
        where $\ds\zeta_m = -\tilde u_{m,t}- O \cd
        \nabla \tilde u_m dt + a_1  \tilde u_m $. By
        Theorem \ref{ch-1-well-mild}, the system
        \eqref{system2n-cha7} admits a unique mild
        solution $ \tilde y_m\in
        C_\dbF([0,T];L^2(\Om;L^2(G)))$.

        \ss

        Let $y_m=\tilde y_m+\tilde u_m$. For any $\tau\in (0,T]$ and $z_{\tau}\in
        L^2_{\cF_\tau}(\Om;L^2(G))$, by  It\^o's
        formula and integration by parts, we have that
        \begin{equation}\label{8.31-eq5-chap7}
            \begin{array}{ll}
                \ds \mE\langle y_{m}(\tau),
                z_{\tau}\rangle_{L^2(G)} - \langle y_0,
                z(0)\rangle_{L^2(G)} \\
                \ns\ds  = \mE\int_0^\tau \langle
                v,Z\rangle_{L^2(G)} dt-\mE\int_0^\tau
                \int_{\G^-}O\cdot\nu u_{m}zd\G^- dt,
            \end{array}
        \end{equation}
        where $(z,Z)$ is the mild solution to \eqref{ch-7-csystem2}.
        Consequently, for any $m_1,m_2\in\dbN$, it holds that
        \begin{equation}\label{8.31-eq2-chap7}
            \mE\langle y_{m_1}(\tau)-y_{m_2}(\tau),z_\tau
            \rangle_{L^2(G)} =-\mE\int_0^\tau
            \int_{\G^-}O\cdot\nu (u_{m_1}-u_{m_2})zd\G^- dt.
        \end{equation}

        By Proposition \ref{2.1-cor16},  we can choose $ z_\tau\in
        L^2_{\cF_\tau}(\Om;L^2(G))$  so that
        $
        |z_\tau|_{L^2_{\cF_\tau}(\Om;L^2(G))}=1
        $
        and
        \begin{equation}\label{8.31-eq3-chap7}
            \mE\langle y_{m_1}(\tau)-y_{m_2}(\tau),
            z_\tau \rangle_{L^2(G)} \geq
            \frac{1}{2}
            |y_{m_1}(\tau)-y_{m_2}(\tau)|_{L^2_{\cF_\tau}(\Om;L^2(G))}.
        \end{equation}
        It follows from \eqref{8.31-eq2-chap7}--\eqref{8.31-eq3-chap7}  and Proposition
        \ref{hidden} that
        $$
        \begin{array}{ll}
            \ds
            |y_{m_1}(\tau)-y_{m_2}(\tau)|_{L^2_{\cF_\tau}(\Om;L^2(G))}
            \leq  2\Big| \mE\int_0^\tau
            \int_{\G^-}O\cdot\nu (u_{m_1}-u_{m_2})zd\G^- ds
            \Big|\\
            \ns\ds \leq \cC
            |u_{m_1}-u_{m_2}|_{L^2_{\dbF}(0,T;L^2_O(\G^-))}
            |z_\tau|_{L^2_{\cF_\tau}(\Om;L^2(G))} \leq \cC
            |u_{m_1}-u_{m_2}|_{L^2_{\dbF}(0,T;L^2_O(\G^-))},
        \end{array}
        $$
        where the constant $\cC$ is independent of
        $\tau$. Consequently, it holds that
        $$
        |y_{m_1}-y_{m_2}|_{C_\dbF([0,T];L^2(\Om;L^2(G)))}
        \leq \cC
        |u_{m_1}-u_{m_2}|_{L^2_{\dbF}(0,T;L^2_O(\G^-))}.
        $$
        Hence, $\{y_{m}\}_{m=1}^\infty$ is a Cauchy
        sequence in $C_\dbF([0,T];L^2(\Om;L^2(G)))$.
        Denote by $y$ the limit of
        $\{y_m\}_{m=1}^\infty$. Letting $m\to \infty$ in
        \eqref{8.31-eq5-chap7}, we see that $y$
        satisfies \eqref{deftraoid}. Thus, $y$ is a
        transposition solution to \eqref{ch-7-csystem1}.

        By Proposition \ref{2.1-cor16}, we see that there exists $z_\tau\!\in\!
        L^2_{\cF_\tau}\!(\Om;L^2(G))$  such that
        $ |z_\tau|_{L^2_{\cF_\tau}\!(\Om;L^2(G))}\!=\!1 $
        and
        \begin{equation}\label{8.31-eq7-cha7}
            \mE\langle y(\tau), z_\tau \rangle_{L^2(G)}\geq
            \frac{1}{2}|y(\tau)|_{L^2_{\cF_\tau}(\Om;L^2(G))}
            .
        \end{equation}
        Combining \eqref{deftraoid},
        \eqref{8.31-eq7-cha7} and using Proposition
        \ref{hidden} again, we obtain that
        $$
        \begin{array}{ll}
            \ds  |y(\tau)|_{L^2_{\cF_\tau}(\Om;L^2(G))} \\
            \ns\ds \leq 2\(\big|\langle y_0,
            z(0)\rangle_{L^2(G)}\big| + \Big|\mE\int_0^\tau
            \langle v,Z\rangle_{L^2(G)} dt\Big| +
            \Big|\mE\int_0^\tau \int_{\G^-}O\cdot\nu uzd\G^-
            dt
            \Big|\) \\
            \ns\ds \leq \cC\big(|y_0|_{L^2(G)}+
            |u|_{L^2_{\dbF}(0,T;L^2_O(\G^-))} +
            |v|_{L^2_{\dbF}(0,T; L^2(G))}\big),
        \end{array}
        $$
        where the constant $\cC$ is independent of
        $\tau$. Therefore, we obtain the desired
        estimate \eqref{well posed est-chapt7}.

        \ss

        Now we prove the uniqueness of
        solutions.  Assume that $y_1$ and $y_2$
        are two  transposition solution to
        \eqref{ch-7-csystem1}. From
        \eqref{deftraoid}, we deduce that, for
        every $\tau\in [0,T]$,
        $$
        \mE\langle y_1(\tau)-y_2(\tau),
        z_{\tau}\rangle_{L^2(G)}=0,\qq \forall\;
        z_\tau \in L^2_{\cF_\tau}(\Om;L^2(G)),
        $$
        which implies that $y_1=y_2$ in
        $C_\dbF([0,T];L^2(\Om;L^2(G)))$. This
        completes the proof of Proposition
        \ref{well posed}.
    \endpf

    Next, we present a weighted identity (for stochastic transport operators), which will play a key role in
    establishing the global Carleman estimate for
    \eqref{ch-7-csystem2}.

    \begin{proposition}\label{prop identity}
        Let $\ell\in C^1([0,T]\times G)$ and
        $\th=e^\ell$. Let $\mathbf{u}$ be an
        $H^1(\dbR^n)$-valued It\^o process, and
        $\mathbf{v}=\theta \mathbf{u}$. Then\footnote{See Remark \ref{12.24-rmk1} for the term $(d\mathbf{v})^2$},
        \begin{equation}\label{identity}
            \begin{array}{lll}
                \ds  -\theta\big(\ell_t + O\cd\nabla \ell\big)\mathbf{v}\big( d\mathbf{u} + O\cd\nabla \mathbf{u} dt \big)\\
                \ns\ds = -\frac{1}{2}d\big[\big(\ell_t +
                O\cd\nabla \ell\big)\mathbf{v}^2 \big] -
                \frac{1}{2} O\cd\nabla\big[\big(\ell_t  +
                O\cd\nabla \ell\big)\mathbf{v}^2\big]dt  \\
                \ns\ds \q + \frac{1}{2}\big[\ell_{tt} +
                O\cd\nabla (O\cd\nabla\ell)+  2O\cd\nabla
                \ell_{t}  \big] \mathbf{v}^2dt+
                \frac{1}{2}\big(\ell_t + O\cd\nabla
                \ell\big)(d\mathbf{v})^2\\
                \ns\ds\q + \big(\ell_t + O\cd\nabla \ell\big)^2
                \mathbf{v}^2dt.
            \end{array}
        \end{equation}
    \end{proposition}

    {\it Proof}.  Clearly,
        $$
        \begin{array}{ll}\ds
            \theta\big(d\mathbf{u}+O\cd\nabla
            \mathbf{u}dt\big)=\theta d\big(\theta^{-1} \mathbf{v}\big)
            + \theta O\cd\nabla\big(\theta^{-1}
            \mathbf{v}\big)dt\\
            \ns\ds = d\mathbf{v} + O\cd\nabla \mathbf{v}dt -
            \big(\ell_t + O\cd\nabla \ell\big)\mathbf{v}dt.
        \end{array}
        $$
        Thus,
        \begin{equation}\label{identity.1}
            \begin{array}{ll}
                -\theta\big(\ell_t + O\cd\nabla \ell\big)\mathbf{v}\big( d\mathbf{u} + O\cd\nabla \mathbf{u}dt\big)\\
                \ns\ds = -\big(\ell_t +O\cd\nabla \ell \big)\mathbf{v} \big[ d\mathbf{v} + O\cd\nabla \mathbf{v}dt - \big(\ell_t + O\cd\nabla \ell\big)\mathbf{v}dt \big]\\
                \ns\ds = -\big(\ell_t + O\cd\nabla
                \ell\big)\mathbf{v} \big(d\mathbf{v} +
                O\cd\nabla \mathbf{v}dt\big) + \big(\ell_t +
                O\cd\nabla \ell\big)^2\mathbf{v}^2dt.
            \end{array}
        \end{equation}
        Direct computations imply that
        $$
        \left\{
        \begin{array}{lll}\ds
            -\ell_t \mathbf{v}d\mathbf{v} = -\frac{1}{2}d(\ell_t \mathbf{v}^2) + \frac{1}{2}\ell_{tt}\mathbf{v}^2dt + \frac{1}{2}\ell_t (d\mathbf{v})^2,\\
            \ns\ds -O\cd\nabla \ell \mathbf{v}d\mathbf{v} = -\frac{1}{2}d(O\cd\nabla \ell \mathbf{v}^2) + \frac{1}{2}(O\cd\nabla \ell)_t \mathbf{v}^2dt + \frac{1}{2}O\cd\nabla \ell (d\mathbf{v})^2,\\
            \ns\ds -\ell_{t} \mathbf{v} O\cd\nabla \mathbf{v}dt = -\frac{1}{2}O\cd\nabla(\ell_t \mathbf{v}^2)dt + \frac{1}{2}O\cd\nabla \ell_{t}\mathbf{v}^2dt,\\
            \ns\ds -O\cd\nabla \ell \mathbf{v}O\cd\nabla
            \mathbf{v}dt= -\frac{1}{2}O\cd\nabla(O\cd\nabla
            \ell \mathbf{v}^2)dt +
            \frac{1}{2}O\cd\nabla(O\cd\nabla
            \ell)\mathbf{v}^2dt.
        \end{array}
        \right.
        $$
        This, together with  \eqref{identity.1}, implies
        \eqref{identity}. This completes the proof of Proposition \ref{prop identity}.
    \endpf

    %%%%%%%%%%%%%%%%%%%%%%%%%%%%%%%%%%%%%%%%%%%%%%

    \subsection{Observability estimate for backward
        stochastic transport equations}\label{ch-7-se3}

    %%%%%%%%%%%%%%%%%%%%%%%%%%%%%%%%%%%%%%%%

    In this subsection, we shall show the following
    observability estimate for the equation
    \eqref{ch-7-csystem2}, which implies Theorem \ref{exact th}.

    \begin{theorem}\label{th obser}
        If $T> 2\max_{x\in \G} |x|_{\mR^n}$, then solutions to the
        equation \eqref{ch-7-csystem2} with
        $\tau=T$ satisfy that
        \begin{equation}\label{exact ob est}
            \ba{ll} \ds|z_T|_{L^2_{\cF_T}(\Omega;L^2(G))}
            \leq \cC \big(|z|_{L_\dbF^2(0,T;L^2_O(\G^-))} +
            |a_3z+Z|_{L^2_\dbF(0,T;L^2(G))}\big),\\
            \ns\ds\qq\qq\qq\qq\qq\qq\qq\qq\forall\;z_T\in
            L^2_{\cF_T}(\Omega;L^2(G)). \ea
        \end{equation}
    \end{theorem}

    {\it Proof}.
        Let $\l>0$, and let $c\in(0,1)$ be such
        that $cT> 2\max_{x\in \G} |x|_{\mR^n}$. Put
        \begin{equation}\label{weight1}
            \ell = \l \[|x|_{\dbR^n}^2  - c
            \(t-\frac{T}{2}\)^2 \].
        \end{equation}
        Applying Proposition \ref{prop
            identity} to the equation
        \eqref{ch-7-csystem2} with
        $\mathbf{u}=z$ and $\ell$ given by
        \eqref{weight1}, integrating
        \eqref{identity} on $Q$, using
        integration by parts and taking expectation,
        we obtain that
        \begin{equation}\label{bobeq1}
            \begin{array}{lll}\ds
                -2\mE \int_Q \theta^2(\ell_t +O\cd\nabla \ell)z
                (dz + O\cd\nabla z dt)dx
                \\
                \ns\ds =\l \mE\int_G  (cT-2O\cd x) \theta^2(T)
                z^2(T) dx +\l \mE\int_G  (cT  + 2O\cd
                x)\theta^2(0)
                z^2(0) dx\\
                \ns\ds \q + \l\mE  \int_{\Si_-}O\cd \nu \big[c(T
                - 2t)-2 O\cd x\big]\theta^2 z^2 d\Si_-+ 2(1 -
                c)\l\mE\int_Q \theta^2 z^2 dxdt
                \\
                \ns\ds \q + \mE\int_Q\theta^2 (\ell_t  + O\cd
                \nabla \ell) (dz)^2 dx+ 2\mE\int_Q\theta^2
                (\ell_t +O\cd \nabla \ell)^2 z^2 dxdt.
            \end{array}
        \end{equation}
        Noting that $(z,Z)$ solves
        \eqref{ch-7-csystem2} with $\tau=T$,  we obtain that
        \begin{equation}\label{bobeq2}
            \begin{array}{ll}\ds
                -2\mE \int_Q \theta^2(\ell_t + O\cd\nabla \ell)
                z (dz + O\cd\nabla zdt) dx
                \\
                \ns\ds = 2\mE \int_Q \theta^2(\ell_t + O \cd
                \nabla \ell) z \big( a_1z +  a_2 Z \big) dxdt
                \\
                \ns\ds\leq \frac{1}{2}\mE\int_Q\theta^2 (\ell_t
                + O\cd\nabla \ell)^2 z^2 dxdt + \cC\mE\int_Q
                \theta^2 ( z^2 + |Z|^2) dxdt\\
                \ns\ds\leq \frac{1}{2}\mE\int_Q\theta^2 (\ell_t
                + O\cd\nabla \ell)^2 z^2 dxdt + \cC\mE\int_Q
                \theta^2 ( z^2 + |a_3z+Z|^2) dxdt
            \end{array}
        \end{equation}
        and
        \begin{equation}\label{bobeq2-chp7}
            \begin{array}{ll}\ds
                \mE\int_Q\theta^2 (\ell_t  + O\cd \nabla \ell)
                (dz)^2 dx \\
                \ns\ds= \mE\int_Q \theta^2 (\ell_t + O\cd\nabla
                \ell) |(a_3z+Z)-a_3z|^2 dxdt
                \\
                \ns\ds \leq \frac{1}{2}\mE\int_Q \theta^2
                (\ell_t  + O\cd\nabla \ell)^2 z^2 dxdt +
                \cC\mE\int_Q \theta^2 z^2 dxdt\\
                \ns\ds\q +  4\mE\int_Q \theta^2 |\ell_t +
                O\cd\nabla \ell||a_3z+Z|^2 dxdt.
            \end{array}
        \end{equation}
        From \eqref{bobeq1}--\eqref{bobeq2-chp7} and
        $z(T)=z_T$ in $G$ a.s., it follows that
        \begin{eqnarray}\label{bobeq3}
            \begin{array}{lll}\ds
                \l \mE\int_G (cT -2O\cd x) \theta^2(T) z_T^2 dx
                +\l \int_G (cT+ 2O\cd x)\theta^2(0) z^2(0) dx\\
                \ns\ds\q  + 2(1-c)\l\mE \int_Q \theta^2 z^2 dxdt
                +\mE\int_Q  \theta^2 (\ell_t + O\cd\nabla
                \ell)^2 z^2 dxdt
                \\
                \ns\ds \leq \cC\mE \! \int_Q\! \theta^2 ( z^2\! +\!\l|a_3z\!+\!Z|^2) dxdt - \l \mE \!\int_{\Si_-}\! O \cd \nu
                \big[c(T\!-\!2t)\!-\!2 O \cd x\big]\theta^2 z^2 d\Si_-,
            \end{array}
        \end{eqnarray}
        where the constant $\cC>0$ is independent of
        $\l$.

        Finally, recalling that $c\in(0,1)$ satisfies
        $cT> 2\max_{x\in \G} |x|_{\mR^n}$, by choosing $\l$  large enough in
        \eqref{bobeq3}, we obtain the desired estimate
        \eqref{exact ob est}.
        This completes the proof of Theorem \ref{th
            obser}.
    \endpf
    %%%%%%%%%%%%%%%%%%%%%%%%%%%%%%%%%%%%%%%%%%%%%%%%%%%%%%%%%%%%%%%%%%%%

    %%%%%%%%%%%%%%%%%%%%%%%%%%%%%%%%%%%%%%%%%%%%%%

    \subsection{Exact controllability of the
        stochastic transport equation}\label{ch-7-se4}

    %%%%%%%%%%%%%%%%%%%%%%%%%%%%%%%%%%%%%%%%

    This subsection is addressed to a proof of the exact controllability of \eqref{ch-7-csystem1}, i.e., Theorem
    \ref{exact th}.

    {\it Proof}.[Proof of  Theorem \ref{exact th}]
        Fix  any
        $y_0\in L^2(G)$ and $y_T\in
        L^2_{\cF_T}(\Om;L^2(G))$. Let us
        introduce a linear subspace of
        $L_{\dbF}^2(0,T;L^2_{O}(\G^-))\times
        L_{\dbF}^2(0,T;L^2(G))$ as follows:
        $$
        \begin{array}{ll}\ds
            \cY\deq\Big\{\big(-z|_{\G^-},
            a_3z+Z\big)\;\Big|\;(z,Z)\hb{ solves the equation
            }\eqref{ch-7-csystem2}\mbox{ with some } \\
            \ns\ds \hspace{3.9cm} z_T\in
            L^2_{\cF_T}(\Om;L^2(G))\Big\},
        \end{array}
        $$
        and define a linear functional $F$ on $\cY$ as
        follows:
        $$
        F(-z|_{\G_{S}^-},
        a_3z+Z)=\mathbb{E}\int_G
        y_Tz_T dx -  \int_G
        y_0 z(0) dx.
        $$
        By Theorem \ref{th obser}, we see that $F$
        is a bounded linear functional on $\cY$. By
        means of Theorem \ref{Hahn-Banach}, $F$ can be
        extended to be a bounded linear functional on
        the space
        $L_{\dbF}^2(0,T;L^2_{O}(\G^-))\times
        L_{\dbF}^2(0,T;L^2(G))$. For
        simplicity, we still use $F$ to denote this
        extension. By Theorem \ref{1t10s},  there exists $
        (u,
        v)\in L_{\dbF}^2(0,T;L^2_{O}(\G^-))\times
        L_{\dbF}^2(0,T;L^2(G))
        $
        such that
        \begin{equation}\label{con eq0}
            \begin{array}{ll}\ds
                \mE\int_G y_Tz_T dx -  \int_G
                y_0 z(0) dx= F(-z|_{\G_{S}^-},
                a_3z+Z)\\
                \ns\ds = -\mE\int_0^T\int_{\G^-} O\cd\nu zu
                d\G^- dt + \mE\int_0^T\int_G
                v(a_3z+ Z) dxdt.
            \end{array}
        \end{equation}
        We claim that $u$ and $
        v$ are the desired controls. Indeed, by
        the definition of the transposition solution to \eqref{ch-7-csystem1}, we have
        \begin{equation}\label{con eq2}
            \begin{array}{ll}\ds
                \mE\int_G y(T,\cd)z_T dx-  \int_G
                y_0 z(0) dx\\
                \ns\ds
                = -\mE\int_0^T\int_{\G^-}O\cd\nu zu
                d\G^- dt + \mE\int_0^T\int_G v
                (a_3z+ Z) dxdt.
            \end{array}
        \end{equation}
        From   \eqref{con eq0} and  \eqref{con eq2}, we
        see that
        $\ds\mE\int_G y_T z_T dx =
        \mE\int_G y(T)z_T dx
        $
        for any $z_T\in L^2_{\cF_T}(\Om;L^2(G))$. Hence,
        $y(T)=y_T$ in $L^2(G)$, a.s. This
        completes the proof of Theorem \ref{exact th}.
    \endpf
    \begin{remark}
        In the proof of Theorem \ref{exact th}, we only give the existence of the controls which drive the state to the desired destination. Some further characterization for these controls can be found in \cite[Section 7.4]{LZ3.1}.
    \end{remark}
    %

    %%%%%%%%%%%%%%%%%%%%%%%%%%%%%%%%%%%%%%%%%%%%%%%%%%%%

    \section{Controllability of stochastic partial differential equations II:\linebreak stochastic parabolic
        equations}\label{s-p}

    %%%%%%%%%%%%%%%%%%%%%%%%%%%%%%%%%%%%%%%%%%%%%%%%%%%%

    This section, based mainly on \cite{Tang-Zhang1}, is devoted to studying the null/approximate controllability and
    observability of stochastic parabolic equations, one class of typical SPDEs.

    %%%%%%%%%%%%%%%%%%%%%%%%%%%%%%%%%%%%%%%%%%%%%%%%%%%%

    \subsection{Formulation of the problem and the main result}\label{sto heat se1}

    %%%%%%%%%%%%%%%%%%%%%%%%%%%%%%%%%%%%%%%%%%%%%%%%%%%%

    Let $G$ be given as in Subsection \ref{sec-mpr-6.1}, and $G_0$ be a given nonempty open
    subset of $G$.  For any $T > 0$, put $Q$ and $\Si$ as in \eqref{2021-1-16e1}, and
    $$
    Q_0\deq (0,T)\times G_0.
    $$

    We begin with the following controlled
    deterministic parabolic equation:
    \be\label{he}\left\{\ba{ll}\ds
    y_{t}-\D y =a_1^0  y+\chi_{G_0} u
    & \hbox{in }Q,\\
    \ns\ds y=0  & \hbox{on }\Si,
    \\
    \ns\ds  y(0)=y_{0}   & \hbox{in }G, \ea\right.
\end{equation}
where $a_1^0  \in L^{\infty}(0,T;L^{\infty}(G))$. In (\ref{he}), $y$ and $u$ are the {\it state } and {\it
    control} variables, while the {\it state} and {\it control} spaces are
chosen to be $L^{2}(G)$ and $L^{2}(Q_0)$, respectively.

\begin{definition}\label{def-dp-non-app}
    The equation (\ref{he}) is said to be {\it null controllable at time $T$}
    (\resp {\it approximately controllable at time $T$}) if for any given
    $y_0\in L^2(G)$ (\resp for any given $\e>0$ and $y_0, y_1\in L^2(G)$),
    one can find a control $u\in L^2(Q_0)$ such that the
    solution $y(\cd)\in C([0,T];L^2(G))\cap
    L^2(0,T;H_0^1(G))$ to
    (\ref{he}) satisfies $y(T)=0$ (\resp $|y(T)-y_1|_{L^2(G)}\le \e$).
\end{definition}
\begin{remark}
    Due
    to the smoothing effect of solutions to  parabolic equations, the exact
    controllability for (\ref{he}) is impossible, i.e., the $\e$ in Definition \ref{def-dp-non-app}
    cannot be zero.
\end{remark}

The dual
equation of (\ref{he}) is
\begin{eqnarray}\label{0719olwsystem1}
    \left\{
    \begin{array}{lll}
        \ds  z_t+ \D z  = - a_1^0z\quad
        &\mbox{ in } Q,\\
        \ns\ds z = 0 & \mbox{ on } \Si, \\
        \ns\ds z(T) = z_T&\mbox{ in } G.
    \end{array}
    \right.
\end{eqnarray}

By means of the standard duality argument, it is easy to show the following result.
\begin{proposition}
    {\rm i)} The equation (\ref{he}) is null controllable at time $T$ if and only if solutions to the equation
    (\ref{0719olwsystem1}) satisfy the following observability estimate:
    \bel{20160719e1}
    |z_0|_{L^2(G)}\le \cC|z|_{L^2(Q_0)},\qq\forall\;w_0\in
    L^2(G);
    \ee

    {\rm ii)} The equation (\ref{he}) is approximately
    controllable at time $T$ if and only if any solution to the equation
    (\ref{0719olwsystem1}) satisfy the following unique continuation property:
    $$
    z=0\;\hbox{ in }Q_0\Longrightarrow z_T=0.
    $$
\end{proposition}

Controllability of
deterministic parabolic equations is now
well-understood. One can use the global
Carleman estimate to prove the observability
inequality \eqref{20160719e1} and the unique continuation property of (\ref{0719olwsystem1}) (c.f. \cite{FLZ1, FI}).
%In the rest of this section, we shall see a
%quite different picture in the stochastic
%setting.

\ss

We now consider the following controlled stochastic
parabolic equation:
\begin{equation}\label{heat 1.1}
    \left\{
    \begin{array}{ll}
        \ds dy-\D y dt =
        \big(a_1 y+\chi_{G_0}u+a_3v\big)dt+(a_2y+v)\, dW(t)&\hb{ in }Q,\\
        \ns
        \ds y=0&\hb{ on }\Si,\\
        \ns y(0)=y_0&\hb{ in } G,
    \end{array}
    \right.
\end{equation}
where $a_{1},  a_2, a_3\in L^{\infty}_{\dbF}(0,T;
L^{\infty}(G))$.
In the system \eqref{heat 1.1}, the initial
state $y_0\in L^2(G)$, $y$ is the {\it state
    variable}, and the {\it control variable} consists of
$(u, v)\in L^2_{\dbF}(0,T;$ $L^2(G_0))\times
L^2_{\dbF}(0,T;L^2(G))$.

\begin{remark}
    Similarly to the control system \eqref{ch-7-csystem1}, the control $v$ in the diffusion term
    affects the drift
    term in the form $a_3v$.
\end{remark}

By Theorem
\ref{ch-1-well-mild2},
for any $y_0 \in L^2(G)$ and $(u, v)\in
L^2_{\dbF}(0,T;L^2(G_0))\times
L^2_{\dbF}(0,T;$ $L^2(G))$, the system \eqref{heat
    1.1} admits a unique weak solution $y\in
L_{\dbF}^2(\Omega;C([0,T];L^2(G)))$ $\cap
L_{\dbF}^2(0,T; H_0^1(G))$. Moreover,
\begin{equation}\label{20160616e3} \ba{ll}\ds |y|_{
        L_{\dbF}^2(\Omega;C([0,T];L^2(G)))\cap
        L_{\dbF}^2(0,T;H_0^1(G))}\\
    \ns\ds \le \cC \big(|y_0|_{L^2(G)}+|(u,
    v)|_{L^2_{\dbF}(0,T;L^2(G_0))\times
        L^2_{\dbF}(0,T;L^2(G))}\big). \ea
\end{equation}
\begin{definition}
    The system \eqref{heat 1.1} is said to be {\it null
        controllable} (\resp  {\it approximately controllable}) {\it  at time $T$} if for any $y_0\in
    L^2(G)$ (\resp for any given $\e>0$, $y_0 \in L^2(G)$ and $ y_1\in L^2_{\cF_T}(\Om;L^2(G))$), there exists a pair of $(u,v)\in
    L^2_{\dbF}(0,T;L^2(G_0))\times
    L^2_{\dbF}(0,T;L^2(G))$ such that the
    corresponding solution to \eqref{heat 1.1}
    fulfills that $y(T)=0$,  a.s. (\resp $|y(T)-y_1|_{L^2_{\cF_T}(\Om;L^2(G))}\le \e$).
\end{definition}

The main result of this section can be stated as follows:

\begin{theorem}\label{f heat control1}
    System \eqref{heat 1.1} is null and approximately controllable at any
    time $T>0$.
\end{theorem}

The rest of this section is addressed to proving Theorem \ref{f heat control1}. For this, we need some preliminaries.

%%%%%%%%%%%%%%%%%%%%%%%%%%%%%%%%%%%%%%
\subsection{A weighted identity for stochastic parabolic operators}
%%%%%%%%%%%%%%%%%%%%%%%%%%%%%%%%%%%%%%

In this subsection, we
derive a weighted identity for the stochastic parabolic operator ``$dh+\Delta hdt$".

In
the rest of this notes, for simplicity, we use
the notation $ y_{ x_j}=\frac{\pa
    y}{\pa x_j}$, $j=1,\cdots,n$ for the partial
derivative of a function $y$ w.r.t.
$x_j$, where $x_j$ is the $j$-th coordinate of a
generic point $x=(x_1,\cdots,x_n)$ in $\dbR^n$.

Assume that $\ell\in C^{1,3}(Q)$ and $\Psi\in
C^{1,2}(Q)$. Write
\begin{equation}\label{c1e15}
    \left\{
    \begin{array}{ll}
        \ds \cA =|\nabla\ell|^2
        -\D\ell-\Psi-\ell_t,\\
        \ns
        \ds
        \cB=2\big[\cA \Psi+
        \div\big(\cA  \nabla\ell \big) \big] -\cA _t+\D\Psi,\\
        \ns
        \ds
        c^{jk}=2\ell_{x_{j}x_{k}}
        - \d^{jk}\D\ell  -\Psi \d^{jk},\q j,k=1,\cds,n,
    \end{array}
    \right.
\end{equation}
where
$$
\d^{jk}=\left\{
\ba{ll}
1,&\hbox{ if }j=k,\\
0,&\hbox{ if }j\not=k.
\ea
\right.
$$
We have the following fundamental weighted
identity.

\begin{theorem}\label{c1t1}
    Let $h$ be an $H^2(G)$-valued It\^o process. Set $\th=e^{\ell }$ and $w=\th
    h$. Then, for any $t\in [0,T]$ and a.e.
    $(x,\om)\in G\times\Omega$\footnote{See Remark \ref{12.24-rmk1} for the term $(dh)^2$},
    \begin{eqnarray}\label{c1e2a}
        \begin{array}{ll}
            \ds 2 \th\big(\D w +\cA w\big)\big(dh+\D h
            dt\big)-2 \div
            (\nabla w dw) \\
            \ns\ds\q + 2 \div\[ 2\big(\nabla\ell\cd \nabla
            w\big) \nabla w - |\nabla w |^2\nabla\ell
            \!-\!\Psi
            w\nabla w\!+ \(\cA \nabla\ell +\frac{\nabla\Psi}{2}\)w^2\] dt \\
            \ns\ds =2 \!\sum_{j,k=1}^n
            c^{jk}w_{x_j}w_{x_k}dt + \cB w^2dt-d\big(
            |\nabla w|^2-\cA w^2\big)   +2 \big( \D w +\cA w\big)^2dt \\
            \ns\ds\q + \th^2|d\nabla h  + \nabla\ell
            dh|^2- \th^2\cA (dh)^2.
        \end{array}
    \end{eqnarray}
\end{theorem}

{\it Proof}.  Recalling that $\th=e^\ell$ and
    $w=\th h$, one has
    $
    dh=\th^{-1}(dw-\ell_twdt)
    $
    and
    $
    h_{x_j}=\th^{-1}(w_{x_j}-\ell_{x_j}w) $ for $j=1,2,\cds,n$.
    Hence,
    \begin{equation}\label{t2}
        \begin{array}{ll}
            \ds \th \D h \3n&\ds
            =\th\sum_{j=1}^n[\th^{-1}(w_{x_j}-\ell_{x_j}
            w)]_{x_j}
            \\ \ns &\ds  =\sum_{j=1}^n[(w_{x_j}-\ell_{x_j} w)]_{x_j}-\sum_{j=1}^n (w_{x_j}-\ell_{x_j} w)\ell_{x_j}\\
            \ns& \ds = \D w - 2\nabla\ell\cd\nabla w
            +(|\nabla\ell|^2 -\D\ell)w.
        \end{array}
    \end{equation}

    Put
    \begin{equation}\label{t4}
        \begin{array}{ll}\ds I\deq \D w +\cA w,\qq I_1\deq \big(\D w +\cA w\big)dt,\\
            \ns\ds I_2\deq dw-2\nabla\ell\cd \nabla wdt,
            \qq I_3\deq \Psi wdt.
        \end{array}
    \end{equation}
    From (\ref{t2}) and (\ref{t4}), it follows that
    \begin{equation}\label{heat 1.1 c1e5}
        \begin{array}{ll}
            \displaystyle 2 \th\big(\D w +\cA
            w\big)\big(dh+\D h dt\big) =2 I(I_1+I_2+I_3).
        \end{array}
    \end{equation}

    Now,  we compute the right hand side of \eqref{heat 1.1 c1e5}.
    Let us first deal with $2 II_2$. We have
    \begin{eqnarray}\label{2c2t8}
        2 \sum_{j,k=1}^n
        \ell_{x_{j}}w_{x_k}w_{x_jx_{k}}  =
        \sum_{j,k=1}^n \ell_{x_{j}}(w_{x_k}^2)_{x_{j}} =
        \sum_{j,k=1}^n \big[(\ell_{x_{j}}w_{x_k}^2)_{x_{j}}-
        \ell_{x_{j}x_{j}}w_{x_k}^2\big].
    \end{eqnarray}
    From \eqref{2c2t8}, we obtain that
    \begin{equation}\label{heat1.1 1c1.3}
        \begin{array}{ll} \displaystyle
            -4\big(\D w +\cA w\big)\nabla\ell\cd\nabla w \\
            \ns \ds =- 4 \div[(\nabla\ell\cd\nabla w)\nabla
            w] +4 \sum_{j,k=1}^n  \ell_{x_{j}x_k}
            w_{x_j}w_{x_{k}} \\
            \ns\ds\q +4 \sum_{j,k=1}^n
            \ell_{x_{j}}w_{x_k}w_{x_jx_{k}} -
            2\cA \nabla\ell\cd\nabla (w^2) \\
            \ns \ds = - 2\div\big[ 2\big(\nabla \ell\cd\nabla
            w\big)\nabla w -  |\nabla w|^2\nabla \ell +\cA
            \nabla\ell  w^2\big]  \\
            \ns\ds\q +2 \sum_{j,k=1}^n \big(2
            \ell_{x_{j}x_{k}} - \d^{jk}\D\ell
            \big)w_{x_j}w_{x_k}+ 2\div(\cA \nabla\ell) w^2.
        \end{array}
    \end{equation}
    Using It\^o's formula, we have
    \begin{equation}\label{2c1}
        \begin{array}{ll}
            \ds
            2 \big(\D w +\cA w\big)dw= 2 \div( \nabla w dw) -2 \nabla w \cd
            d\nabla w +2 \cA wdw
            \\
            \ns \ds =2 \div( \nabla w dw) + d\big(-|\nabla
            w|^2+\cA w^2\big)  + |d\nabla w|^2- \cA _tw^2dt
            - \cA (dw)^2.
        \end{array}
    \end{equation}
    Combining (\ref{t4}), (\ref{heat1.1 1c1.3}) and
    (\ref{2c1}), we arrive at
    \begin{equation}\label{2c11}
        \begin{array}{ll}
            \ds 2 II_2\3n&\ds= - 2 \div\big[ 2
            \big(\nabla\ell \cd\nabla w\big) \nabla w
            -  |\nabla w|^2\nabla \ell +\cA  w^2 \nabla\ell\big] dt +2\div( \nabla w dw)\\
            \ns&\ds \q + d\big(-|\nabla w|^2+\cA w^2\big) +2
            \sum_{j,k=1}^n \big(2 \ell_{x_{j}x_{k}} -
            \d^{jk}\D\ell
            \big)w_{x_j}w_{x_k}dt\\
            \ns &\ds\q- \big[\cA_t- 2\div(\cA
            \nabla\ell)\big] w^2dt+ |d\nabla w|^2- \cA
            (dw)^2.
        \end{array}
    \end{equation}

    \ss

    Next, we compute $2 II_3$. By
    (\ref{t4}), we get
    \begin{equation}\label{heat 1.1 2c2t11}
        \begin{array}{ll}
            \displaystyle 2 II_3\3n&\ds =2 \big(\D w +\cA w\big)\Psi wdt\\
            \ns&\ds= \big[ 2\div \big(\Psi w \nabla w \big) -
            2\Psi |\nabla w|^2 -\nabla\Psi\cd\nabla(w^2) + 2\cA \Psi w^2\big]dt\\
            \ns&\ds=  \big[\div\big(2\Psi w \nabla w -
            w^2\nabla\Psi \big)  - 2\Psi |\nabla w|^2
            +\big(\D\Psi +2\cA \Psi\big) w^2\big]dt.
        \end{array}
    \end{equation}
    \ss

    Finally, combining the equalities
    (\ref{heat 1.1 c1e5}), (\ref{2c11}) and
    (\ref{heat 1.1 2c2t11}), and noting that
    $$
    \begin{array}{ll}
        \ds- |d\nabla w|^2+ \cA (dw)^2 = -\th^2 |\nabla
        h +dh\nabla\ell|^2+ \th^2\cA (dh)^2,
    \end{array}
    $$
    we obtain the desired equality (\ref{c1e2a})
    immediately. \endpf

Next, for any fixed nonnegative and nonzero function
$\psi\in C^4(\cl{G})$, and parameters $\l>1$ and
$\mu>1$, we choose
\be
\label{alphad}
\th=e^{\ell }, \q\ell=\l\a,\q\a(t,x)=\frac{e^{\mu\psi(x)}-e^{2\mu|\psi|_{C(\cl{G })}}}{ t(T-t)},\q
\varphi(t,x)=\frac{e^{\mu\psi(x)}}{t(T-t)},
\end{equation}
and
\begin{equation}\label{h5}
\Psi=-2\D\ell.
\end{equation}
In what follows, for a nonnegative integer $m$, we
denote by $O(\mu^m)$ a function of order $\mu^m$
for large $\mu$ (which is independent of $\l$);
by $O_{\mu}(\l^m)$ a function of order $\l^m$
for fixed $\mu$ and for large $\l$. Likewise, we use the notation
$O(e^{\mu|\psi|_{C(\cl{G })}})$ and so on. For
$j,k=1,2,\cdots,n$, it is easy to check that
\begin{equation}\label{heat 1.1 h2}
\ell_t=\l\a_t,\q \ell_{x_j}=\l\mu\f\psi_{x_j},\q
\ell_{x_jx_k}=\l\mu^2\f\psi_{x_j}\psi_{x_k}+\l\mu\f\psi_{x_jx_k}
\end{equation}
and that
\begin{equation}\label{h4}
\begin{array}{ll}\ds
    \a_t=\f^2O(e^{2\mu |\psi|_{C(\cl{G})}}),\qq \a_{tt} =\f^3O(e^{2\mu |\psi|_{C(\cl{G})}}),\\
    \ns\ds \f_t
    =\f^2O(1),\qq\qq\q\;\, \f_{tt} =\f^3O(1).
\end{array}
\end{equation}

We shall need the following estimates on the $\cA$, $\cB$ and
$c^{jk}$ appearing in the  weighted identity \eqref{c1e2a} (See \eqref{c1e15} for their definitions, with $\ell$ and $\Psi$ given in \eqref{alphad}--\eqref{h5}):

\begin{proposition}\label{c1t2}
When $\l$
and $\mu$ are large enough, it holds that
\begin{eqnarray}\label{h6}
    \left\{
    \begin{array}{ll}\ds
        \cA = \l^2\mu^2\f^2|\nabla\psi|^2
        +\l\f^2O(e^{2\mu
            |\psi|_{C(\cl{G})}}), \\
        \ns
        \ds \cB\ge \ds 2 \l^3\mu^4\f ^3|\n \psi|^4
        \! +\!\l^3\f^3O(\mu^3)
        \!   +\!\l^2\f^3 O(\mu^2e^{2\mu
            |\psi|_{C(\cl{G})}})\!+\!\l^2\f^2 O(\mu^4), \\
        \ns
        \ds
        \sum_{j,k=1}^n c^{jk}\xi_j\xi_k\ge [ \l\mu^2\f |\n\psi|^2+\l\f  O(\mu)]|\xi
        |^2,\q\forall\;\xi=(\xi_1,\cdots,\xi_n)\in\dbR^n,
    \end{array}\right.
\end{eqnarray}
for any $t\in [0,T]$ and $x\in G$ with $|\n \psi(x)|>0$.
\end{proposition}

{\it Proof}. Noting (\ref{h5})--(\ref{heat 1.1
h2}), from (\ref{c1e15}), we have that
$$
\ell_{x_jx_k}=\l\mu^2\f
\psi_{x_j}\psi_{x_k}+\l\f O(\mu)
$$
and that
\begin{eqnarray*}
    \begin{array}{ll}\ds
        \sum_{j,k=1}^n c^{jk}\xi_j\xi_k  =\sum_{j,k=1}^n\big(2 \ell_{x_{j}x_{k}}  +
        \d^{jk}\D\ell \big)\xi_j\xi_k\\
        \ns\ds  =\sum_{j,k=1}^n \big[\big(2\l\mu^2\f
        \psi_{x_{j}}\psi_{x_{k}} + \l\mu^2\f
        \d^{jk}|\nabla\psi|^2  + \l\f  O(\mu)\big)\big]\xi_j\xi_k\\
        \ns\ds =2\l\mu^2\f \big(\nabla\psi\cd\xi
        \big)^2+\l\mu^2\f |\nabla\psi|^2|\xi|^2+\l\f |\n w|^2 O(\mu)\\
        \ns\ds \ge \big[ \l\mu^2\f |\n\psi|^2+\l\f O(\mu)\big]|\xi|^2,
    \end{array}
\end{eqnarray*}
which gives the last inequality in (\ref{h6}).

\ss

Similarly, by the definition of $\cA$ in
\eqref{c1e15}, and noting \eqref{h4},  we see
that
$$
\begin{array}{ll}
    \ds \cA = |\nabla\ell|^2 +\D\ell -\ell_t\\
    \ns \ds \q = \l\mu \big( \l\mu\f^2|\nabla\psi|^2
    +  \mu\f|\nabla\psi|^2 +\f\D\psi \big)
    +\l\f^2O(e^{2\mu
        |\psi|_{C(\cl{G})}})\\
    \ns \ds \q =
    \l^2\mu^2\f^2|\nabla\psi|^2+\l\f^2O(e^{2\mu
        |\psi|_{C(\cl{G})}}).
\end{array}
$$
Hence, we get the first estimate in \eqref{h6}.

Now, let us estimate $\cB$ (recall \eqref{c1e15}
for the definition of $\cB$). By
\eqref{heat 1.1 h2}, and recalling the
definition of $\Psi$ (in \eqref{h5}), we see
that
$$
\begin{array}{ll}\ds
    \ds\Psi\ds =-2\l\mu (\mu\f |\nabla\psi|^2 +
    \f
    \D\psi )=-2\l\mu^2\f |\nabla\psi|^2 +\l\f O(\mu),
    \\
    \ns\ds
    \Psi_{x_k}=-2 \D\ell_{x_k} = -2\l\mu^3\f |\nabla\psi|^2\psi_{x_k}+\l\f O(\mu^2),\\
    \ns\ds
    \Psi_{x_jx_k}=-2 \D\ell_{x_jx_k}
    =-2\l\mu^4\f|\nabla\psi|^2\psi_{x_j}\psi_{x_k}+\l\f
    O(\mu^3),
    \\
    \ns\ds
    \D\Psi\ds= -2\l\mu^4\f|\nabla\psi|^4 +\l\f
    O(\mu^3).
\end{array}
$$
Hence, recalling  the definition of $\cA$ (in
\eqref{c1e15}), and using \eqref{heat 1.1 h2}
and \eqref{h4}, we have that
$$
\begin{array}{ll}
    \cA\Psi=-2\l^3\mu^4\f^3|\nabla\psi|^4+\l^3\f^3O(\mu^3)+\l^2\mu^2 \f^3O\big(e^{2\mu
        |\psi|_{C(\cl{G})}}\big),
    \\
    \ns\ds \cA_{x_k}= \sum_{j=1}^n\big(2
    \ell_{x_{j}}\ell_{x_{j}x_k}
    +  \ell_{x_{j}x_{j}x_k}\big)-\ell_{tx_k} \\
    \ns \ds \q\;\,\;\; =
    2\l^2\mu^3\f^2|\nabla\psi|^2\psi_{x_k}+\l^2\f^2
    O(\mu^2)+\l \f
    O(\mu^3)+\l\mu \f^2O\big(e^{2\mu
        |\psi|_{C(\cl{G})}}\big),
    \\
    \ns\ds
    -\sum_{j=1}^n\cA_{x_j}
    \ell_{x_j}\!=\!-2\l^3\mu^4\f^3 |\nabla\psi|^4
    \!+\!\l^3\f^3 O(\mu^3)\!+\!\l^2\f^2
    O(\mu^4)\!+\!\l^2\mu^2\f^3O\big(e^{2\mu
        |\psi|_{C(\cl{G})}}\big),
    \\
    \ns\ds -\sum_{j=1}^n\!\big(\cA
    \ell_{x_j}\big)_{x_j}\!
    =-\sum_{j=1}^n\cA_{x_j} \ell_{x_j}-\cA \D\ell \\
    \ns\ds\qq\qq\qq\;\, = -3\l^3\mu^4\f^3 |\nabla\psi|^4
    + \l^2\f^2
    O(\mu^4) + \l^3\f^3 O(\mu^3) \\
    \ns\ds\qq\qq\qq\;\, \q + \l^2\mu^2\f^3O\big(e^{2\mu
        |\psi|_{C(\cl{G})}}\big),
    \\
    \ns\ds\cA_t=
    \big(|\nabla\ell|^2
    +\D\ell +\ell_t\big)_t=\l^2\mu^2\f^3 O(e^{2\mu
        |\psi|_{C(\cl{G})}}).
\end{array}
$$
From the definition of $\cB$ (see
\eqref{c1e15}), we have that
$$
\begin{array}{ll}
    \cB \3n&=\ds -4\l^3\mu^4\f^3|\nabla\psi|^4+\l^3\f^3O(\mu^3)+\l^2\mu^2\f^3O(e^{2\mu
        |\psi|_{C(\cl{G})}})\\
    \ns
    &\ds\q +6\l^3\mu^4\f^3|\nabla\psi|^4+\l^3\f^3
    O(\mu^3)+\l^2\mu^2\f^3O(e^{2\mu
        |\psi|_{C(\cl{G})}})\\
    \ns&\ds\q+\l^2\mu^2\f^3 O(e^{2\mu
        |\psi|_{C(\cl{G})}})
    +\l^2\f^2 O(\mu^4)\\
    \ns
    &\geq\ds 2\l^3\mu^4\f ^3|\nabla\psi|^4+\l^3\f^3O(\mu^3)
    +\l^2\f^3 O(\mu^2e^{2\mu
        |\psi|_{C(\cl{G})}})+\l^2\f^2 O(\mu^4),
\end{array}
$$
which leads to the second estimate in
(\ref{h6}).\endpf

%%%%%%%%%%%%%%%%%%%%%%%%%%%%%%%%%%%%%%%%%%%%%%%%%%%%%

\subsection{Global Carleman estimate for
backward stochastic parabolic equations}

%%%%%%%%%%%%%%%%%%%%%%%%%%%%%%%%%%%%%%%%%%%%%%%%%%%%%

To study the null and approximate controllability problems for \eqref{heat 1.1}, we introduce the
following backward stochastic parabolic
equation:
\begin{equation}\label{dual heat 1.1}
\left\{
\begin{array}{ll}
    \ds dz+\D z dt =-\big(a_1 z + a_2 Z\big)dt+ Z
    dW(t) \q&\hb{ in }Q,\\
    \ns
    z=0&\hb{ on }\Si,\\
    \ns z(T)=z_T&\hb{ in }G.
\end{array}
\right.
\end{equation}
By Theorem \ref{ch-1-bwell-mild2}, for any $z_T \in L_{\cF_T}^2(\Omega;L^2(G))$,
the system \eqref{dual heat 1.1} admits a unique
weak solution $(z,Z)\in
\big(L^2_{\dbF}(\Omega;C([0,T];L^2(G))) \cap
L_{\dbF}^2(0,T;H_0^1(G))\big) \times
L_{\dbF}^2(0,T; L^2(G))$. Moreover, for any
$t\in [0,T]$,
\begin{equation}\label{20160619e1}
\begin{array}{ll}
    \ds|(z(\cd),Z(\cd))|_{
        \left(L^2_{\dbF}(\Omega;C([0,t];L^2(G)))\cap
        L_{\dbF}^2(0,t;H_0^1(G))\right)\times
        L_{\dbF}^2(0,t;L^2(G))}\le \cC
    |z(t)|_{L_{\cF_t}^2(\Omega;L^2(G))}.
\end{array}
\end{equation}

We also need the following known
result (\cite[p. 4, Lemma 1.1]{FI}).

\begin{lemma}\label{hl1}
For any nonempty open subset $G_1$ of $G$, there
is a  $\psi\in C^4({\overline{G}})$ such
that $\psi>0$ in $G$, $\psi=0$ on $\Gamma$, and
$| \nabla\psi(x)|>0$ for all $x\in \cl{G
    \setminus G _1}$.
\end{lemma}

Choose $\th$ and $\ell$ as the ones in
\eqref{alphad}, and $\psi$ given by Lemma
\ref{hl1} with $G_1$ being any fixed nonempty
open subset of $G$ such that $\cl{G_1}\subset
G_0$. We have the following global Carleman estimate for
\eqref{dual heat 1.1}:

\begin{theorem}\label{c1t4}
There is a constant $\mu_0=\mu_0(G,
G_0,T)>0$ such that for all
$\mu\ge \mu_0$, one can find two constants
$\cC=\cC(\mu)>0$ and $\l_0=\l_0(\mu)>0$ such
that for all $\l\ge \l_0$ and $z_T\in
L^2_{\cF_T}(\Omega;L^2(G))$, the solution
$(z,Z)$  to \eqref{dual heat 1.1}
satisfies that
\begin{equation}\label{h5.2}
    \begin{array}{ll}\ds
        \ds
        \l^3\mu^4\mathbb{E}\int_Q\th^2\f^3z^2dxdt+\l\mu^2\mathbb{E}\int_Q\th^2\f|\n
        z|^2dxdt\\
        \ns \le\ds  \cC \( \l^3\mu^4\mathbb{E}
        \int_{Q_0} \th^2\f^3 z^2dxdt  + \l^2\mu^2\mathbb{E} \int_Q
        \th^2\f^2(a_3z+Z)^2dxdt \).
    \end{array}
\end{equation}
\end{theorem}

{\it Proof}. We shall use Theorem \ref{c1t1} and
Proposition \ref{c1t2}. Integrating the equality
\eqref{c1e2a} on $Q$, taking expectation in both
sides, and noting \eqref{h6}, we conclude that
\begin{eqnarray}\label{6h10}
    && \3n\3n 2\mathbb{E}\int_Q\th\big(\D w + \cA
    w\big)\big(dz+\D z dt \big)dx-2\mathbb{E}\int_Q
    \div (dw\nabla w ) dx\nonumber\\
    &&\3n\3n \q+ 2\mathbb{E}\int_Q\div\[
    2\big(\nabla\ell\cd \nabla w\big) \nabla w -
    |\nabla w |^2\nabla\ell  - \Psi
    w\nabla w + \(\cA \nabla\ell +\frac{\nabla\Psi}{2}\)w^2\]dxdt\nonumber\\
    &&\3n\3n  \geq 2 \mathbb{E}\int_Q\big[\f
    \big(\l\mu^2|\n\psi|^2+\l O(\mu)\big)|\n
    w|^2+\f^3\big(\l^3\mu^4\f^3|\n
    \psi|^4+\l^3\f^3O(\mu^3)\\
    && \qq\q+\l^2\f^2O(\mu^4)\! +\!\l^2\mu^2\f^3O(e^{2\mu
        |\psi|_{C(\cl{G})}})\big)w^2\big] dxdt
    \!+\!2\mathbb{E}\! \int_Q |\D w\! +\! \cA
    w |^2dxdt\nonumber\\
    &&\3n\3n \q +\mathbb{E}\int_Q\th^2 |d\nabla z
    +dz\nabla\ell|^2 dx- \mathbb{E} \int_Q
    \th^2\cA(dz)^2dx.\nonumber
\end{eqnarray}

Recall that $\nu\left(x\right)=(\nu^{1}(x),
\cdots, $ $\nu^{n}(x))$ stands for the unit outward normal vector (of $G$) at $x \in
\G $. Noting that $z=0$ on $\Si$, we have that
$$
w =0 \; \mbox{ and } \;w_{x_j} =\th v_{x_j} =\th \frac{\pa
    z}{\pa\nu}\nu^{j} \q\mbox{  on } \Si.
$$
Similarly, by Lemma \ref{hl1}, we get
$$
\ell_{x_j}=\l\mu\f\psi_{x_j}=\l\mu\f\frac{\pa\psi}{\pa\nu}\nu^j\q\hb{ and }\q \frac{\pa\psi}{\pa\nu}< 0\q \hb{ on }\Si.
$$
Therefore, utilizing integration by parts, we obtain that
\begin{equation}\label{h12}
    \ds
    -2\mathbb{E}\int_Q
    \div (dw\nabla w ) dx=-2\mathbb{E}\int_{\Si}\sum_{j=1}^n
    w_{x_j}\nu^{j}dwdx =-2\mathbb{E}\int_{\Si}\th \frac{\pa
        z}{\pa\nu}dwdx=0
\end{equation}
and that
\begin{eqnarray}\label{h12.1}
    &&\3n\3n\3n\ds
    2\mathbb{E}\int_Q\div\[
    2\big(\nabla\ell\cd \nabla w\big) \nabla w -
    |\nabla w |^2\nabla\ell  - \Psi
    w\nabla w + \(\cA \nabla\ell +\frac{\nabla\Psi}{2}\)w^2\]dxdt\nonumber\\
    &&\3n\3n\3n
    =2\mathbb{E}\int_{\Si}\[2\(\sum_{j=1}^n \l\mu\f\frac{\pa\psi}{\pa\nu}\nu^j \l\mu\f\frac{\pa z}{\pa\nu}\nu^j\)\th \frac{\pa
        z}{\pa\nu}\nu^{j} - \(\sum_{j=1}^n \th^2 \|\frac{\pa
        z}{\pa\nu}\|^2|\nu^{j}|^2\)\l\mu\f\frac{\pa\psi}{\pa\nu} \]dxdt \nonumber\\
    &&
    \3n\3n\3n=2\l\mu\mathbb{E}\int_{\Si}\th^2\f\frac{\pa\psi}{\pa\nu}\(\frac{\pa
        z}{\pa\nu}\)^2
    d\G
    dt\le0.
\end{eqnarray}

It follows from \eqref{dual heat 1.1} that
\begin{eqnarray}\label{6h11}
    && 2\mathbb{E}\int_Q\th\big(\D w + \cA
    w\big)\big(dz+\D z dt \big)dx\nonumber\\
    &&=2\mathbb{E}\int_Q\th\big(\D w + \cA
    w\big)\big[-\big(a_1 z + a_2
    Z\big) dt+ZdW(t)\big]dx
    \nonumber\\&& =-2\mathbb{E}\int_Q\th\big(\D w + \cA
    w\big)(a_1 z + a_2
    Z\big)dxdt\\
    &&\le \mathbb{E}\int_Q\big(\D w + \cA
    w\big)^2dtdx+\mathbb{E}\int_Q\th^2
    (a_1 z + a_2
    Z\big)^2dxdt.\nonumber
\end{eqnarray}
By \eqref{6h10}--\eqref{6h11},  we have
that
\begin{equation}\label{6h16}
    \begin{array}{ll}
        \displaystyle 2\mathbb{E}\int_Q\big[\f
        \big(\l\mu^2|\n\psi|^2+\l O(\mu)\big)|\n
        w|^2\\
        \ns \displaystyle\qq\q
        +\f^3\big(\l^3\mu^4|\n
        \psi|^4+\l^3O(\mu^3)+\l^2\mu^2O(e^{2\mu
            |\psi|_{C(\cl{G})}})+\l^2\f^{-1}O(\mu^4) \big)w^2\big]dxdt\\
        \ns\ds \leq \mathbb{E}\int_Q\th^2\big[\big(a_1 z + a_2
        Z\big)^2+2\cA (Z+a_3z)^2+2a_3^2\cA z^2\big]
        dxdt.
    \end{array}
\end{equation}

Choose a cut-off function $\zeta\in
C_0^{\infty}(G_0;[0,1])$ so that $\zeta\equiv 1$
in $G_1$. By
$$ d(\th^2\f h^2)=h^2(\th^2\f)_tdt+2\th^2\f h
dh+\th^2\f(dh)^2,$$
recalling $$\lim_{t\to0^+}\f(t,\cd) =\lim_{t\to
    T^-}\f(t,\cd)\equiv 0$$ and using \eqref{dual heat 1.1},
we find that
$$
\begin{array}{ll} 0\3n
    &\ds=\mathbb{E}\int_{Q_0}\th^2\big[\zeta^2z^2(\f_t+2\l\f\eta_t)
    +2\zeta^2\f|\nabla z|^2+2\mu \zeta^2\f(1
    +2\l\f)z\nabla z\cd\nabla\psi \\
    \ns &\ds\qq\qq\q+4\zeta\f z\nabla
    z\cd\nabla\zeta +2\zeta^2\f fz+\zeta^2\f
    Z^2\big]dxdt.
\end{array}
$$
Therefore, for any $\e>0$, one has
\begin{equation*}\label{nh21}
    \begin{array}{ll}
        \ds 2\mathbb{E}\int_{Q_0}\th^2\zeta^2\f|\nabla z|^2dxdt+\mathbb{E}\int_{Q_0}\th^2\zeta^2\f Z^2dxdt\\
        \ns \ds\le
        \e\mathbb{E}\int_{Q_0}\th^2\zeta^2\f|\n
        z|^2dxdt+ \frac{\cC}{\e}\mathbb{E}\int_{Q_0}
        \th^2\[\frac{1}{\l^2\mu^2}\big(a_1 z + a_2
        Z\big)^2+\l^2\mu^2\f^3z^2\]dxdt.
    \end{array}
\end{equation*}
This yields that
\begin{equation}\label{nh22}
    \begin{array}{ll}\ds
        \mathbb{E}\int_0^T\int_{G_1}\th^2\f|\n
        z|^2dxdt
        \le \cC\mathbb{E}\int_{Q_0}
        \th^2\[\frac{1}{\l^2\mu^2}\big(a_1 z + a_2
        Z\big)^2+\l^2\mu^2\f^3z^2\]dxdt.
    \end{array}
\end{equation}

From \eqref{6h16} and \eqref{nh22}, we conclude
that there is a $\mu_0>0$ such that for all
$\mu\ge \mu_0$, one can find a constant
$\l_0=\l_0(\mu)$ so that for any $\l\ge \l_0$,
the desired estimate \eqref{h5.2} holds. This
completes the proof of Theorem \ref{c1t4}.\endpf

As a consequence of Theorem \ref{c1t4}, we have the following
observability estimate and unique continuation for \eqref{dual heat 1.1}:
\begin{corollary}\label{b heat ob1}
{\rm 1)} Solutions $(z,Z)$ to the system \eqref{dual heat
    1.1} satisfy
\begin{equation}\label{b heat obser1}
    \begin{array}{ll}\ds
        |z(0)|_{L_{\cF_0}^2(\Omega;L^2(G))}\le \ds \cC
        \big(|\chi_{G_0}z|_{L^2_{\dbF}(0,T;L^2(G_0))}+|a_3z+Z|_{L^2_{\dbF}(0,T;L^2(G))}\big),\\
        \ns\ds\qq\qq\qq\qq\qq\qq\qq\qq\qq\q\forall\;z_T\in L^2_{\cF_T}(\Omega;L^2(G)).
    \end{array}
\end{equation}

{\rm 2)} If for some $z_T\in L^2_{\cF_T}(\Omega;L^2(G))$,
the corresponding solution  $(z,Z)$ satisfies that $z=0$ in $Q_0$ and $a_3z+Z=0$ in $Q$, a.s., then $z_T=0$ in $G$, a.s.
\end{corollary}

{\it Proof}. We only prove the conclusion 1).
Applying Theorem \ref{c1t4} to the equation
\eqref{dual heat 1.1}, and choosing $\mu=\mu_0$
and $\l=\cC$, from \eqref{h5.2}, we deduce that,
\begin{equation}\label{zt7}
    \mathbb{E}\int_Q\th^2\f^3 z^2dxdt\le \cC
    \[\mathbb{E}\int_{Q_0}\th^2\f^3 z^2dxdt+
    \mathbb{E}\int_Q\th^2\f^2(a_3z+Z)^2dxdt\].
\end{equation}
Recalling \eqref{alphad}, it follows from
\eqref{zt7} that
\begin{eqnarray}\label{zhoppo07}
    &&\3n\3n\3n
    \mathbb{E}\int_{T/4}^{3T/4}\int_G z^2dxdt\nonumber\\
    &&\3n\3n\3n\le \cC \frac{\ds\sup_{(t,x)\in
            Q}\(\th^2(t,x)\f^3(t,x)\! +\!\th^2(t,x)\f^2(t,x)
        \)}{\ds\inf_{x\in G}\(\th^2(T/4,x)\f^3(T/2,x)
        \)}\[\mathbb{E}\!\int_{Q_0}\!
    z^2dxdt + \mathbb{E}\!\int_Q\!(a_3z\!+\!Z)^2dxdt\]\nonumber\\
    &&\3n\3n\3n\le \cC \[\mathbb{E}\int_{Q_0} z^2dxdt+
    \mathbb{E}\int_Q(a_3z+Z)^2dxdt\].
\end{eqnarray}
From \eqref{20160619e1}, it follows that
\begin{equation}\label{lk---j12}
    \mathbb{E}\int_G z^2(0)dx\le \cC\mathbb{E}\int_G
    z^2(t)dx,\qq\forall\;t\in [0,T].
\end{equation}
Combining  (\ref{zhoppo07}) and
(\ref{lk---j12}), we conclude that, the solution
$(z,Z)$ to the equation \eqref{dual heat 1.1}
satisfies \eqref{b heat obser1}.
\endpf

%%%%%%%%%%%%%%%%%%%%%%%%%%%%%%%%%%%%%%%%%%%%%%%%%%%%%

\subsection{Null and approximate controllability of stochastic parabolic equations}

%%%%%%%%%%%%%%%%%%%%%%%%%%%%%%%%%%%%%%%%%%%%%%%%%%%%%

This subsection is addressed to a proof of Theorem
\ref{f heat control1}.

{\it Proof}.[Proof of  Theorem \ref{f heat
control1}] Similarly to the proof of  Theorem
\ref{exact th}], we use the classical duality
argument.

We first prove that \eqref{heat 1.1} is null controllable at any time $T>0$.
Let us
introduce a linear subspace of
$L_{\dbF}^2(0,T;L^2(G_0))\times
L_{\dbF}^2(0,T;L^2(G))$ as follows:
$$
\begin{array}{ll}\ds
    \cY\deq\Big\{\big(\chi_{G_0} z,
    a_3z+Z\big)\;\Big|\;(z,Z)\hb{ solves the equation
    }\eqref{dual heat 1.1}\mbox{ with some } \\
    \ns\ds \hspace{3.9cm} z_T\in
    L^2_{\cF_T}(\Om;L^2(G))\Big\}
\end{array}
$$
and define a linear functional $F$ on $\cY$ as
follows:
$$
F(\chi_{G_0} z,
a_3z+Z)= -  \int_G
y_0 z(0) dx.
$$
By the conclusion 1) in Corollary \ref{b heat ob1}, we see that $F$
is a bounded linear functional on $\cY$. By
means of Theorem \ref{Hahn-Banach}, $F$ can be
extended to be a bounded linear functional on
the space
$L_{\dbF}^2(0,T;L^2(G_0))\times
L_{\dbF}^2(0,T;L^2(G))$. For
simplicity, we still use $F$ to denote this
extension. By Theorem \ref{1t10s},  there exists $
(u,
v)\in L_{\dbF}^2(0,T;L^2(G_0))\times
L_{\dbF}^2(0,T;L^2(G))
$
such that
\begin{equation}\label{con eq0-1}
    \begin{array}{ll}\ds
        -  \int_G
        y_0 z(0) dx = F(\chi_{G_0} z,
        a_3z+Z)\\
        \ns\ds =  \mE\int_0^T\int_{G_0} zu
        dxdt + \mE\int_0^T\int_G
        v(a_3z+ Z) dxdt.
    \end{array}
\end{equation}
We claim that the above obtained $u$ and
$v$ are the desired controls. Indeed, by
It\^o's formula and
integration by parts, we have
\begin{equation}\label{con eq2-1}
    \begin{array}{ll}\ds
        \mE\int_G y(T)z_T dx-  \int_G
        y_0 z(0) dx\\
        \ns\ds
        = -\mE\int_0^T\int_{G_0} zu
        dx dt + \mE\int_0^T\int_G v
        (a_3z+ Z) dxdt,
    \end{array}
\end{equation}
where $z$ is the solution to \eqref{dual heat 1.1} with $z_T=\eta$ and $y$ is the state of \eqref{heat 1.1}.
From   \eqref{con eq0-1} and  \eqref{con eq2-1}, we
see that
\begin{equation}\label{con eq3-1}
    \mE\int_G y(T)z_T dx=0.
\end{equation}
Since $z_T$ can be an arbitrary element in
$L^2_{\cF_T}(\Om;L^2(G))$, from the
equality \eqref{con eq3-1}, we conclude that
$y(T)=0$ in $L^2(G)$, a.s. This
concludes the null controllability of  \eqref{heat 1.1}.

\ms

Next, we  prove that \eqref{heat 1.1} is approximately controllable at time $T$.  It suffices to show that the
set
$$
\begin{array}{ll}\ds
    A_T\deq \big\{y(T)\;\big|\; y \mbox{ is the state of \eqref{heat 1.1} with some controls }\\
    \ns\ds \hspace{2cm}(u,v)\in L_{\dbF}^2(0,T;L^2(G_0))\times
    L_{\dbF}^2(0,T;L^2(G))\big\}
\end{array}
$$
is dense in
$L^2_{\cF_T}(\Om;L^2(G))$. Let us prove this by the contradiction argument.  Assume that there was a nonzero $\eta\in L^2_{\cF_T}(\Om;L^2(G))$ such that
$$\ds\mathbb{E}\int_G
y(T)\eta dx=0,\qq \forall\,y(T)\in A_T.
$$
Then, by It\^o's formula and
integration by parts,
we would obtain that
\begin{equation}\label{app1}
    \mathbb{E}\int_G y(T)\eta dx-\int_G y_0z(0) dx = \mathbb{E}\int_0^T\int_{G_0} zudxdt + \mathbb{E}\int_0^T\int_{G} (a_3z+Z)vdxdt,
\end{equation}
where $z$ is the solution to \eqref{dual heat 1.1} with $z_T=\eta$ and $y$ is the state of \eqref{heat 1.1}.
Hence,
\begin{equation}
    \label{12.30-eq1}
    \begin{array}{ll}\ds
        -\int_G y_0z(0) dx =\mathbb{E}\int_0^T\int_{G_0} zudxdt + \mathbb{E}\int_0^T\int_{G} (a_3z+Z)vdxdt,\\
        \ns\ds\qq\qq\qq\q \forall\, (u,v)\in L_{\dbF}^2(0,T;L^2(G_0))\times
        L_{\dbF}^2(0,T;L^2(G)).
    \end{array}
\end{equation}
Noting that the left hand side of \eqref{12.30-eq1} is independent of $(u,v)$. By choosing $(u,v)=(0,0)$, we find that $\int_G y_0z(0) dx=0$. Thus,
$$
\begin{array}{ll}\ds
    0 =\mathbb{E}\int_0^T\int_{G_0} zudxdt + \mathbb{E}\int_0^T\int_{G} (a_3z+Z)vdxdt,\\
    \ns\ds\qq\qq\qq\q \forall\, (u,v)\in L_{\dbF}^2(0,T;L^2(G_0))\times
    L_{\dbF}^2(0,T;L^2(G)),
\end{array}
$$
which yields that $z=0$  in $G_0\t (0,T)$,
a.s. and $a_3z+Z$ in $G\t (0,T)$,
a.s. By the conclusion 2) in Corollary \ref{b heat ob1}, we then have $\eta=0$, a
contradiction.
\endpf
%

%%%%%%%%%%%%%%%%%%%%%%%%%%%%%%%%%%%%%%%%%%%%%%%%%%%%%

\section{Controllability of stochastic partial differential equations III:  stochastic hyperbolic equations}
\label{s-h}

%%%%%%%%%%%%%%%%%%%%%%%%%%%%%%%%%%%%%%%%%%%%%%%%%%%%

Similarly to parabolic equations, hyperbolic
equations (including particularly the wave
equations) are another class of typical PDEs. In
the deterministic setting, the controllability
problems for hyperbolic equations is extensively
studied. As we shall see in this section,  the
usual stochastic hyperbolic equation, i.e., the
classic hyperbolic equation perturbed by a term
of It\^o's integral, is not exactly controllable
even if the controls are effective everywhere in
both drift and diffusion terms, which differs
dramatically from its deterministic counterpart.
This section is based on \cite{LZ6}.

%%%%%%%%%%%%%%%%%%%%%%%%%%%%%%%%%%%%%%%%%%%%%%%%%%%%

\subsection{Formulation of the problem}

%%%%%%%%%%%%%%%%%%%%%%%%%%%%%%%%%%%%%%%%%%%%%%%%%%%%

Let $T > 0$, $G$ be given as in Subsection \ref{sec-mpr-6.1} and $\G_0$ be a subset of $\G$.  Put $Q$ and $\Si$ as in \eqref{2021-1-16e1}, and
$\Si_0 \deq  (0,T) \t \G_0$.

Similarly to the previous two sections, we begin with the following controlled (deterministic)
hyperbolic equation:
\begin{equation}\label{system0001}
\left\{
\begin{array}{ll}
    \ds y_{tt}-
    \D y =a_1^0 y&\mbox{ in }Q,\\
    \ns\ds y=\chi_{\Si_0}h&\mbox{ on }\Si,\\
    \ns\ds y(0)=y_0,\q y_t(0)=y_1&\mbox{ in }G.
\end{array}
\right.
\end{equation}
Here $(y_0,y_1)\in
L^2(G)\times H^{-1}(G)$,  $a_1^0\in
L^\infty(Q)$, $(y,y_t)$ is the {\it state variable}, and
$h\in L^2(\Si_0)$ is the
{\it control variable}. It is well-known that the system (\ref{system0001}) admits a unique transposition solution $y\in C([0,T];L^2(G))\cap C^1([0,T];H^{-1}(G))$ (e.g., \cite{Lions3}).

\begin{definition}
We say the system \eqref{system0001} is {\it exactly controllable at time $T$} if for any given $(y_0,y_1),$ $(\tilde y_0,\tilde y_1)\in
L^2(G)\times H^{-1}(G)$, one can find a control $h\in L^2(\Si_0)$ such that the solution to (\ref{system0001}) satisfying $y(T)=\tilde y_0$ and $y_t(T)=\tilde y_1$.
\end{definition}

To show the exact controllability of \eqref{system0001}, people introduce its dual
equation as follows
\begin{eqnarray}\label{system0001.1}
\left\{
\begin{array}{lll}
    \ds  z_{tt}- \D z  =  a_1^0z\quad
    &\mbox{ in } Q,\\
    \ns\ds z = 0 & \mbox{ on } \Si, \\
    \ns\ds z(T) = z_0,\; z_t(T)=z_1&\mbox{ in } G,
\end{array}
\right.
\end{eqnarray}
where $(z_0,z_1)\in
H_0^1(G)\times L^2(G)$.

By means of the standard duality argument, it is easy to show the following result (e.g., \cite{Lions0}).
\begin{proposition}\label{20160719pro1.1}
The equation \eqref{system0001} is exactly controllable at time $T$ if and only if solutions to the equation
\eqref{system0001.1} satisfy the following observability estimate:
\begin{equation}\label{20160719e1.1}
    |(z_0,z_1)|_{H_0^1(G)\times L^2(G)}\le \cC\Big|\frac{\pa z}{\pa\nu}\Big|_{L^2(\Si_0)},\qq\forall\;(z_0,z_1)\in
    H_0^1(G)\times L^2(G).
\end{equation}
\end{proposition}

It is known that (e.g., \cite{Zhangxu2, Zu1}), under some assumptions on $(T, G,\G_0)$(for
example,   $\G_0=\{ x\in \G \;|\;
(x-x_0)\cdot \nu(x)>0\}$ for some
$x_0\in\mathbb{R}^n$ and
$T>2\max_{x\in \cl G}|x-x_0|$), the inequality \eqref{20160719e1.1} holds. As a result,  the system (\ref{system0001}) is exactly controllable.
In the rest of this section, we shall see that the controllability property of stochastic
hyperbolic equation is quite different from its deterministic counterpart.

Now, let us consider the following controlled stochastic
hyperbolic equation:
\begin{equation}\label{system1}
\left\{
\begin{array}{ll}
    \ds dy_t- \D y dt=(a_1 y+g_1)dt + (a_2y+g_2)
    dW(t)&\mbox{ in }Q,\\
    \ns\ds y= h&\mbox{ on }\Si,\\
    \ns\ds y(0)=y_0,\q y_t(0)=y_1&\mbox{ in }G,
\end{array}
\right.
\end{equation}
where $(y_0,y_1)\in
L^2(G)\times H^{-1}(G)$, $a_1,a_2 \in
L^\infty_\dbF(0,T;L^\infty(G))$, $(y,y_t)$ are the
{\it state variable},  $g_1,g_2\in
L^\infty_\dbF(0,T;H^{-1}(G))$ and $h\in
L^2_\dbF(0,T;L^2(\G))$ are the {\it control variables}. As we shall see in Subsection \ref{ssec-well}, the equation \eqref{system1} admits a unique transposition solution $y\in
C_{\dbF}([0,T];L^2(\Om; L^2(G)))\cap
C_{\dbF}^1([0,T];L^2(\Om;H^{-1}(G))) $.

\begin{definition}\label{exact def}
The system \eqref{system1} is called {\it exactly
    controllable at time $T$} if for any $(y_0,y_1)\in
L^2(G)\times H^{-1}(G)$ and $(\tilde y_0,\tilde y_1)\in
L^2_{\cF_T}(\Om;L^2(G))\times
L^2_{\cF_T}(\Om;H^{-1}(G))$, there is  a triple
of controls $(g_1,g_2,h)\in L^2_{\dbF}(0,T;
H^{-1}(G))\times L^2_{\dbF}(0,T;H^{-1}(G))\times
L^2_\dbF(0,T;L^2(\G))$ such that the
corresponding solution $y$ to the system
\eqref{system1} satisfies that $(y(T),y_t(T)) =
(\tilde y_0,\tilde y_1)$, a.s.
\end{definition}

Since three controls are introduced in
\eqref{system1}, one may guess that the desired
exact controllability should be trivially
correct. Surprisingly, we have the following
negative result.

\begin{theorem}\label{th-non-con}
The system \eqref{system1} is not exactly
controllable  at any time $T>0$.
\end{theorem}

Clearly, the controls in \eqref{system1} are
the strongest ones people can introduce. Nevertheless, the result in Theorem \ref{th-non-con} differs significantly from the above mentioned
controllability property of deterministic hyperbolic
equations. Since \eqref{system1} is a
generalization of the classical hyperbolic equation to
the stochastic setting, from the viewpoint of
control theory, we believe that some key feature
has been ignored in the derivation of the
equation \eqref{system1}. A proof of Theorem \ref{th-non-con} and some further discussions will be given in Subsection \ref{20210117subs6.3}.

%%%%%%%%%%%%%%%%%%%%%%%%%%%%%%%%%%%%%%%%%%%%%%%%%%%

\subsection{Backward stochastic hyperbolic equation and hidden regularity}

%%%%%%%%%%%%%%%%%%%%%%%%%%%%%%%%%%%%%%%%%%%%%%%%%%%

System  \eqref{system1} is a nonhomogeneous
boundary value problem. Similarly to the control system \eqref{ch-7-csystem1}, its solution  is
understood in the sense of transposition. To
this end, we introduce the following
``reference" equation:
\begin{equation}\label{bsystem1}
\left\{
\begin{array}{ll}
    \ds dz=\hat zdt   +Z dW(t) &\mbox{ in } Q_\tau\deq(0,\tau)\times
    G,\\
    \ns\ds d\hat z -  \D z dt
    = (a_1  z + a_2Z )dt  +  \widehat Z dW(t)  &\mbox{ in } Q_\tau,\\
    \ns\ds z = 0 &\mbox{ on }  \Si_\tau\deq(0,\tau)\times \G,\\
    \ns\ds z(\tau) = z^{\tau},\q \hat z(\tau) = \hat
    z^{\tau} &\mbox{ in } G,
\end{array}
\right.
\end{equation}
where $\tau\in (0,T]$ and $(z^{\tau},
\hat z^{\tau}) \in L^2_{\cF_\tau}(\Om;H_0^1(G))\times
L^2_{\cF_\tau}(\Om; L^2(G))$.
By Theorem \ref{ch-1-well-bmild}, the system
\eqref{bsystem1} admits a unique solution
$
(z,\hat z,Z,\widehat Z)\in
L^2_\dbF(\Om;C([0,\tau];H_0^1(G)))\times L^2_\dbF(\Om;C([0,\tau];L^2(G)))\times
L^2_\dbF(0,\tau;H^1_0(G))\times
L^2_\dbF(0,\tau;L^2(G)).
$
Moreover,
\begin{eqnarray}\label{best1.1}
\begin{array}{ll}\ds
    |z|_{L^2_\dbF(\Om;C([0,\tau]; H_0^1(G)))}+|\hat
    z|_{L^2_\dbF(\Om;C([0,\tau];L^2(G)))}   \\
    \ns\ds\q+
    |Z|_{L^2_\dbF(0,\tau;H^1_0(G))}  + |\widehat
    Z|_{L^2_\dbF(0,\tau;L^2(G))}\\
    \ns\ds  \leq \cC
    \big(|z^{\tau}|_{L^2_{\cF_\tau}(\Om;H_0^1(G))}+|\hat
    z^{\tau}|_{L^2_{\cF_\tau}(\Om;L^2(G))}\big).
\end{array}
\end{eqnarray}

Similarly to Section \ref{s-t}, we need to
establish a hidden regularity for solutions to
\eqref{bsystem1}.

\begin{proposition}\label{ch-h-prop-hid}
For any $(z^{\tau},\hat z^{\tau})\in
L^2_{\cF_\tau}(\Om;H_0^1(G))\times
L^2_{\cF_\tau}(\Om;L^2(G))$, the solution
$(z,\hat z,Z,\widehat Z)$ to \eqref{bsystem1}
satisfies $\frac{\pa z}{\pa\nu}\big|_{\G}\in
L^2_{\dbF}(0,\tau;L^2(\G))$. Furthermore,
\begin{equation}\label{hid-eq1}
    \Big| \frac{\pa z}{\pa\nu}
    \Big|_{L^2_\dbF(0,\tau;L^2(\G))} \leq \cC
    \big(|z^{\tau}|_{L^2_{\cF_\tau}(\Om;H_0^1(G))}+|\hat
    z^{\tau}|_{L^2_{\cF_\tau}(\Om;L^2(G))}\big),
\end{equation}
where the constant $\cC$ is independent of
$\tau$.
\end{proposition}

{\it Proof}. One can find a vector field
$\Xi=(\Xi^1,\cdots,\Xi^n)\in
C^1(\mathbb{R}^n;\mathbb{R}^n)$ such that
$\Xi(x)=\nu(x)$ for $x\in\G$ (See \cite[p.
29]{Lions0})). By It\^o's formula and the first
equation of \eqref{bsystem1}, we have
\begin{equation}\label{8.6-eq3.1}
    \ba{ll}\ds
    d(\hat z\Xi\cdot\nabla z)=d\hat z\Xi\cdot\nabla z +\hat z \Xi\cdot \nabla dz+d\hat z\Xi\cdot \nabla dz\\
    \ns
    \ds=d\hat z\Xi\cdot\nabla z +\hat z \Xi\cdot \nabla \big(\hat zdt + ZdW(t)\big)+d\hat z\Xi\cdot \nabla dz\\
    \ns
    \ds=d\hat z\Xi\cdot\nabla z
    +\frac{1}{2}\big[\div(\hat z^2\Xi)-(\div
    \Xi)\hat z^2\big]dt +\hat z \Xi\cdot \nabla
    ZdW(t)+d\hat z\Xi\cdot \nabla dz. \ea
\end{equation}
It follows from a direct computation that
\begin{equation}\label{8.6-eq3.2}
    \begin{array}{ll}
        \ds \div\big[ 2(\Xi\cd\nabla
        z) \nabla z -  \Xi|\nabla z|^2\big] =  \ds 2 \Big(
        \D z \Xi \cd \nabla z +\sum_{j,k=1}^n\Xi^k_{x_j} z_{x_j} z_{x_k}\Big) -
        (\div
        \Xi)|\nabla z|^2.
    \end{array}
\end{equation}
Combining \eqref{8.6-eq3.1} and \eqref{8.6-eq3.2}, we obtain that
\begin{equation}\label{8.6-eq3}
    \begin{array}{ll}
        \ds - \div\big[ 2(\Xi\cd\nabla z) \nabla z  +
        \Xi \big(\hat z^2 -
        |\nabla z|^2\big) \big] dt\\
        \ns =  \ds 2 \Big[ - d(\hat z \Xi\cd \nabla z)
        +\big(d\hat z - \D z dt\big) \Xi \cd \nabla z -
        \sum_{j,k=1}^n\Xi^k_{x_j} z_{x_j} z_{x_k}
        dt\Big]\\
        \ns\ds\q  - (\div \Xi) \hat z^2dt +  (\div
        \Xi)|\nabla z|^2 dt + 2d\hat z \Xi\cd\nabla dz
        +2 \hat z \Xi\cd\nabla ZdW(t).
    \end{array}
\end{equation}
Integrating \eqref{8.6-eq3} in $Q$,
taking expectation on $\Om$, using the second
equation of \eqref{bsystem1} and noting that $z=0$
on $(0,\tau)\times \G$, we get that
$$
\begin{array}{ll}
    \ds -\mE\int_{\Si_\tau}\Big|\frac{\pa z}{\pa\nu}\Big|^2d\Si_\tau\\
    \ns\ds =  \ds -2 \mE\int_G \hat z^T \Xi\cd
    \nabla z^T dx + 2 \mE\int_G \hat z(0)
    \Xi\cd \nabla z(0)dx + 2
    \int_{Q_\tau}\!\Big[\big(a_1 z + a_2Z \big) \Xi
    \cd \nabla z  \\
    \ns\ds\q - \sum_{j,k=1}^n \Xi^k_{x_j} z_{x_j}
    z_{x_k}  - (\div \Xi) \hat z^2 + (\div \Xi) |\nabla z|^2
    + 2 \widehat Z \Xi\cd\nabla Z \Big]dxdt \deq \cI.
\end{array}
$$
By \eqref{best1.1}, we have
$
|\cI|\leq \cC
\big(|z^{\tau}|_{L^2_{\cF_\tau}(\Om;H_0^1(G))}+|\hat
z^{\tau}|_{L^2_{\cF_\tau}(\Om;L^2(G))}\big),
$
which, together with the above equality, implies
\eqref{hid-eq1}.
\endpf

%%%%%%%%%%%%%%%%%%%%%%%%%%%%%%%%%%%%%%%%%%%%%%%%%%%

\subsection{Well-posedness of the control system}\label{ssec-well}

%%%%%%%%%%%%%%%%%%%%%%%%%%%%%%%%%%%%%%%%%%%%%%%%%%%

We begin with the following notion.

\begin{definition}\label{1-def1}
A stochastic process $y\in
C_{\dbF}([0,T];L^2(\Om; L^2(G)))\cap
C_{\dbF}^1([0,T];L^2(\Om;$ $H^{-1}(G))) $ is called a
{\it transposition solution} to \eqref{system1} if for
any $\tau\!\in\! (0,T]$, $(z^{\tau},\hat
z^{\tau})\!\in\!  L^2_{\cF_\tau}(\Om;H_0^1(G))\times L^2_{\cF_\tau}(\Om;L^2(G))$ and the corresponding solution $(z,\hat z,Z,\widehat Z)$ to \eqref{bsystem1},
it holds that
\begin{eqnarray}\label{def id}
    &&\3n\3n\mE \langle
    y_t(\tau),z^{\tau}\rangle_{H^{-1}(G),H^1_0(G)} -
    \langle
    y_1,z(0)\rangle_{H^{-1}(G),H^1_0(G)}\nonumber\\
    &&  - \mE\langle y(\tau),\hat
    z^{\tau}\rangle_{L^2(G)} + \langle y_0,\hat
    z(0)\rangle_{L^2(G)}
    \\ &&\3n\3n = \mE\int_0^\tau
    \langle g_1,z\rangle_{H^{-1}(G),H_0^1(G)} dt +
    \mE\int_0^\tau \langle
    g_2,Z\rangle_{H^{-1}(G),H_0^1(G)} dt
    -\mE\int_{\Si_\tau}\frac{\pa
        z}{\pa\nu}hd\Si_\tau.
    \nonumber
\end{eqnarray}
\end{definition}
\begin{remark}
Note that, by Proposition \ref{ch-h-prop-hid}, the solution
$(z,\hat z, Z,$ $\widehat Z)$ to \eqref{bsystem1}
satisfies $\frac{\pa z}{\pa\nu}\big|_{\G}\in
L^2_{\dbF}(0,\tau;L^2(\G))$. Hence, the term ``$\mE\int_{\Si_\tau}\frac{\pa
    z}{\pa\nu}hd\Si_\tau$'' in \eqref{def id}  makes sense.
\end{remark}

We have the following well-posedness result for
\eqref{system1}.

\begin{proposition}\label{well posed1-1}
For each $(y_0,y_1)\in L^2(G)\times
H^{-1}(G)$ and $(g_1,g_2,h)\in
L^2_{\dbF}(0,T; H^{-1}(G)) \times
L^2_{\dbF}(0,T;$ $H^{-1}(G))\times
L^2_\dbF(0,T;L^2(\G_{0}))$, the system \eqref{system1} admits a
unique transposition solution $y\in C_{\dbF}([0,T];$ $L^2(\Om;L^2(G)))\cap C_{\dbF}^{1}([0,T];L^2(\Om;H^{-1}(G)))$, and
\begin{equation}\label{well posed est1}
    \begin{array}{ll}\ds
        |y|_{C_{\dbF}([0,T];L^2(\Om;L^2(G)))\cap C_{\dbF}^{1}([0,T];L^2(\Om;H^{-1}(G)))}\\
        \ns\ds \leq \cC \big( |y_0|_{L^2(G)} + |y_1|_{H^{-1}(G)} +
        |g_1|_{L^2_{\dbF}(0,T;H^{-1}(G))}\\
        \ns\ds\qq + |g_2|_{L^2_{\dbF}(0,T;H^{-1}(G))}+
        |h|_{L^2_{\dbF}(0,T;L^2(\G))}\big).
    \end{array}
\end{equation}
\end{proposition}

{\it Proof}. {\bf Uniqueness}. Assume that $y$
and $\tilde y$ are two transposition solutions
of \eqref{system1}. It follows from Definition
\ref{1-def1} that for any $\tau\in (0,T]$ and
$(z^{\tau},\hat z^{\tau})\in
L^2_{\cF_\tau}(\Om;$ $H_0^1(G))\times
L^2_{\cF_\tau}(\Om; L^2(G))$,
\begin{equation}\label{def id1-1}
    \begin{array}{ll}
        \ds \mE \langle
        y_t(\tau),z^{\tau}\rangle_{H^{-1}(G),H^1_0(G)}
        -\! \mE\langle y(\tau),\hat z^{\tau}
        \rangle_{L^2(G)}\\
        \ns\ds = \mE \langle \tilde{
            y}_t(\tau),z^{\tau}\rangle_{H^{-1}(G),H^1_0(G)}
        - \mE\langle \tilde y(\tau),\hat z^{\tau}
        \rangle_{L^2(G)},
    \end{array}
\end{equation}
which implies that
$
\big(y(\tau), y_t(\tau)\big)=\big(\tilde{
    y}(\tau), \tilde y_t(\tau)\big)$,
a.s.,  $\forall\;\tau\in (0,T]$.
Hence, $y=\tilde y\; \mbox{ in
}C_{\dbF}([0,T];L^2(\Om;$ $L^2(G)))\cap
C_{\dbF}^{1}([0,T];L^2(\Om;H^{-1}(G))).$

\vspace{0.2cm}

{\bf Existence}. Since $h\in
L^2_\dbF(0,T;L^2(\G))$, there exists a sequence
$\{h_m\}_{m=1}^\infty\subset C_\dbF^2([0,T];$
$H^{3/2}(\G))$ with $h_m(0)=0$ for all
$m\in\dbN$ such that
\begin{equation}\label{8.31-eq1.1}
    \lim_{m\to\infty} h_m = h \q \mbox{ in
    }L^2_\dbF(0,T;L^2(\G)).
\end{equation}
For each $m\in\dbN$, we can find an $\tilde h_m
\in C_\dbF^2([0,T];H^2(G))$ such that $\tilde
h_m|_{\G}=h_m$ and $\tilde h_m(0)=0$.

Consider the following equation:
\begin{equation}\label{system2n}
    \left\{\!
    \begin{array}{ll}
        \ds d\tilde { y}_{m,t} - \D\tilde y_{m} dt=(a_1
        \tilde y_m  + \zeta_m)dt + [a_2(\tilde y_m +
        \tilde h_{m}) + g_2]
        dW(t)&\mbox{ in }Q, \\
        \ns\ds \tilde y_m =0&\mbox{ on }\Si,\\
        \ns\ds \tilde y_m(0)= y_0,\q \tilde {
            y}_{m,t}(0)=y_1&\mbox{ in }G,
    \end{array}
    \right.
\end{equation}
where $\ds\zeta_m = g_1+\D\tilde h_{m} + a_1
\tilde h_m$. By Theorem \ref{ch-1-well-mild},
the system \eqref{system2n} admits a unique mild
(also weak) solution $\tilde y_m\in
C_\dbF([0,T];L^2(\Om;L^2(G)))\cap
C_\dbF^{1}([0,T];L^2(\Om;H^{-1}(G)))$.

Let $y_m=\tilde y_m+\tilde h_m$.  By It\^o's formula and
integration by parts, we have that
\begin{eqnarray}\label{8.31-eq5}
    && \mE \langle
    y_{m,t}(\tau),z^\tau\rangle_{H^{-1}(G),H^1_0(G)}
    - \langle y_1,z(0)\rangle_{H^{-1}(G),H^1_0(G)}\nonumber\\
    &&\q  - \mE\langle y_{m}(\tau),\hat z^\tau
    \rangle_{L^2(G)} + \langle y_0,\hat
    z(0)\rangle_{L^2(G)}
    \\ && = -
    \mE\int_0^\tau\!\!\langle g_1,z\rangle_{H^{-1}(G),H_0^1(G)}dt +
    \mE\int_0^\tau\!\! \langle
    g_2,Z\rangle_{H^{-1}(G),H_0^1(G)}
    dt -\mE\int_{\Si_\tau}\!\!\frac{\pa
        z}{\pa\nu}h_{m}d\Si_\tau.\nonumber
\end{eqnarray}
Consequently, for any
$m_1,m_2\in\dbN$,
\begin{eqnarray}\label{8.31-eq2}
    &&\mE \langle y_{m_1,t}(\tau) -
    y_{m_2,t}(\tau),z^\tau\rangle_{H^{-1}(G),H^1_0(G)}
    -  \mE\langle y_{m_1}(\tau) - y_{m_2}(\tau),\hat
    z^\tau
    \rangle_{L^2(G)} \nonumber\\
    &&=  -\mE \int_{\Si_\tau}\frac{\pa
        z}{\pa\nu}(h_{m_1} - h_{m_2})d\Si_\tau.
\end{eqnarray}

By Proposition \ref{2.1-cor16}, there is $(z^\tau,\hat z^\tau)\in
L^2_{\cF_\tau}(\Om;H_0^1(G))\times
L^2_{\cF_\tau}(\Om;L^2(G))$  so that
$$
|z^\tau|_{L^2_{\cF_\tau}(\Om;H_0^1(G))}=1,\qq
|\hat z^\tau|_{L^2_{\cF_\tau}(\Om;L^2(G))}=1
$$
and that
\begin{eqnarray}\label{8.31-eq3}
    \begin{array}{ll}\ds
        \mE \langle  y_{m_1,t}(\tau)-
        y_{m_2,t}(\tau),z^\tau\rangle_{H^{-1}(G),H^1_0(G)}
        - \mE\langle y_{m_1}(\tau)-y_{m_2}(\tau),\hat
        z^\tau
        \rangle_{L^2(G)}\\
        \ns\ds \geq
        \frac{1}{2}\big(|y_{m_1}(\tau)-y_{m_2}(\tau)|_{L^2_{\cF_\tau}(\Om;L^2(G))}
        + |y_{m_1,t}(\tau)-
        y_{m_2,t}(\tau)|_{L^2_{\cF_\tau}(\Om;H^{-1}(G))}\big).
    \end{array}
\end{eqnarray}
It follows from \eqref{8.31-eq2},
\eqref{8.31-eq3}  and Proposition
\ref{ch-h-prop-hid} that
$$
\begin{array}{ll}
    \ds
    |y_{m_1}(\tau)-y_{m_2}(\tau)|_{L^2_{\cF_\tau}(\Om;L^2(G))}
    + |y_{m_1,t}(\tau)-
    y_{m_2,t}(\tau)|_{L^2_{\cF_T}(\Om;H^{-1}(G))}
    \\
    \ns\ds \leq  2\Big| \mE\int_{\Si_\tau}\frac{\pa
        z}{\pa\nu}(h_{m_1}-h_{m_2})d\Si_\tau
    \Big|\\
    \ns\ds \leq \cC
    |h_{m_1}-h_{m_2}|_{L^2_\dbF(0,T;L^2(\G))}
    |(z^\tau,\hat
    z^\tau)|_{L^2_{\cF_\tau}(\Omega;H^1_0(G))\times
        L^2_{\cF_\tau}(\Omega;L^2(G))}\\
    \ns\ds \leq \cC
    |h_{m_1}-h_{m_2}|_{L^2_\dbF(0,T;L^2(\G))},
\end{array}
$$
where the constant $\cC$ is independent of
$\tau$. Consequently, it holds that
$$
\begin{array}{ll}
    |y_{m_1}-y_{m_2}|_{C_\dbF([0,T];L^2(\Om;L^2(G)))}
    + |y_{m_1,t}-
    y_{m_2}|_{C_\dbF([0,T];L^2(\Om;H^{-1}(G)))}\\
    \ns\ds \leq \cC
    |h_{m_1,t}-h_{m_2}|_{L^2_\dbF(0,T;L^2(\G))}.
\end{array}
$$
This concludes that $\{y_{m}\}_{m=1}^\infty$ is
a Cauchy sequence in
$C_\dbF([0,T];L^2(\Om;L^2(G)))\cap
C_\dbF^{1}([0,T];L^2(\Om;$ $H^{-1}(G)))$. Denote
by $y$ the limit of $\{y_m\}_{m=1}^\infty$.
Letting $m\to \infty$ in \eqref{8.31-eq5}, we
see that $y$ satisfies \eqref{def id}.
Thus, $y$ is a transposition solution to
\eqref{system1}.

By Proposition \ref{2.1-cor16}, there is $(z^\tau,\hat z^\tau)\in
L^2_{\cF_\tau}(\Om;H_0^1(G))\times
L^2_{\cF_\tau}(\Om;L^2(G))$  such that
$
|z^\tau|_{L^2_{\cF_\tau}(\Om;H_0^1(G))}$ $=1$,
$|\hat z^\tau|_{L^2_{\cF_\tau}(\Om;L^2(G))}=1
$
and
\begin{equation}\label{8.31-eq7}
    \begin{array}{ll}\ds \mE \langle
        y_t(\tau),z^\tau\rangle_{H^{-1}(G),H^1_0(G)}   -
        \mE\langle y(\tau),\hat z^\tau
        \rangle_{L^2(G)}\\
        \ns\ds\geq
        \frac{1}{2}\big(|y(\tau)|_{L^2_{\cF_\tau}(\Om;L^2(G))}
        + |\hat
        y(\tau)|_{L^2_{\cF_\tau}(\Om;H^{-1}(G))}\big).
    \end{array}
\end{equation}
Combining \eqref{def id}, \eqref{8.31-eq7} and
Proposition \ref{ch-h-prop-hid}, we obtain that, for any $\tau\in(0,T]$,
$$
\begin{array}{ll}
    \ds  |y(\tau)|_{L^2_{\cF_\tau}(\Om;L^2(G))} +
    |y_t(\tau)|_{L^2_{\cF_\tau}(\Om;H^{-1}(G))}
    \\
    \ns\ds \leq 2\(\big|\langle y_1,z(0)\rangle_{H^{-1}(G),H^1_0(G)}\big| +
    \big|\langle y_0,\hat z(0)\rangle_{L^2(G)}\big|+
    \Big|\mE\int_0^\tau \langle g_1,z\rangle_{H^{-1}(G),H_0^1(G)}dt\Big|\\
    \ns\ds\qq + \Big|\mE\int_0^\tau \langle
    g_2,Z\rangle_{H^{-1}(G),H_0^1(G)} dt\Big|+ \Big|
    \mE\int_{\Si_\tau}\frac{\pa z}{\pa\nu}hd\Si_\tau
    \Big|\) \\
    \ns\ds \leq \cC \big( |y_0|_{L^2(G)} + |y_1|_{H^{-1}(G)} +
    |g_1|_{L^2_{\dbF}(0,T;H^{-1}(G))} +
    |g_2|_{L^2_{\dbF}(0,T;H^{-1}(G))} \\
    \ns\ds\qq + |h|_{L^2_{\dbF}(0,T;L^2(\G))}\big)
    \times |(z^\tau,\hat
    z^\tau)|_{L^2_{\cF_\tau}(\Omega;H^1_0(G))\times
        L^2_{\cF_\tau}(\Omega;L^2(G))}\\
    \ns\ds \leq \cC \big( |y_0|_{L^2(G)} + |y_1|_{H^{-1}(G)} +
    |g_1|_{L^2_{\dbF}(0,T;H^{-1}(G))} +
    |g_2|_{L^2_{\dbF}(0,T;H^{-1}(G))} \\
    \ns\ds\qq + |h|_{L^2_{\dbF}(0,T;L^2(\G))}\big),
\end{array}
$$
where the constant $\cC$ is independent of
$\tau$. Therefore, we obtain the estimate \eqref{well posed est1}.
This completes the proof of Proposition
\ref{well posed1-1}.
\endpf

\subsection{Lack of exact controllability and a refined stochastic wave equation}\label{20210117subs6.3}

First of all, let us prove the negative controllability result in Theorem \ref{th-non-con}:

{\it Proof}.[Proof of Theorem \ref{th-non-con}]
We use the contradiction argument.  Choose
$\psi\in H_0^1(G)$ satisfying
$|\psi|_{L^2(G)}=1$ and let $\tilde y_0=\xi
\psi$, where $\xi$ is given in Lemma
\ref{6.18-cor1}. Assume that \eqref{system1} was
exactly controllable for some time $T>0$. Then,
for any $y_0\in L^2(G)$, we would find a triple
of controls $(g_1,g_2,h)\in
L^2_{\dbF}(0,T;H^{-1}(G))\times
L^2_{\dbF}(0,T;H^{-1}(G))\times
L^2_\dbF(0,T;L^2(\G))$ such that the
corresponding solution $ y\in C_{\dbF}([0,T];$
$L^2(\Om;L^2(G)))\cap C_{\dbF}^1([0,T];L^2(\Om;$
$H^{-1}(G))) $ to the equation \eqref{system1}
satisfies that $y(T)=\tilde y_0$. Clearly,
$$
\int_G \tilde y_0 \psi dx-\int_G  y_0 \psi dx =
\int_0^T\langle
y_t,\psi\rangle_{H^{-1}(G),H_0^1(G)}dt,
$$
which leads to
$$
\xi = \int_G  y_0 \psi dx + \int_0^T\langle
y_t,\psi\rangle_{H^{-1}(G),H_0^1(G)}dt.
$$
This contradicts Lemma \ref{6.18-cor1}.
\endpf

Then, as we did in \cite{LZ6}, motivated by Theorem \ref{th-non-con}, we
propose below a refined model to describe the
DNA molecule considered in Example \ref{ch-4-ex2}. For this purpose, we partially
employ a dynamical theory of Brownian motions,
developed in \cite{Nelson}, to describe the
motion of a particle perturbed by random forces.

According to \cite[Chapter 11]{Nelson}, we may
suppose that
\begin{equation}\label{8.6-eq1.1}
y(t,x)=\int_0^t \tilde y(s,x)ds + \int_0^t
F(s,x,y(s))dW(s).
\end{equation}
Here $\tilde y(\cd,\cd)$ is the expected
velocity, $F(\cd,\cd,\cd)$ is the random
perturbation from the fluid molecule. When $y$
is small, one can assume that $F(\cd,\cd,\cd)$
is linear in the third argument, i.e., for a suitable $b_1(\cd,\cd)$,
\begin{equation}\label{3.18-eq2}
F(s,x,y(t,x))=b_1 (t,x)y(t,x).
\end{equation}

Formally, the acceleration at position $x$ along the
string at time $t$ is $\tilde y_{t}(t,x)$.  By
Newton's second law, it follows that
$$
\tilde y_{t}(t,x) = F_1(t,x) + F_2(t,x) +
F_3(t,x) + F_4(t,x).
$$
Similar to \eqref{1.17-eq2},
we obtain then that
\begin{equation}\label{1.17-eq2.1}
d\tilde y (t,x) = y_{xx}(t,x)dt + F_2(t,x)dt +
F_3(t,x)dt + a_4 (t,x)y(t,x)dW(t).
\end{equation}
Combining \eqref{8.6-eq1.1}, \eqref{3.18-eq2}
and \eqref{1.17-eq2.1}, we arrive at the following modified version of \eqref{1.17-eq2}:
\bel{8.600eq2.1}\left\{
\ba{ll} dy = \tilde y dt + b_1 y dW(t)
&\mbox{ in }(0,T)\times
(0,L),\\
\ns\ds d\tilde y= (y_{xx}+a_1y_x+a_2y_t+a_3y ) dt
+ a_4 y dW(t) &\mbox{ in }(0,T)\times (0,L).
\ea\right.
\ee
Then, similarly to \eqref{ch-4-ex2-eq1-con}, we
obtain the following control system:
$$
\left\{
\begin{array}{ll}\ds
dy = \tilde y dt + (b_1 y+u) dW(t) &\mbox{ in
}(0,T)\times (0,L),
\\ \ns\ds
d\tilde y   = (y_{xx}+a_1y_x+a_2y_t+a_3y) dt + (a_4y+v) dW(t)
&\mbox{ in } (0,T)\times (0,L),\\
\ns\ds y =f_1 &\mbox{ on } (0,T)\times \{0\},\\
\ns\ds y =f_2 &\mbox{ on } (0,T)\times \{L\},
\\
\ns\ds y(0) = y_0,\q \tilde y(0)=y_1 &\mbox{ in }
(0,L),
\end{array}
\right.
$$
where $(f_1,f_2,u,v)$ are
controls which belong to some suitable spaces. Under some assumptions, one can show the exact controllability of
the above control system (See \cite[Theorem 10.12 in Chapter 10]{LZ3.1}). This, in turn, justifies our modification \eqref{8.600eq2.1}.

Some further result related to
the exact controllability of the above refined stochastic wave equations (even in several space dimensions) can be found in \cite{LZ6}.

%

%%%%%%%%%%%%%%%%%%%%%%%%%%%%%%%%%%%%%%%%%%%%%%%%%%%%%%%%

\section{Stochastic linear quadratic optimal control
problems I: open-loop optimal
control}\label{s-lq-o}

%%%%%%%%%%%%%%%%%%%%%%%%%%%%%%%%%%%%%%%%%%%%%%%%%%%%%%%%%%%

From this section to Section \ref{sec-slq-d}, we deal with linear quadratic optimal control problems (LQ problems for short) for SEEs. The content of these sections are taken from \cite{Luqi9, LZ3.1}.

%%%%%%%%%%%%%%%%%%%%%%%%%%%%%%%%%%%%%%%%%%%%%%%%%%%%%%%%%%%%

\subsection{Formulation of the problem}

%%%%%%%%%%%%%%%%%%%%%%%%%%%%%%%%%%%%%%%%%%%%%%%%%%%%%%%%%%%%

Fix any $T>0$, let both $H$ and $A$ (generating a $C_0$-semigroup $S(\cd)$ on $H$) be given as in Subsection \ref{Ch-SEE}, and let $U$ be another Hilbert space.

In what follows, for a Hilbert space $\wh H$, denote by $\dbS(\wh H)$ the set  of all
self-adjoint operators on $\wh H$.  For $M, N\in \dbS(\wh H)$, we use the notation $M \geq N$ (\resp $M >
N$) to indicate that $M-N$ is positive semi-definite (resp. positive definite).
For any $\dbS(\wh H)$-valued stochastic process $F$ on $[0, T]$, we write $F\geq 0$ (\resp
$F>0$, $F\gg 0$) if $F(t,\om)\geq 0$ (\resp $F(t,\om)>0$, $F(t,\om) \geq \d I$ for some $\d>0$ and the identity operator $I$ on $\wh H$)
for a.e. $(t,\om) \in [0, T] \times\Om$.

Consider a control system governed by the following linear SEE:
\begin{equation}\label{LQsystem1}
\left\{
\begin{array}{ll}\ds
    dX(t)=\big(A X(t) + A_1(t)X(t)+B(t)u(t)\big)dt \\
    \ns\ds\qq\qq + \big(C(t)X(t)+D(t)u(t)\big)dW(t) &\mbox{ in }(0,T],\\
    \ns\ds X(0)=\eta.
\end{array}
\right.
\end{equation}
Here $\eta\in H$, $A_1, C\in
L^\infty_\dbF(0,T;\cL(H))$ and $ B, D \in
L^\infty_\dbF(0,T;\cL(U;H))$.  In \eqref{LQsystem1}, $X(\cd)$ is
the {\it state variable} and $u(\cd)\in\cU[0,T]\deq L^2_\dbF(0, T; U))$ is the {\it control variable}.  For any
$u\in\cU[0,T]$, by Theorem \ref{ch-1-well-mild}, the system
\eqref{LQsystem1} admits a unique solution
$X(\cd;\eta,u)\in C_\dbF([0,T]; L^2(\Om;H))$
such that
\begin{equation}\label{est1}
|X(\cd;\eta,u)|_{C_\dbF([0,T];L^2(\Om;H))}\leq
\cC\big(|\eta|_{H} + |u|_{L^2_\dbF(0,T;U)}\big).
\end{equation}
When there is no confusion, we simply denote the
solution by $X(\cd)$.

Associated with the system \eqref{LQsystem1}, we consider the following quadratic
cost functional
\begin{equation}\label{LQcost}
\begin{array}{ll}\ds
    \cJ(\eta;u(\cd))\3n&\ds=\frac{1}{2}\mE\Big[
    \int_0^T\big(\langle
    M(t)X(t),X(t)\rangle_H +\langle
    R(t)u(t),u(t)\rangle_U\big)dt
    \\
    \ns&\ds\qq\;\,+ \langle
    GX(T),X(T)\rangle_H\Big],
\end{array}
\end{equation}
Here $M(\cd)\in L^\infty_{\dbF}(0,T;\dbS(H))$,
$R(\cd)\in L^\infty_{\dbF}(0,T;\dbS(U))$  and $G(\cd)\in
L_{\cF_T}^\infty(\Om;\dbS(H))$. In what follows,
to simplify the notations, the time variable $t$ will be suppressed in $X(\cd)$, $B(\cd)$, $C(\cd)$, $D(\cd)$, $M(\cd)$ and $R(\cd)$, and therefore we shall simply
write them as $B$, $C$, $D$, $M$ and $R$, respectively (if there is no confusion).

Let us now state our stochastic LQ problem as follows:

\ms

\no\bf Problem (SLQ). \rm For each $ \eta \in
H$, find a $\bar u(\cd)\in\cU[0,T]$ such that
\begin{equation}\label{SLQP}
\cJ(\eta;\bar u(\cd))=
\inf_{u(\cd)\in\cU[0,T]}\cJ(\eta;u(\cd)).
\end{equation}

\ms

Problem (SLQ) is now extensively studied for the
case of finite dimensions (i.e., $\dim
H<\infty$) and natural filtration. A nice treatise for this topic is \cite{Sun} and the readers can find rich references therein. To handle the infinite dimensional case, we borrow some ideas, such as the optimal feedback control operator,  from \cite{Sun}. At first glance, one might think that the study
of  Problem (SLQ) for $\dim
H=\infty$ is simply a routine extension
of that for $\dim
H<\infty$. However, the infinite
dimensional setting leads to significantly new
difficulties.

We begin with the following notions.

\begin{definition}\label{def1}
1)  Problem (SLQ) is said to be {\it
    finite at $\eta\in H$} if
$$\inf_{u(\cd)\in\cU[0,T]}\cJ(\eta;u(\cd))>-\infty.
$$

\ss

2) Problem (SLQ) is said to be {\it
    solvable at $ \eta \in H$} if there exists a
control $\bar u(\cd)\in\cU[0,T]$, such that
$$
\cJ(\eta;\bar
u(\cd))=\inf_{u(\cd)\in\cU[0,T]}\cJ(\eta;u(\cd)).
$$
In this case, $\bar u(\cd)$ is called an {\it
    optimal control}. The corresponding state $\cl
X(\cd)$ and $(\cl X(\cd),\bar u(\cd))$ are
called an {\it optimal state} and an
{\it optimal pair}, respectively.

\ss

3) Problem (SLQ) is said to be {\it
    finite} (\resp  {\it solvable})   if it is
finite (\resp solvable) at any
$\eta\in H$.

\end{definition}
%

%%%%%%%%%%%%%%%%%%%%%%%%%%%%%%%%%%%%%%%%%%%%%%%%

\subsection{Finiteness and solvability of Problem (SLQ)}

%%%%%%%%%%%%%%%%%%%%%%%%%%%%%%%%%%%%%%%%%%%%%%%%

In this subsection, we are concerned with the finiteness and solvability of
Problem (SLQ).

By the variation of constants formula, the solution to \eqref{LQsystem1} can be explicitly
expressed in terms of the initial state and the control in a linear form. Substituting that into the cost functional (which is quadratic in the state and
control), one can obtain a quadratic functional w.r.t.  the
control. Thus, Problem (SLQ) can be transformed to a quadratic optimization
problem on $\cU[0,T]$. This leads to some necessary and
sufficient conditions for the finiteness and solvability of Problem (SLQ). Let
us present the details below.

To begin with, we define four linear
operators $\Xi:\cU[0,T]\to L^2_\dbF (0,T;H)$, $\h
\Xi:\cU[0,T]\to L^2_{\cF_T}(\Omega;H)$,
$\G :H\to L^2_\dbF (0,T;H)$ and $\h\G:H\to L^2_{\cF_T}(\Omega;H)$ as follows:
\begin{equation}\label{LQoperator}
\left\{
\begin{array}{ll}\ds
    (\Xi u(\cd))(\cd)\deq X(\cd;0,u),\q \h
    \Xi u(\cd)\deq X(T; 0,u), \q\forall\,  u(\cd)
    \in\cU[0,T],\\
    \ns\ds (\G  \eta)(\cd)\deq X(\cd; \eta,0),\qq\,
    \h \G  \eta\deq X(T; \eta,0)\qq\q \forall\,  \eta\in H.
\end{array}
\right.
\end{equation}
where $X(\cd)\equiv X(\cd; \eta, u(\cd))$ solves
the equation \eqref{LQsystem1}. By \eqref{est1},
the  above operators are all bounded.

By \eqref{LQoperator}, for any $ \eta \in  H$ and
$u(\cd)\in\cU[0,T]$, the corresponding state
process $X(\cd)$ and its terminal value $X(T)$
satisfy
\begin{equation}\label{LQsystem4}
X(\cd)=(\G  \eta)(\cd)+(\Xi u(\cd))(\cd),\q X(T)=\h\G  \eta +\h \Xi u(\cd).
\end{equation}

We need to compute
the adjoint operators of $\Xi$, $\h
\Xi$,
$\G$ and $\h\G$.  To this end, let us introduce
the following BSEE\footnote{Here and henceforth, for any operator-valued process (\resp
random variable) $R$, we denote by $R^*$ its
pointwise dual operator-valued process (\resp
random variable). For example, if $R\in
L^{r_1}_\dbF(0,T; L^{r_2}(\Omega; \cL(H)))$,
then $R^*\in L^{r_1}_\dbF(0,T; L^{r_2}(\Omega;
\cL(H)))$, and $|R|_{L^{r_1}_\dbF(0,T;
    L^{r_2}(\Omega;
    \cL(H)))}=|R^*|_{L^{r_1}_\dbF(0,T;
    L^{r_2}(\Omega; \cL(H)))}$.}:
\begin{equation}\label{bLQsystem1}
\left\{\2n
\begin{array}{ll}\ds
    dY(t)=\!-\big(A^* Y(t)\! +\! A_1^* Y(t)\!+\!C^*
    Z(t)\!+\!\xi(t)\big)dt \!+\!Z(t)dW(t) & \mbox{in }[0,T),\\
    \ns\ds Y(T)=Y_T,
\end{array}
\right.
\end{equation}
where $Y_T\in L^2_{\cF_T}(\Omega;H)$ and $\xi(\cd)\in L^2_\dbF (0,T;H)$. Since we do not assume that the filtration $\mathbf{F}$ is the natural one generated by $W(\cd)$, the equation \eqref{bLQsystem1} may not have a weak/mild solution. As we have explained in Subsection \ref{sec-BSEE}, its solution is  understood in the sense of transposition.  For any $Y_T\in L^2_{\cF_T}(\Omega;H)$ and
$\xi(\cd)\in L^2_\dbF (0,T;H)$, by Theorem \ref{vw-th1},
there exists a unique transposition solution
$(Y(\cd),Z(\cd))\in
D_\dbF([0,T];L^2(\Om;H))\times L^2_\dbF (0,T;H)$ to
\eqref{bLQsystem1}. Moreover,
\begin{equation}\label{12.4-eq1}
\sup_{0\leq t\leq T}\mE|Y(t)|_H^2+
\mE\int_0^T|Z(t)|_H^2dt\leq \cC\mE\Big(|Y_T|_H^2
+\int_0^T|\xi(t)|_H^2dt\Big),
\end{equation}
for some constant $\cC>0$.

We have the following result.
\begin{proposition}\label{LQprop1}
For any $\xi(\cd)\in L^2_\dbF (0,T;H)$,
\begin{equation}\label{LsT}
    \left\{
    \begin{array}{ll}\ds
        (\Xi^*\xi)(t)= B^* Y_0(t)+D^* Z_0(t),\qq
        t\in[0,T],\\
        \ns\ds \G^*\xi =  Y_0(0),\\
    \end{array}
    \right.
\end{equation}
where $(Y_0(\cd),Z_0(\cd))$ is the transposition
solution to \eqref{bLQsystem1} with $Y_T=0$.
Similarly, for any $Y_T\in L^2_{\cF_T}(\Omega;H)$,
\begin{equation}\label{LssT}
    \left\{
    \begin{array}{ll}\ds
        (\widehat \Xi^*Y_T)(t)=  B^* Y_1(t)+D^*
        Z_1(t),\qq t\in[0,T],\\
        \ns\ds \h\G^*Y_T=  Y_1(0),
    \end{array}
    \right.
\end{equation}
where $(Y_1(\cd),Z_1(\cd))$
is the transposition solution to
\eqref{bLQsystem1} with $\xi(\cd)=0$.
\end{proposition}

{\it Proof}. For any $\eta\in H$ and
$u(\cd)\in\cU[0,T]$, let $X(\cd)$ be the
solution of \eqref{LQsystem1}. From the
definition of the transposition solution to the
equation \eqref{bLQsystem1}, we find that
\begin{equation*}
    \begin{array}{ll}\ds
        \mE\langle X(T),Y_T\rangle_H-\mE\langle
        \eta,Y(0)\rangle_H \\
        \ns\ds =\mE\int_0^T\big(\langle
        u(t),B^* Y(t)+D^* Z(t)\rangle_U
        -\langle X(t),\xi(t)\rangle_H \big)dt.
    \end{array}
\end{equation*}
This implies that
\begin{eqnarray}\label{ex eq1}
    \begin{array}{ll}\ds
        \mE\big(\langle\h\G  \eta +\h \Xi
        u(\cd),Y_T\rangle_H -\langle
        \eta,Y(0)\rangle_H\big)\\
        \ns\ds=\mE\int_0^T \big(\langle u(t),B^* Y(t) + D^*
        Z(t)\rangle_U - \langle(\G \eta)(t)+(\Xi
        u(\cd))(t),\xi(t)\rangle_H\big)dt.
    \end{array}
\end{eqnarray}

Choosing $Y_T=0$ and $\eta=0$ in \eqref{ex eq1},
we have
$$
\mE \int_0^T \langle(\Xi
u(\cd))(t),\xi(t)\rangle_H dt =  \mE\int_0^T
\langle u(t),B^* Y_0(t) + D^* Z_0(t)\rangle_U
dt,
$$
which yields the first equality in \eqref{LsT}.
Choosing $u(\cd)=0$ and $Y_T=0$ in \eqref{ex
    eq1}, we obtain that
$$
\mE\int_0^T\langle(\G \eta)(t),\xi(t)\rangle_H
dt=\mE\langle \eta,y_0(0)\rangle_H.
$$
This proves the second equality in \eqref{LsT}.

Letting $\eta=0$ and $\xi(\cd)=0$ in \eqref{ex eq1}, we find that
$$
\mE\langle\h \Xi u(\cd),Y_T\rangle_H=
\mE\int_0^T\langle u(t), B^* Y_1(t)+D^*  Z_1(t)
\rangle_H dt.
$$
This proves the first equality in \eqref{LssT}.
Letting $u(\cd)=0$ and $\xi(\cd)=0$ in \eqref{ex eq1}, we see
that $ \mE\langle \h\G \eta,Y_T\rangle_H =
\mE\langle \eta,Y_1(0)\rangle_H$, which gives
the second equality in \eqref{LssT}.
\endpf

From Proposition \ref{LQprop1}, we
get immediately the following result, which is a
representation for the cost functional
\eqref{LQcost}, and will play an important role in
the sequel.

\begin{proposition}\label{LQprop2}
The cost functional given by \eqref{LQcost} can
be presented as follows:
\begin{equation}\label{LQcost1}
    \cJ(\eta;u(\cd))=\frac{1}{2}\[\mE\int_0^T\big(\langle
    \cN u,u\rangle_U +2\langle \cH(\eta),u \rangle_U
    \big)dt +  \cM(\eta)\],
\end{equation}
where
\begin{equation}\label{cost1.1}
    \left\{
    \begin{array}{ll}\ds
        \cN =R+\Xi^*M\Xi  +\h \Xi^*G\h
        \Xi,\\
        \ns\ds \cH(\eta)= \big(\Xi^*M \G \eta\big)(\cd)
        +\big(\h
        \Xi^*G \h\G \eta\big)(\cd), \\
        \ns\ds \cM(\eta)=\big\langle M \G\eta,\G\eta
        \big\rangle_{ L^2_\dbF(0,T;H)} +\big\langle G\h
        \G\eta,\h \G\eta\big\rangle_{
            L^2_{\cF_T}(\Omega;H)}.
    \end{array}
    \right.
\end{equation}
\end{proposition}

As an application of Proposition \ref{LQprop2}, we have the following result for the finiteness
and solvability of  Problem (SLQ).

\begin{theorem}\label{th ext}
{\rm 1)}  If Problem (SLQ) is finite at some
$\eta\in H$, then
\begin{equation}\label{LQth1eq1}
    \cN \geq 0.
\end{equation}

\ss

{\rm 2)} Problem (SLQ) is solvable at
$\eta\in H$ if and only if $\cN \geq 0$ and
there exists a control $\bar u(\cd)\in\cU[0,T]$  such
that
\begin{equation}\label{LQth1eq2}
    \cN \bar u(\cd)+\cH (\eta)=0.
\end{equation}
In this case, $\bar u(\cd)$ is an  optimal
control.

\ss

{\rm 3)}  If $\cN \gg 0$, then for any $\eta\in H$,
$\cJ(\eta;\cd)$ admits a  unique minimizer $\bar
u(\cd)$ given by
\begin{equation}\label{LQth1eq3}
    \bar u(\cd)=-\cN^{-1} \cH(\eta).
\end{equation}
In this case, it holds
\begin{equation*}
    \begin{array}{ll}\ds
        \inf_{u(\cd)\in\cU[0,T]}\cJ(\eta;u(\cd))=\cJ(\eta;\cl
        u(\cd)) = \frac{1}{2} \(\cM(\eta) - \langle
        \cN^{-1}\cH(\eta),\cH(\eta)\rangle_H \).
    \end{array}
\end{equation*}
\end{theorem}

{\it Proof}. We prove the assertion 1) by
contradiction. Suppose that \eqref{LQth1eq1} was
not true. Then there would exist $u_0\in
\cU[0,T]$ such that
$$
\mE\int_0^T\langle \cN  u_0, u_0\rangle_U ds <
0.
$$
Define a sequence
$\{u_k\}_{k=1}^\infty\subset \cU[0,T]$ as $
u_k(\cd)= ku_0(\cd)$ for $k\in\dbN$. For $u_k(\cd)$, we have
$$
\begin{array}{ll}\ds
    \cJ(\eta;u_k(\cd))  = \frac{k^2}{2}\mE\(
    \int_0^T\langle \cN u_0, u_0 \rangle_U dt +
    \frac{2}{k} \int_0^T\langle
    \cH(\eta,u_0)\rangle_U dt + \frac{1}{k^2} \cM
    (\eta)\).
\end{array}
$$
Since
$
\ds\mE\int_0^T\langle \cH(\eta,u_0)\rangle_U dt +
\mE \cM (\eta)<+\infty,
$
there is a $k_0>0$ such that for
all $k\geq k_0$,
$$
\cJ(\eta;u_k)\leq \frac{k^2}{4}
\mE\int_0^T\langle \cN u_0, u_0 \rangle_U dt.
$$
By letting $k\to \infty$, we find that
$\cJ(\eta;u_k(\cd))\to-\infty$, which
contradicts that  Problem (SLQ) is
finite at $\eta$.

\medskip

Now we prove the assertion 2). For the ``if" part, let $\bar
u(\cd)\in \cU[0,T]$ be an optimal control of
Problem (SLQ) for $\eta\in H$. From the
optimality of $\bar u(\cd)$, for any $u(\cd)\in
\cU[0,T]$, we see that
$$
\begin{array}{ll}\ds
    0 \leq \liminf_{\l\to 0 } \frac{1}{\l} \(
    \cJ(\eta;\bar u+\l u) - \cJ(\eta;\bar u) \) =
    \int_0^T\langle \cN \bar u + \cH (\eta), u
    \rangle_U dt.
\end{array}
$$
Consequently,
$
\cN \bar u+\cH (\eta)=0.
$

For the ``only if" part, let
$(\eta,\bar u(\cd))\in   H \times \cU[0,T]$
satisfy \eqref{LQth1eq2}. For any
$u\in\cU[0,T]$, from \eqref{LQth1eq1}, we have
$$
\begin{array}{ll}\ds
    \cJ(\eta;u(\cd)) - \cJ(\eta;\bar u(\cd)) = \cJ\big(\eta;\bar u(\cd) + u(\cd) - \bar u(\cd)\big) - \cJ(\eta;\bar u(\cd))\\
    \ns\ds = \mE\int_0^T\big\langle \cN \bar u  +
    \cH (\eta), u -\bar u \big\rangle_Udt
    + \frac{1}{2}\mE\int_0^T\big\langle \cN\big(u -\bar u \big),u -\bar u \big\rangle_U dt\\
    \ns\ds = \frac{1}{2}\mE\int_0^T\big\langle \cN
    \big(u -\bar u \big),u -\bar u \big\rangle_Udt
    \geq 0.
\end{array}
$$
This concludes that $\bar u(\cd)$ is an optimal
control.

\medskip

Finally, we prove the assertion 3). Since all the optimal controls should
satisfy \eqref{LQth1eq2} and $\cN$ is
invertible, we get the assertion 3)
immediately.
\endpf

When
\begin{equation}\label{4.7-eq3}
R\gg 0,\;M \geq 0,\;G\geq 0,
\end{equation}
Problem (SLQ) is called a
{\it standard SLQ problem}. In such a case, $\cN \gg
0$. Then, by Theorem \ref{th ext}, it is
uniquely solvable and the optimal control is given by \eqref{LQth1eq3}.

\begin{remark}
The  main drawback of the formula \eqref{LQth1eq3} is that  it is very difficult
to compute the inverse of the operator $\cN$.
\end{remark}
%

%%%%%%%%%%%%%%%%%%%%%%%%%%%%%%%%%%%%%%%%%%%%%%%%%%%%%%%%%%%%%%

\subsection{Pontryagin type maximum principle for
Problem (SLQ)}

%%%%%%%%%%%%%%%%%%%%%%%%%%%%%%%%%%%%%%%%%%%%%%%%%%%%%%%%%%%%%%

In this subsection, we  derive the following Pontryagin type maximum principle  for
Problem (SLQ).

\begin{theorem}\label{LQth max}
Let Problem (SLQ) be solvable at
$\eta\in H$ with $(\cl X(\cd), \bar u(\cd))$
being an optimal pair. Then for the
transposition solution $(Y(\cd), Z(\cd))$ to
\begin{equation}\label{bsystem2.1}
    \left\{\2n
    \begin{array}{ll}\ds
        d Y =-\big(A^*Y  + A_1^*Y + C^*
        Z -  M\cl X \big)dt + Z dW(t) &
        \mbox{in }[0,T),\\
        \ns\ds Y(T)=-G\cl X(T),
    \end{array}
    \right.
\end{equation}
we have that
\begin{equation}\label{max pr}
    \begin{array}{ll}\ds
        R\bar u -B^* Y -D^*Z  =0,
    \end{array}\q\ae
    (t,\om) \in[0,T]\times \Om.
\end{equation}
\end{theorem}
To prove Theorem \ref{LQth max}, we need the
following preliminary result.
\begin{lemma}\label{lemma4}
Let $\wt U$ be a convex subset of $U$. If $F(\cdot),\bar
u(\cdot)\in L^2_\dbF(0,T;\wt U)$ satisfy
\begin{equation}\label{lemma4 ine1}
    \dbE \int_0^T \big\langle F(t,\cd), u(t,\cd) -
    \bar u(t,\cd) \big\rangle_{U} dt \leq 0,\q \forall\; u(\cdot)\in L^2_\dbF(0,T;\wt U),
\end{equation}
then,
\begin{equation}\label{lemma4 ine2}
    \big\langle F(t,\omega), \rho- \bar
    u(t,\omega)\big\rangle_{U} \leq 0, \,\ \ae
    (t,\omega)\in [0,T]\times\Omega,\q \forall\; \rho\in \wt U.
\end{equation}
\end{lemma}

{\it Proof}. We prove \eqref{lemma4 ine2} by the
contradiction argument. Suppose that
\eqref{lemma4 ine2} did not hold. Then, there
would exist  $u_0\in \wt U$ and $\e>0$ such that
$$
\a_\e\deq \int_\Omega\int_0^T
\chi_{\L_\e}(t,\omega)dtd\dbP
>0,
$$
where $ \L_\e \deq \big\{ (t,\omega)\in
[0,T]\times \Omega\;\big|\; \big\langle
F(t,\omega), u_0 - \bar
u(t,\omega)\big\rangle_{U} \geq \e \big\}$. For any $m\in \dbN$, let
$$
\L_{\e,m}\deq\L_\e\cap \big\{(t,\omega)\in
[0,T]\times\Omega\;\big|\;|\bar
u(t,\omega)|_{U}\le m\big\}.
$$
It is clear
that $\ds\lim_{m\to\infty}\L_{\e,m}=\L_\e$.
Hence, there is $m_\e\in\dbN$ such that
$$
\int_\Omega\int_0^T
\chi_{\L_{\e,m}}(t,\omega)dtd\dbP
>\frac{\a_\e}{2}>0, \qq \forall\;m\ge m_\e.
$$
Since $\big\langle F(\cdot), u_0 - \bar
u(\cdot)\big\rangle_{U}$ is
$\dbF$-measurable, so is the process
$\chi_{\L_{\e,m}}(\cdot)$. Set
$$
\hat u_{\e,m}(t,\omega) \deq u_0
\chi_{\L_{\e,m}}(t,\omega)+ \bar
u(t,\omega)\chi_{\L_{\e,m}^c}(t,\omega),\q
(t,\omega)\in [0,T]\times \Omega.
$$
Since $|\bar u(\cd)|_{U}\le m$ on
$\L_{\e,m}$, we see that $\hat
u_{\e,m}(\cdot)\in L^2_\dbF(0,T;\wt U)$ and
$\hat u_{\e,m}(\cd)-\bar u(\cd)\in
L^2_\dbF(0,T; U)$. Hence, for any $m\ge m_\e$,
we have that
$$
\begin{array}{ll}\ds
    \dbE\int_0^T \big\langle F(t), \hat u_{\e,m}(t)
    -\bar u(t) \big\rangle_{U} dt\\
    \ns\ds= \int_\Omega\int_0^T
    \chi_{\L_{\e,m}}(t,\omega) \big\langle
    F(t,\omega), u_0 -
    \bar u(t,\omega) \big\rangle_{U}\, dtd\dbP  \\
    \ns \ds \geq \e\int_\Omega\int_0^T
    \chi_{\L_{\e,m}}(t,\omega)dtd\dbP  \geq
    \frac{\e\a_\e}{2}
    >0,
\end{array}
$$
which contradicts \eqref{lemma4 ine1}. This
completes the proof of Lemma \ref{lemma4}.\endpf

Now, let us prove Theorem \ref{LQth max}.

{\it Proof}.[Proof of Theorem \ref{LQth max}]
For the optimal pair $(\cl X(\cdot),\bar
u(\cdot))$ and a control $u(\cdot)\in \cU[0,T]$,
we see that
$$
u^\e(\cdot) \deq  \bar u(\cdot) + \e [u(\cdot) -
\bar u(\cdot)] = (1-\e)\bar u(\cdot) + \e
u(\cdot) \in \cU[0,T], \q\forall\;\e \in [0,1].
$$
Denote by $X^\e(\cdot)$ the state process of
\eqref{LQsystem1} corresponding to the control
$u^\e(\cdot)$. Write
$$
X_1^\e(\cd) \deq
\frac{X^\e(\cd)-\cl
    X(\cd)}{\e},\qq\d u(\cd) \deq u(\cd) - \bar
u(\cd).
$$
It is easy to see that $X_1^\e(\cdot)$
solves the following SEE:
\begin{eqnarray}\label{LQfsystem3x}
    \3n\3n\!\!\!\!\left\{
    \begin{array}{lll}\ds
        dX_1^\e = \big(AX_1^\e\! + \! A_1 X^\e_1 \!+ B\d
        u \big)dt + \big(C X^\e_1\! + D\d u \big)dW(t)
        &\mbox{ in
        }(0,T],\\
        \ns\ds X_1^\e(0)=0.
    \end{array}
    \right.
\end{eqnarray}
Since $(\cl X(\cdot),\bar u(\cdot))$ is an
optimal pair of Problem (SLQ), we have
\begin{eqnarray}\label{LQvar 1}
    \begin{array}{ll}\ds
        0\leq \lim_{\e\to 0}\frac{\cJ(\eta;u^\e(\cdot)) - \cJ(\eta;\bar u(\cdot))}{\e} \\
        \ns\ds\;\;\;= \dbE\int_0^T
        \big(\big\langle M \cl
        X(t),X_1^\e(t)\big\rangle_H    +
        \big\langle R \bar u(t),\d u(t)
        \big\rangle_{U}\big) dt +
        \dbE\big\langle G\cl
        X(T),X_1^\e(T)\big\rangle_H.
    \end{array}
\end{eqnarray}

By the definition of the transposition solution
to \eqref{bsystem2.1}, we obtain that
\begin{equation}\label{LQmax eq1}
    \begin{array}{ll}\ds
        -\dbE \big\langle G\cl
        X(T),X_1^\e(T)\big\rangle_H +
        \dbE\int_0^T \big\langle  M \cl X(t), X_1^\e(t)\big\rangle_H dt \\
        \ns\ds = \dbE \int_0^T \big(\big\langle
        B \d u(t), Y(t)\big\rangle_H +
        \big\langle D \d u(t),
        Z(t)\big\rangle_H\big)dt.
    \end{array}
\end{equation}
Combining \eqref{LQvar 1} and \eqref{LQmax eq1},
we arrive at
\begin{equation}\label{LQmax ine2}
    \dbE\int_0^T\big\langle R \bar u(t) -
    B^*Y(t) - D^*Z(t), u(t)-\bar u(t)
    \big\rangle_{U}dt\geq 0, \q \forall\;
    u(\cdot)\in \cU[0,T].
\end{equation}
Hence, by Lemma \ref{lemma4} (where we choose $\wt U=U$), we conclude that
\begin{equation}\label{LQmax ine3}
    \begin{array}{ll}\ds
        \big\langle R \bar u - B^*Y -
        D^*Z, \rho-\bar u \big\rangle_{U}
        \geq 0,\q\ae [0,T]\times
        \Omega,\;\forall\; \rho \in U.
    \end{array}
\end{equation}
This implies  \eqref{max pr}.
\endpf

\ms

Next, we introduce the following decoupled
forward-backward stochastic evolution equation
(FBSEE for short):
\begin{equation}\label{fbLQsystem1}
\left\{
\begin{array}{ll}\ds
    d X =\big(AX +A_1X +Bu \big)dt +\big(CX +Du \big)dW(t) &\mbox{ in }(0,T],\\
    \ns\ds dY = -\big(A^* Y +
    A_1^*Y - MX + C^*Z \big)dt +Z dW(t)&\mbox{ in }[0,T),\\
    \ns\ds X(0)=\eta,\qq Y(T)=-GX(T).
\end{array}
\right.
\end{equation}
We call $(X(\cd),Y(\cd),Z(\cd))$ a {\it transposition
solution} to the equation \eqref{fbLQsystem1} if
$X(\cd)$ is the mild solution to the forward SEE
and $(Y(\cd),Z(\cd))$ is the transposition
solution to the backward one.

Since the equation \eqref{fbLQsystem1} is decoupled, its well-posedness is easy to be
obtained. Given $\eta\in H$ and $u(\cd)\in\cU[0,T]$, one
can first solve the forward SEE to get $X(\cd)$,
and then find the transposition solution
$(Y(\cd),Z(\cd))$ to the backward one. Hence,  \eqref{fbLQsystem1} admits a unique transposition
solution $(X(\cd), Y(\cd),Z(\cd))$ corresponding
to $\eta$ and $u(\cd)$. The following result is
a further consequence of Proposition
\ref{LQprop1}.

\begin{proposition}\label{LQprop3}
For any
$(\eta,u(\cd))\in H\times\cU[0,T]$, let
$(X(\cd),Y(\cd),Z(\cd))$ be the transposition
solution to \eqref{fbLQsystem1}. Then
\begin{equation}\label{prop3eq1}
    \big(\cN u +\cH(\eta)\big)(t) = Ru(t)
    -B^* Y(t)-D^* Z(t),\q \ae
    (t,\om)\in[0,T]\times\Om.
\end{equation}
\end{proposition}

\medskip

{\it Proof}. \rm Let $(X(\cd),Y(\cd),Z(\cd))$ be
the transposition solution of
\eqref{fbLQsystem1}. From \eqref{cost1.1}, we
obtain that
\begin{equation}\label{prop3eq3}
    \begin{array}{ll}\ds
        \cN u  +
        \cH (\eta)\3n&\ds  =(R+\Xi^*M\Xi +\h \Xi^*G\h \Xi)u  +
        (\Xi^*M \G \eta)  + \h
        \Xi^*G \h\G \eta \\
        \ns&\ds  =Ru  +\Xi^*M\big(\G \eta  +\Xi u  \big)
        +\h \Xi^*G(\h\G \eta + \h
        \Xi u) \\
        \ns&\ds =Ru  +\Xi^* MX  +\h \Xi^*GX(T).
    \end{array}
\end{equation}
By \eqref{LsT} and \eqref{LssT}, we find
$
\Xi^* MX  + \h \Xi^*GX(T)= -B^* Y -D^* Z.
$
This, together with \eqref{prop3eq3}, implies
\eqref{prop3eq1} immediately.
\endpf

As a variant of Theorem
\ref{th ext}, we have the following result, in which the conditions are given in
terms of an FBSEE.
\begin{theorem}\label{prop4}
Problem (SLQ) is solvable at $\eta\in H$ with an
optimal pair $(\cl X(\cd),\bar u(\cd))$ if and
only if there exists a unique $(\cl X(\cd),\bar
u(\cd), Y(\cd),$ $Z(\cd))$ satisfying the FBSEE
\begin{equation}\label{fbsystem2}
    \left\{\!
    \begin{array}{ll}\ds
        d\cl X =\big(A\cl X + A_1\cl
        X + B\bar u \big)dt +\big(C\cl
        X + D\bar
        u \big)dW(t) &\mbox{in }(0,T],\\
        \ns\ds d Y =-\big(A^* Y + A_1^*
        Y +C^*Z - M\cl X \big)dt
        +Z dW(t)&\mbox{in }[0,T),\\
        \ns\ds\cl X(0)=\eta,\qq y(T)=-G\cl X(T)
    \end{array}
    \right.
\end{equation}
with the condition
\begin{equation}\label{prop4eq1}
    R \bar u -B^* Y -D^* Z =0,\q
    \ae
    (t,\om)\in [0,T]\times\Om,
\end{equation}
and for any $u(\cd)\in\cU[0,T]$, the unique
transposition solution $(X_0(\cd), Y_0(\cd),Z_0(\cd))$
to \eqref{fbLQsystem1} with $\eta=0$ satisfies
\begin{equation}\label{prop4eq2}
    \mE\int_0^T \big\langle R u(t) -B^*(t)
    Y_0(t)-D^* Z_0(t),u(t)\big\rangle_U dt \geq 0.
\end{equation}
\end{theorem}

{\it Proof}.  {\bf The ``only if" part}.
Clearly, \eqref{prop4eq1} and \eqref{prop4eq2}
follow  from Theorem \ref{LQth max} immediately.

\medskip

{\bf The ``if" part}. From Proposition \ref{LQprop3},
the inequality \eqref{prop4eq2} is equivalent to
$\cN\geq 0$. Now, let $(\cl
X(\cd),Y(\cd),Z(\cd))$ be a transposition
solution to \eqref{fbsystem2} such that
\eqref{prop4eq1} holds. Then, by Proposition
\ref{LQprop3}, we see that \eqref{prop4eq1} is
the same as \eqref{LQth1eq2}. Hence by Theorem
\ref{LQth max},  Problem (SLQ) is
solvable.
\endpf

\ms

Theorem \ref{prop4} is nothing but a restatement of Theorem \ref{th ext}. However,
when   $R(t)$ is
invertible for a.e. $(t,\om)\in [0,T]\times\Om$ and
\begin{equation}\label{11.13-eq2}
R(\cd)^{-1}\in L^\infty_\dbF(0,T;\cL(H)),
\end{equation}
it gives a way to find the optimal control by solving the following coupled
FBSEE:
\begin{equation}\label{LQfbsystem3}
\left\{
\begin{array}{ll}\ds
    d\cl X =\big[(A + A_1)\cl X +
    BR^{-1}B^*Y +
    BR^{-1}D^*Z \big]dt\cr
    \ns\ds\hspace{1.25cm} +\big(C\cl X
    + DR^{-1} B^*Y  +  D R^{-1} D^*
    Z \big)dW(t) & \mbox{ in } (0,T],
    \\
    \ns\ds dY =- \big[ M\cl X -(A +
    A_1)^* Y - C^* Z \big]dt+  Z
    dW(t) & \mbox{ in } [0,T),\cr
    \ns\ds \cl X(0)=\eta,\qq Y(T)=-G\cl X(T).
\end{array}
\right.
\end{equation}
Clearly, \eqref{LQfbsystem3} is a {\it coupled}
linear FBSEE. As consequences of Theorem \ref{prop4}, we have the following results.

\begin{corollary}\label{LQcor1}
Let \eqref{11.13-eq2} hold and $\cN \ge0$. Then
Problem (SLQ) is uniquely solvable at
$\eta\in H$ if and only if the FBSEE
\eqref{LQfbsystem3} admits a unique
transposition solution $(\cl X(\cd), Y(\cd),
Z(\cd))$. In this case,
\begin{equation}\label{cor1eq1}
    \bar u(t)= R^{-1}\big(B^* Y(t)+D^*
    Z(t)\big), \q \ae
    (t,\om)\in[0,T]\times\Om
\end{equation}
is an
optimal control.
\end{corollary}
\begin{corollary}\label{cor2}
Let \eqref{4.7-eq3} hold. Then FBSEE
\eqref{LQfbsystem3} admits a unique transposition
solution $(\cl X(\cd),$ $Y(\cd),Z(\cd))$ and
Problem (SLQ) is uniquely solvable  with the optimal control given by \eqref{cor1eq1}.
\end{corollary}
%

%%%%%%%%%%%%%%%%%%%%%%%%%%%%%%%%%%%%%%%%%%%%%%%%%%%%%%%%

\section{Stochastic linear quadratic optimal control
problems II:  optimal feedback
control}\label{s-lq-c}

%%%%%%%%%%%%%%%%%%%%%%%%%%%%%%%%%%%%%%%%%%%%%%%%%%%%%%%%%%%

The
optimal control given in the previous section (say \eqref{cor1eq1} in Corollary \ref{LQcor1}) is not in a feedback form. In Control Theory, one of the
fundamental issues is to find optimal feedback controls,
which is particularly important in practical
applications. The main advantage of
feedback controls is that they keep the
corresponding control strategy to be robust with
respect to  (small) perturbations/disturbances,
which are usually unavoidable in realistic
background. Unfortunately, it is actually very
difficult to find optimal feedback controls for a
control problem.

In this section, we study the optimal feedback
control for  Problem (SLQ) introduced in Section \ref{s-lq-o}. To simplify
notations, without loss of generality,
hereafter, we assume that $A_1=0$.

To begin with,  we introduce some notations.
Let $H_1$ and $H_2$ be Hilbert spaces. For $1\leq p_1, p_2,q_1,q_2\leq\infty$, let
\begin{equation}\label{12.13-eq1}
\begin{array}{ll}\ds
    \cL_{pd}\big(L_{\dbF}^{p_1}(0,T; L^{q_1}(\Omega;H_1)); L^{p_2}(0,T; L_{\dbF}^{q_1}(\Omega;H_2))\big)\\
    \ns\ds \deq \Big\{L\in \cL\big(L_{\dbF}^{p_1}(0,T; L^{q_1}(\Omega;H_1)); L^{p_2}(0,T; L_{\dbF}^{q_1}(\Omega;H_2))\big) \;\Big| \\
    \ns\ds\q \mbox{ for }\ae (t,\omega)\in [0,T]\times\Omega, \mbox{ there
        is }\wt L(t,\omega)\in\cL (H_1;H_2)\mbox{ satisfying  }\\
    \ns\ds\q   \mbox{ that } \big(\cL u(\cd)\big)(t,\omega)
    =\wt L (t,\omega)u(t,\omega),\q
    \forall\, u(\cd)\in
    L_{\dbF}^{p_1}(0,T; L^{q_1}(\Omega;H_1))\Big\}.
\end{array}
\end{equation}
In what follows, if there is no confusion,   we identify the above $L\in \cL_{pd}\big(L_{\dbF}^{p_1}(0,T; L^{q_1}(\Omega;H));$ $L^{p_2}(0,T; L_{\dbF}^{q_1}(\Omega;H))\big)$ with $\wt L(\cd,\cd)$.

Fix $p\geq 1$ and $q\geq 1$. Write
\begin{equation}\label{10.10-eq30}
\begin{array}{ll}\ds
    \Upsilon_{p,q}(H_1;H_2) \deq  \big\{J(\cd,\cd)\in
    \cL_{pd}(L^p_{\dbF}(\Om;L^\infty(0,T;H_1));L_{\dbF}^p(\Om;L^q(0,T;H_2)))\;\big|\\
    \ns \ds \hspace{1.53cm}\qq\qq |J(\cd,\cd)|_{\cL(H_1;H_2)}\in
    L^\infty_\dbF(\Om;L^q(0,T))\big\}.
\end{array}
\end{equation}
In the sequel, we shall simply denote $\Upsilon_{p,p}(H_1;H_2)$ (\resp $\Upsilon_{p,p}(H_1;H_1)$) by  \index{$\Upsilon_{p}(H_1;H_2)$} $\Upsilon_{p}(H_1;$ $H_2)$ (\resp  \index{$\Upsilon_{p}(H_1)$} $\Upsilon_{p}(H_1)$).

For Problem (SLQ), let us introduce the
notion of optimal feedback operator as follows:
\begin{definition}\label{5.7-def1}
An operator $\Th(\cd)\in \Upsilon_2(H;U)$ is
called an {\it optimal feedback  operator} for
Problem (SLQ) if
\begin{equation}\label{5.7-eq2}
    \begin{array}{ll}\ds
        \cJ(\eta;\Th(\cd)\cl X(\cd))\leq
        \cJ(\eta;u(\cd)), \q \forall\;  (\eta,
        u(\cd))\in H\times\cU[0,T],
    \end{array}
\end{equation}
where $\cl X(\cd)=\cl X(\cd;\eta, \Th(\cd)\cl
X(\cd))$ solves the following equation:
\begin{equation}\label{5.2-eq1.1}
    \left\{\begin{array}{ll}\ds d\cl X =\big(A\cl X + B\Th \cl X \big)dt + \big(C \cl X +D \Th
        \cl X \big)dW(t) &\mbox{ in
        }(0,T],\\
        \ns\ds \cl X(0)=\eta.
    \end{array}
    \right.
\end{equation}
\end{definition}
\begin{remark}
In Definition \ref{5.7-def1}, $\Th(\cd)$ is
required to be independent of $ \eta \in  H$.
For a fixed $ \eta \in H$, the inequality
\eqref{5.7-eq2} implies that the control $\bar
u(\cd)\equiv \Th(\cd)\cl X(\cd)\in
\cU[0,T]$
is optimal for Problem (SLQ).
Therefore, for Problem (SLQ), the
existence of an optimal feedback  operator implies the existence of optimal
controls for any $\eta\in H$, but not vice versa.
\end{remark}

Stimulated by
the pioneer work \cite{Bismut1} (for stochastic LQ problems in
finite dimensions), to study $\Th$, we introduce the following {\it
operator-valued, backward stochastic Riccati
equation} for Problem (SLQ):
\begin{equation}\label{5.5-eq6}
\left\{
\begin{array}{ll}\ds
    dP =-\big( PA + A^* P + \L C + C^* \L + C^* PC +
    M - L^* K^{-1} L
    \big)dt \\
    \ns\ds\qq\q + \L dW(t) \qq\mbox{ in }[0,T),\\
    \ns\ds P(T)=G,
\end{array}
\right.
\end{equation}
where
\begin{equation}\label{9.7-eq10}
K\equiv R+D^*PD>0, \qq L= B^* P+D^* (PC+\L).
\end{equation}

There exists a new essential difficulty in the
study of \eqref{5.5-eq6} when $\dim H=\infty$. Indeed, in the
infinite dimensional setting, although $\cL(H)$
is still a Banach space, it is neither reflexive
(needless to say to be a Hilbert space) nor
separable even if the Hilbert space $H$ itself is separable. As
far as we know, there exists no such a
stochastic integration/evolution equation theory
in general Banach spaces that can be employed to
treat the well-posedness of \eqref{5.5-eq6},
especially to handle the (stochastic integral)
term ``$\L dW(t)$". For example, the existing
results on stochastic integration/evolution
equations in UMD Banach spaces (e.g.
\cite{vanNeerven1}) do not fit the present case
because, $\cL(H)$ is not a UMD Banach space.

Because of the above mentioned difficulty, we
have to employ the stochastic transposition method (\cite{LZ1}) to introduce a new type of solutions to \eqref{5.5-eq6}. To this
end, let us first introduce the following two
SEEs:
\begin{equation}\label{op-fsystem1.2}
\left\{
\begin{array}{ll}
    \ds d\f_1 = \big(A \f_1 + u_1\big)d\tau + \big(C \f_1 + v_1\big)dW(\tau) &\mbox{ in } (t,T],\\
    \ns\ds \f_1(t)=\xi_1
\end{array}
\right.
\end{equation}
and
\begin{equation}\label{op-fsystem2.2}
\left\{
\begin{array}{ll}
    \ds d\f_2 =\big(A \f_2 + u_2\big)d\tau +\big( C \f_2 + v_2 \big)dW(\tau) &\mbox{ in } (t,T],\\
    \ns\ds \f_2(t)=\xi_2.
\end{array}
\right.
\end{equation}
Here $t\in [0,T)$, $\xi_1,\xi_2$ are suitable
random variables and $u_1,u_2,v_1,v_2$ are
suitable stochastic processes.

Also, we need to give the solution spaces for
\eqref{5.5-eq6}.
Let $V$ be a Hilbert space such that $H\subset V$ and the embedding from $H$ to $V$ is a Hilbert-Schmidt operator. Denote by $V'$ the dual space of $V$ with the pivot space $H$.
Put
$$
\ba{ll} \ds  D_{\dbF,w}([0,T];
L^{\infty}(\Om;\cL(H)))\\
\ns \ds\deq\big\{P\in D_{\dbF}([0,T];
L^{\infty}(\Om;\cL_2(H;V))) \;\big|\;  P(t,\omega) \in
\dbS(H),\,\ae (t,\omega)\in
[0,T] \times \Omega, \\
\ns \ds\q  \hbox{and
}\chi_{[t,T]}P(\cd)\zeta\in D_\dbF([t,T];L^2(\Om;H)),\
\forall\;\zeta\in L^2_{\cF_t}(\Om;H)\big\} \ea
$$
and
$$
\begin{array}{ll}\ds
L^2_{\dbF,w}(0,T;\cL(H))\\
\ns\ds \deq\big\{ \L \in
L^2_{\dbF}(0,T;\cL_2(H;V))\;
\big|\; D^*\L\in \cL_{pd}\big(L_{\dbF}^{\infty}(0,T; L^{2}(\Omega;H)); L^{2}_{\dbF}(0,T;U) \big) \big\}.
\end{array}
$$

Now, we introduce the notion of {\it transposition
solutions} to \eqref{5.5-eq6}:
\begin{definition}\label{4.8-def2}
We call $\big(P(\cd),\L(\cd)\big)\!\in\!
D_{\dbF,w}([0,T];
L^{\infty}(\Om;\cL(H)))\times
L^2_{\dbF,w}(0,T;\cL(H))$ a {\it transposition solution}
to \eqref{5.5-eq6} if the following three
conditions hold:

\ss

{\rm 1)}  $K(t,\om)\big(\equiv R(t,\om) +
D(t,\om)^*P(t,\om)D(t,\om)\big)> 0$ and its left
inverse $K(t,\om)^{-1}$ is a densely defined
closed operator for a.e. $(t,\om)\in
[0,T]\times\Om$;

\ss

{\rm 2)}  For any $t\in [0,T)$, $\xi_1,\xi_2\in
L^4_{\cF_t}(\Om;H)$, $u_1(\cd), u_2(\cd) \in L^4_\dbF(\Om; L^2(t,T;H))$ and  $v_1(\cd), v_2(\cd)\in L^4_\dbF(\Om;L^2(t,T;V'))$, it holds that
\bel{6.18-eq1}
\ba{ll}\ds
\mE\langle
G\f_{1}(T),\f_{2}(T)\rangle_{H} +\mE
\int_t^T \big\langle M(\tau) \f_{1}(\tau),
\f_{2}(\tau) \big\rangle_{H}d\tau\\[2mm]\ds
\q - \mE \int_t^T
\big\langle K(\tau)^{-1} L(\tau)\f_{1}(\tau),
L(\tau)\f_{2}(\tau) \big\rangle_{H}d\tau
\\[2mm]\ds = \mE\big\langle P(t)
\xi_{1},\xi_{2} \big\rangle_{H} + \mE \int_t^T
\big\langle P(\tau)u_{1}(\tau),
\f_{2}(\tau)\big\rangle_{H}d\tau \\[2mm]\ds \q + \mE
\int_t^T \big\langle P(\tau)\f_{1}(\tau),
u_{2}(\tau)\big\rangle_{H}d\tau  + \mE
\int_t^T\big\langle
P(\tau)C(\tau)\f_{1}(\tau),
v_{2}(\tau)\big\rangle_{H}d\tau  \\[2mm]\ds \q + \mE
\int_t^T \big\langle  P(\tau)v_{1}(\tau), C(\tau)\f_{2}(\tau)+v_{2}(\tau)\big\rangle_{H}d\tau\\[2mm]\ds  \q + \mE \int_t^T \big\langle
v_{1}(\tau),
\L(\tau)\f_2(\tau)\big\rangle_{V',V}d\tau+ \mE
\int_t^T \big\langle \L(\tau)\f_1(\tau),
v_{2}(\tau) \big\rangle_{V,V'}d\tau,
\ea\ee
where $\f_1(\cd)$ and $\f_2(\cd)$ solve
\eqref{op-fsystem1.2} and \eqref{op-fsystem2.2},
respectively\footnote{By Theorem
    \ref{ch-1-well-mild}, one has $\f_1(\cd),
    \f_2(\cd)\in C_\dbF([0,T];L^4(\Om;H))$.}.

\ss

{\rm 3)}  For any $t\in [0,T)$, $\xi_1,\xi_2\in
L^2_{\cF_t}(\Om;H)$, $u_1(\cd), u_2(\cd)\in
L^2_\dbF(t,T;H)$  and $v_1(\cd), v_2(\cd)$ $\in
L^2_\dbF(t,T;U)$, it holds that
\begin{eqnarray}\label{10.10-eq10}
    && \mE\langle G\f_{1}(T),\f_{2}(T)\rangle_{H}
    +\mE \int_t^T \big\langle M(\tau) \f_{1}(\tau),
    \f_{2}(\tau)
    \big\rangle_{H}d\tau\nonumber\\
    &&\q - \mE \int_t^T \big\langle K(\tau)^{-1}
    L(\tau) \f_{1}(\tau), L(\tau)\f_{2}(\tau)
    \big\rangle_{H}d\tau\nonumber
    \\
    && = \mE\big\langle P(t) \xi_{1},\xi_{2}
    \big\rangle_{H} + \mE \int_t^T \big\langle
    P(\tau)u_{1}(\tau),
    \f_{2}(\tau)\big\rangle_{H}d\tau \\
    &&  \q + \mE \int_t^T \big\langle
    P(\tau)\f_{1}(\tau),
    u_{2}(\tau)\big\rangle_{H}d\tau  + \mE
    \int_t^T\big\langle
    P(\tau)C(\tau)\f_{1}(\tau), D(\tau)
    v_{2}(\tau)\big\rangle_{H}d\tau \nonumber\\
    && \q+ \mE
    \int_t^T \big\langle  P(\tau)D(\tau) v_{1}(\tau), C(\tau)\f_{2}(\tau)+D (\tau)v_{2}(\tau)\big\rangle_{H}d\tau\nonumber\\
    && \q + \mE \int_t^T \big\langle
    v_{1}(\tau),
    D(\tau)^* \L(\tau)\f_2(\tau)\big\rangle_Ud\tau+ \mE
    \int_t^T \big\langle
    D(\tau)^* \L(\tau)\f_1(\tau), v_{2}(\tau)
    \big\rangle_Ud\tau.\nonumber
\end{eqnarray}
Here, $\f_1(\cd)$ and $\f_2(\cd)$ solve
\eqref{op-fsystem1.2} and \eqref{op-fsystem2.2}
with $v_1$ and $v_2$ replaced by $Dv_1$ and
$Dv_2$, respectively\footnote{By Theorem
    \ref{ch-1-well-mild}, one has $\f_1(\cd),
    \f_2(\cd)\in C_\dbF([0,T];L^2(\Om;H))$.}.
\end{definition}
\begin{theorem}\label{5.7-th1}
If the Riccati equation \eqref{5.5-eq6} admits a
transposition solution
$\big(P(\cd), \L(\cd)\big)\!\in\! D_{\dbF,w}([0,\!T];$ $
L^{\infty}(\Om;\cL(H)))\times
L^2_{\dbF,w}(0,T;\cL(H))$ such that
\begin{equation}\label{5.7-eq5}
    \begin{array}{ll}\ds
        K(\cd)^{-1}\big[B(\cd)^* P(\cd) +D(\cd)^*
        P(\cd)C(\cd) + D(\cd)^*\L(\cd)\big] \in
        \Upsilon_2(H;U),
    \end{array}
\end{equation}
then Problem (SLQ) has an optimal
feedback operator $\Th(\cd)\in \Upsilon_2(H;U)$.
In this case, the optimal feedback operator
$\Th(\cd)$ is given by
\begin{equation}\label{5.10}
    \begin{array}{ll}
        \ns\ds\Th(\cd)=-K(\cd)^{-1}[B(\cd)^* P(\cd)
        +D(\cd)^* P(\cd)C(\cd) + D(\cd)^*\L(\cd)].
    \end{array}
\end{equation}
Furthermore,
\begin{equation}\label{Value}
    \inf_{u\in
        \cU[0,T])}\cJ(\eta;u)=\frac{1}{2} \langle
    P(0)\eta,\eta\rangle_H.
\end{equation}
\end{theorem}

{\it Proof}. For any $\eta\in H$ and $u(\cd)\in
\cU[0,T]$, let $X(\cd)\equiv
X(\cd\,;\eta,u(\cd))$ be the corresponding state
for \eqref{LQsystem1}. Choose
$\xi_1=\xi_2=\eta$, $u_1=u_2=Bu$ and
$v_1=v_2=Du$ in
\eqref{op-fsystem1.2}--\eqref{op-fsystem2.2}.
From \eqref{9.7-eq10}, \eqref{10.10-eq10} and
the pointwise symmetry of $K(\cd)$, we obtain
that
\begin{eqnarray}\label{6.8-eq19}
    &&\mE\langle G X(T),X(T)\rangle_H + \mE
    \int_0^T \big\langle M(r)X(r), X(r)
    \big\rangle_{H}dr\nonumber  \\
    &&\q - \mE \int_0^T \big\langle
    \Th(r)^* K(r)\Th(r) X(r), X(r) \big\rangle_{H}dr
    \\
    && =\mE\big\langle P(0) \eta,\eta
    \big\rangle_{H}+\mE\int_0^T\big\langle
    P(r)B(r)u(r), X(r)\big\rangle_{H}dr \nonumber  \\
    &&\q +\mE \int_0^T\big\langle
    P(r)X(r), B(r)u(r)\big\rangle_{H}dr  +\mE
    \int_0^T  \big\langle P(r)C(r)X(r),
    D(r)u(r)\big\rangle_{H}dr\nonumber  \\
    &&\q+ \mE
    \int_0^T \big\langle  P(r)D(r)u(r), C(r)X(r) + D(r)u(r)\big\rangle_{H}dr\nonumber\\
    && \q + \mE \int_0^T \big\langle u(r),
    D(r)^*\L(r)X(r)\big\rangle_Udr+ \mE \int_0^T
    \big\langle D(r)^*\L(r)X(r), u(r)
    \big\rangle_Udr.\nonumber
\end{eqnarray}
Then, by  \eqref{6.8-eq19}, and recalling the
definition of $L(\cd)$ and $K(\cd)$, we arrive
at
\begin{eqnarray*}
    && 2\cJ(\eta;u(\cd))\\
    &&= \dbE\Big[\int_0^T\big(\big\langle
    MX(r),X(r)\big\rangle_H
    +\big\langle Ru(r),u(r)\big\rangle_U \big)dr+\langle GX(T),X(T)\rangle_H\Big]\\
    &&=\mE\big\langle P(0) \eta,\eta
    \big\rangle_{H}+\mE\int_0^T\big\langle PBu(r),
    X(r)\big\rangle_{H}dr+\mE \int_0^T\big\langle
    PX(r),
    Bu(r)\big\rangle_{H}dr \\
    && \q  +\mE\int_0^T \big\langle PCX(r),
    Du(r)\big\rangle_{H}dr + \mE
    \int_0^T\big\langle  PDu(r), CX(r)+Du(r)\big\rangle_{H}dr\\
    && \q + \mE \int_0^T \big\langle u(r),
    D^*\L(r)X(r)\big\rangle_Udr+ \mE \int_0^T
    \big\langle D^*\L(r)x(r), u(r)
    \big\rangle_Udr\\
    && \q  + \mE \int_0^T \big\langle \Th^*
    K \Th X(r), X(r) \big\rangle_{H}dr  +
    \mE\int_0^T \big\langle Ru(r),u(r)\big\rangle_U
    dr
    \\
    &&= \dbE\[\big\langle
    P(0)\eta,\eta\big\rangle_H\\
    &&\qq + \int_0^T \big( \big\langle \Th^*
    K\Th X(r),X(r)\big\rangle_H + 2\big\langle
    Lx(r),u(r)\big\rangle_U+\big\langle K
    u(r),u(r)\big\rangle_U\big)dr\].
\end{eqnarray*}
This, together with \eqref{5.10}, implies that
\begin{eqnarray*}\label{5.31-eq1}
    &&\cJ(\eta;u(\cd))  \\
    &&=\frac{1}{2}\dbE\Big[\big\langle
    P(0)\eta,\eta\big\rangle_H+\int_0^T\big(\big\langle
    K\Th X,\Th X\big\rangle_U-2\big\langle K\Th X,u\big\rangle_U+\langle Ku,u\rangle_U\big)dr\Big]\\
    &&=\frac{1}{2}\dbE\Big[\big\langle
    P(0)\eta,\eta\big\rangle_H+\int_0^T\big\langle K(u-\Th X),u-\Th X\big\rangle_Udr\Big]\\
    &&=\cJ\big(\eta;\Th X\big)
    +\frac{1}{2}\dbE\int_0^T\big\langle K(u-\Th
    X),u-\Th X\big\rangle_U dr.
\end{eqnarray*}
Hence,
$$
\cJ(\eta;\Th X)\leq \cJ(\eta;u),\q\forall\;
u(\cd)\in \cU[0,T].
$$
Consequently, $\Th(\cd)$ is an optimal feedback
operator for Problem (SLQ), and \eqref{Value}
holds. This completes the proof of Theorem \ref{5.7-th1}.
\endpf

\begin{remark}
In Definition \ref{4.8-def2}, we only ask that
$K(t,\om)$ has left inverse for a.e. $(t,\om)\in
(0,T)\times\Om$, and therefore $K(t,\om)^{-1}$
may be unbounded. Nevertheless, this result
cannot be improved. An example can be found in
\cite{LZ5}.
\end{remark}

By Theorem \ref{5.7-th1}, the existence of an optimal feedback operator is reduced to the existence of suitable transposition solutions to the Riccati equation \eqref{5.5-eq6}.
One may expect that \eqref{5.5-eq6} would admit
a transposition solution $(P,\L) $ without
further assumptions. Unfortunately, this is
incorrect even in finite dimensions, i.e.,
$H=\dbR^n$ (e.g., \cite[Example 6.2]{LWZ}).
So far, there is only some partial answer to this challenging problem (e.g., \cite{LZ5}).

%%%%%%%%%%%%%%%%%%%%%%%%%%%%%%%%%%%%%%%%%%%%%%%%%%%%%%%%%

%%%%%%%%%%%%%%%%%%%%%%%%%%%%%%%%%%%%%%%%%%%%%%%%%%%%%%%%

\section{Stochastic linear quadratic optimal control
problems III: deterministic coefficients}\label{sec-slq-d}

%%%%%%%%%%%%%%%%%%%%%%%%%%%%%%%%%%%%%%%%%%%%%%%%%%%%%%%%%%%

In this section, we study a special case of
Problem (SLQ), i.e., all coefficients appeared in the state equation \eqref{LQsystem1} and the cost functional \eqref{LQcost} are deterministic, and
$$\left\{\2n\begin{array}{ll}
\ns\ds A_1(\cd)=0,\q B(\cd), D(\cd)\in L^\infty(0,T;\cL(U;H)), \q
C(\cd)\in L^\infty(0,T;\cL(H)), \\
\ns\ds G\in\dbS(H),\q M(\cd)\in
L^1(0,T;\dbS(H)), \q R(\cd)\in
L^\infty(0,T;\dbS(U)).
\end{array}\right.$$
One can consider a more general case such as $B(\cd)\in L^2(0,T;\cL(U;H))$ and $C(\cd)\in L^2(0,T;\cL(H))$ (e.g. \cite{Luqi9}).

\subsection{Formulation of the problem and the main results}

Denote by  $C_\cS([0,T];$ $\dbS(H))$ the set of all
strongly continuous mappings $F:[0,T]\to
\dbS(H)$, that is, $F(\cd)\xi$ is continuous on
$[0,T]$ for each $\xi\in H$. A sequence
$\{F_n\}_{n=1}^\infty$ $\subset
C_\cS([0,T];\dbS(H))$ is said
to converge strongly to
$F\in C_\cS([0,T];\dbS(H))$ if
$$
\lim_{n\to\infty}F_n(\cd)\xi=F(\cd)\xi,\qq
\forall\; \xi\in H.
$$
In this case, we write
$
\ds\lim_{n\to\infty}F_n =F$  in
$C_\cS([0,T];\dbS(H))$.
If $F\in C_\cS([0,T];\dbS(H))$, then, by the
Uniform Boundedness Theorem (i.e., Theorem \ref{2.1ubth10}), the quantity
\bel{20210119e1}
|F|_{C_\cS([0,T];\dbS(H))}\deq \sup_{t\in
[0,T]}|F(t)|_{\cL(H)}
\ee
is finite.

Let $\cX$ and $\cY$ be two Banach spaces. For $1\leq p \leq\infty$, let
$$
\begin{array}{ll}\ds
L^{p,\cS}(0,T;\cL(\cX;\cY))\deq \{F:[0,T]\to \cL(\cX;\cY)\big| F\eta\in L^p(0,T;\cY),\;\forall \eta\in \cX, \\
\ns\ds\hspace{6.63cm} |F|_{\cL(\cX;\cY)}\in L^p(0,T)\}.
\end{array}
$$
Particularly,
$$
\begin{array}{ll}\ds
L^{p,\cS}(0,T;\dbS(H))\deq \{F:[0,T]\to \dbS(H)\big| F \in L^{p,\cS}(0,T;\cL(H))\}.
\end{array}
$$

As a special case of \eqref{5.5-eq6}, we introduce the following operator-valued Riccati equation:
\begin{equation}\label{Riccati}
\left\{\2n
\begin{array}{ll}
    \ns\ds\frac{dP}{dt} +P A + A^* P +C^* PC
    +M-L^* K^{-1} L=0 &\mbox{ in }[0,T),\\
    \ns\ds P(T)=G,
\end{array}
\right.
\end{equation}
where
\bel{20210119e2}
L(\cd)=B(\cd)^*P(\cd)+D(\cd)^* P(\cd)C(\cd),\qq
K(\cd)=R(\cd)+D(\cd)^* P(\cd)D(\cd).
\ee

Solutions to \eqref{Riccati} are understood in the following sense.
\begin{definition}\label{def3}
We call $P\in C_\cS([0,T];\dbS(H))$ a {\it strongly regular  mild
    solution} to \eqref{Riccati} if for any $\eta\in
H$ and $s\in [0,T]$,
\begin{equation}\label{8.20-eq23}
    \begin{array}{ll}\ds
        P(s)\eta\3n&\ds=S(T-s)^* GS(T-s) \eta \\
        \ns&\ds\q +
        \int_s^T S(\tau-s)^* \big(C^* PC + M -
        L^*K^{-1} L\big)S(\tau-s) \eta d\tau
    \end{array}
\end{equation}
and
\begin{equation}\label{strong-regular}
    K(s)\geq c_0 I,\qq\ae~s\in[0,T],
\end{equation}
for some $c_0>0$.
\end{definition}
\begin{remark}
To define the mild solution to \eqref{Riccati}, we only need the equality \eqref{8.20-eq23} make sense. To this end,  \eqref{strong-regular} is unnecessary (e.g. \cite{Luqi9}). Nevertheless, on order to find the desired optimal feedback control for the corresponding Problem (SLQ), it is more convenient to consider more ``regular" solutions to \eqref{Riccati}. A natural one is the so called ``regular mild solutions" (e.g. \cite{Luqi9}). In this notes, to avoid too many technical details, we consider only strongly regular solutions introduced in Definition \ref{def3}.
\end{remark}

Similarly to the proof of Theorem \ref{5.7-th1}, one can show the following
result (and hence we omit its proof).

\begin{theorem}\label{cor4.8}
If the equation \eqref{Riccati} admits a  strongly regular mild solution $P\in C_\cS([0,T];\dbS(H))$, then,  Problem {\rm(SLQ)}
has an optimal feedback control
\begin{equation}\label{opti-biaoshi}
    \bar u(\cd)=-K(\cd)^{-1}L(\cd)\cl X(\cd),
\end{equation}
where $\cl X(\cd)$ solves the
following closed-loop system (See \eqref{20210119e2} for $K$ and $L$):
\begin{equation}\label{closed-loop-state}
    \left\{\2n
    \begin{array}{ll}
        \ds d\cl X = \big(A -BK^{-1}L\big)\cl Xdt +
        \big(C-DK^{-1}L\big)\cl X dW(t) &\mbox{ in }(0,T], \\
        \ns\ds \cl X(0)=\eta.
    \end{array}
    \right.
\end{equation}
\end{theorem}

Further, we have the following two results.

\begin{theorem}\label{12.6-th1}
The Riccati equation \eqref{Riccati} admits at most one strongly regular  mild solution.
\end{theorem}
\begin{theorem}\label{th4.6}
The Riccati equation
{\rm(\ref{Riccati})} admits a strongly regular  solution if and only if the map $u(\cd)\mapsto \cJ(0;u(\cd))$
is uniformly convex, i.e., for some constant $c_1>0$,
\begin{equation}\label{J>l*}
    \cJ(0;u(\cd))\geq c_1\,\dbE\1n\int_0^T|u(s)|_U^2ds,\qq\forall\;
    u(\cd)\in\cU[0,T].
\end{equation}
\end{theorem}
\begin{remark}
Clearly, if \eqref{4.7-eq3} holds, then the map
$u(\cd)\mapsto \cJ(0;u(\cd))$ is uniformly
convex. On the other hand, there are some
interesting cases for which such a map is uniformly convex but \eqref{4.7-eq3}
does not hold (e.g., \cite{Luqi9}).
\end{remark}

The proofs of Theorems \ref{12.6-th1}--\ref{th4.6} will be given in Subsection \ref{sec-pr-main2}.

\br
By \eqref{LQoperator}, for any  $(\eta,u(\cd))\in H\times\cU[0,T]$, the
cost functional can be written as
$$
\ds \cJ(\eta;u(\cd))=\lan G\big(\widehat\G
\eta+\widehat \Xi u\big),\widehat\G
\eta+\widehat \Xi u\rangle_H +\lan M(\G \eta+\Xi
u),\G \eta+\Xi u\rangle_H +\langle Ru,u
\rangle_U.
$$
Hence, $u(\cd)\mapsto \cJ(0;u(\cd))$ is
uniformly convex if and only if for some $c_1>0$,
\begin{equation}\label{M_2>l}
\widehat \Xi^*G\widehat \Xi+\Xi^*M\Xi+R\geq c_1
I.
\end{equation}

If  $R \geq \d
I$ for some $\d>0$, then \eqref{M_2>l} holds for $c_0=\d$.
On the other hand, even if $R \geq \d I$ does not
hold, \eqref{M_2>l} may still be true when  $G$ is
large enough.  Such an example can be found in \cite{Luqi9}.
\er

%%%%%%%%%%%%%%%%%%%%%%%%%%%%%%%%%%%%%%%%%%%%%%%

\subsection{Some preliminaries}\label{sec-pre}

%%%%%%%%%%%%%%%%%%%%%%%%%%%%%%%%%%%%%%%%%%%%%%%
In order to prove Theorems \ref{12.6-th1}--\ref{th4.6}, we need some preliminary results.

Let us consider the following (operator-valued) Lyapunov equation:
\begin{equation}\label{lm2.4-eq1}
\left\{ \2n\begin{array}{ll}
    \ns\ds \frac{d \wt P}{dt} + \wt P (A+\wt A )+(A+\wt A )^*  \wt P
    +\wt C^*  \wt P\wt C+\wt M=0 &\mbox{ in }[0,T),\\
    \ns\ds  \wt P(T)=\wt G,
\end{array}
\right.
\end{equation}
where $\wt A(\cd)\in L^{1,\cS}(0,T;\cL(H))$, $\wt
C(\cd)\in L^{2,\cS}(0,T;\cL(H))$, $\wt M(\cd)\in L^{1,\cS}(0,T;\dbS(H))$ and $\wt
G\in\dbS(H)$.
We call $ \wt P\in C_\cS([0,T];\dbS(H))$ a {\it mild
solution} to \eqref{lm2.4-eq1} if for any $s\in
[0,T]$ and $\eta\in H$,
\begin{equation}\label{8.20-eq22}
\begin{array}{ll}\ds
    \wt P(s)\eta\3n&\ds=S(T-s)^* \wt GS(T-s)\eta\\
    \ns&\ds\q + \int_s^T
    S(\tau-s)^*\big( \wt P \wt A +\wt A^*
    \wt P +\wt C^*  \wt P\wt C+\wt M\big)S(\tau-s)\eta d\tau.
\end{array}
\end{equation}

We have the following well-posedness result for the equation \eqref{lm2.4-eq1}.

\begin{lemma}\label{12.7-lm2.4}
The equation \eqref{lm2.4-eq1} admits a unique
mild solution $P(\cd)\in C_\cS([0,T];$ $\dbS(H))$.  Moreover,
\begin{equation}\label{12.7-lm2.4-eq1}
    \begin{array}{ll}\ds
        | \wt P|_{C_\cS([0,T];\dbS(H))} \leq  \cC e^{\int_0^T (2|\wt A|_{\cL(H)}
            +|\wt C|_{\cL(H)}^2 ) ds}\(| \wt G|_{\cL(H)} +
        \int_0^T
        |\wt M|_{\cL(H)} ds\).
    \end{array}
\end{equation}
\end{lemma}

{\it Proof}. Let
$\cM_T\deq\sup_{t\in[0,T]}|S(t)|_{\cL(H)}$. Let
${T_0}\in [0,T)$ such that
$$
\int_{T_0}^T \big(2 |\wt A|_{\cL(H)}
+|\wt C|_{\cL(H)}^2\big) ds
<\frac{1}{2\cM_T^2}.
$$
Define a map $\cG:C_\cS([{T_0},T];\dbS(H))\to
C_\cS([{T_0},T];\dbS(H))$ as follows:
$$
\begin{array}{ll}\ds
    \cG( \wt P)(r)\zeta\3n&\ds\deq S(T-r)^*GS(T-r) \eta\\
    \ns&\ds\q +
    \int_r^T S(s-r)^* \big( \wt P\wt A+\wt A^*
    \wt P +\wt C^*  \wt P\wt C+\wt M\big)S(s-r) \zeta ds,\q  \forall\,\zeta\in H.
\end{array}
$$

Let $ \wt P_1, \wt P_2\in C_\cS([0,T];\dbS(H))$. For
each $\zeta\in H$,
$$
\begin{array}{ll}  \ds
    \sup_{r\in[{T_0},T]}\big|\big(\cG( \wt P_1)-\cG( \wt P_2)\big)(r)\zeta\big|_H\\
    \ns\ds= \sup_{r\in[{T_0},T]} \Big| \int_r^T
    S(s-r)^* [( \wt P_1- \wt P_2) \wt A + \wt A^*
    ( \wt P_1- \wt P_2) \\
    \ns\ds\qq\qq\q+ \wt C^*
    ( \wt P_1 -  \wt P_2)\wt C]S(s - r) \zeta ds\Big|_H\\
    \ns\ds \leq \int_{T_0}^T
    \big|S(s-r)^*\big|_{\cL(H)}\big[\big(\big|\wt A\big|_{\cL(H)}+\big|\wt A^*\big|_{\cL(H)}\big)\big|( \wt P_1- \wt P_2)\big|_{\cL(H)}
    \\
    \ns\ds\qq\q\;+\big|\wt C^*\big|_{\cL(H)}
    \big| \wt P_1 - \wt P_2\big|_{\cL(H)}\big|\wt C\big|_{\cL(H)}\big]\big|S(s-r)\big|_{\cL(H)}\big|\zeta\big|_H
    ds\\
    \ns\ds\leq \cM_T^2\int_{T_0}^T
    \big(2\big|\wt A\big|_{\cL(H)}
    +\big|\wt C\big|_{\cL(H)}^2\big) ds
    \sup_{r\in[{T_0},T]}\big|( \wt P_1- \wt P_2)(r)\big|_{\cL(H)}
    \big|\zeta\big|_H\\
    \ns\ds \leq
    \frac{1}{2}\sup_{r\in[{T_0},T]}\big|( \wt P_1- \wt P_2)(r)\big|_{\cL(H)}
    \big|\zeta\big|_H.
\end{array}
$$
Therefore,
\begin{equation}\label{2.9-eq41}
    \begin{array}{ll}\ds
        \sup_{r\in[{T_0},T]}\big|\big(\cG( \wt P_1)-\cG( \wt P_2)\big)(r)\big|_{\cL(H)}\leq
        \frac{1}{2}\sup_{r\in[{T_0},T]}\big|( \wt P_1- \wt P_2)(r)\big|_{\cL(H)}.
    \end{array}
\end{equation}
This implies that $\cG$ is contractive.
Consequently, there is a unique fixed point of
$\cG$, which is the mild solution to
\eqref{lm2.4-eq1} on $[{T_0},T]$. Repeating this
process gives us the unique $ \wt P\in
C_\cS([0,T];\dbS(H))$ which satisfies
\eqref{8.20-eq22}. The uniqueness of the
solution is obvious.

\vspace{0.1cm} From \eqref{8.20-eq22}, we see
that for any $r\in[0,T]$ and $\zeta\in H$,
$$
\begin{array}{ll}\ds
    \big| \wt P(r)\zeta\big|_H\leq
    \big|S(T-r)\big|_{\cL(H)}^2\big| \wt G\big|_{\cL(H)}
    \big|\zeta\big|_H
    + \int_r^T
    \big|S(s-r)\big|_{\cL(H)}^2
    \big[\big(2\big|\wt A\big|_{\cL(H)}\\
    \ns \ds\qq\qq\q\,+\big|\wt C \big|_{\cL(H)}^2
    \big)\big| \wt P\big|_{\cL(H)}
    +\big|\wt M\big|_{\cL(H)}\big] \big|\zeta\big|_H ds.
\end{array}
$$
Consequently,
$$
\big| \wt P(r)\big|_H\leq
\cC\Big\{\big| \wt G\big|_{\cL(H)}
+ \int_r^T
\big[\big(2\big|\wt A\big|_{\cL(H)}+\big|\wt C \big|_{\cL(H)}^2
\big)\big| \wt P\big|_{\cL(H)}
+\big|\wt M\big|_{\cL(H)}\big] ds\Big\}.$$
This, together with Gronwall's inequality,
implies \eqref{12.7-lm2.4-eq1}.
\endpf

\ms

The following result illustrates the
differentiability of mild solutions to  \eqref{lm2.4-eq1}.
\begin{proposition}\label{prop1}
Let $ \wt P$ be a mild solution to \eqref{lm2.4-eq1}.
Then for any $\eta,\zeta\in D(A)$, $\lan
\wt P(\cd)\eta, \zeta\ran_H$ is differentiable in
$[0,T]$ and
\begin{equation}\label{8.20-eq20}
    \begin{array}{ll}\ds
        \frac{d}{dt}\lan  \wt P\eta, \zeta\rangle_H
        \3n&\ds=-\lan  \wt P\eta, (A +\wt A)\zeta\rangle_H -
        \lan  \wt P(A +\wt A)\eta, \zeta\rangle_H\\
        \ns&\ds\q - \lan  \wt P\wt C\eta,  \wt C\zeta
        \rangle_H  - \lan \wt M \eta, \zeta \rangle_H.
    \end{array}
\end{equation}
\end{proposition}
{\it Proof}. For any $\eta,\zeta\in H$, we have
that
\begin{equation}\label{12.7-eq6}
    \begin{array}{ll}\ds
        \lan  \wt P(r)\eta,\zeta\rangle_H\3n&
        \ds=\lan  \wt GS(T-r) \eta, S(T-r) \zeta\rangle_H
        \\
        \ns&\ds\q+ \int_r^T \lan \big( \wt P\wt A+\wt A^*
        \wt P +\wt C^*  \wt P\wt C+\wt M\big)S(s-r) \eta, S(s-r)\zeta\rangle_H ds.
    \end{array}
\end{equation}
If $\eta,\zeta\in D(A)$, it follows \eqref{12.7-eq6} that $\lan
\wt P(r)\eta,\zeta\rangle_H$ is differentiable w.r.t. $r$. A simple computation gives
\eqref{8.20-eq20}.
\endpf

Now, let us show the following result:

\begin{lemma}\label{lm2.4}
If
\begin{equation}\label{lm2.4-eq2}
    \wt G\geq0,\qq\wt M(t)\geq0,\qq \ae~t\in[0,T],
\end{equation}
then the mild solution $ \wt P(\cd)$ to \eqref{lm2.4-eq1} satisfies that $ \wt P(t)\geq 0$ for all $t\in [0,T]$.
\end{lemma}

{\it Proof}. Let $t\in [0,T)$. Consider the
following equation:
$$
\left\{
\begin{array}{ll}\ds
    d\wt X =(A+\wt A )\wt Xds+\wt C \wt X
    dW(s) &\mbox{ in }(t,T],\cr
    \ns\ds \wt X(t)=\eta\in H.
\end{array}
\right.
$$
Clearly this equation admits a unique solution $\wt X(\cd)$.

For each
$\l\in\rho(A)$, the resolvent of $A$, write
$\wt X_\l(\cd)\deq R(\l) \wt X(\cd)$, where $R(\l)=  \l I(\l I-A)^{-1}$. Then $\wt X_\l(\cd)$ solves
the following equation:
$$
\left\{\2n
\begin{array}{ll}
    \ds d\wt X_\l =\big(A\wt X_\l + R(\l)\wt A \wt X\big)ds +R(\l)\wt C\wt X dW(s) &\mbox{in }(t,T], \\
    \ns\ds \wt X_\l(t)=R(\l)\eta.
\end{array}
\right.
$$
By It\^o's formula and Proposition \ref{prop1},
we have
\begin{eqnarray}\label{12.7-eq9}
    \begin{array}{ll}\ds
        \dbE \lan \wt G
        \wt X_\l(T),\wt X_\l(T)\ran_H-\dbE \lan \wt P(t) R(\l)\eta,R(\l)\eta\ran_H\\
        \ns\ds = -\dbE \int_t^T \big[\lan  \wt P\wt X_\l, (A +\wt A)\wt X_\l\ran_H +
        \lan  \wt P(A +\wt A)\wt X_\l, \wt X_\l\ran_H\\
        \ns \ds\qq\qq\ + \lan  \wt P\wt C\wt X_\l, \wt  P\wt C\wt X_\l
        \ran_H + \lan
        \wt M(s)\wt X_\l,\wt X_\l \ran_H\big] ds \\
        \ns \ds\q  +  \dbE \int_t^T \big[\lan  \wt P\wt X_\l, (A\wt X_\l +R(\l)\wt A\wt X)\rangle_H +
        \lan  \wt P(A\wt X_\l +R(\l)\wt A\wt X), \wt X_\l\ran_H\\
        \ns \ds\qq\qq\ + \lan  \wt PR(\l)\wt C\wt X, R(\l)\wt C\wt X
        \rangle_H \big] ds\\
        \ns\ds = -\dbE \int_t^T \big[\lan  \wt P\wt X_\l,  \wt A \wt X_\l-R(\l)\wt A\wt X\rangle_H +
        \lan  \wt P \big(\wt A \wt X_\l-R(\l)\wt A\wt X\big), \wt X_\l\ran_H\\
        \ns \ds\qq\qq\ + \lan  \wt P\wt C\wt X_\l,  \wt P\wt C\wt X_\l
        \rangle_H - \lan  \wt PR(\l)\wt C\wt X, R(\l)\wt C\wt X
        \ran_H + \lan
        \wt M \wt X_\l,\wt X_\l \rangle_H\big] ds.
    \end{array}
\end{eqnarray}

By Theorem \ref{ch-2-app1}, we have that
\begin{equation}\label{8.20-eq26}
    \lim_{\l\to\infty} \wt X_\l=\wt X \q\mbox{ in
    }\;C_\dbF([t,T];L^2(\Om;H)).
\end{equation}
Noting that for any $\zeta\in H$,
\begin{equation}\label{8.20-eq26.1}
    \lim_{\l\to\infty}R(\l)\zeta = \zeta \q\mbox{
        in }H,
\end{equation}
we
get from \eqref{8.20-eq26} that for a.e. $s\in [t,T]$,
\begin{equation}\label{8.20-eq27}
    \begin{array}{ll}\ds
        \lim_{\l\to\infty}\big(\lan  \wt P(s)\wt X_\l(s),  \wt A(s) \wt X_\l(s)-R(\l)\wt A(s)\wt X(s)\ran_H\big)=0,\q \as
    \end{array}
\end{equation}
It follows from the definition of $R(\l)$
that
$$
\begin{array}{ll}\ds
    \big|\lan  \wt P(s)\wt X_\l(s),  \wt A(s) \wt X_\l(s)-R(\l)\wt A(s)\wt X(s)\ran_H\big|\\
    \ns\ds\leq \cC | \wt P(s)|_{\cL(H)}
    |\wt A(s)|_{\cL(H)} \big(|\wt X(s)|_H^2+1\big),\qq \as
\end{array}
$$
This, together with \eqref{8.20-eq27} and
Dominated Convergence Theorem (Theorem \ref{Ap-th4}),
implies that
$$
\lim_{\l\to\infty}\dbE \int_t^T \lan  \wt P\wt X_\l,  \wt A \wt X_\l-R(\l)\wt A\wt X\ran_H ds = 0.
$$
By a similar argument, letting $\l\to\infty$ in both sides of
\eqref{12.7-eq9}, we get
$$
\lan \wt P(t) \eta,\eta\ran_H =\dbE \lan \wt G
\wt X(T),\wt X(T)\ran_H + \dbE \int_t^T \lan
\wt M \wt X,\wt X \ran_H ds,\qq
\forall\; t\in [0,T].
$$
This, together with \eqref{lm2.4-eq2}, implies
that $\wt P(t)\geq 0$ for all $t\in [0,T]$.
\endpf

\begin{remark}
Since $\wt X(\cd)$ may not be $D(A)$-valued, in the
proof of Lemma \ref{lm2.4}, we introduce a
family of $\{\wt X_\l(\cd)\}_{\l\in\rho(A)}$ as Theorem \ref{ch-2-app1} to apply It\^o's formula and
Proposition \ref{prop1}. In the rest of this
section, we omit such procedures to save space. The readers are encouraged to give the omitted part themselves.
\end{remark}

Similarly to the proof of Proposition \ref{prop1}, we can obtain the following result.
\begin{proposition}\label{prop2}
Let $P$ be a strongly regular mild solution to \eqref{Riccati}.
Then for any $\eta,\zeta\in D(A)$, $\lan
P(\cd)\eta, \zeta\rangle_H$ is differentiable in
$[0,T]$ and
\begin{equation}\label{8.20-eq24}
    \begin{array}{ll}\ds
        \frac{d}{dt}\lan P\eta,
        \zeta\ran_H\3n&\ds=-\lan P\eta,  A
        \zeta\ran_H - \lan P A \eta, \zeta\ran -
        \lan P C \eta, C \zeta \ran_H \\
        \ns&\ds\q - \lan M \eta, \zeta \ran_H+\lan
        K^{-1}L \eta, L \zeta \ran_H.
    \end{array}
\end{equation}
\end{proposition}

Now, for any $\Th(\cd)\in L^{2,\cS}(0,T;\cL(H;U))$, let us consider
the following (operator-valued) Lyapunov equation:
\begin{equation}\label{eq-Lya}
\left\{\2n\begin{array}{ll}
    \ds\frac{dP_\Th}{dt}+P_\Th(A + B\Th)+(A + B\Th)^* P_\Th \\
    \ns\ds \q+(C + D\Th)^* P_\Th(C+D\Th)+ \Th^* R\Th  + M =0 &\mbox{ in }[0,T),\\
    \ns\ds P_\Th(T)=G.
\end{array}
\right.
\end{equation}

As an immediate consequence of Lemma \ref{12.7-lm2.4}, we have the following result.
\begin{corollary}\label{prop3}
There is a unique mild solution to
\eqref{eq-Lya}. Moreover,
\begin{equation}\label{prop3-eq1}
    \begin{array}{ll}\ds
        |P_\Th|_{C_\cS([0,T];\dbS(H))}\\
        \ns\ds \leq  \cC e^{\int_0^T (2|B\Th|_{\cL(H)}
            +|C + D\Th|_{\cL(H)}^2 ) dt}\[|G|_{\cL(H)}\! +
        \int_0^T\! \big( |\Th|_{\cL(H;U)}^2 |R|_{\cL(U)} \!+
        |M|_{\cL(H)}\big)dt\].
    \end{array}
\end{equation}
\end{corollary}

For any $r\in [0,T)$ and
$\xi\in L^2_{\cF_r}(\Om;H)$, we consider
the following control system:
\begin{equation}\label{state-r}
\left\{\2n
\begin{array}{ll}
    \ds dX =\big(AX   +B u
    \big)dt +\big(CX
    +D u \big)dW(t) &\mbox{ in }(r,T],\\
    \ns\ds X(r)=\xi
\end{array}
\right.
\end{equation}
with  the cost functional
\begin{equation}\label{cost-r}
\ds \cJ(r,\xi;u(\cd))
\deq\frac{1}{2}\dbE\Big[\int_r^T\big(\lan M X, X\ran_H +
\lan
Ru, u \ran_U \big)ds+\lan
GX(T),X(T)\ran_H\],
\end{equation}
where $u\in \cU[r,T]\deq L^2_\dbF(r,T;U)$.
We need to introduce the following  optimal control problem (parameterized by $r\in [0,T)$):

\ms

\bf Problem (SLQ-$r$). \rm Find a control $\bar u(\cd)\in\cU[r,T]$ such
that
\begin{equation}\label{optim-r}
\cJ(r,\xi;\bar
u(\cd))=\inf_{u(\cd)\in\cU[r,T]}\cJ(r,\xi;u(\cd)).
\end{equation}

Clearly, if the map $u(\cd)\mapsto \cJ(0;u(\cd))$
is uniformly convex, then so is the map $u(\cd)\mapsto \cJ(r,0;u(\cd))$, i.e., for some $c_0>0$,
\begin{equation}\label{J>l*-r}
\cJ(r,0;u(\cd))\geq c_0 \dbE \int_r^T|u(s)|_U^2ds,\qq\forall\;
u(\cd)\in\cU[r,T].
\end{equation}

The  result below gives a relation between the
cost functional \eqref{cost-r} and the Lyapunov equation
\eqref{eq-Lya}.

\begin{lemma}\label{lm2.3}
Let $P_\Th(\cd)$ solve \eqref{eq-Lya}. Then,
\begin{equation}\label{8.20-eq29}
    \begin{array}{ll}
        \ns\ds \cJ(r,\xi;\Th(\cd)X(\cd)+u(\cd))\\
        \ns\ds =\dbE\int_r^T\big[2\lan\big(L_\Th+K_\Th\Th\big)X,u\ran_U
        +\lan K_\Th u,u\ran_U \big]dt+\dbE\lan P_\Th(r)\xi,\xi\ran_H,
    \end{array}
\end{equation}
where
\bel{20210120e1}
K_\Th(\cd)\deq R(\cd)+D(\cd)P_\Th(\cd) D(\cd),\qq L_\Th(\cd)\deq B(\cd)^*P_\Th(\cd) + D(\cd)^*P_\Th(\cd) C(\cd).
\ee
\end{lemma}

{\it Proof}. For any $\eta\in H$ and
$u(\cd)\in\cU[r,T]$, let $X(\cd)$ be the
solution to
$$
\left\{
\begin{array}{ll}
    \ns\ds dX =\big[(A +B\Th)X+Bu\big]dt
    +\big[(C+D\Th)X+Du\big]dW(t) &\mbox{ in }(r,T], \\
    \ns\ds X(r)=\xi.
\end{array}
\right.
$$
By It\^o's formula and Proposition \ref{prop1}, proceeding as \eqref{12.7-eq9} in the proof of Lemma \ref{lm2.4},
we can obtain that
$$
\ba{ll}
\ds \mE\[\int_r^T\big(\lan MX,X\ran_H + \lan
R(\Th X+u),\Th X+u\ran_U\big)
dt+\lan GX(T),X(T)\ran_H\]\\\ns\ds=  \dbE\int_r^T\big[2\lan\big(L_\Th+K_\Th\Th\big)X,u\ran_U
+\lan K_\Th u,u\ran_U\big]dt+\dbE\lan P_\Th(r)\xi,\xi\rangle_H,
\ea
$$
which implies \eqref{8.20-eq29}.
\endpf

Next, we give a result for the existence and uniqueness
of optimal controls for Problem (SLQ-$r$).
\begin{proposition}\label{prop4.1}
Suppose the map $u(\cd)\mapsto \cJ(r,0;u(\cd))$ is
uniformly convex. Then Problem {\rm(SLQ-$r$)} admits
a unique optimal control, and for some
constant $\a\in\dbR$,
\begin{equation}\label{uni-convex-prop0}
    \inf_{u\in\cU[r,T]}\cJ(r,\xi;u)\geq\a\mE|\xi|_H^2,\qq\forall\;
    \xi\in L^2_{\cF_r}(\Om;H).
\end{equation}
\end{proposition}

{\it Proof}. Clearly, for some $c_0>0$,
\begin{equation}\label{J>l*.1}
    \cJ(r,0;u(\cd))\geq c_0
    \dbE\1n\int_r^T|u(t)|_U^2dt,\qq\forall\;
    u(\cd)\in\cU[r,T].
\end{equation}
Denote by $X_0$ (\resp $X_1$) the solution to \eqref{state-r}  with $u=0$ (\resp $\xi=0$). Then, $X=X_0+X_1$, and
\begin{equation}\label{J-rep}
    \begin{array}{ll}
        \ns\ds \cJ(r,\xi;u(\cd)) \\
        \ns \ds = \cJ(r,\xi;0)+\cJ(r,0;
        u(\cd))+\int_r^T\lan MX_0,X_1\ran_U
        dt+ \lan G X_0(T),X_1(T)\ran_H\\
        \ns\ds\geq \cJ(r,\xi;0)+ c_0\dbE\int_r^T|u(t)|_U^2dt- \frac{c_0}{2}\dbE
        \int_r^T|u(t)|_U^2dt - \cC\dbE \int_r^T
        |X_0|_H^2dt\\
        \ns\ds\geq \cJ(r,\xi;0)+\frac{c_0}{2}\dbE\int_0^T|u(t)|_U^2dt
        -\cC\mE|\xi|_H^2dt,\q\forall\;
        \xi\in L^2_{\cF_r}(\Om;H).
    \end{array}
\end{equation}
This implies that $u(\cd)\mapsto
\cJ(r,\xi;u(\cd))$ is coercive. Clearly,
$u(\cd)\mapsto \cJ(r,\xi;u(\cd))$ is continuous and
convex. Consequently, by a standard argument involving a minimizing
sequence and locally weak compactness of Hilbert
spaces, this functional has a unique minimizer.
Moreover,  \eqref{J-rep} implies that
\begin{equation}\label{uni-convex-prop2}
    \inf_{u\in\cU[0,T]}\cJ(r,\xi;u)\geq
    \cJ(r,\xi;0)-\cC\mE|\xi|_H^2dt.
\end{equation}
Since the terms in the right-hand side of
(\ref{uni-convex-prop2}) are quadratic in $\xi$,
we get (\ref{uni-convex-prop0}).
\endpf

The next result shows that the solution to
\eqref{eq-Lya} is bounded from below.
\begin{proposition}\label{prop4.5}
Let  {\rm(\ref{J>l*})} hold. Then for any
$\Th(\cd)\in L^{2,\cS}(0,T;\cL(U;H))$, the mild solution
$P_\Th(\cd)$ to
\eqref{eq-Lya} and the process $K_\Theta(\cd)$ defined by \eqref{20210120e1} satisfy
\begin{equation}\label{Convex-prop-1}
    K_\Th(t)\geq c_0 I, \q\ae~t\in[0,T]
\end{equation}
and
\begin{equation}\label{Convex-prop-1.1}
    P_\Th(t)\geq \a I,\q\forall\; t\in[0,T],
\end{equation}
where $c_0>0$ and $\a\in\dbR$ are the constants appearing in \eqref{J>l*} and
\eqref{uni-convex-prop0}, respectively.
\end{proposition}

{\it Proof}.   For any $u(\cd)\in\cU[0,T]$, let
$X_0(\cd)$ be the solution to
$$
\left\{\2n
\begin{array}{ll}
    \ns\ds dX_0 =\big[(A +B\Th) X_0+Bu\big]dt+\big[(C+D\Th)X_0+Du\big]dW(t) & \mbox{ in }(0,T], \\
    \ns\ds X_0(0)=0.
\end{array}
\right.
$$
It follows from {\rm(\ref{J>l*})} and Lemma
\ref{lm2.3} that
$$
\begin{array}{ll}\ds
    c_0\dbE \int_0^T |\Th X_0 + u |_U^2dt \\
    \ns\ds \leq \cJ(0;\Th(\cd)X_0(\cd) + u(\cd)) = \dbE \int_0^T \big[2\lan\big(L_\Th +
    K_\Th\Th\big)X_0,u\rangle_U + \lan
    K_\Th u,u\rangle_U\big]dt,
\end{array}
$$
which yields that, for any $u(\cd)\in\cU[0,T]$,
\begin{equation}\label{P>LI}
    \begin{array}{ll}
        \ds\dbE\int_0^T\big\{2\lan\big[L_\Th+(K_\Th-c_0
        I)\Th\big]X_0,u\ran_U+\lan (K_\Th-c_0
        I)u,u\ran_U\big\}dt\\
        \ns\ds=c_0\dbE\int_0^T|\Th(t)X_0(t)|_U^2dt\geq0.
    \end{array}
\end{equation}
Fix any $u_0\in U$ and $t_0\in (0,T)$, and choose any $h>0$ such that $t_0+h\leq T$. Take
$u(\cd)=u_0\chi_{[t_0,t_0+h]}(\cd)$. Then
\begin{equation}\label{8.20-eq35}
    |X_0|_{C_\dbF([0,T];L^2(\Om;H))}\leq
    \cC|u|_{L^2_\dbF(0,T;U)}\leq
    \cC\sqrt{h}|u_0|_{U}.
\end{equation}
Dividing both sides of \eqref{P>LI} by $h$ and
then letting $h\to 0$, noting \eqref{8.20-eq35},
we obtain
$$\lan\big(K_\Th(t_0)-c_0 I\big)u_0,u_0\rangle_U\geq 0,\qq\ae~t_0\in (0,T),\q \forall\; u_0\in U.$$
This gives
(\ref{Convex-prop-1}).

\ss

Now we prove (\ref{Convex-prop-1.1}).
For any $\Th(\cd)\in L^{2,\cS}(0,T;\cL(H;U))$ and $\xi\in
L^2_{\cF_r}(\Om; H)$, let $X(\cd)$ be
the solution of the following closed-loop
system:
\begin{equation}\label{Feb12-01-r}
    \left\{\2n\begin{array}{ll}
        \ds dX= \big(A +B\Th \big)X  dt + \big(C +D \Th \big)X dW(t)&\mbox{ in }(r,T],\\
        \ns\ds X(r)=\xi.
    \end{array}
    \right.
\end{equation}
It follows from Proposition \ref{prop4.1} and
Lemma \ref{lm2.3} that
$$
\begin{array}{ll}
    \ds \a\mE|\xi|^2 \3n&\ds\leq
    \cJ(r,\xi;\Th(\cd)X(\cd)) =\dbE\lan P_\Th(r)\xi,\xi
    \ran_H,\q\forall\;\xi\in
    L^2_{\cF_r}(\Om; H).
\end{array}
$$
In particular,  $ \lan P_\Th(r)\eta,\eta\ran_H \geq\a |\eta|^2$ for all  $\eta\in H$,
which yields (\ref{Convex-prop-1.1}).
\endpf

We shall also need the following result:
\begin{lemma}\label{lm2.5}
For any $u(\cd)\in\cU[0,T]$, let $X$ be the
corresponding solution to \eqref{LQsystem1} with
$\eta = 0$.  Then for every $\Th(\cd)\in
L^{2,\cS}(0,T;\cL(H;U))$, there exists a constant
$c_\Th>0$ such that
\begin{equation}\label{lem-2.6}
    \dbE \int_0^T\big|u(s)-\Th(s) X(s)\big|_U^2ds
    \geq c_\Th\dbE\int_0^T|u(s)|_U^2ds,\q\forall\;
    u(\cd)\in\cU[0,T].
\end{equation}
\end{lemma}

{\it Proof}. Define a bounded linear operator
$\mathfrak{L}:\cU[0,T]\to\cU[0,T]$ by
$\mathfrak{L}u=u-\Th X$. Then $\mathfrak{L}$ is
bijective and its inverse $\mathfrak{L}^{-1}$ is
given by $\mathfrak{L}^{-1}u=u+\Th\wt X$,
where $\wt X(\cd)$ is the solution to
$$
\left\{\2n
\begin{array}{ll}
    \ns\ds d\wt X =\big[\big(A  +B \Th \big)\wt X +B u \big]ds+ \big[\big(C +D \Th \big)\wt X +D u \big]dW(s) &\mbox{ in }(0,T], \\
    \ns\ds\wt X(0)=0.
\end{array}
\right.
$$
By the Inverse Mapping Theorem (i.e., Theorem \ref{2.1-th10}),
$\mathfrak{L}^{-1}$ is bounded with the norm
$|\mathfrak{L}^{-1}|_{\cL(\cU[0,T])}>0$. Thus,
$$\begin{array}{ll}
    \ns\ds\dbE\int_0^T|u(s)|_U^2ds
    \3n&\ds=\dbE\int_0^T|(\mathfrak{L}^{-1}\mathfrak{L}u)(s)|_U^2ds
    \leq|\mathfrak{L}^{-1}|_{\cL(\cU[t,T])}\dbE\int_0^T|(\mathfrak{L}u)(s)|_U^2ds\\
    \ns&\ds
    =|\mathfrak{L}^{-1}|_{\cL(\cU[t,T])}\dbE\int_0^T\big|u(s)-\Th(s)
    X(s)\big|_U^2ds, \q \forall\;
    u(\cd)\in\cU[0,T],\end{array}$$
which implies (\ref{lem-2.6}) with
$c_\Th=|\mathfrak{L}^{-1}|_{\cL(\cU[0,T])}^{-1}$.
\endpf

%%%%%%%%%%%%%%%%%%%%%%%%%%%%%%%%%%%%%%%%%%%%%%%%%%%%%%%%%%%%

\subsection{Uniqueness and existence of solutions to the deterministic Riccati equation}\label{sec-pr-main2}

%%%%%%%%%%%%%%%%%%%%%%%%%%%%%%%%%%%%%%%%%%%%%%%%%%%%%%%%%%%%

The purpose of this subsection is to prove Theorems \ref{12.6-th1}--\ref{th4.6}, addressed to the uniqueness and existence of solutions to the equation \eqref{Riccati}, respectively.

{\it Proof}.[Proof of Theorem \ref{12.6-th1}]
Let $P_1, P_2\in C_\cS([0,T];\dbS(H))$ be two
strongly regular mild solutions to
\eqref{Riccati}. Then, it follows from
\eqref{8.20-eq23} that for $j=1,2$, for any
$\eta\in H$ and $t\in [0,T]$,
\begin{equation}\label{12.7-eq1}
    \begin{array}{ll}\ds
        P_j(t)\eta\3n&\ds=S(T-t)^*Ge^{A(T-t)}\eta\\
        \ns&\ds\q +
        \int_t^T S(\tau-t)^*\big(C^* P_jC + M  -
        L_j^*K_j^{-1} L_j\big)S(\tau-t)\eta d\tau,
    \end{array}
\end{equation}
where
$$
L_j(\cd)=B(\cd)^*P_j(\cd)+D(\cd)^* P_j(\cd)C(\cd),\qq
K_j(\cd)=R(\cd)+D(\cd)^* P_j(\cd)D(\cd).
$$
From \eqref{12.7-eq1} and \eqref{strong-regular}, we have that
\begin{eqnarray}\label{12.7-eq2}
    &&\3n\3n
    \big|\big(P_1(t)-P_2(t)\big)\eta\big|_H\nonumber\\
    &&\3n\3n=
    \Big|\int_t^T S(\tau-t)^*\big[C^* \big(P_1-P_2\big)C -
    L_1^*K_1^{-1}\big(L_1-L_2\big)\\
    && \q -  \big(L_1^*-L_2^*\big)K_1^{-1}L_2 -  L_2^*\big(K_1^{-1}-K_2^{-1}\big)L_2 \big]S(\tau-t)\eta d\tau \Big|_H \nonumber\\
    &&\3n\3n \leq  \cC(|P_1|_{ C_\cS([0,T];\dbS(H))},|P_2|_{ C_\cS([0,T];\dbS(H))},c_0,R,B,C,D)|\eta|_H\int_t^T|P_1-P_2|_{\cL(H)}d\tau.\nonumber
\end{eqnarray}
By \eqref{12.7-eq2} and the arbitrariness of $\eta\in H$, we get that
\begin{equation*}\label{12.7-eq3}
    \ba{ll}\ds
    \big|P_1(t)-P_2(t)\big|_H \\
    \ns\ds\leq  \cC(|P_1|_{ C_\cS([0,T];\dbS(H))},|P_2|_{ C_\cS([0,T];\dbS(H))},c_0,R,B,C,D) \int_t^T|P_1-P_2|_{\cL(H)}d\tau.
    \ea
\end{equation*}
This, together with Gronwall's inequality, implies that $P_1(t)=P_2(t)$, $\forall\;t\in [0,T]$.
\endpf

\ms

{\it Proof}.[Proof of Theorem \ref{th4.6}]
``{\bf The ``if" part}". We divide the proof
into three steps.

\ss

{\bf Step 1}. In this step, we introduce a
sequence of operator-valued functions
$\{P_j\}_{j=1}^N$.

\ss

Let $P_0$ be the solution to
\begin{equation}\label{8.20-eq30}
    \left\{\2n
    \begin{array}{ll}
        \ds\frac{dP_0}{dt}+P_0 A + A^* P_0+C^* P_0C+M=0 &\mbox{ in }[0,T),\\
        \ns\ds P_0(T)=G.
    \end{array}
    \right.
\end{equation}
Applying Proposition \ref{prop4.5} to
\eqref{8.20-eq30} with $\Th=0$, we obtain that
\begin{equation}\label{8.20-eq30.1}
    R(t)+D(t)^* P_0(t)D(t)\geq c_0 I,\qq P_0(t)\geq\a
    I,\q\;\ae~t\in[0,T].
\end{equation}
Inductively, for $j = 0,1,2, \cdots$, we set
\begin{equation}\label{Iteration-i}
    \begin{array}{ll}\ds
        K_j\deq R+D^* P_jD,\qq L_j\deq B^* P_j+D^* P_jC,
        \\
        \ns\ds \Th_j\deq -K_j^{-1}L_j, \qq\cA_j\deq
        B\Th_j,\qq C_j\deq C+D\Th_j,
    \end{array}
\end{equation}
and let $P_{j+1}$ be the solution to
\begin{equation}\label{}
    \left\{\2n
    \begin{array}{ll}
        \ds\frac{dP_{j+1}}{dt}+P_{j+1}(A+\cA_j)+(A+\cA_j)^*
        P_{j+1}\\
        \ns\ds \q + C_j^* P_{j+1}C_j
        +\Th_j^* R\Th_j+M=0 &\mbox{ in }[0,T),\\
        \ns\ds P_{j+1}(T)=G.
    \end{array}
    \right.
\end{equation}

\ss

{\bf Step 2}. In this step, we show the uniform
boundedness of the sequence
$\{P_j\}_{j=1}^\infty$.

\ss

From \eqref{8.20-eq30.1}, we have that
\begin{equation}\label{11.27-eq1}
    K_0(t)\geq c_0 I,\q P_0(t)\geq\a
    I,\qq\ae~t\in[0,T].
\end{equation}
This implies that
$$
\Th_0=-K_0^{-1}L_0\in
L^{2,\cS}(0,T;\cL(H;U)).
$$
It follows from Proposition
\ref{prop4.5} (with $P_\Th$ and $\Th$ in
\eqref{eq-Lya} replaced by $P_{1}$ and $\Th_0$,
respectively) that
$$
K_{1}(t)\geq c_0 I,\q P_{1}(t)\geq\a I,\qq
\ae~t\in[0,T].
$$
Inductively, we have that
\begin{equation}\label{R+Pi-lowerbound}
    K_{j+1}(t)\geq c_0 I,\q P_{j+1}(t)\geq\a I,\qq
    \ae~t\in[0,T],\qq j=0,1,2,\cdots
\end{equation}

Let
$$\D_j\deq P_j-P_{j+1},\qq \Upsilon_j\deq\Th_{j-1}-\Th_j,\qq j\geq1.$$
Then for $j\geq1$ and $\zeta\in H$, we have
\begin{eqnarray}\label{Di-equa1}
    &&
    -\D_j(t)\zeta=P_{j+1}(t)\zeta-P_j(t)\zeta \nonumber\\
    &&= \int_t^T
    S(r-t)^*\big[\big(P_j -P_{j+1} \big) \cA_j +
    \cA_j^* \big(P_j -P_{j+1} \big)+ C^*_j
    \big(P_j -P_{j+1} )C_j \nonumber\\
    &&\qq\qq\qq\q +P_j(\cA_{j-1}-\cA_j)+
    (\cA_{j-1}-\cA_j)^*P_j+C_{j-1}^*
    P_jC_{j-1}\\
    &&\qq\qq\qq\q-C_j^* P_jC_j+\Th_{j-1}^*
    R\Th_{j-1}-\Th_j^* R\Th_j\big]S(r-t) \zeta dr \nonumber\\
    &&= \int_t^T S(r-t)^* \big[\D_j \cA_j +
    \cA_j^* \D_j + C^*_j
    \D_j C_j+P_j(\cA_{j-1}-\cA_j)+ (\cA_{j-1}-\cA_j)^*P_j \nonumber\\
    &&\qq\qq\qq\q +
    C_{j-1}^* P_jC_{j-1}-C_j^* P_jC_j+\Th_{j-1}^*
    R\Th_{j-1}  -\Th_j^*
    R\Th_j\big]S(r-t) \zeta dr.\nonumber
\end{eqnarray}
From \eqref{Iteration-i}, we have that
$$
\cA_{j-1}-\cA_j= B\Th_{j-1}  - B\Th_j =
B(\Th_{j-1} - \Th_j) =B\Upsilon_j,\qq\qq
$$
$$
C_{j-1}-C_j = C+D\Th_{j-1} - C -
D\Th_j=D(\Th_{j-1} -  \Th_j) =D\Upsilon_j,
$$
and
$$
\begin{array}{ll}\ds
    C_{j-1}^* P_jC_{j-1}-C_j^* P_jC_j \\
    \ns\ds= (C+D\Th_{j-1})^*P_j(C+D\Th_{j-1}) -
    (C+D\Th_{j})^*P_j(C+D\Th_{j})\\
    \ns\ds = \Th_{j-1}^*D^* P_jD\Th_{j-1} -
    \Th_{j}^*D^*P_jD\Th_{j} + \Th_{j-1}^*D^*\!P_jC +
    C^* P_jD\Th_{j-1}\\
    \ns\ds\q  - \Th_{j}^* D^* P_jC -
    C^*P_jD\Th_{j}\\
    \ns\ds =
    (\Th_{j-1}-\Th_{j})^*D^*P_jD(\Th_{j-1}-\Th_{j})+
    (C+D\Th_{j})^* P_jD(\Th_{j-1}-\Th_{j})\\
    \ns\ds\q +(\Th_{j-1}-\Th_{j})^* D^*
    P_j(C+D\Th_{j})
    \\ \ns\ds =\Upsilon_j^* D^* P_jD\Upsilon_j+C_j^*
    P_jD\Upsilon_j+\Upsilon_j^* D^* P_jC_j.
\end{array}
$$
Similarly,
\begin{equation*}\label{Di-equa2}
    \left\{\2n
    \begin{array}{ll}
        \ds  \Th_{j-1}^* R\Th_{j-1}-\Th_j^*
        R\Th_j=\Upsilon_j^* R\Upsilon_j+\Upsilon_j^*
        R\Th_j+\Th_j^* R\Upsilon_j,\\
        \ns\ds B^* P_j+D^* P_jC_j+R\Th_j=B^* P_j+D^*
        P_jC+(R+D^* P_jD)\Th_j=0.
    \end{array}
    \right.
\end{equation*}
These, together with \eqref{Di-equa1}, yields
that for any $\zeta\in H$,
\begin{eqnarray}\label{Di-equa3}
    &&  \D_j(t)\zeta-\int_t^T S(r-t)^*\big(\Delta_j\cA_j+\cA_j^*\D_j+C_j^*\D_jC_j\big)S(r-t)\zeta dr \nonumber\\
    &&=\int_t^T S(r-t)^*\big(P_j B\Upsilon_j+\Upsilon_j^* B^* P_j+\Upsilon_j^* D^* P_jD\Upsilon_j+C_j^* P_jD\Upsilon_j\\
    &&\qq\qq\qq\q+\Upsilon_j^* D^* P_jC_j+\Upsilon_j^* R\Upsilon_j+\Upsilon_j^* R\Th_j+\Th_j^* R\Upsilon_j\big)S(r-t)\zeta dr \nonumber\\
    &&= \int_t^T
    S(r-t)^*\big[\Upsilon_j^*K_j\Upsilon_j + (P_jB
    +C_j^* P_jD + \Th_j^* R)\Upsilon_j\nonumber\\
    &&\qq\qq\qq\q + \Upsilon_j^*(B^* P_j\!+\!D^*
    P_jC_j\!+\!R\Th_j)\big]S(r-t)\zeta dr\nonumber
    \\
    &&=\int_t^T S(r-t)^* \Upsilon_j^*K_j\Upsilon_j
    S(r-t)\zeta dr.\nonumber
\end{eqnarray}
By \eqref{Di-equa3},  $\D_j(\cd)$
solves \eqref{lm2.4-eq1} with $\wt
G=0$, $\wt A = \cA_j$, $\wt C = \cC_j$ and $\wt
M = \Upsilon_j^*K_j\Upsilon_j\geq 0$. Using
Lemma \ref{lm2.4}, we have $\D_j(\cd)\geq
0$, namely, $P_{j-1}(\cd)-P_j(\cd)\geq 0$ for
$j\geq1$. By (\ref{R+Pi-lowerbound}), for each $j\in\dbN$ and $ t\in [0,T]$,
$P_1(t)\geq P_j(t)\geq P_{j+1}(t)\geq\a I.$
Hence, the sequence $\{P_j\}_{j=1}^\infty$
is uniformly bounded. Consequently, there exists
a constant $\cC>0$ such that (noting
(\ref{R+Pi-lowerbound})) for all $j\geq0$ and
a.e. $t\in [0,T]$,
\begin{equation}\label{Di-equa4}
    \left\{\2n
    \begin{array}{ll}
        \ns\ds|P_j(t)|_{\cL(H)}\leq \cC,\q
        |K_j(t)|_{\cL(U)}\leq \cC,
        \\
        \ns\ds |\Th_j(t)|_{\cL(H;U)}\leq \cC\big(|B(t)|_{\cL(U;H)}+|C(t)|_{\cL(H)}\big),\\
        \ns\ds|\cA_j(t)|_{\cL(H)}\leq
        \cC|B(t)|_{\cL(U;H)}\big(|B(t)|_{\cL(U;H)}+|C(t)|_{\cL(H)}\big),\\
        \ns\ds |C_j(t)|_{\cL(H)}\leq
        \cC\big(|B(t)|_{\cL(U;H)}+|C(t)|_{\cL(H)}\big).
    \end{array}
    \right.
\end{equation}

\ss

{\bf Step 3}. In this step, we prove the
convergence of the sequence
$\{P_j\}_{j=1}^\infty$.

\ss

Noting
$
\L_j=\Th_{j-1}-\Th_j
=K_j^{-1}D^*\D_{j-1}DK_{j-1}^{-1}L_j-K_{j-1}^{-1}\big(B^*\D_{j-1}+D^*\D_{j-1}C\big)
$
and (\ref{Di-equa4}), one has
\begin{equation}\label{3.22}
    \begin{array}{ll}
        \ds|\Upsilon_j(t)^*
        K_j(t)\Upsilon_j(t)|_{\cL(H)}\\
        \ns\ds\leq\big(|\Th_j(t)|_{\cL(U)}+|\Th_{j-1}(t)|_{\cL(U)}\big)
        |K_j(t)|_{\cL(U)}
        |\Th_{j-1}(t)-\Th_j(t)|_{\cL(U)}\\
        \ns\ds \leq
        \cC\big(|B(t)|_{\cL(U;H)}+|C(t)|_{\cL(H)}\big)^2|\D_{j-1}(t)|_{\cL(H)}.
    \end{array}
\end{equation}
By (\ref{Di-equa3}), it follows that for any $\zeta\in H$,
$$
\D_j(t)\zeta=-\int^T_t S(r-t)^* \big(\D_j\cA_j+
\cA_j^*\D_j+C_j^*\D_iC_j+\Upsilon_j^*
K_j\Upsilon_j\big)S(r-t)\zeta dr.
$$
Making use of (\ref{3.22}) and noting
(\ref{Di-equa4}), by the arbitrariness of $\zeta\in H$, we get
$$
|\D_j(t)|_{\cL(H)}\leq
\int^T_t\f(r)\big(|\D_j(r)|_{\cL(H)}
+|\D_{j-1}(r)|_{\cL(H)}\big)dr,\qq\forall\;
t\in[0,T],$$
where $\f(\cd)$ is a nonnegative integrable
function independent of $\D_j(\cd)$ ($j\in\dbN$). Hence,
$$
|\D_j(t)|_{\cL(H)} \leq
a\int^T_t\f(r)|\D_{j-1}(r)|_{\cL(H)}dr,
$$
where $a=e^{\int_0^T\f(r)dr}$. Set
$b\deq\max_{0\leq r\leq T}|\D_0(r)|_{\cL(H)}$.
By induction, we deduce that
$$
|\D_j(t)|_{\cL(H)}\leq
b\frac{a^j}{j!}\(\int_t^T\f(r)dr\)^j,\qq\forall\;
t\in[0,T],
$$
which implies the uniform convergence of
$\{P_j\}_{j=1}^\infty$. Denote by $P$ the limit
of $\{P_j\}_{j=1}^\infty$, then (noting
(\ref{R+Pi-lowerbound}))
$\ds
K(t)=\lim_{j\to\infty}K_j(t) \geq c_0 I$
for a.e. $t\in[0,T]$,
and
$$\left\{\2n\begin{array}{ll}
    \ds \lim_{j\to\infty}\Th_j=-K^{-1}L
    \equiv\Th\q\hb{in } L^{2,\cS}(0,T;\cL(H;U)),\\
    \ns\ds \lim_{j\to\infty}\cA_j= B\Th\q\hb{in }
    L^{1,\cS}(0,T;\cL(H)),\\
    \ns\ds \lim_{j\to\infty}C_j= C+D\Th\q\hb{in
    }L^{2,\cS}(0,T;\cL(H)).\end{array}\right.$$
Therefore, $P(\cd)$ solves the following
equation (in the sense of strongly regular mild solution):
$$\left\{\2n
\begin{array}{ll}
    \ds\dot P+P(A+B\Th)+(A+B\Th)^* P\\
    \ns\ds \q+(C+D\Th)^* P(C+D\Th)+\,\Th^* R\Th+M=0 &\mbox{ in }[0,T),\\
    \ns\ds P(T)=G,
\end{array}
\right.$$
which is equivalent to (\ref{Riccati}).

\ms

``{\bf The ``only if" part}". Let $P(\cd)$ be the strongly
regular solution to (\ref{Riccati}). Then \eqref{strong-regular} holds for some $c_0>0$.
Put
$
\Th\deq-K^{-1}L(\in L^{2,\cS}(0,T;\cL(H;U))).
$
For any $u(\cd)\in\cU[0,T]$, let
$X(\cd)=X(\cd;0,u)$ be the solution to
\eqref{LQsystem1} with $\eta=0$. Applying
It\^o's formula to $t\mapsto\langle
P(t)X(t),X(t)\rangle_H$, we have
$$
\begin{array}{ll}
    \ds \cJ(0;u(\cd))=\dbE\[\int_0^T\big(\lan MX,X\rangle_H + \lan Ru,u\rangle_U\big)dt+\lan GX(T),X(T)\rangle_H\]\\
    \ns\ds=\dbE\int_0^T\Big\{\lan-\big(P A + A^* P
    +C^* PC+M-L^* K^{-1} L\big)X,X\rangle_H\\
    \ns\ds\q+\lan P\big(AX+Bu\big),X\rangle_H+\lan PX,AX+Bu\rangle_H \\
    \ns\ds\q+\lan P\big(CX+Du\big),CX+Du\rangle_H+\lan MX,X\rangle_H+\lan Ru,u\rangle_U\Big\}dt  \\
    \ns\ds=\dbE\int_0^T\big(\lan\Th^*K\Th
    X,X\rangle_H
    -2\lan K\Th X,u\rangle_U +\lan Ku,u\rangle_U\big)dt\\
    \ns\ds=\dbE\int_0^T\lan K\big(u-\Th X\big),u-\Th
    X\rangle_U dt.
\end{array}$$
Hence, by (\ref{strong-regular}) and  Lemma
\ref{lm2.5}, it holds that, for some $c_1>0$ and
all $u(\cd)\in\cU[0,T]$,
$$
\begin{array}{ll}\ds
    \cJ(0;u(\cd)) \ds =\dbE\int_0^T\lan K\big(u-\Th
    x\big),u-\Th x\rangle_U ds \geq c_1
    \dbE\int_0^T|u(s)|_U^2ds,
\end{array}
$$
which completes the proof of Theorem
\ref{th4.6}. \endpf

%%%%%%%%%%%%%%%%%%%%%%%%%%%%%%%%%%%%%%%%%%%%%%%%%%%%%%%%%%

\section{Pontryagin-type maximum
principle for controlled stochastic evolution equations}\label{ch-Pon}

%%%%%%%%%%%%%%%%%%%%%%%%%%%%%%%%%%%%%%%%%%%%%%%%%%%%%%%%%%%

In this section, we shall present the first order
necessary optimality condition, or more precisely, the
Pontryagin-type maximum principle, for optimal
control problems of nonlinear stochastic
evolution equations in infinite dimensions, in
which both the drift and diffusion terms may contain
the control variables, and the control domain is
allowed to be nonconvex. For the second order
necessary optimality conditions, we refer the
readers to \cite{FL, LZZ}.  The results in this
section are taken from \cite{LZ1, LZ2} (See also \cite{LZ-2019, LZ3.1}). Some related
results can be found in
\cite{Du, Fuhrman1}.

%%%%%%%%%%%%%%%%%%%%%%%%%%%%%%%%%%%%%%%%%%%%%%%%%%
\subsection{Formulation of the problem}

Throughout this section,  both $H$ and $A$ (generating a $C_0$-semigroup $S(\cd)$ on $H$) are the same as before,
$(\Omega,\cF,\mathbf{F},\dbP)$ (with
$\mathbf{F}\deq \{\cF_t\}_{t\in [0,T]}$ for some $T>0$) is a
fixed filtered probability space (satisfying
the usual condition), on which a $1$-dimensional standard Brownian motion $W(\cdot)$ is defined. Denote by $\dbF$ the
progressive $\si$-field w.r.t. $\mathbf{F}$.
Let $U$ be a separable
metric space with a metric $\mathbf{d}(\cd,\cd)$. Put
\bel{200619e1}
\cU[0,T] \deq \big\{u(\cdot):\,[0,T]\times\Omega\to U\;\big|\; u(\cd) \mbox{ is $\mathbf{F}$-adapted} \big\}
\ee

For any given functions $a(\cd,\cd,\cd,\cd):[0,T]\times\Omega\times H\times U\to H$ and $b(\cd,\cd,\cd,\cd):[0,T]\times\Omega \times H\times U\to H$, let us consider the following controlled SEE:
\begin{equation}\label{ch-10-fsystem1}
\left\{
\begin{array}{lll}\ds
    dX = \big(AX +a(t,X,u)\big)dt + b(t,X,u)dW(t) &\mbox{ in }(0,T],\\
    \ns\ds X(0)=X_0,
\end{array}
\right.
\end{equation}
where $u(\cd)\in \cU[0,T]$ is the {\it control variable}, $X(\cd)$ is the {\it state variable}, and the (given) initial state $X_0\in
H$.
Throughout this section, For $\psi=a(\cd,\cd,\cd,\cd),b(\cd,\cd,\cd,\cd)$, we assume the following
condition:

\ms

\no{\bf (AS1)} {\it i) For any $(x,u)\in H\times U$, the
function $\psi(\cd,\cd,x,u):[0,T]\times\Om\to H$ is $\dbF$-measurable; ii) For any $(t,x)\in [0,T]\times
H$, the function $\psi(t,x,\cd):U\to H$ is
continuous; and iii) For any $(x_1,x_2,u)\in  H\times
H\times U$ and a.e. $(t,\omega)\in (0, T)\times\Omega$,
\begin{equation}\label{ab0}
    \left\{
    \begin{array}{ll}\ds
        |\psi(t,x_1,u) - \psi(t,x_2,u)|_H  \leq
        \cC |x_1-x_2|_H,\\
        \ns\ds |\psi(t,0,u)|_H \leq \cC.
    \end{array}
    \right.
\end{equation}}

By Theorem \ref{ch-1-well-mild-eq1}, it is clear that
\begin{proposition}\label{12.9-prop1}
Let the assumption (AS1) hold. Then, for any $p_0\geq 2$,
$X_0\in H$ and
$u(\cd)\in \cU[0,T]$, the
equation \eqref{ch-10-fsystem1} admits a unique
mild solution $X\in C_\dbF([0,T];L^{p_0}(\Om;H))$. Furthermore,
$$
|X(\cd)|_{C_\dbF([0,T];L^{p_0}(\Omega;H))} \leq
\cC\big(1+|X_0|_{H}\big).
$$
\end{proposition}

Also, we need the following condition:

\ms

\no{\bf (AS2)} {\it Suppose that
$g(\cd,\cd,\cd,\cd):[0,T]\times \Omega\times H\times U\to \dbR$
and $h(\cd,\cd):\Omega\times H\to \dbR$ are two functions
satisfying: i) For any $(x,u)\in H\times U$, the
function $g(\cd,\cd,x,u):[0,T]\times\Om\to \dbR$ is $\dbF$-measurable and the
function $h(\cd,x):\Omega\to \dbR$ is $\cF_T$-measurable; ii) For a.e. $(t,\om)\in [0,T]\times\Om$ and any $x\in
H$, the function $g(t,\om,x,\cd):U\to \dbR$ is
continuous; and iii) For all $(t,x_1,x_2,u)\in
[0,T]\times H\times H\times U$,
\begin{equation}\label{ch-10-gh} \left\{
    \begin{array}{ll}\ds
        |g(t,x_1,u) - g(t,x_2,u)| +|h(x_1) - h(x_2)|
        \leq \cC|x_1-x_2|_H,\q \as\!\!,\\
        \ns\ds |g(t,0,u)| +|h(0)| \leq \cC,\q \as
    \end{array}
    \right.
\end{equation}}

Define a cost functional $\cJ(\cdot)$ (for the
controlled system \eqref{ch-10-fsystem1}) as
follows:
\begin{equation}\label{jk1xu}
\cJ(u(\cdot))\deq
\dbE\Big(\int_0^T g(t,X(t),u(t))dt +
h(X(T))\Big),\q\forall\; u(\cdot)\in
\cU[0,T],
\end{equation}
where $X(\cd)$ is the state of
\eqref{ch-10-fsystem1}  corresponding to $u$.

Consider the following optimal control problem:

\no {\bf Problem (OP)} {\it Find a $\bar
u(\cdot)\in \cU[0,T]$ such that
\begin{equation}\label{jk2}
    \ds\cJ (\bar u(\cdot)) = \inf_{u(\cdot)\in
        \cU[0,T]} \cJ (u(\cdot)).
\end{equation}
Any $\bar u(\cdot)$ satisfying (\ref{jk2}) is
called an {\it optimal control}. The
corresponding state $\cl X(\cdot)$ is
called an {\it optimal state}, and
$(\cl X(\cdot),\bar u(\cdot))$ is called an
{\it optimal pair}.}

\begin{remark}
We do not put state/end point constraints for
our optimal control problem. Readers who are
interested in this are referred to \cite{FL} for
some recent results.
\end{remark}

Problem (OP) is now well-understood for the
case of finite dimensions (i.e., $\dim
H<\infty$) and natural filtration. In this case,
a Pontryagin-type maximum principle was obtained
in \cite{Peng} for general stochastic control
systems with control-dependent diffusion
coefficients and possibly non-convex control
regions, and it was found there that the corresponding
result differs significantly from its
deterministic counterpart. The main purpose in this section is to see what happens when $\dim
H=\infty$.

%

%%%%%%%%%%%%%%%%%%%%%%%%%%%%%%%%%%%%%%%%%%%%%%%%%%%%%

\subsection{Pontryagin-type maximum principle for
convex control domain}\label{subsec20210120}

%%%%%%%%%%%%%%%%%%%%%%%%%%%%%%%%%%%%%%%%%%%%%%%%%%%%%%

In this subsection,  we shall present a necessary condition
for optimal controls of  Problem (OP) under the
following conditions:

\ms

\no{\bf(AS3)} {\it The control region $U$ is a convex
subset of a separable Hilbert space  $H_1$ and
the metric of $U$ is induced by the norm of
$H_1$, i.e., $\mathbf{d}(u_1,u_2)=|u_1-u_2|_{H_1}$.}

\ms

\no{\bf (AS4)} {\it For a.e. $(t,\omega)\in (0, T)\times\Omega$, the functions $a(t,\cd,\cd):H\times U\to H$ and
$b(t,\cd,\cd):H\times U\to H$,  $g(t,\cd,\cd):H\times U\to \dbR$ and
$h(\cd): H\to \dbR$ are $C^1$. Moreover,
for any $(x,u)\in H\times U$ and a.e. $(t,\omega)\in (0, T)\times\Omega$,
$$
\left\{
\begin{array}{ll}\ds
    |a_x(t,x,u)|_{\cL(H)}+|b_x(t,x,u)|_{\cL(H)} + |g_x(t,x,u)|_H+|h_x(x)|_{H}\leq \cC,\\
    \ns\ds|a_u(t,x,u)|_{\cL(H_1;H)}+
    |b_u(t,x,u)|_{\cL(H_1;H)} +
    |g_u(t,x,u)|_{H_1} \leq \cC.
\end{array}
\right.
$$
}

Further, we need to introduce the following $H$-valued BSEE:
\begin{eqnarray}\label{ch-10-bsystem1}
\left\{
\begin{array}{lll}
    \ds dY(t) = -  A^* Y(t) dt + f(t,Y(t),Z(t))dt + Z(t) dW(t) &\mbox{ in }[0,T),\\
    \ns\ds y(T) = Y_T.
\end{array}
\right.
\end{eqnarray}
Here $Y_T \in L_{\cF_T}^{2}(\Omega;H)$) and
$f(\cd,\cd,\cd,\cd):[0,T]\times \Om\times H\times H \to H$
satisfies
\begin{equation}\label{Lm1.1.1}
\left\{
\begin{array}{ll}\ds
    f(\cd,\cd,0,0)\in L^1_{\dbF}(0,T;
    L^2(\Omega;H)),\\\ns\ds
    |f(t,\om,y_1,z_1)-f(t,\om,y_2,z_2)|_H\leq
    \cC\big(|y_1-y_2|_H+|z_1-z_2|_{H} \big),\\
    \ns\ds\hspace{3cm} \ae (t,\omega)\in
    [0,T]\times\Omega,\;\;
    \forall\;y_1,y_2, z_1,z_2 \in H.
\end{array}
\right.
\end{equation}
By Theorem \ref{vw-th1}, the equation
\eqref{ch-10-bsystem1} admits a unique
transposition solution $(Y(\cdot), $ $Z(\cdot)) \in
D_{\dbF}([0,T];L^{2}(\Omega; H)) \times
L^{2}_{\dbF}(0,T;H)$. Furthermore,
\begin{equation}\label{s4eq1z}
\begin{array}{ll}\ds
    |(Y(\cdot), Z(\cdot))|_{
        D_{\dbF}([0,T];L^2(\Omega;H)) \times
        L^2_{\dbF}(0,T;H)}\\
    \ns\ds \leq \cC\big(|Y_T|_{
        L^2_{\cF_T}(\Omega;H)}+ |f(\cdot,0,0)|_{
        L^1_{\dbF}(0,T;L^2(\Omega;H))}\big).
\end{array}
\end{equation}

We have the following necessary condition for optimal controls of Problem
(OP) with a convex control domain.

\begin{theorem}\label{th max}
Let
the assumptions (AS1)--(AS4)
hold, and let $(\cl X(\cd),\bar u(\cd))$ be an
optimal pair of Problem (OP). Let
$(Y(\cdot),$ $Z(\cdot))$ be the transposition
solution to the equation \eqref{ch-10-bsystem1}
with $Y_T$ and $f(\cd,\cd,\cd)$ given by
\begin{equation}\label{zv1}
    \left\{
    \begin{array}{ll}
        \ds Y_T =
        -h_x\big(\cl X(T)\big),\\
        \ns \ds f(t,y_1,y_2)=-a_x(t,\cl
        X(t),\bar u(t))^*y_1 - b_x\big(t,\cl
        X(t),\bar u(t)\big)^*y_2 +
        g_x\big(t,\cl X(t),\bar u(t)\big).
    \end{array}
    \right.
\end{equation}
Then,
\begin{equation}\label{1.25-eq1}
    \begin{array}{ll}\ds
        \big\langle a_u(t,\cl X(t),\bar
        u(t))^* Y(t) + b_u(t,\cl X(t),\bar
        u(t))^*Z(t)- g_u(t,\bar u(t),\cl
        X(t)), \rho-\bar u(t) \big\rangle_{H_1}
        \leq 0,\\
        \ns\ds \qq\qq\qq\qq\qq\qq\qq\qq\qq \forall\,\rho \in U,\q \ae (t,\omega)\in
        [0,T]\times \Omega.
    \end{array}
\end{equation}
\end{theorem}

%%%%%%%%%%%%%%%%%%%%%%%%%%%%%%%%%%%%%%%%%%%%%%%%%%%%%%%%%%%%%%%%

{\it Proof}. We use the convex perturbation
technique and divide the proof into three steps.

\ms

{\bf Step 1}. We fix any
control $u(\cdot)\in \cU[0,T]$. Since $U$ is
convex, for any $\e \in [0,1]$,
$$
u^\e(\cdot) \deq \bar u(\cdot) + \e \big(u(\cdot) -
\bar u(\cdot)\big) = (1-\e)\bar u(\cdot) + \e
u(\cdot) \in \cU[0,T].
$$
Denote by $X^\e(\cdot)$ the state
of \eqref{ch-10-fsystem1}
corresponding to the control
$u^\e(\cdot)$. It follows from Proposition \ref{12.9-prop1} that
\begin{equation}\label{th max eq0.0xz}
    |X^\e(\cd)|_{C_\dbF([0,T];L^2(\Omega;H))}\leq C
    \big(1+|X_0|_{H}\big),\q
    \forall\; \e \in [0,1].
\end{equation}
Write
$$
\ds X_1^\e(\cd)
= \frac{1}{\e}\big(X^\e(\cd)-\cl
X(\cd)\big),\qq \d u(\cd) = u(\cd) -
\bar u(\cd).
$$
It is easy to see that $X_1^\e(\cdot)$
solves the following SEE:
\begin{equation}\label{fsystem3x}
    \left\{
    \begin{array}{lll}\ds
        dX_1^\e = \big(AX_1^\e +  a_1^\e X^\e_1
        + a_2^\e\d u  \big)dt + \big( b_1^\e
        X^\e_1 + b_2^\e\d u \big)dW(t) &\mbox{
            in
        }(0,T],\\
        \ns\ds X_1^\e(0)=0,
    \end{array}
    \right.
\end{equation}
where for $\psi=a,b$,
\begin{equation}\label{tatb}
    \begin{cases}
        \ds   \psi_1^\e (t)  = \int_0^1
        \psi_x(t,\cl X(t) +
        \si\e X_1^\e(t), u^\e(t))d\si, \\
        \ns\ds \psi_2^\e (t) = \int_0^1
        \psi_u(t,\cl X(t), \bar u(t)+\si\e\d
        u(t))d\si.
    \end{cases}
\end{equation}
Consider the following SEE:
\begin{equation}\label{fsystem3.1}
    \left\{
    \begin{array}{lll}\ds
        dX_2 = \big(AX_2 + a_1 X_2 +  a_2 \d u \big)dt +
        \big( b_1 X_2 + b_2 \d u \big) dW(t) &\mbox{ in
        }(0,T],\\
        \ns\ds X_2(0)=0,
    \end{array}
    \right.
\end{equation}
where for $\psi=a,b$,
\begin{equation}\label{barab}
    \psi_1(t) = \psi_x(t,\cl X(t),\bar u(t)),\q
    \psi_2(t) = \psi_u(t,\cl X(t),\bar u(t)).
\end{equation}

\ms

{\bf Step 2}. In this step, we shall show that
\begin{equation}\label{th max eq0.2}
    \lim_{\e\to 0+} \big|X_1^\e -
    X_2\big|_{L_\dbF^\infty(0,T;L^2(\Omega;H))}=0.
\end{equation}

First, using Burkholder-Davis-Gundy's inequality
(i.e., Proposition \ref{BDGQ1}) and by the assumption
(AS4), we find that
\begin{equation}\label{th max eq0}
    \begin{array}{ll}\ds
        \mE|X_1^\e(t)|^2_H \\
        \ns\ds  = \mE\Big|
        \int_0^t S(t-s) a_1^\e(s) X_1^\e(s) ds
        + \int_0^t S(t-s)
        a_2^\e(s)\d u(s) ds  \\
        \ns \ds \qq +
        \int_0^t S(t-s) b_1^\e(s) X_1^\e(s) dW(s) +  \int_0^t S(t-s)b_2^\e(s)\d u(s) dW(s)\Big|_H^2\\
        \ns \ds  \leq \cC \mE\(\Big| \int_0^t\!
        S(t\!-\!s) a_1^\e(s) X_1^\e(s) ds \Big|_H^2
        +  \Big|\int_0^t\! S(t\!-\!s) b_1^\e(s) X_1^\e(s) dW(s) \Big|_H^2 \\
        \ns \ds \qq  + \Big| \int_0^t\! S(t\!-\!s)a_2^\e(s)\d u(s)ds\Big|_H^2 + \Big| \int_0^t\! S(t\!-\!s)b_2^\e(s)\d u(s)dW(s)\Big|_H^2\) \\
        \ns \ds   \leq \cC\( \int_0^t
        \mE|X_1^\e(s)|_H^2ds + \int_0^T\mE |\d
        u(s)|_{H_1}^2 dt\).
    \end{array}
\end{equation}
It follows from \eqref{th max eq0} and
Gronwall's inequality that
\begin{equation}\label{th max eq0.1}
    \begin{array}{ll}\ds
        \mE|X_1^\e(t)|^2_H \leq  \cC|\bar u -
        u|^2_{L^2_\dbF(0,T;H_1)},\q\forall\;t\in
        [0,T].
    \end{array}
\end{equation}
Similarly,
\begin{equation}\label{th max eq2}
    \begin{array}{ll}\ds
        \mE|X_2(t)|^2_H \leq  \cC |\bar u -
        u|^2_{L^2_\dbF(0,T;H_1)},\q\forall\;t\in [0,T].
    \end{array}
\end{equation}

Put $X_3^\e \deq X_1^\e - X_2$. Then,
$X_3^\e$ solves the following equation:
\begin{equation}\label{fsystem3x2}
    \left\{
    \begin{array}{lll}\ds
        dX_3^\e = \!\big[AX_3^\e \!+\!  a_1^\e(t)
        X^\e_3 \!+\! \big(a_1^\e(t)\!- \!a_1(t)\big)X_2
        +
        \big( a_2^\e(t)\!-\!a_2(t)\big)\d u \big]dt \\
        \ns\ds \hspace{1.2cm} + \big[
        b_1^\e(t) X^\e_3 + \big(b^\e_1(t) -
        b_1(t)\big) X_2 + \big(
        b_2^\e(t)-b_2(t)\big)\d u \big] dW(t)
        &\mbox{in
        }(0,T],\\
        \ns\ds X_3^\e(0)=0.
    \end{array}
    \right.
\end{equation}
It follows from \eqref{th max
    eq2}--\eqref{fsystem3x2} that
$$
\ba{ll}
\ds\mE|X_3^\e(t)|^2_{H} \\
\ns\ds = \mE\Big| \int_0^t S(t-s)
a_1^\e(s)X_3^\e(s)ds
+ \int_0^t S(t-s) b_1^\e(s)X_3^\e(s) dW(s) \\
\ns\ds \q  + \int_0^t S(t - s)\big[
a^\e_1(s) -a_1(s)\big] X_2(s) ds +
\int_0^t S(t-s)\big[ b_1^\e(s) -
b_1(s)\big] X_2 (s) dW(s) \nonumber
\\
\ns\ds \q + \int_0^t S(t - s)\big[ a^\e_2
-a_2 \big]\d u(s) ds
+ \int_0^t S(t - s)\big[ b_2^\e -  b_2 \big]\d u dW(s)\Big|_H^2\\
\ns\ds \leq \cC\[\mE\int_0^t |X_3^\e
|_H^2 ds  + |X_2
|^2_{L^\infty_\dbF(0,T;L^2(\Omega;H))}
\int_0^T \mE\big(|a^\e_1 -a_1 |_{\cL(H)}^2
+ |b^\e_1 - b_1 |_{\cL(H)}^2 \big) dt
\\
\ns\ds  \q +  |u - \bar
u|^2_{L^2_\dbF(0,T;L^2(\Omega;H_1))} \int_0^T
\mE\big(|a_2^\e - a_2 |_{\cL(H_1;H)}^2 + |b_2^\e
- b_2 |_{\cL(H_1;H)}^2 \big)dt
\]\\
\ns\ds \leq  \cC(1 + |u- \bar
u|^2_{L^2_\dbF(0,T;L^2(\Omega;H_1))})\Big[\mE
\int_0^t |X_3^\e(s)|_H^2 ds +\int_0^T
\mE\big(|a^\e_1 - a_1 |_{\cL(H)}^2   \\
\ns\ds \q +  |b^\e_1 - b_1 |_{\cL(H)}^2+
|a_2^\e - a_2 |_{\cL(H_1;H)}^2 + |b_2^\e - b_2
|_{\cL(H_1;H)}^2
\big)dt\Big].
\ea
$$
This, together with Gronwall's inequality,
implies that
\begin{equation}\label{th max eq3}
    \begin{array}{ll}\ds
        \mE|X_3^\e(t)|^2_{H} \3n&\ds \leq \cC
        e^{\cC|u-\bar
            u|_{L^2_\dbF(0,T;L^2(\Omega;H_1))}}\int_0^T
        \mE\big(|a^\e_1 -a_1 |_{\cL(H)}^2
        + |b^\e_1 -b_1 |_{\cL(H)}^2\\
        \ns&\ds \q + |a_2^\e - a_2 |_{\cL(H_1;H)}^2 +
        |b_2^\e - b_2 |_{\cL(H_1;H)}^2 \big)ds,\q
        \forall\; t\in [0,T].
    \end{array}
\end{equation}
Note that \eqref{th max eq0.1} implies
$X^\e(\cd)\to \cl X(\cd)$ (in $H$) in
probability, as $\e\to0$. Hence, by
\eqref{tatb}, \eqref{barab} and the
continuity of $a_x(t,\cd,\cd)$,
$b_x(t,\cd,\cd)$, $a_u(t,\cd,\cd)$ and
$b_u(t,\cd,\cd)$, we deduce that
$$
\ba{ll}\ds \lim_{\e\to 0}\int_0^T\mE\big(
|a^\e_1(s)-a_1(s)|_{\cL(H)}^2 +
|b^\e_1(s)-b_1(s)|_{\cL(H)}^2
\\\ns\ds\qq\qq+ |a_2^\e(s) - a_2(s)|_{\cL(H_1;H)}^2 +|b_2^\e(s) -
b_2(s)|_{\cL(H_1;H)}^2 \big)ds=0. \ea
$$
This, together with (\ref{th max eq3}), gives
\eqref{th max eq0.2}.

\ms

{\bf Step 3}. Since $(\cl X(\cdot),\bar
u(\cdot))$ is an optimal pair of Problem
(OP), from \eqref{th max eq0.2}, we find that
\begin{eqnarray}\label{var 1}
    \begin{array}{ll}\ds
        0\leq \lim_{\e\to 0}\frac{\cJ(u^\e(\cdot)) - \cJ(\bar u(\cdot))}{\e} \\
        \ns\ds\;\;\;=  \dbE\int_0^T
        \Big(\big\langle g_1(t),X_2(t)\big\rangle_H
        \!+\! \big\langle g_2(t),\d u(t)
        \big\rangle_{H_1}\Big) dt + \dbE\big\langle
        h_x(\cl X(T)),X_2(T)\big\rangle_H,
    \end{array}
\end{eqnarray}
where $ g_1(t) = g_x(t,\cl X(t),\bar
u(t))$ and $g_2(t) = g_u(t,\cl
X(t),\bar u(t))$.

By the definition of transposition solution
to \eqref{ch-10-bsystem1}, we have that
\begin{equation}\label{max eq1}
    \begin{array}{ll}\ds
        -\dbE \big\langle h_x(\cl
        X(T)),X_2(T)\big\rangle_H - \dbE\int_0^T
        \big\langle g_1(t),X_2(t)\big\rangle_H dt \\
        \ns\ds= \dbE \int_0^T \Big(\big\langle
        a_2(t)\d u(t), Y(t)\big\rangle_H +
        \big\langle b_2(t)\d u(t),
        Z(t)\big\rangle_H\Big)dt.
    \end{array}
\end{equation}
Combining \eqref{var 1} and \eqref{max eq1}, we
find that
\begin{eqnarray}\label{max ine2}
    \dbE\int_0^T\big\langle a_2(t)^* Y(t) +
    b_2(t)^*Z(t) - g_2(t), u(t)-\bar u(t)
    \big\rangle_{H_1}dt\leq 0
\end{eqnarray}
holds for any $u(\cdot)\in \cU[0,T]$ satisfying
$u(\cd)-\bar u(\cd)\in L^2_\dbF(0,T;H_1)$.
Hence, by Lemma \ref{lemma4} and \eqref{max
ine2}, we conclude \eqref{1.25-eq1}. This
completes the proof of Theorem \ref{th
max}.\endpf

%%%%%%%%%%%%%%%%%%%%%%%%%%%%%%%%%%%%%%%%%%%%%%%%%%%%%

\subsection{Pontryagin-type maximum principle for the
general control domain}

%%%%%%%%%%%%%%%%%%%%%%%%%%%%%%%%%%%%%%%%%%%%%%%%%%%%%%

In this subsection, we give a necessary
condition for optimal controls of Problem
(OP) for a general control domain.

As we shall see later, the main difficulty in
dealing with the non-convex control
domain $U$ is that one needs   to solve the
following $\cL(H)$-valued BSEE:
\begin{equation}\label{op-bsystem3}
\left\{
\begin{array}{ll}
    \ds dP  =  - (A^*  + J^* )P dt  -  P(A + J )dt
    -K^*PKdt \\
    \ns\ds \hspace{1cm} - (K^* Q +  Q K)dt
    +   Fdt  +  Q dW(t) \qq\q\mbox{ in } [0,T),\\
    \ns\ds P(T) = P_T.
\end{array}
\right.
\end{equation}
Here and henceforth, $F\in
L^1_\dbF(0,T;L^2(\Omega;\cL(H)))$, $P_T\in
L^2_{\cF_T}(\Omega;\cL(H))$ and $J,K\in
L^4_\dbF(0,T; L^\infty(\Omega; $ $\cL(H)))$.
As we have explained in Section
\ref{s-lq-c}, in the previous literature there
exists no such a stochastic integration theory
in general Banach spaces that can be employed to
treat the well-posedness of \eqref{op-bsystem3}.
In order to overcome this difficulty,
similarly to Section
\ref{s-lq-c}, we employ
the stochastic transposition method introduced in our previous works
\cite{LZ, LZ1}.

%%%%%%%%%%%%%%%%%%%%%%%%%%%%%%%%%%%%%%%%%%%%%%%%%%%%%

\subsubsection{Relaxed transposition solutions}

%%%%%%%%%%%%%%%%%%%%%%%%%%%%%%%%%%%%%%%%%%%%%%%%%%%%%%
To define solutions to \eqref{op-bsystem3}, for $t\in [0,T)$, $ \xi_i \in L^{4}_{\cF_t}(\Omega;H)$ and
$u_i,v_i\in L^2_\dbF(t,T;L^{4}(\Omega;H))$ ($i=1,2)$,
we introduce two SEEs:
\begin{equation}\label{op-fsystem1}
\left\{
\begin{array}{ll}
    \ds d\f_1 = (A+J)\f_1ds + u_1ds + K\f_1 dW(s) + v_1 dW(s) &\mbox{ in } (t,T],\\
    \ns\ds \f_1(t)=\xi_1
\end{array}
\right.
\end{equation}
and
\begin{equation}\label{op-fsystem2}
\left\{
\begin{array}{ll}
    \ds d\f_2 = (A+J)\f_2ds + u_2ds + K\f_2 dW(s) + v_2 dW(s) &\mbox{ in } (t,T],\\
    \ns\ds \f_2(t)=\xi_2.
\end{array}
\right.
\end{equation}

Also, we shall introduce the solution space for
\eqref{op-bsystem3}.   Write
$$
\begin{array}{ll}\ds
\cP[0,T] \deq \Big\{P(\cd,\cd)\;\Big|\; P(\cd,\cd)\in
\cL_{pd}\big(L_{\dbF}^{2}(0,T;$ $L^{4}(\Omega;H)); L^2(0,T;L_{\dbF}^{\frac{4}{3}}(\Omega;H))\big) \\
\ns\ds\qq\qq \;\;\mbox{ and
    for every }
t\in[0,T]\hb{ and }\xi\in
L^{4}_{\cF_t}(\Omega;H), \\
\ns\ds\qq\qq \;\;\; P(\cd,\cd)\xi\in
D_{\dbF}([t,T];L^{\frac{4}{3}}(\Omega;H)) \mbox{
    and }\\
\ns\ds\qq\qq\q
|P(\cd,\cd)\xi|_{D_{\dbF}([t,T];L^{\frac{4}{3}}(\Omega;H))}
\leq \cC|\xi|_{L^{4}_{\cF_t}(\Omega;H)} \Big\}
\end{array}
$$
and
$$
\begin{array}{ll}\ds \cQ[0,T]\deq\!
\Big\{\big(Q^{(\cd)},\!\widehat
Q^{(\cd)}\big)\;\Big|\;\mbox{For any}\, t\!\in\!
[0,\!T], Q^{(t)}\,\mbox{and}\,\widehat
Q^{(t)}\;\mbox{are bounded
    linear}\\
\ns\ds\hspace{1.57cm}\mbox{ operators from
}L^{4}_{\cF_t}(\Omega;H)\!\times\!
L^2_\dbF(t,T;L^{4}(\Omega;H))\!\times\!
L^2_\dbF(t,T;L^{4}(\Omega;H)) \mbox{ to }
\\
\ns\ds
\hspace{1.67cm}L^{2}_\dbF(t,T;L^{\frac{4}{3}}(\Omega;H))
\mbox{ and }Q^{(t)}(0,0,\cd)^*=\widehat
Q^{(t)}(0,0,\cd)\Big\}.
\end{array}
$$
The norms on $\cP[0,T]$ and $\cQ[0,T]$ are given respectively by
$$
|P(\cd,\cd)|_{\cP[0,T]}\deq |P(\cd,\cd)|_{\cL\big(L_{\dbF}^{2}(0,T;L^{4}(\Omega;H)); L^2(0,T;L_{\dbF}^{\frac{4}{3}}(\Omega;H))\big)}
$$
and
$$
\begin{array}{ll}\ds
|\big(Q^{(\cd)}, \widehat
Q^{(\cd)}\big)|_{\cQ[0,T]}\\
\ns\ds\deq \sup_{t\in[0,T]}|\big(Q^{(t)},\!\widehat
Q^{(t)}\big)|_{\cL(L^{4}_{\cF_t}(\Omega;H)\!\times\!
    L^2_\dbF(t,T;L^{4}(\Omega;H))\!\times\!
    L^2_\dbF(t,T;L^{4}(\Omega;H)); L^{2}_\dbF(t,T;L^{\frac{4}{3}}(\Omega;H)))^2}
\end{array}
$$

The following notion will play a fundamental role in the sequel.
\begin{definition}\label{op-definition2x}
We call
$ \big(P(\cd),Q^{(\cd)},\widehat
Q^{(\cd)}\big)\in \cP[0,T]\times \cQ[0,T] $
a {\it relaxed transposition solution} to the equation
\eqref{op-bsystem3} if for any $t\in [0,T]$,
$\xi_1,\xi_2\in L^{4}_{\cF_t}(\Omega;H)$ and
$u_1(\cd),u_2(\cd), v_1(\cd), v_2(\cd)\in
L^2_{\dbF}(t,T;L^{4}(\Omega;H))$, it holds that
$$
\begin{array}{ll}
    \ds \mE\big\langle P_T \f_1(T), \f_2(T)
    \big\rangle_{H} - \mE \int_t^T \big\langle
    F(s) \f_1(s), \f_2(s) \big\rangle_{H}ds\\
    \ns\ds =\mE\big\langle P(t) \xi_1,\xi_2
    \big\rangle_{H} + \mE \int_t^T \big\langle
    P(s)u_1(s), \f_2(s)\big\rangle_{H}ds + \mE
    \int_t^T \big\langle P(s)\f_1(s),
    u_2(s)\big\rangle_{H}ds \\
    \ns\ds \q \! + \mE \int_t^T \big\langle P(s)K
    (s)\f_1 (s), v_2 (s)\big\rangle_{H}ds + \mE
    \int_t^T \big\langle  P(s)v_1 (s), K (s)\f_2 (s) + v_2(s)\big\rangle_{H}ds\\
    \ns\ds \q \! + \mE \int_t^T  \big\langle v_1(s),
    \widehat
    Q^{(t)}(\xi_2,u_2,v_2)(s)\big\rangle_{H}ds + \mE
    \int_t^T \big\langle Q^{(t)}(\xi_1,u_1,v_1)(s),
    v_2(s) \big\rangle_{H}ds.
\end{array}
$$
Here, $\f_1(\cd)$ and $\f_2(\cd)$ solve
\eqref{op-fsystem1} and \eqref{op-fsystem2},
respectively.
\end{definition}

We have the following well-posedness result for
the equation (\ref{op-bsystem3}) (See \cite{LZ1}
for its proof).

\begin{theorem}\label{OP-th2}
Suppose that  $L^p_{\cF_T}(\Omega)$ ($1\leq p <
\infty$) is separable. Then the equation
\eqref{op-bsystem3} has a unique relaxed
transposition solution
$\big(P(\cd),Q^{(\cd)},\widehat
Q^{(\cd)}\big)\in \cP[0,T]\times \cQ[0,T]$.
Furthermore,
$$
\begin{array}{ll}\ds
    \q |P|_{\cP[0,T]} + \big|\big(Q^{(\cd)},\widehat
    Q^{(\cd)}\big)\big|_{\cQ[0,T]}\leq \cC\big(
    |F|_{L^1_\dbF(0,T; L^2(\Omega;\cL(H)))} +
    |P_T|_{L^2_{\cF_T}(\Omega; \cL(H))}\big).
\end{array}
$$
\end{theorem}
\begin{remark}
It is well known that $L^p_{\cF_T}(\Om)$ is
separable if the probability space $(\Om,\cF_T,\dbP)$ is separable,
i.e.,  there exists a countable family
$\cE\subset\cF_T$ such that, for any $\e>0$ and
$B\in \cF_T$ one can find $B_1\in\cE$ with
$\dbP\big((B\setminus B_1)\cup (B_1\setminus
B)\big)<\e$ (e.g., \cite[Section 13.4]{BBT}).
Except some artificial
examples, almost all frequently used probability
spaces are separable (e.g.,
\cite{Ito}).
\end{remark}

Next, we give a regularity result for
relaxed transposition solutions to \eqref{op-bsystem3}. To this
end, we recall two preliminary results
(See \cite{LZ2} for their proofs). The first one is as follows:
\begin{lemma}\label{lemma5}
For each $t\in[0,T]$, if $u_2=v_2=0$ in the
equation \eqref{op-fsystem2}, then there exists
an operator $U(\cd,t)\in
\cL\big(L^{4}_{\cF_t}(\Omega;H);
C_\dbF([t,T];L^{4}(\Omega;H))\big)$ such that
the solution to \eqref{op-fsystem2} can be
represented as $\f_2(\cd) = U(\cd,t)\xi_2$.
\end{lemma}

Let $\{\D_n\}_{n=1}^\infty$ be a sequence of
partitions of $[0,T]$, that is,
$$
\D_n\deq  \Big\{t_i^n\;\Big|\;i=0,1,\cdots,n, \hb{
and }0=t_0^n < t_1^n < \cds < t_{n}^n =T\Big\},
$$
such that $\D_n\subset \D_{n+1}$ and
$\max_{0\leq i\leq n-1} (t_{i+1}^n -
t_{i}^n)\to 0 \;\mbox{ as }\;n\to\infty.$ We
introduce the following subspaces of
$L^2_\dbF(0,T;L^{4}(\Om;H))$:
\begin{equation}\label{cH}
\cH_n\deq\Big\{\sum_{i=0}^{n-1}
\chi_{[t_i^n,t_{i+1}^n)}(\cd)U(\cd,t_i^n)h_i\;\Big|\;
h_i\in L^{4}_{\cF_{t_i^n}}(\Om;H),\;t_i^n\in \D_n,\; i=0,\cds, n-1\Big\}.
\end{equation}
Here $U(\cd,\cd)$ is the operator introduced in
Lemma \ref{lemma5}. We have the following
result.
\begin{lemma}\label{5.13-prop1}
The set $\bigcup_{n=1}^\infty \cH_n$ is dense in
$L^2_\dbF(0,T;L^{4}(\Om;H))$.
\end{lemma}

The regularity result for solutions to
\eqref{op-bsystem3} is stated as follows (See
\cite{LZ2, LZ3.1} for its proof).

\begin{lemma}\label{10.1th}
Suppose that the assumptions in Theorem
\ref{OP-th2} hold and let $(P(\cd),$
$Q^{(\cd)},\widehat Q^{(\cd)})$ be the relaxed
transposition solution to the equation
\eqref{op-bsystem3}. Then, for each$n\in\dbN$, there exist two pointwise defined, linear
operators\footnote{One can define pointwise defined, linear
    operators similarly to that in \eqref{12.13-eq1}.} $Q^n$ and $\widehat Q^n$, both of
which are from $\cH_n$ to
$L^{2}_\dbF(0,T;L^{\frac{4}{3}}(\Om;H))$, such
that for any $\xi_1,\xi_2\in
L^{4}_{\cF_0}(\Om;H)$, $u_1(\cd), u_2(\cd)\in
L^{4}_\dbF(\Om;L^2(0,T;H))$ and
$v_1(\cd),v_2(\cd)\in \cH_{n}$, it holds that
\begin{eqnarray}\label{10.9eq2}
    &&\mE \int_{0}^T \big\langle v_1(s), \widehat
    Q^{(0)}(\xi_2,u_2,v_2)(s) \big\rangle_{H}ds +
    \mE \int_{0}^T \big\langle
    Q^{(0)}(\xi_1,u_1,v_1) (s),
    v_2(s) \big\rangle_{H}ds \nonumber\\
    && =\mE \int_{0}^T
    \[\big\langle
    Q^n(s) v_1(s), \f_2 (s)
    \big\rangle_{H}+\big\langle \f_1(s),
    \widehat Q^n(s) v_2(s)
    \big\rangle_{H}\]ds,
\end{eqnarray}
where $\f_1(\cd)$ and $\f_2(\cd)$ solve
\eqref{op-fsystem1} and \eqref{op-fsystem2} with
$t=0$, respectively.
\end{lemma}

%%%%%%%%%%%%%%%%%%%%%%%%%%%%%%%%%%%%%%%%%%%%%%%%%

\subsubsection{Pontryagin-type
maximum principle for nonconvex control domain}

%%%%%%%%%%%%%%%%%%%%%%%%%%%%%%%%%%%%%%%%%%%%%%%%%

We first assume the following further conditions for
Problem (OP).

\ms

\medskip

\no{\bf (AS5)} {\it The functions  $a(t,\cd,u)$,
$b(t,\cd,u)$,  $g(t,\cd,u)$ and $h(\cd)$ are $C^2$
w.r.t.  $x\in H$, such that for
$\f(t,x,u)=a(t,x,u), b(t,x,u)$, $\psi(t,x,u)=
g(t,x,u),h(x)$, it holds that
$\f_x(t,x,\cd)$, $\psi_x(t,x,\cd)$,\linebreak
$\f_{xx}(t,x,\cd)$ and $\psi_{xx}(t,x,\cd)$  are
continuous. Moreover, for
all $(x,u)\in  H\times U$ and a.e. $(t,\omega)\in [0,T]\times\Omega$,
\begin{equation}\label{ab1}
    \left\{
    \begin{array}{ll}\ds
        |\f_x(t,x,u)|_{\cL(H)} +|\psi_x(t,x,u) |_H
        \leq \cC,\\
        \ns\ds |\f_{xx}(t,x,u)|_{\cL(H, H;H)}
        +|\psi_{xx}(t,x,u)|_{\cL(H)} \leq \cC.
    \end{array}
    \right.
\end{equation}}

\ms

Let
\begin{equation}\label{H}
\begin{array}{ll}\ds
    \dbH(t,x,u,k_1,k_2) \deq  \big\langle k_1,a(t,x,u)  \big\rangle_H + \big\langle k_2, b(t,x,u)  \big\rangle_{H} - g(t,x,u),\\
    \ns\ds \hspace{4cm} (t,x,u,k_1,k_2)\in
    [0,T]\times H \times U\times H\times H.
\end{array}
\end{equation}

We have the following result.
\begin{theorem}\label{maximum p2}
Suppose that $L^p_{\cF_T}(\Omega)$ ($1\leq p <
\infty$) is a separable Banach space. Let the assumptions
(AS1), (AS2) and
(AS5) hold, and let $(\cl X(\cd),\bar
u(\cd))$ be an optimal pair of Problem
(OP). Let $\big(Y(\cdot),Z(\cdot)\big)$ be the
transposition solution to \eqref{ch-10-bsystem1}
with $Y_T$ and $f(\cd,\cd,\cd)$ given by
\eqref{zv1}. Assume that
$(P(\cd),Q^{(\cd)},\widehat Q^{(\cd)})$ is the
relaxed transposition solution to the equation
\eqref{op-bsystem3} in which $P_T$, $J(\cd)$,
$K(\cd)$ and $F(\cd)$ are given by
\begin{equation}\label{MP2-eq9}
    \left\{
    \begin{array}{ll} \ds P_T =
        -h_{xx}\big(\cl X(T)\big),\q J(t) =
        a_x(t,\cl X(t),\bar
        u(t)), \\
        \ns \ds K(t) =b_x(t,\cl X(t),\bar
        u(t)), \q F(t)= -\dbH_{xx}\big(t,\cl
        X(t),\bar u(t),Y(t),Z(t)\big).
    \end{array}
    \right.
\end{equation}
Then,   for $\ae$ $(t,\omega)\in [0,T]\times
\Omega$ and for all $u \in U$,
\begin{equation}\label{MP2-eq1}
    \begin{array}{ll}\ds
        \dbH\big(t,\cl X(t),\bar u(t),Y(t),Z(t)\big) - \dbH\big(t,\cl X(t),u,Y(t),Z(t)\big) \\
        \ns\ds \q - \frac{1}{2}\big\langle
        P(t)\big[ b\big(t,\cl X(t),\bar
        u(t)\big)-b\big(t,\cl X(t),u\big)
        \big], b\big(t,\cl X(t),\bar
        u(t)\big)-b\big(t,\cl X(t),u\big)
        \big\rangle_{H}\\
        \ns\ds \geq 0.
    \end{array}
\end{equation}
\end{theorem}

{\it Proof}. We divide the proof into several
steps.

\ms

{\bf Step 1}. For each $\e\in (0,T)$, let $E_\e\subset
[0,T]$ be a Lebesgue measurable set with the measure $\e$.
Put
$
u^\e(\cd) = \bar u(\cd)\chi_{[0,T]\setminus E_\e}(\cd)+u(\cd)\chi_{E_\e}(\cd),
$
where  $u(\cdot)$ is an arbitrarily given
element in $\cU[0,T]$.

Let us introduce some notations which will be used later. For $\psi=a,b,g$,  let
\begin{equation}\label{s7tatb1}
    \left\{
    \begin{array}{ll}
        \ds \psi_1(t) = \psi_x(t,\cl X(t),\bar
        u(t)), \q
        \psi_{11}(t) = \psi_{xx}(t,\cl X(t),\bar u(t)),\\
        \ns\ds \tilde \psi_1^\e (t)  = \int_0^1 \psi_x\big(t,\cl X(t) + \si (X^\e (t) - \cl X(t)), u^\e (t)\big)d\si, \\
        \ns\ds \tilde \psi_{11}^\e (t)  =
        2\int_0^1 (1-\si) \psi_{xx}\big(t,\cl
        X(t) + \si (X^\e(t) - \cl X(t)), u^\e
        (t)\big)d\si,
    \end{array}
    \right.
\end{equation}
and
\begin{equation}\label{s7tatb2}
    \left\{
    \begin{array}{ll}
        \ds \d \psi(t)  = \psi(t,\cl X(t), u(t)) -
        \psi(t,\cl X(t),\bar u(t)),
        \\
        \ns\ds \d \psi_1(t) = \psi_x(t,\cl X(t),
        u(t)) -
        \psi_x(t,\cl X(t),\bar u(t)),  \\
        \ns\ds \d \psi_{11}(t) = \psi_{xx}(t,\cl
        X(t), u(t)) - \psi_{xx}(t,\cl X(t),\bar
        u(t)).
    \end{array}
    \right.
\end{equation}

Consider the following equation:
\begin{eqnarray}\label{s7fsystem2}
    \left\{
    \begin{array}{lll}\ds
        d X^\e = \big(AX^\e +a(t,X^\e,u^\e)\big)dt + b(t,X^\e,u^\e)dW(t) &\mbox{ in }(0,T],\\
        \ns\ds X^\e(0)=X_0.
    \end{array}
    \right.
\end{eqnarray}
Let $X_1^\e(\cd) \deq X^\e(\cd)-\cl
X(\cd)$. Then, we see that
$X_1^\e(\cdot)$ solves the following
SEE:
\begin{equation}\label{fsystem3zx}
    \left\{
    \begin{array}{lll}\ds
        dX_1^\e = \big(AX_1^\e + \tilde
        a_1^\e(t) X^\e_1
        + \chi_{E_\e} (t)\d a(t) \big)dt\\
        \ns\ds\qq \q + \big(\tilde b_1^\e(t)
        X^\e_1 + \chi_{E_\e} (t)\d b(t)
        \big)dW(t) &\mbox{ in
        }(0,T],\\
        \ns\ds X_1^\e(0)=0.
    \end{array}
    \right.
\end{equation}

The equation \eqref{fsystem3zx} reminds us to introduce the following first order ``variational" equation:
\begin{equation}\label{s7fsystem3.1}
    \left\{
    \begin{array}{lll}\ds
        dX_2^\e = \big(AX_2^\e + a_1(t)X_2^\e \big)dt +
        \big(b_1(t) X_2^\e + \chi_{E_\e} (t)\d b(t)
        \big)dW(t) \q\mbox{ in
        }(0,T],\\
        \ns\ds X_2^\e(0)=0.
    \end{array}
    \right.
\end{equation}
We also need to introduce the following second order ``variational" equation\footnote{Recall that, by Subsection \ref{subsec0120}, for any
    $C^2$-function $f(\cd)$ defined on a
    Banach space $X$ and $X_0\in X$,
    $f_{xx}(X_0)\in \cL (X,X;X)$. This
    means that, for any $X_1,X_2\in X$,
    $f_{xx}(X_0)(X_1,X_2)\in X$. Hence, by
    (\ref{s7tatb1}),
    $a_{11}(t)\big(X_2^\e,X_2^\e\big)$ (in
    (\ref{s7fsystem3.2})) stands for
    $a_{xx}(t,\cl X(t),\bar
    u(t))\big(X_2^\e(t),X_2^\e(t)\big)$.
    One has a similar meaning for
    $b_{11}(t)\big(X_2^\e,X_2^\e\big)$ and
    so on.}
\begin{equation}\label{s7fsystem3.2}
    \left\{
    \begin{array}{lll}\ds
        dX_3^\e = \[AX_3^\e + a_1(t)X_3^\e +
        \chi_{E_\e}(t)\d a(t) +
        \frac{1}{2}a_{11}(t)\big(X_2^\e,X_2^\e\big) \]dt \\
        \ns\ds \hspace{1cm} + \[ b_1(t) X_3^\e +
        \chi_{E_\e} (t)\d b_1(t)X_2^\e +
        \frac{1}{2}b_{11}(t)\big(X_2^\e,X_2^\e\big)\]
        dW(t) &\mbox{ in
        }(0,T],\\
        \ns\ds X_3^\e(0)=0.
    \end{array}
    \right.
\end{equation}

Denote by $\cC(X_0)$ a generic
constant (depending on $X_0$, $T$ and $A$), which may change from line to
line. One can show that
\begin{equation}\label{s7th max eq2}
    |X_2^\e(\cd)|_{C_\dbF([0,T];L^8(\Om;H))}
    \leq \cC(X_0)\sqrt{\e},
\end{equation}
\begin{equation}\label{s7th max eq2.1b}
    |X_3^\e(\cd)|_{C_\dbF([0,T];L^2(\Om;H))}
    \leq \cC(X_0) \e.\;\;\,
\end{equation}

\ms

{\bf Step 2}. We now compute $\cJ(u^\e(\cd)) -
\cJ(\bar u(\cd))$.

Using Taylor's formula and the continuity
of both $h_{xx}(x)$ and $g_{xx}(t,x,u)$ with respect
to $x$, noting \eqref{s7th max eq2}--\eqref{s7th max eq2.1b}, after a careful computation, we can show that
\begin{equation}\label{s7eq5}
    \begin{array}{ll}
        \ds
        \cJ(u^\e(\cd)) -  \cJ(\bar u(\cd))\\
        \ns\ds = \mE \int_0^T\(\big\langle
        g_1(t),X_2^\e(t) + X_3^\e(t) \big\rangle_H +
        \frac{1}{2}\big\langle
        g_{11}(t)X_2^\e(t),X_2^\e(t) \big\rangle_H \\
        \ns\ds\qq\qq+
        \chi_{E_\e}(t)\d g(t)
        \)dt  + \mE \big\langle h_x\big(\bar
        x(T)\big), X_2^\e(T)+X_3^\e(T) \big\rangle_H \\
        \ns\ds\q+
        \frac{1}{2}\mE\big\langle h_{xx}\big(\bar
        x(T)\big)X_2^\e(t),X_2^\e(t) \big\rangle_H +
        o(\e),\qq\hb{as
        }\e\to0.
    \end{array}
\end{equation}

Next, we shall get rid of $X_2^\e(\cd)$ and
$X_3^\e(\cd)$ in \eqref{s7eq5} by solutions to
the equations \eqref{ch-10-bsystem1} and
\eqref{op-bsystem3}. It follows from the
definition of the transposition solution to the
equation \eqref{ch-10-bsystem1} (with $Y_T$ and
$f(\cd,\cd,\cd)$ given by \eqref{zv1}) that
\begin{equation}\label{5.26-eq6}
    \ba{ll}\ds
    -\mE\big\langle h_x(\bar
    x(T))),X_2^\e(T)\big\rangle_H -\mE
    \int_0^T \big\langle
    g_1(t),X_2^\e(t)\big\rangle_H dt\\
    \ns\ds = \mE
    \int_0^T\big\langle Z(t), \d
    b(t)\big\rangle_H\chi_{E_\e}(t) dt
    \ea
\end{equation}
and
\begin{eqnarray}\label{5.26-eq7}
    &&-\mE\big\langle h_x(\bar
    x(T))),X_3^\e(T)\big\rangle_H - \mE \int_0^T
    \big\langle g_1(t),X_3^\e(t)\big\rangle_H dt\nonumber
    \\
    && = \mE \int_0^T
    \Big[\frac{1}{2}\big( \big\langle
    Y(t),a_{11}(t)\big(X_2^\e(t),
    X_2^\e(t)\big) \big\rangle_H +
    \big\langle Z(t), b_{11}(t)\big(
    X_2^\e(t), X_2^\e(t)\big) \big\rangle_H \big) \nonumber\\
    && \hspace{1.8cm} + \chi_{E_\e}(t)
    \big( \big\langle Y(t),\d a(t)
    \big\rangle_H + \big\langle Y,\d
    b_1(t)X_2^\e(t) \big\rangle_H
    \big)\Big]dt.
\end{eqnarray}
According to \eqref{s7eq5}--\eqref{5.26-eq7}, we
conclude that
\begin{equation}\label{5.26-eq8}
    \begin{array}{ll}\ds
        \cJ(u^\e(\cd)) -  \cJ(\bar u(\cd)) \\
        \ns\ds  =
        \frac{1}{2}\mE\int_0^T\big(\big\langle
        g_{11}(t)X_2^\e(t),
        X_2^\e(t)\big\rangle_H - \big\langle
        Y(t),a_{11}(t)\big(X_2^\e(t),
        X_2^\e(t)\big) \big\rangle_H\\
        \ns\ds \qq\qq\q - \big\langle Y,
        b_{11}(t)\big(X_2^\e(t),
        X_2^\e(t)\big)\big\rangle_H \big)dt \\
        \ns\ds \q  + \mE\int_0^T
        \chi_{E_\e}(t)\big(\d g(t) -
        \big\langle Y(t),\d a(t)\big\rangle_H
        -\big\langle Z(t),\d b(t) \big\rangle_H
        \big)dt \\
        \ns\ds \q + \frac{1}{2}\mE \big\langle
        h_{xx}\big(\bar x(T)\big)X_2^\e(T), X_2^\e(T)
        \big\rangle_H + o(\e),\qq\hb{as
        }\e\to0.
    \end{array}
\end{equation}
By the definition of the relaxed transposition
solution to the equation \eqref{op-bsystem3}
(with $P_T$, $J(\cd)$, $K(\cd)$ and $F(\cd)$
given by \eqref{MP2-eq9}), we obtain that
\begin{eqnarray}\label{5.26-eq9}
    &&\3n\3n-\mE\big\langle
    h_{xx}\big(\bar x(T)\big) X_2^\e(T), X_2^\e(T)
    \big\rangle_H\nonumber\\
    && +  \mE \int_0^T \big\langle
    \dbH_{xx}\big(t,\bar
    X(t),\bar u(t),Y(t),Z(t)\big) X_2^\e(t), X_2^\e(t) \big\rangle_H dt\nonumber\\
    &&\3n\3n = \dbE\!\!\int_0^T\!\!\!
    \chi_{E_\e}(t)\big\langle b_1(t)X_2^\e(t),
    P(t)^*\d b(t)\big\rangle_{H} dt\! +
    \!\dbE\!\!\int_0^T\!\!\! \chi_{E_\e}(t)\big\langle
    P(t)\d b(t),
    b_1(t)X_2^\e(t)\big\rangle_{H} dt\nonumber\\
    &&  + \dbE\!\int_0^T\!\!
    \chi_{E_\e}(t)\big\langle P(t)\d b(t), \d
    b(t)\big\rangle_{H} dt \!+\! \dbE\!\int_0^T\!\!
    \chi_{E_\e}(t)\big\langle
    \d b(t),\widehat Q^{(0)}(0,0,\chi_{E_\e}\d b)(t)\big\rangle_{H} dt \nonumber\\
    &&   + \dbE\int_0^T
    \chi_{E_\e}(t)\big\langle Q^{(0)}(0,0,\d
    b)(t),\d b(t)\big\rangle_{H} dt.
\end{eqnarray}

Now, we estimate the terms in the right hand
side of \eqref{5.26-eq9}. By \eqref{s7th max
    eq2}, we have
\begin{equation}\label{s7eq9.1}
    \begin{array}{ll}\ds
        \Big|\dbE\int_0^T \chi_{E_\e}(t)\big\langle
        b_1(t)X_2^\e(t), P(t)^*\d b(t)\big\rangle_{H} dt
        \\
        \ns\ds\q+\dbE\int_0^T \chi_{E_\e}(t)\big\langle
        P(t)\d b(t), b_1(t)X_2^\e(t)\big\rangle_{H}
        dt\Big|= o(\e).
    \end{array}
\end{equation}

In what follows, for any $\tau\in [0,T)$, we
choose $E_{\e}=[\tau,\tau+\e]\subset [0,T]$.

By Lemma \ref{5.13-prop1}, we can find a
sequence $\{\beta_n\}_{n=1}^\infty$ such that
$\b_n\in\cH_n$ (Recall \eqref{cH} for the
definition of $\cH_n$) and
$\ds
\lim_{n\to\infty}\beta_n = \d b \;\mbox{ in }\;
L^2_\dbF(0,T;L^4(\Om;H)).$
Hence,
\begin{equation}\label{10.9qqq3}
    |\beta_n|_{L^2_\dbF(0,T;L^4(\Om;H))}\le
    \cC(X_0)<\infty,\qq\forall\;n\in\dbN,
\end{equation}
and there is a subsequence
$\{n_k\}_{k=1}^\infty\subset \dbN$ such that
\begin{equation}\label{s7eq9.2-11}
    \lim_{k\to\infty} |\b_ {n_k}(t)-\d
    b(t)|_{L^4_{\cF_t}(\Om;H)} = 0\q\mbox{ for }\ae
    t\in [0,T].
\end{equation}
Denote by $Q^{n_k}$ and $\widehat Q^{n_k}$ the
corresponding pointwise defined, linear operators
from $\cH_{n_k}$ to $L^2_\dbF(0,T;$
$L^{\frac{4}{3}}(\Om;H))$, given in Lemma
\ref{10.1th}.

Consider the following SEE:
\begin{equation}\label{s7fsystem3.1x}
    \left\{\2n\!
    \begin{array}{lll}\ds
        dX_{2,n_k}^{\e} \!=\! \big(AX_{2,n_k}^{\e}\!\! +
        a_1(t)X_{2,n_k}^{\e} \big)dt\! + \big( b_1(t)
        X_{2,n_k}^{\e} \!\!+\! \chi_{E_{\e}} (t)\b_ {n_k}(t)
        \big) dW(t) \q\mbox{in
        }(0,T],\\
        \ns\ds X_{2,n_k}^{\e}(0)=0.
    \end{array}
    \right.
\end{equation}
Applying Theorem \ref{ch-1-well-mild} to
\eqref{s7fsystem3.1x} and by \eqref{10.9qqq3},
we obtain that
\begin{equation}\label{s7th max eq2x}
    \begin{array}{ll}\ds
        |x_{2,n_k}^{\e}(\cd)|^4_{C_\dbF([0,T];L^4(\Om;H))} \\
        \ns\ds   \leq \cC \mE\Big( \int_0^T\chi_{E_{\e}}(s)\big|\b_ {k}(s)\big|_{H}^2ds\Big)^2   \leq \cC \e
        \int_{E_{\e}}\mE|\beta_{k}(s)|_{H}^4ds\leq
        \cC(x_0,k)\e^2.
    \end{array}
\end{equation}
Here and henceforth, $\cC(X_0,k)$ is a generic
constant (depending on $X_0$,
$k$, $T$ and $A$), which may
be different from line to line.

For any fixed $
k\in\dbN$, since $Q^{n_k}\b_ {n_k}\in
L^2_{\dbF}(0,T;L^{\frac{4}{3}}(\Om;H))$, by
\eqref{s7th max eq2x}, we find that
\begin{equation}\label{s7eq9.3}
    \begin{array}{ll}\ds
        \Big|\dbE\int_0^T \chi_{E_{\e}}(t)\big\langle
        \big(Q^{n_k} \b_
        {n_k}\big)(t),X_{2,n_k}^{\e}(t)\big\rangle_{H}
        dt \Big|\\
        \ns\ds \leq
        |X_{2,n_k}^{\e}(\cd)|_{L^\infty_\dbF(0,T;L^4(\Om;H))}
        \int_{E_{\e}}\big|\big(Q^{n_k} \b_ {n_k}\big)(t)\big|_{L^{\frac{4}{3}}_{\cF_t}(\Om;H)}dt\\
        \ns\ds \leq \cC(X_0,k)\sqrt{{\e}}
        \int_{E_{\e}}\big|\big(Q^{n_k} \b_
        {n_k}\big)(t)\big|_{L^{\frac{4}{3}}_{\cF_t}(\Om;H)}dt=
        o({\e}), \qq\hbox{as }\e\to0.
    \end{array}
\end{equation}
Similarly,
\begin{equation}\label{s7eq9.3x}
    \Big|\dbE\int_0^T \chi_{E_{\e}}(t)\big\langle
    X_{2,n_k}^{\e}(t),\big(\widehat Q^{n_k}\b_
    {n_k}\big)(t)\big\rangle_{H} dt \Big| = o({\e}),
    \qq\hbox{as }\e\to0.
\end{equation}

From \eqref{10.9eq2} in Lemma \ref{10.1th},
and noting that both $Q^{n_k}$ and $\widehat
Q^{n_k}$ are pointwise defined, we get
\begin{equation}\label{wws1}
    \begin{array}{ll}\ds
        \mE \int_{0}^T \big\langle \chi_{E_{\e}}(t)\b_
        {n_k}(t), \widehat Q^{(0)}(0,0,\chi_{E_{\e}}
        \b_{n_k})(t) \big\rangle_{H}dt\\
        \ns\ds\q +  \mE \int_{0}^T \big\langle
        Q^{(0)}(0,0,\chi_{E_{\e}} \b_ {n_k}) (t), \chi_{E_{\e}}(t)\b_ {n_k}(t) \big\rangle_{H}dt \\
        \ns\ds=\mE \int_{0}^T \chi_{E_{\e}}(t)\[\big\langle
        \big(Q^{n_k} \b_ {n_k}\big)(t), X_{2,n_k}^{\e}
        (t) \big\rangle_{H}+\big\langle X_{2,n_k}^{\e}
        (t), \big(\widehat Q^{n_k} \b_ {n_k}\big)(t)
        \big\rangle_{H}\]dt.
    \end{array}
\end{equation}

From \eqref{s7eq9.2-11} and the density of the
Lebesgue points, we find that for $\ae\tau\in
[0,T)$, it holds that
\begin{eqnarray}\label{s7eq9.3-1}
    &&\3n\3n\3n\lim_{k\to\infty}\lim_{\e\to
        0}\frac{1}{{\e}}\Big|\mE \int_{0}^T \big\langle
    \chi_{E_{\e}}(t)\d b(t),\widehat
    Q^{(0)}(0,0,\chi_{E_{\e}} \d b)(t)
    \big\rangle_{H}dt\nonumber\\
    &&\qq - \mE \int_{0}^T
    \big\langle
    \chi_{E_{\e}}(t)\d b(t),\widehat Q^{(0)}(0,0,\chi_{E_{\e}} \b_ {n_k})(t) \big\rangle_{H}dt \Big|\nonumber\\
    &&\3n\3n\3n\leq \lim_{k\to\infty}\!\lim_{\e\to
        0}\frac{1}{{\e}}\!\[\!\int_0^T\!
    \chi_{E_{\e}}(t) \(\!\mE|\d
    b(t)|^4_{H}\!\)^{\frac{1}{2}} dt\!
    \Big]^{\frac{1}{2}}|\widehat
    Q^{(0)}(0,0,\chi_{E_{\e}} (\d b-\b_ {n_k}))
    |_{L^2_\dbF(0,T;L^{\frac{4}{3}}(\Om;H))}\nonumber\\
    &&\3n\3n\3n\leq \cC\lim_{k\to\infty}\lim_{\e\to
        0}\frac{1}{{\e}}
    \[\int_0^T \chi_{E_{\e}}(t) \(\mE|\d b(t) |^4_{H}\)^{\frac{1}{2}} dt
    \Big]^{\frac{1}{2}}\big|\chi_{E_{\e}} (\d b\!-\!\b_ {n_k})\big|_{L^2_\dbF(0,T;L^4(\Om;H))}\nonumber\\
    &&\3n\3n\3n\leq \cC\lim_{k\to\infty}\lim_{\e\to
        0}\frac{|\d
        b(\tau)|_{L^4_{\cF_\tau}(\Om;H)}}{\sqrt{{\e}}}\[\int_0^T
    \chi_{E_{\e}}(t) \(\mE|\d b(t) - \b_
    {n_k}(t)|^4_{H}\)^{\frac{1}{2}} dt
    \Big]^{\frac{1}{2}}\\
    &&\3n\3n\3n= \cC\lim_{k\to\infty}\lim_{\e\to
        0}|\d
    b(\tau)|_{L^4_{\cF_\tau}(\Om;H)}\[\frac{1}{{\e}}\int_\tau^{\tau+{\e}}
    |\d b(t) - \b_ {n_k}(t)|_{L^4_{\cF_t}(\Om;H)}^2
    dt
    \Big]^{\frac{1}{2}}\nonumber\\
    &&\3n\3n\3n= \cC\lim_{k\to\infty}|\d
    b(\tau)|_{L^4_{\cF_\tau}(\Om;H)}|\d b(\tau) -
    \b_
    {n_k}(\tau)|_{L^4_{\cF_\tau}(\Om;H)}\nonumber\\
    &&\3n\3n\3n= 0.\nonumber
\end{eqnarray}
Similarly,
\begin{eqnarray}\label{s7eq9.3-1x}
    && \3n\3n\lim_{k\to\infty}\lim_{\e\to
        0}\frac{1}{{\e}}\Big|\mE \int_{0}^T \big\langle
    \chi_{E_{\e}}(t)\d b(t), \widehat
    Q^{(0)}(0,0,\chi_{E_{\e}} \b_ {n_k})(t)
    \big\rangle_{H}dt \nonumber\\
    && \3n\3n \qq\qq\q - \mE  \int_{0}^T
    \big\langle \chi_{E_{\e}}(t)\b_ {n_k}(t),
    \widehat Q^{(0)}(0,0,\chi_{E_{\e}} \b_ {n_k})(t)
    \big\rangle_{H}dt \Big|\nonumber\\
    && \3n\3n \leq
    \lim_{k\to\infty}\lim_{\e\to
        0}\frac{1}{{\e}}\big|\widehat
    Q^{(0)}(0,0,\chi_{E_{\e}} \b_
    {n_k})\big|_{L^2_\dbF(0,T;L^{\frac{4}{3}}(\Om;H))}\nonumber\\
    &&\3n\3n\qq\qq\qq\times
    \[\int_0^T \chi_{E_{\e}}(t) \(\mE|\d b(t) - \b_ {n_k}(t)|^4_{H}\)^{\frac{1}{2}} dt
    \Big]^{\frac{1}{2}}\nonumber\\
    && \3n\3n \leq
    \cC\lim_{k\to\infty}\lim_{\e\to
        0}\frac{1}{{\e}}\big|\chi_{E_{\e}} \b_
    {n_k}\big|_{L^2_\dbF(0,T;L^4(\Om;H))}
    \[\int_0^T \chi_{E_{\e}}(t) \(\mE|\d b(t) - \b_ {n_k}(t)|^4_{H}\)^{\frac{1}{2}} dt
    \Big]^{\frac{1}{2}}\nonumber\\
    && \3n\3n \leq
    \cC\lim_{k\to\infty}\lim_{\e\to
        0}\frac{1}{{\e}}\Big\{\big|\chi_{E_{\e}} \d
    b\big|_{L^2_\dbF(0,T;L^4(\Om;H))}
    \[\int_0^T \chi_{E_{\e}}(t) \(\mE|\d b(t) - \b_ {n_k}(t)|^4_{H}\)^{\frac{1}{2}} dt
    \Big]^{\frac{1}{2}}\nonumber\\
    &&\3n\3n \qq\qq\qq\qq\q+\int_0^T \chi_{E_{\e}}(t) \(\mE|\d b(t) - \b_ {n_k}(t)|^4_{H}\)^{\frac{1}{2}} dt\Big\}\\
    && \3n\3n \leq
    \cC\lim_{k\to\infty}\lim_{\e\to
        0}\bigg\{\frac{|\d
        b(\tau)|_{L^4_{\cF_\tau}(\Om;H)}}{\sqrt{{\e}}}\[\int_0^T
    \chi_{E_{\e}}(t) \(\mE|\d b(t) - \b_
    {n_k}(t)|^4_{H}\)^{\frac{1}{2}} dt
    \Big]^{\frac{1}{2}}\nonumber\\&&\3n\3n
    \qq\qq\qq\q+\frac{1}{{\e}}\int_0^T
    \chi_{E_{\e}}(t) \(\mE|\d b(t) - \b_
    {n_k}(t)|^4_{H}\)^{\frac{1}{2}} dt
    \bigg\}\nonumber\\
    && \3n\3n =
    \cC\lim_{k\to\infty}\lim_{\e\to 0}\Big[|\d
    b(\tau)|_{L^4_{\cF_\tau}(\Om;H)}\(\frac{1}{{\e}}\int_\tau^{\tau+{\e}}
    |\d b(t) - \b_ {n_k}(t)|_{L^4_{\cF_t}(\Om;H)}^2
    dt
    \Big)^{\frac{1}{2}}\nonumber\\
    &&\3n\3n
    \qq\qq\qq\q+\frac{1}{{\e}}\int_\tau^{\tau+{\e}}
    |\d b(t) - \b_ {n_k}(t)|_{L^4_{\cF_t}(\Om;H)}^2
    dt
    \Big]\nonumber\\
    && \3n\3n =
    \cC\lim_{k\to\infty}\big[|\d
    b(\tau)|_{L^4_{\cF_\tau}(\Om;H)}|\d b(\tau) -
    \b_ {n_k}(\tau)|_{L^4_{\cF_\tau}(\Om;H)}+|\d
    b(\tau) - \b_
    {n_k}(\tau)|_{L^4_{\cF_\tau}(\Om;H)}^2\big]\nonumber
    \\
    && = 0.\nonumber
\end{eqnarray}
From \eqref{s7eq9.3-1}--\eqref{s7eq9.3-1x}, we
find that
\begin{eqnarray}\label{s7eq9.3-1xx}
    &&
    \lim_{k\to\infty} \lim_{\e\to
        0}\frac{1}{{\e}}\Big|\mE \int_{0}^T \big\langle
    \chi_{E_{\e}}(t)\d b(t), \widehat
    Q^{(0)}(0,0,\chi_{E_{\e}} \d b)(t)
    \big\rangle_{H}dt\\
    && \qq\qq\q - \mE  \int_{0}^T \big\langle
    \chi_{E_{\e}}(t)\b_ {n_k}(t), \widehat
    Q^{(0)}(0,0,\chi_{E_{\e}} \b_ {n_k})(t)
    \big\rangle_{H}dt \Big| = 0.\nonumber
\end{eqnarray}
By a similar argument, we get that
\begin{equation}\label{s7eq9.3-1xxx}
    \begin{array}{ll}\ds
        \lim_{k\to\infty} \lim_{\e\to
            0}\frac{1}{{\e}}\Big|\mE \int_{0}^T \big\langle
        Q^{(0)}(0,0,\chi_{E_{\e}} \d
        b)(t),\chi_{E_{\e}}(t)\d b(t)
        \big\rangle_{H}dt \\
        \ns\ds \qq\qq\q - \mE  \int_{0}^T \big\langle
        Q^{(0)}(0,0,\chi_{E_{\e}} \b_
        {n_k})(t),\chi_{E_{\e}}(t)\b_ {n_k}(t)
        \big\rangle_{H}dt \Big| = 0.
    \end{array}
\end{equation}

From \eqref{s7eq9.3}--\eqref{s7eq9.3x} and
\eqref{s7eq9.3-1xx}--\eqref{s7eq9.3-1xxx}, we
obtain that
\begin{equation}\label{s7eq9.2-111}
    \ba{ll}\ds \Big|\dbE\int_0^T
    \chi_{E_{\e}}(t)\big\langle \d b(t),\widehat
    Q^{(0)}(0,0,\chi_{E_{\e}}\d b)(t)\big\rangle_{H}
    dt\\
    \ns\ds   \q+ \dbE\int_0^T
    \chi_{E_{\e}}(t)\big\langle Q^{(0)}(0,0,\d
    b)(t),\d b(t)\big\rangle_{H} dt\Big| =o({\e}),
    \qq\hbox{as }\e\to0.
    \ea
\end{equation}

Combining \eqref{5.26-eq8}--\eqref{s7eq9.1} and
\eqref{s7eq9.2-111}, we end up with
$$
\ba{ll}\ds
\cJ(u^{\e}(\cd)) -  \cJ(\bar u(\cd))\\
\ns\ds = \mE\int_0^T \( \d g(t) -
\big\langle Y(t),\d a(t)\big\rangle_H
-\big\langle Z(t),\d b(t)
\big\rangle_H \\
\ns\ds\qq\qq - \frac{1}{2}\big\langle P(t)\d
b(t), \d b(t) \big\rangle_H \)\chi_{E_{\e}}(t)dt
+ o({\e}),\qq\hbox{as }\e\to0. \ea
$$
Since $\bar u(\cd)$ is the optimal control,
$\cJ(u^{\e}(\cd)) - \cJ(\bar u(\cd))\geq 0$.
Thus,
\begin{equation}\label{s7eq11}
    \begin{array}{ll}\ds
        \mE\int_0^T \chi_{E_{\e}}(t)\(
        \big\langle Y(t),\d a(t)\big\rangle_H
        +\big\langle Z(t),\d
        b(t) \big\rangle_H -\d g(t) \\
        \ns\ds\qq\qq\qq + \frac{1}{2}\big\langle P(t)\d
        b(t), \d b(t) \big\rangle_H \)dt \leq o({\e}),
    \end{array}
\end{equation}
as $\e\to0$.

\ms

{\bf Step 3}. We are now in a position to
complete the proof.

Since $L^2_{\cF_T}(\Om)$ is separable, for any
$t\in [0,T]$, $\cF_t$ is countably generated by
a sequence $\{M_k\}_{k=1}^\infty \subset \cF_t$,
that is, for any $M\subset \cF_t$, there exists
a subsequence $\{M_{k_j}\}_{j=1}^\infty \subset
\{M_k\}_{k=1}^\infty$ such that
$$
\lim_{j\to\infty} \dbP\big((M\setminus M_{k_j})
\cup(M_{k_j}\setminus M)\big) = 0.
$$

Denote by $\{t_i\}_{i=1}^\infty$ the sequence of
rational numbers in $[0, T)$, and by
$\{v_i\}_{i=1}^\infty$ a dense subset of $U$.
For each $i\in\dbN$, we choose
$\{M_{ij}\}_{j=1}^\infty\subset \cF_{t_i}$ to be
a sequence which generates  $\cF_{t_i}$.

Fix $i, j, k \in \dbN$ arbitrarily. For any
$\tau\in [t_i, T)$ and $\th\in (0, T-\tau)$,
write $E_\th^i = [\tau, \tau+\th)$. Put
$$
u_{ij}^{k,\th}=\begin{cases} v^k, &\mbox{ if
    }(t,\omega)\in E_\th^i\times M_{ij},\\
    \ns\ds \bar u(t,\om), &\mbox{ if }(t,\omega)\in
    ([0,T]\times\Om)\setminus(E_\th^i\times M_{ij}).
\end{cases}
$$
Clearly, $u_{ij}^{k,\th}\in \cU[0,T]$  and
$$
u_{ij}^{k,\th}(t,\om)-\bar u(t,\om)=(v^k-\bar
u(t,\om))\chi_{M_{ij}}(\om)\chi_{E_\th^i}(t),\qq
(t,\om)\in [0,T]\times\Om.
$$

From \eqref{s7eq11}, we know that
\begin{eqnarray}\label{s7eq12}
    &&\mE\int_\th^{\th+\tau} \(\dbH\big(t,\cl X(t),\bar u(t),Y(t),Z(t)\big)- \dbH\big(t,\cl X(t),u_{ij}^{k,\th}(t),Y(t),Z(t)\big) \nonumber\\
    &&\hspace{1.21cm}  -
    \frac{1}{2}\big\langle P(t)\big(
    b\big(t,\cl X(t),\bar
    u(t)\big)-b\big(t,\cl X(t),
    u_{ij}^{k,\th}(t)\big) \big), \\
    &&\hspace{1.81cm}  b\big(t,\bar
    X(t),\bar u(t)\big)-b\big(t,\bar
    X(t),u_{ij}^{k,\th}(t)\big)
    \big\rangle_{H}
    \)dt\nonumber\\
    && \leq o({\e}).\nonumber
\end{eqnarray}
Dividing both sides of \eqref{s7eq12} by $\tau$
and letting $\tau\to 0^+$, by the property of
Lebesgue point, we conclude that for any $i, j, k
\in \dbN$, there exists a zero measure set
$E_{i,j}^k \subset [t_i, T)$  such that for all
$t\in [t_i, T)\setminus E_{i,j}^k$,
$$
\begin{array}{ll}\ds
    \mE\[\(\dbH\big(t,\cl X(t),\bar u(t),Y(t),Z(t)\big) - \dbH\big(t,\cl X(t),v^{k}(t),Y(t),Z(t)\big) \\
    \ns\ds \q  - \frac{1}{2}\big\langle
    P(t)\big( b\big(t,\cl X(t),\bar
    u(t)\big)-b\big(t,\cl
    X(t), v^{k}(t)\big) \big),\\
    \ns\ds\qq\q b\big(t,\cl X(t),\bar
    u(t)\big)-b\big(t,\cl
    X(t),v^{k}(t)\big) \big\rangle_{H}
    \)\chi_{M_{ij}}\] \leq 0.
\end{array}
$$
Let $$\ds E_0 = \bigcup_{i,j,k\in\dbN}
E_{i,j}^k.$$ Then its Lebesgue measure is zero,
and for any $i,j,k\in\dbN$,
\begin{equation*}\label{s7eq14}
    \begin{array}{ll}\ds \mE\[\(\dbH\big(t,\cl X(t),\bar u(t),Y(t),Z(t)\big) - \dbH\big(t,\cl X(t),v^{k}(t),Y(t),Z(t)\big) \\
        \ns\ds \hspace{0.7cm}  -
        \frac{1}{2}\big\langle P(t)\big(
        b\big(t,\cl X(t),\bar
        u(t)\big)-b\big(t,\cl X(t),
        v^{k}(t)\big)
        \big),\\
        \ns\ds  \hspace{1.19cm}  b\big(t,\cl
        X(t),\bar u(t)\big)-b\big(t,\cl
        X(t),v^{k}(t)\big)
        \big\rangle_{H} \)\chi_{M_{ij}}\] \leq 0, \\
        \ns\ds \hspace{6.7cm} \forall\; t\in [t_i,
        T)\setminus E_{0}.
    \end{array}
\end{equation*}
This, together with the construction of
$\{M_{ij}\}_{j=1}^\infty$, the right continuity
of the filtration $\dbF$ and the density of
$\{v_k\}_{k=1}^\infty$, implies that
\eqref{MP2-eq1} holds. This completes the proof
of Theorem \ref{maximum p2}.
\endpf

%%%%%%%%%%%%%%%%%%%%%%%%%%%%%%%%%%%%%%%%%%%%%%%%%%%%%%%%%%

\section{A sufficient condition  for an optimal control}\label{ch-Pon2}

%%%%%%%%%%%%%%%%%%%%%%%%%%%%%%%%%%%%%%%%%%%%%%%%%%%%%%%%%%%

Both Theorems \ref{th max} and \ref{maximum p2} are about first odder necessary
conditions for optimal controls of Problem
(OP).  In this
section,  we shall consider very quickly the related first order sufficient
(optimality) condition. For simplicity, we assume that the control region $U$ is a convex
subset of the separable Hilbert space  $H_1$ appeared in (AS3) (in Subsection \ref{subsec20210120}), and $L^p_{\cF_T}(\Om)$ ($p\geq 1$) is
separable. The finite version is obtained in \cite{Zhou1}. We borrow some idea from that paper.

\subsection{Clarke's generalized gradient}

In this subsection, as preliminaries,  we recall the definition and
some basic properties of Clarke's generalized
gradient.  A detailed introduction to this topic
can be found in \cite{Clarke}.

In the rest of this subsection, $O$ is a domain in $H$ and $\f:O\to
\dbR$ is a locally Lipschitz continuous
function.

\begin{definition}\label{3.30-def1}
For any $x\in O$, the {\it (Clarke)
    generalized gradient} of $\f$ at $x$ is defined by
\begin{equation}\label{3.1-eq20}
    \pa\f(x)\deq\Big\{\xi\in
    H\;\Big|\;\langle\xi,y\rangle_H\leq\limsup_{{z\to
            x, z\in
            O}\atop{r\to 0^+}} \frac{\f(z+ry)-\f(z)}{r}\Big\}.
\end{equation}
\end{definition}
The following lemma provides some properties of the Clarke
generalized gradient.
\begin{lemma}\label{3.1-lm3}
The following results hold:

\ms

{\rm 1)} $\pa\f(x)$ is a nonempty, bounded, convex subset of $H$;

\ms

{\rm 2)} $\pa(-\f)(x)=-\pa\f(x)$;

\ms

{\rm 3)} $0\in\pa\f(x)$ if $\f$ attains a local
minimum or maximum at $x$.
\end{lemma}

{\it Proof}. 1) The boundedness of $\pa\f(x)$
follow immediately from \eqref{3.1-eq20}. We
only need to show the convexity and
non-emptiness. Since $\f:O\to \dbR$ is locally
Lipschitz continuous, for the given $x\in O$ and
a small $\d>0$ with
$$
B_\d(x)\deq\big\{z\in
O\bigm||z-x|\le\d\big\}\subseteq O,
$$
there
exists a constant $\cC>0$ such that
$$|\f(z)-\f(\h z)|\leq \cC|z-\h z|_H,\qq\forall\; z,\h z\in B_\d(x).$$
Thus, for any fixed $y\in H$ and $r>0$  small enough,
$$
\frac{|\f(z+ry)-\f(z)|}{r}\le \cC |y|_H,\qq\forall\; z\in B_{\d/2}(x).
$$
Hence, the following functional is well-defined:
\begin{equation}\label{4.11-eq4}
    \f^0(x;y)\deq\limsup_{{z\to
            x, z\in
            O}\atop{r\to 0^+}}\frac{\f(z+ry)-\f(z)}{r}.
\end{equation}
It is easy to check that
\begin{equation}\label{3.1-eq23}
    \f^0(x;\l y)=\l\f^0(x;y),\q\forall\; y\in
    H,~\l\geq0
\end{equation}
and that
\begin{equation}\label{3.1-eq23.1}
    \f^0(x;y+z)\le\f^0(x;y)+\f^0(x;z),\q\forall\;
    y,z\in H.
\end{equation}
Thus, $y\mapsto\f^0(x;y)$ is a convex function.
Also,  \eqref{3.1-eq23.1}
implies that
\begin{equation}\label{3.1-eq24}
    -\f^0(x;y)\le\f^0(x;-y),\q\forall\; y\in H.
\end{equation}

Next, we fix a $z\in H$ and define $F:\{\l
z\bigm|\l\in\dbR\}\to\dbR$ by
$$
F(\l z)=\l\f^0(x;z),\q\forall\;\l\in\dbR.
$$
Then for $\l\ge0$ (noting \eqref{3.1-eq23}--\eqref{3.1-eq24}),
$$
\left\{
\begin{array}{ll}\ds
    F(\l z)\equiv\l\f^0(x;z)=\f^0(x;\l z),\\
    \ns\ds F(-\l z)\equiv-\l\f^0(x;z)=-\f^0(x;\l
    z)\le\f^0(x,-\l z),
\end{array}
\right.
$$
which implies
\begin{equation}\label{3.1-eq25}
    F(\l z)\le\f^0(x;\l z),\q\forall\;\l\in\dbR.
\end{equation}
Therefore, $F$ is a linear functional defined on
the linear space spanned by $z$, and it is
dominated by the convex function $\f^0(x;\cd)$.
By the Hahn-Banach theorem (i.e., Theorem \ref{Hahn-Banach}), $F$ can be extended
to be a bounded linear functional on $H$. Then,
by the classical Riesz representation theorem (i.e., Theorem \ref{1t10s}), there exists
$\xi\in H$, such that
\begin{equation}\label{3.1-eq26}
    \left\{
    \begin{array}{ll}\ds
        \lan\xi,\l z\rangle=F(\l
        z)\equiv\l\f^0(x;z),\q\forall\;\l\in\dbR,\\
        \ns\ds\lan\xi,y\rangle\le\f^0(x;y),\q\forall\;
        y\in H.
    \end{array}
    \right.
\end{equation}
This implies $\xi\in\pa\f(x)$. Consequently, $\pa\f(x)$ is nonempty.

\ms

2) It follows from \eqref{4.11-eq4} that
$$
(-\f)^0(x;y) =\limsup_{z\to x\atop{r\to 0^+}}\frac{-\f(z+ry)+\f(z)}{r} =\limsup_{z'\to x\atop{r\to 0^+}}\frac{-\f(z')+\f(z'-ry)} {r}=\f^0(x;-y).
$$
Thus, $\xi\in\pa(-\f)(x)$ if and only if
$$
\begin{array}{ll}\ds
    \lan-\xi,y\rangle_H
    =\lan\xi,-y\rangle_H\le(-\f)^0(x;-y)=\f^0(x;y),\q\forall\;
    y\in H,
\end{array}
$$
which is equivalent to $-\xi\in\pa\f(x)$.

\ms

3) Suppose $\f$ attains a local minimum at $x$.
Then
$$
\begin{array}{ll}
    \ds\f^0(x;y)\3n&\ds=\limsup_{{z\to x\atop{r\to 0^+}}}\frac{\f(z+ry)-\f(z)}{r}\\
    \ns&\ds  \ge\limsup_{r\to
        0^+}\frac{\f(x+ry)-\f(x)}{r}\geq0=\lan0,y\rangle_H,\q\forall\; y\in H,
\end{array}
$$
This implies $0\in\pa\f(x)$.

If $\f$ attains a local maximum at $x$, then the
conclusion follows from the fact that $-\f$
attains a local minimum at $x$.
\endpf

By fixing some arguments of a function, one may
naturally define its {\it partial generalized
gradient}. For example,   if $\psi: H\times H_1\to \dbR$ is locally Lipschitz,
by $\pa_x\psi(x, u)$ (\resp $\pa_u\psi(x, u)$), we mean the partial generalized gradient of $\psi$
in $x$ (\resp in $u$) at $(x, u) \in H \times U$.

Next, we present a technical lemma.

\begin{lemma}\label{3.1-lm4}
Let $\psi$ be a convex or concave function on
$H\times H_1$. Assume
that $\psi(\cd,u)$ is differentiable and
$\psi_x(\cd,\cd)$ is continuous. Then
\begin{equation}\label{3.1-eq36}
    \begin{array}{ll}\ds
        \big\{(\psi_x(\hat x,\hat
        u),r)\bigm|r\in\pa_u\psi(\hat x,\hat
        u)\big\}\subseteq\pa_{x,u} \psi(\hat x,\hat
        u),\q\forall\;(\hat x,\hat u)\in H\times U.
    \end{array}
\end{equation}
\end{lemma}

{\it Proof}.  We first handle the case that
$\psi$ is convex. For any $\xi\in H$ and $u\in
H_1$, we choose a sequence
$\{(x_j,\d_j)\}_{j=1}^\infty\subset H\times
\dbR$ as follows:
$$
\begin{array}{ll}\ds
    (x_j,\hat u)\in H\times U,\q
    (x_j+\d_j\xi,\hat u+\d_ju)\in H\times U, \\
    \ns\ds \d_j\to 0^+,\mbox{ as }j\to\infty,\hb{
        and }|x_j-\hat x|_H\le \d_j^2.
\end{array}
$$
By the convexity of $\psi$, we have
\begin{equation}\label{4.11-eq5}
    \begin{array}{ll}\ds
        \lim_{j\to\infty}\frac{\psi(x_j+\d_j\xi,\hat
            u+\d_ju)-\psi(\hat x,\hat u+\d_ju)}{\d_j}\cr
        \ns\ds
        \ge\lim_{j\to\infty}\frac{\langle\psi_x(\hat
            x,\hat u+\d_ju),x_i-\hat
            x+\d_j\xi\rangle_H}{\d_j}=\langle\psi_x(\hat
        x,\hat u),\xi\rangle_H.
    \end{array}
\end{equation}
Similarly,
\begin{equation}\label{4.11-eq6}
    \lim_{j\to\infty}\frac{\psi(\hat x,\hat
        u+\d_ju)-\psi(\hat x,\hat u)}{\d_j} \geq \langle
    r,u\rangle_{H_1}.
\end{equation}
Also,
\begin{equation}\label{4.11-eq7}
    \lim_{j\to\infty}\frac{\psi(\hat x,\hat
        u)-\psi(x_j,\hat u)}{\d_j}\geq\lim_{j\to\infty}
    \frac{\langle\psi_x(x_j,\hat u),\hat
        x-x_j\rangle_H}{\d_j}=0.
\end{equation}
It follows from
\eqref{4.11-eq5}--\eqref{4.11-eq7} that
$$
\lim_{j\to\infty}\frac{\psi(x_j+\d_j\xi,\hat
    u+\d_ju)-\psi(x_j,\hat u)}{\d_j}\geq
\langle\psi_x(\hat x,\hat u),\xi\rangle_H+\langle
r,u\rangle_{H_1}.
$$
This, together with the definition of the
generalized gradient, implies that\linebreak $(\psi_x(\hat
x,\hat u),r)\in\pa_{x,u}\psi(\hat x,\hat u)$.

\ss

If $\psi$ is concave, the desired result follows
immediately by noting that $-\psi$ is convex and
the conclusion 2) in Lemma \ref{3.1-lm3}. \endpf

\subsection{A sufficient condition of the optimal control}

Let $u(\cd)$ be an
admissible control and $X(\cd)$ be the
corresponding state. Let
$\big(Y(\cdot),Z(\cdot)\big)$ be the
transposition solution to \eqref{ch-10-bsystem1}
with $Y_T$ and $f(\cd,\cd,\cd)$ given by
\begin{equation}\label{zv1.1.1}
\left\{
\begin{array}{ll}
    \ds Y_T =
    -h_x\big(X(T)\big),\\
    \ns \ds f(t,
    X(t), u(t))=-a_x(t,
    X(t), u(t))^*X(t) - b_x\big(t,
    X(t), u(t)\big)^*u(t)\\
    \ns\ds\hspace{2.65cm} +
    g_x\big(t, X(t), u(t)\big)
\end{array}
\right.
\end{equation}
and
$(P(\cd),Q^{(\cd)},\widehat Q^{(\cd)})$ be the
relaxed transposition solution to the equation
\eqref{op-bsystem3} in which $F(\cd)$, $J(\cd)$,
$K(\cd)$ and $P_T$ are given by
\begin{equation}\label{20160221MP2-eq9.1}
\left\{
\begin{array}{ll} \ds F(t)= -\dbH_{xx}\big(t,
    X(t), u(t),Y(t),Z(t)\big),\q P_T = -h_{xx}\big(x(T)\big), \\
    \ns \ds K(t) =b_x(t, X(t), u(t)), \q J(t) = a_x(t, X(t),
    u(t)).
\end{array}
\right.
\end{equation}
Recall \eqref{H} for the definition of $\dbH$ and let
$$
\3n\3n\begin{array}{ll}\ds
\cH(t, \xi, \eta)\!\3n&\ds\deq \dbH\big(t,
\xi,\eta,Y(t),Z(t)\big) +\frac{1}{2} \big\langle
P(t)b\big(t,\xi, \eta\big), b\big(t,
\xi, \eta\big)\big\rangle_{H}\\
\ns&\ds\q -\big\langle P(t)b\big(t, X(t),
u(t)\big), b\big(t,
\xi, \eta\big)\big\rangle_{H},\qq \forall\;
(t,\xi,\eta)\in [0,T]\times H\times U.
\end{array}
$$
\begin{lemma}\label{3.4-lm1}
Suppose that (AS1)--(AS3) and (AS5) hold.
Then for a.e. $(t,\om)\in[0,T]\times \Om$,
\begin{equation}\label{3-eq1}
    \pa_u \dbH(t, X(t),
    u(t),Y(t),Z(t))=\pa_u\cH(t, X(t),
    u(t)).
\end{equation}
\end{lemma}
\ms

{\it Proof}. \rm Fix a $t\in[0,T]$. Let
$$
\left\{
\begin{array}{ll}\ds
    \dbH(\eta)\deq  \dbH\big(t, X(t),\eta,Y(t),Z(t)\big),  \\
    \ns\ds \cH(\eta)\deq  \cH\big(t, X(t),\eta\big),\q b(\eta)\deq  b\big(t, X(t),\eta\big),\\
    \ns\ds \psi(\eta)\deq  \frac{1}{2} \big\langle
    P(t)b(\eta), b(\eta)\big\rangle_H - \big\langle
    P(t)b(u(t)),b(\eta)\big\rangle_H.
\end{array}
\right.
$$
Then,
$
\cH(\eta)=\dbH(\eta)+\psi(\eta).
$
Note that for any $r\to 0^+$, $\eta,\k\in U$, with
$\eta\to\bar u(t)$,
$$
\begin{array}{ll}\ds
    \psi(\eta+r\k)-\psi(\eta)\3n&\ds =\frac{1}{2}
    \big\langle P(t)\big[b(\eta+r\k)+b(\eta)
    -2b(u(t))\big],
    b(\eta+r\k)-b(\eta)\big\rangle_H\\
    \ns&\ds=o(r).
\end{array}
$$
Thus,
$$
\limsup_{{\eta\to\bar u(t),\eta\in U}\atop{r\to
        0^+}}{\frac{\cH(\eta+r\k)-\cH(\eta)}{r}}=
\limsup_{{\eta\to\bar u(t),\eta\in U}\atop{r\to
        0^+}}{\frac{\dbH(\eta+r\k)-\dbH(\eta)}{r}}.
$$
Consequently, by \eqref{3.1-eq20}, the desired
result \eqref{3-eq1} follows. \endpf

Now, we present the following sufficient
condition of optimality for Problem
(OP).

\begin{theorem}\label{3.4-th2}
Let  (AS1)--(AS3) and (AS5) hold.
Suppose that $h(\cd)$ is convex,
$\dbH(t,\cd\,,\cd\,,Y(t),$ $Z(t))$ is
concave for all $t\in[0,T]$,
a.s.,  and
\begin{equation}\label{3.4-eq2}
    \cH(t, X(t), u(t))=\max_{u\in U}\cH(t,
    X(t),u),\q\ae (t,\om)\in[0,T]\times\Om.
\end{equation}
Then, $(X(\cd), u(\cd))$ is an optimal pair of
Problem (OP).
\end{theorem}

{\it Proof}.  By the maximum condition
\eqref{3.4-eq2}, Lemma \ref{3.4-lm1} and the
assertion 3) in  Lemma \ref{3.1-lm3}, we have
\begin{equation}\label{3.4-eq3}
    0\in\pa_u\cH(t, X(t), u(t))=\pa_u
    \dbH(t, X(t), u(t),Y(t),Z(t)).
\end{equation}
Hence, by  Lemma \ref{3.1-lm4}, we conclude that
\begin{equation}\label{3.4-eq4}
    \Big(\dbH_x(t, X(t),
    u(t),Y(t),Z(t)),0\Big)\in\pa_{x,u}\dbH(t,
    X(t), u(t),Y(t),Z(t)).
\end{equation}
This, together with the concavity of
$\dbH(t,\cd\,,\cd\,,Y(t),Z(t))$, implies that
\begin{equation}\label{3.4-eq5}
    \begin{array}{ll}\ds
        \int_0^T\big(\dbH(t,\wt
        X(t),\tilde u(t),Y(t),Z(t))-\dbH(t,
        X(t), u(t),Y(t),Z(t))
        \big)dt\\
        \ns\ds\le \int_0^T\big\langle \dbH_x(t,
        X(t), u(t),Y(t),Z(t)),\wt X(t)-
        X(t)\big\rangle_H dt,
    \end{array}
\end{equation}
for any admissible pair $(\wt
X(\cd),\tilde u(\cd))$. Let
$\xi(t)\deq \wt X(t)- X(t)$. Then $\xi(\cd)$ solves
\begin{equation}\label{3.4-eq6}
    \left\{
    \begin{array}{ll}\ds
        d\xi(t)=\big[A\xi(t) + a_x\big(t, X(t),
        u(t)\big)\xi(t)+\a(t)\big]dt\\
        \ns\ds \hspace{1.2cm} + \big[b_x\big(t,
        X(t), u(t)\big)\xi(t)+\b(t)\big]dW(t)
        &\mbox{ in }
        (0,T],\\
        \ns\ds\xi(0)=0,
    \end{array}
    \right.
\end{equation}
where
$$
\left\{
\begin{array}{ll}
    \ds\a(t)\deq -a_x\big(t, X(t),
    u(t)\big)\xi(t) +a\big(t,\wt
    X(t),\tilde u(t)\big)-a\big(t,
    X(t), u(t)\big),\\
    \ns\ds\b(t)\deq -b_x\big(t, X(t),
    u(t)\big)\xi(t) +b\big(t,\wt
    X(t),\tilde u(t)\big)-b\big(t, X(t),
    u(t)\big).
\end{array}
\right.
$$
It follows from \eqref{3.4-eq5}, \eqref{3.4-eq6} and the definition of the
transposition solution to \eqref{ch-10-bsystem1}
that
\begin{eqnarray}\label{3.4-eq6.1}
    &&\mE\big\langle h_x(X(T)),\xi(T)\big\rangle_H\nonumber\\
    &&=-\mE\big\langle
    Y(T),\xi(T)\big\rangle_H+\mE\big\langle Y(0),
    \xi(0)\big\rangle_H\nonumber\\
    &&=-\mE\int_0^T\[\big\langle g_x(t,
    X(t),
    u(t)),\xi(t)\big\rangle_H+\big\langle
    Y(t),\a(t) \big\rangle_H+ \big\langle
    Z(t),\b(t)\big\rangle_H\]dt\nonumber\\
    &&=\mE\int_0^T\big\langle \dbH_x(t,
    X(t),
    u(t),Y(t),Z(t)),\xi(t)\big\rangle_H dt\nonumber \\
    &&\q -\mE\int_0^T\(\big\langle
    Y(t),a(t,\wt X(t),\tilde u(t))-a(t,
    X(t),
    u(t))\big\rangle_H\nonumber\\
    &&\qq\qq\q+\big\langle Z(t),b(t, \wt
    X(t),\tilde u(t))-b(t, X(t),
    u(t))\big\rangle_H\)dt\\
    &&\ge \mE\int_0^T\(\dbH(t,\wt
    X(t),\tilde u(t),Y(t),Z(t))-\dbH(t,
    X(t),
    u(t),Y(t),Z(t)) \)dt\nonumber\\
    &&\q-\mE\int_0^T\(\big\langle
    Y(t),a(t,\wt X(t),\tilde u(t))-a(t,
    X(t),
    u(t))\big\rangle_H \nonumber\\
    &&\qq\qq\q +\big\langle Z(t),b(t, \wt
    X(t),\tilde u(t))-b(t, X(t),
    u(t))\big\rangle_H\)dt\nonumber\\
    &&=-\mE\int_0^T\(g(t,\wt X(t),\tilde
    u(t))-g(t, X(t), u(t))\)dt.\nonumber
\end{eqnarray}
On the other hand, the convexity of $h$ implies that
$$
\mE\big\langle h_x(
X(T)),\xi(T)\big\rangle_H\leq \mE h(\wt
X(T))-\mE h(X(T)).
$$
This, together with \eqref{3.4-eq6.1}, yields
$
\cJ( u(\cd))\leq \cJ(\tilde u(\cd)).
$
Since $u(\cd)\in\cU[0,T]$ is arbitrary, the
desired result follows.
\endpf
%

%%%%%%%%%%%%%%%%%%%%%%%%%%%%%%%%%%%%%%%%%%%%%%%%%%%%%%%%%%

\section{Further comments and open problems}

%%%%%%%%%%%%%%%%%%%%%%%%%%%%%%%%%%%%%%%%%%%%%%%%%%%%%%%%%%

\subsection{Further comments}

\begin{itemize}

\item For an earlier version of this notes, we refer to \cite{LZ-2019}. The latter, considering the controllability
and optimal control problems for both stochastic differential equations in finite
dimensions and stochastic evolution equations in infinite dimensions, is a lecture notes for the LIASFMA Hangzhou
Autumn School on ``Control and Inverse Problems for Partial Differential
Equations" which was held during October 17--22, 2016 at Zhejiang University,
Hangzhou, China.

\item Since the seminal paper \cite{Kalman}, controllability and observability became the
basis of Control Theory.  In this notes, we reduce the controllability problem for SEEs to observability problems for BSEEs (by means of the so called ``duality method") without proofs due to the limitation of space. Readers are referred to \cite{LZ3.1} for a detailed introduction to this topic.

\item Besides the ``duality method", there are several other methods to solve controllability problems for deterministic PDEs. Nevertheless, most of them cannot be simply applied to the controllability problems of SPDEs. Let us explain this below briefly for transport equations and parabolic equations.

To establish the exact controllability of deterministic transport equations, there are another two approaches. The first one is to utilize the
explicit formula of their solutions. By this method, for some simple transport
equations, one can explicitly construct a control steering the system from a given
initial state to another given final state, provided that the time is large enough.
It seems that this method cannot be used to solve our stochastic problem since generally we do not have the explicit formula for solutions to the system
\eqref{ch-7-csystem1}. The second one is the
extension method.
However, it seems that it is only valid for time
reversible systems, and therefore, it is not
suitable for the stochastic problems in the framework of It\^o's integral.

In
\cite{FR}, the  null
controllability of the heat equation in one
space dimension is obtained by solving a moment problem,
based on the classical result
on the linear independence of a suitably chosen
families of real exponentials in $L^2(0,T)$.  However, it seems that
it is very hard to employ the same idea to prove
the null controllability of stochastic parabolic equations. Indeed, so far it is
unclear that how to reduce the stochastic null
controllability problem to a suitable moment problem.
For example, let us consider the following
equation
\begin{equation}\label{ch-3-se8-eq1}
    \left\{
    \begin{array}{ll}\ds
        dy(t,x) - y_{xx}(t,x)dt \\
        \ns\ds\q= f(x)u(t)dt+a(x)y(t,x)dW(t)\q
        &\mbox{ in } (0,T]\times (0,1),\\
        \ns\ds y(t,0)=y(t,1)=0 &\mbox{ on } (0,T),\\
        \ns\ds y(0,x)=y_0(x) &\mbox{ in } (0,1).
    \end{array}
    \right.
\end{equation}
Here $y_0\in L^2(0,1)$,
$a\in L^\infty(0,1)$, $f\in L^2(0,1)$, $y$ is the state and $u\in
L_{\dbF}^2(0,T)$ is the control. One can see that it is not easy to reduce the null controllability
problem for the system \eqref{ch-3-se8-eq1} to the usual moment
problem. Under some conditions, it seems possible to reduce it to a stochastic moment problem but this remains to be done.

In \cite{Russell1}, it was
shown that if the wave equation is exactly
controllable for some $T
> 0$ with controls supported in some subdomain $G_0$ of
$G$, then the heat equation is null
controllable for any $T > 0$ with controls
supported in $G_0$.
It seems that one can follow this idea to establish a connection
between the null controllability of stochastic
heat equations and that of stochastic wave equations.
Nevertheless, at this moment the null controllability of
stochastic wave equations is still open, which seems even
more difficult than that of stochastic heat equations.

\item
In this notes, we only give a very brief introduction to
controllability problems for three kinds of
SPDEs. Recently, there are also some works for
controllability problems for other types of SPDEs,
such as \cite{Fu} for stochastic complex
parabolic equations,  \cite{Gao} for
stochastic Kuramoto-Sivashinsky equations, \cite{Liu-Yu1, WuCW} for degenerate stochastic parabolic equations, \cite{Liuxu2} for coupled  stochastic parabolic equations,
\cite{Luqi4} for stochastic Schr\"odinger
equations and \cite{LZ6} for a refined stochastic wave equation.

\item By means of the tools that we developed for solving the stochastic controllability problems, in \cite{Luqi3, Luqi7, LZ2015} we initiated the study of inverse problems for SPDEs (See also \cite{WuCW, Yuan1, Yuan} for further interesting progress), in which the point is that, unlike most of the previous works in this topic, the problem is genuinely stochastic and therefore it cannot be reduced to any deterministic one. Especially, it was found in \cite{LZ2015} that both the formulation of
stochastic inverse problems and the tools to solve them differ
considerably from their deterministic counterpart. Indeed, as the counterexample in \cite[Remark 2.7]{LZ2015} shows, the inverse problem considered therein makes sense only in the stochastic setting!

\item The classical transposition method was introduced by J.-L.
Lions and E. Magenes (\cite{Lions2, Lions3}) to  solve the non-homogeneous boundary value problems (including in particular the boundary control problems) for deterministic
PDEs. The main idea is to interpret
solutions to a less understood equation by means
of that to another well understood one (See also Remark \ref{rem20210121}). Nevertheless, in the stochastic setting the philosophy is quite different. Indeed, the main purpose to introduce the stochastic transposition method is not to solve the non-homogeneous boundary value problems for SPDEs (though it does work for this sort of problems) but to solve the backward stochastic equations (especially the difficult operator-valued BSEEs \eqref{5.5-eq6} and \eqref{op-bsystem3}), in which there exist no boundary conditions explicitly!

\item In this notes, we do not
present the dynamic programming method for
solving optimal control problems of SPDEs. We
refer the readers to  \cite{Fabbri} for a
comprehensive introduction to this topic.

\end{itemize}

\subsection{Open problems}

In our opinion, control theory for
SPDEs is still at its early stage. There are lots of
challenging open problems in this topic. We shall
list below some of them which, in our opinion,
are particularly interesting:

\ss

\begin{itemize}

\item{\bf Controllability of the stochastic parabolic
    equation with one control}

\ss

Note that we introduce two controls $u$ and $v$
in \eqref{heat 1.1}. In view of the
controllability result for deterministic
parabolic equations, it is more natural to use
only one control and consider the following
controlled stochastic parabolic equation (which is
a special case of  \eqref{heat 1.1} with
$v\equiv0$):
\begin{equation}\label{heat---1.1}
    \left\{
    \begin{array}{ll}
        \ds dy-\D y dt =
        \big( a_{1} y+\chi_{G_0}u\big)dt+a_2y dW(t)\q&\hb{ in }Q,\\
        \ns
        \ds y=0&\hb{ on }\Si,\\
        \ns
        y(0)=y_0&\hb{ in } G.
    \end{array}
    \right.
\end{equation}

In order to obtain the null controllability of
\eqref{heat---1.1},  we shall prove that
solutions to the system \eqref{dual heat 1.1}
satisfy the following observability estimate:
\begin{equation}\label{btrulyheat obser1}
    |z(0)|_{L^2(G)}\le \cC
    |z|_{L^2_{\dbF}(0,T;L^2(G_0))},\q
    \forall\;z_T\in L^2_{\cF_T}(\Omega;L^2(G)).
\end{equation}
Unfortunately, at this moment, we do not know how
to obtain
\eqref{btrulyheat obser1}. Indeed, positive results are available only for some special cases
of \eqref{heat---1.1} when the spectral method can be applied (c.f. \cite{Liuxu2, Lu1, Yang}). Even for these special cases, some new phenomenons
appear in the stochastic setting. Inspired by \cite{Liuxu2, Lu1, Yang}, we believe that
one control is enough for the null controllability of a
stochastic parabolic equation. The
main difficulty in establishing
\eqref{btrulyheat obser1}  is that though the
correction   term ``$Z$" plays a ``coercive"
role for the well-posedness of \eqref{dual heat 1.1}, it seems to be a ``bad" (non-homogeneous)
term when one tries to prove \eqref{btrulyheat
    obser1} using the global Carleman estimate. In such case, that term behaves like a nonhomogeneous term and appears  in the right hand side
of the Carleman estimate (See \eqref{h5.2} in Theorem \ref{c1t4}).

\ss

\item {\bf Null/approximate controllability for
    stochastic hyperbolic equations}

\ss

In Section \ref{s-h}, we show that stochastic
hyperbolic equations are not exactly
controllable. It is interesting to study whether
they are null/approximately controllable.
As far as we know, there are no nontrivial results published in this respect.

A possible way to solve these problems is to
establish suitable observability estimate/unique
continuation property for some backward
stochastic hyperbolic equation.   For example,
to get the null controllability of
\eqref{system1}, one should establish the following
observability estimate for the equation \eqref{bsystem1}
with $\tau=T$:
$$
|(z(0),\hat z(0))|_{H_0^1(G)\times L^2(G)}^2
\leq \cC\(\int_{\Si_0}\Big|\frac{\pa
    z}{\pa\nu}\Big|^2d\Si_0 +
|Z|_{L^2_\dbF(0,T;H_0^1(G))}^2\).
$$
However, until now, we can only prove that (see
\cite{LZ6})
$$
\begin{array}{ll}\ds
    |(z(0),\hat z(0))|_{H_0^1(G)\times L^2(G)}^2\\
    \ns\ds
    \leq \cC\(\int_{\Si_0}\Big|\frac{\pa
        z}{\pa\nu}\Big|^2d\Si_0 +
    |Z|_{L^2_\dbF(0,T;H_0^1(G))}^2 + |\wh
    Z|_{L^2_\dbF(0,T;L^2(G))}^2\).
\end{array}
$$
How to get rid of the above extra term $|\wh
Z|_{L^2_\dbF(0,T;L^2(G))}^2$ is an unsolved
problem.

\ss

\item{\bf Controllability problem for one dimensional stochastic
    hyperbolic systems}

\ss

When the system is  linear, it seems possible to
generalize the powerful method in \cite{Li} (which is for deterministic problems) to
get the exact/null controllability results for stochastic
hyperbolic systems in one space dimension.
However, when the system is quasilinear, to our
best knowledge, it is unknown about how to define
the characteristics. Thus, it is unclear how to
study the controllability problem for such sort of nonlinear
systems.

\ss

\item{\bf The existence of an optimal feedback control operator
    for stochastic linear\linebreak quadratic optimal problems with less restriction}

\ss

To guarantee the existence of an optimal
feedback control operator for stochastic linear quadratic optimal problems, we need
very restrictive conditions (e.g., \cite{LZ5}). It is interesting to
relax these conditions. We believe, very likely, the theory for forward-backward stochastic evolution equations and
Malliavin Calculus are helpful to handle this problem.

\ss

\item{\bf Linear quadratic optimal problems for SEEs with unbounded control operators}

\ss

In Sections \ref{s-lq-o}--\ref{sec-slq-d}, we assume that $B$ , $D$ and $Q$ are bounded, linear operator-valued functions. In many concrete control systems, it is very common
that both the control and the observation operators are
unbounded w.r.t. the state spaces.
Typical examples are systems governed by SPDEs in which the actuators
and the sensors act on lower-dimensional
hypersurfaces or on the boundary of a spatial
domain. The unboundedness of the control operators and the
observation operators leads to substantial
technical difficulties even for the formulation
of the state spaces.   To study such systems, one needs to make
some further assumptions, such as the semigroup $\{S(t)\}_{t\geq 0}$ has some smoothing effect. This is studied when $D=0$ (e.g., \cite{Lasiecka 2017}).
In the deterministic
setting, to overcome this sort of difficulties, people introduced the notions of the
admissible control operator and the admissible
observation operator (e.g., \cite{Salamon}).
On the other hand, people
introduced the notion of well-posed linear
system (with possibly unbounded control and observation operators), which satisfies, roughly speaking, that
the map from the input space to the output one
is bounded (e.g., \cite{Salamon}). The
well-posed linear systems form a very general
class whose basic properties are rich enough to
study LQ problems. The concept of
well-posed linear system was generalized to the stochastic setting in
\cite{Luqi8} but it seems that there are many things to be done to develop a satisfactory theory for stochastic well-posed linear control
systems.

\ss

\item{\bf Optimal control problems for SEEs with state/endpoint constraints}

\ss

We do not consider state/endpoint constraints
for the optimal control problem in this notes.
Some results can be found in \cite{FL}. However,
the constraints considered in \cite{FL} are not so strong, and hence many
interesting cases cannot be covered, say, the endpoint constraint such that
$X(T)\in \cS$,  a.s. for some set $\cS$.

\ss

\item{\bf High order necessary conditions for optimal controls}

\ss

Pontryagin type maximum principle is a first
order necessary condition for optimal controls.
In addition to the first-order necessary
conditions, some higher order necessary
conditions should be established to distinguish
optimal controls from the candidates which
satisfy the first order necessary conditions,
especially when the optimal controls are
singular, i.e., optimal control satisfy the
first order necessary conditions trivially. For
instance, when the Hamiltonian corresponding to
optimal controls is equal to a constant in a
subset of the control region or the gradient and
the Hessian (w.r.t. the control
variable $u$) of the corresponding Hamiltonian
vanish/degenerate. In those cases, the first
order necessary conditions are not enough to
determine the optimal controls. When the control
domain $U$ is convex, some pointwise second
order necessary conditions are obtained in
\cite{LZZ}. When $U$ is not convex, some
integral type second order necessary conditions
are obtained in \cite{FL}. However, compared
with the deterministic counterpart, the results in
\cite{FL, LZZ} are not satisfactory.

\ss

\item{\bf Time-inconsistent optimal control problems for SEEs }

\ss

In recent years, there are many works on the time-inconsistent optimal control problems for stochastic differential equations (See \cite{Yong2014} and the references cited therein). It would be quite interesting to study the same problems but for stochastic evolution equations in infinite dimensions. As far as we know, there are only very few publications on this topic (e.g., \cite{DL1}).

\ss

\item{\bf Optimal control problems for mean-field SEEs }

\ss

Mean field approximation is extremely useful when describing macroscopic phenomena from microscopic overviews. Over
the last years, due to various applications
in economics, finance and physics, there has
been a growing interest in control problems for
mean-field stochastic differential equations (e.g., \cite{Bensoussan 2013, Carmona-Delarue 2018} and the rich references therein). However, the publications on control problems for
mean-field SEEs are quite limited. Indeed, to the best of our
knowledge, \cite{Oksendal 2019, Sulem
    2018, Luqi10, Tang2019} are the only papers in dealing
with optimal control problems for this sort of systems.

\ss

\item{\bf Control problems for other types of stochastic control systems}

\ss

The control problems considered in this paper and mentioned above make
sense also for forward-backward stochastic evolution equations, stochastic evolution equations driven by fractional Brownian motions or $G$-Brownian motions, forward-backward
doubly stochastic evolution equations, stochastic Volterra integral equation in infinite dimensions or stochastic evolution equations driven by general martingales (even with incomplete information).  As far as we
know, all of them remain to be done, and very likely people might
need to develop new tools to obtain interesting new results.

\end{itemize}

%\section*{Acknowledgement}

%%%%%%%%%%%%%%%%%%%%%%%%%%%%%%%%%%%%%%%%%%%%%%%%%%%%%%%%%%%%%%%%%%%%%%%%%%%%%%%%%%%%

% {\bf \large
%Acknowledgments.}

%For acknowledgements section, please don't number the section, please begin it with \section*{Acknowledgements}
%\section*{Acknowledgments} We would like to thank you for \textbf{following
%the instructions above} very closely in advance. It will definitely
%save us lot of time and expedite the process of your paper's
%publication.

% You may incorporate your references as follows in your main tex file.
% Using BibTex is not recommended but can be handled.

\end{document}